\documentclass[letterpaper]{amsart}

\usepackage{stmaryrd}
\usepackage{latexsym}
\usepackage{amssymb}
\usepackage{amscd}
\usepackage{amsopn}
\usepackage{amsfonts}
\usepackage{amsmath}
\usepackage{graphicx}
\usepackage{psfrag}
\usepackage{amsthm}
\usepackage{enumerate}

\def \eps {\varepsilon}
\def \bfx {\mathbf{x}}

\def \bfe {\mathbf{e}}
\def \bfn {\mathbf{n}}

\def \bfh {\mathbf{h}}
\def \bfp {\mathbf{p}}

\def \bfv {\mathbf{v}}
\def \bfm {\mathbf{m}}

\def \bfd {\mathbf{d}}
\def \bfu {\mathbf{u}}

\def \bfalpha {\boldsymbol{\alpha}}

\def \bfomega {{\boldsymbol{\omega}}}

\def \bfDelta {\boldsymbol{\Delta}}

\def \cR {\mathbb{R}}
\def \cN {\mathbb{N}}
\def \cQ {\mathbb{Q}}
\def \cC {\mathbb{C}}

\def \cM {\mathbb{M}}
\def \cZ {\mathbb{Z}}

\def \cS {\mathbb{S}}
\def \vD {\mathfrak{D}}
\def \vF {\mathcal{F}}
\def \vI {\mathcal{I}}

\def \Strip {\mathcal{Q}}

\def \vN {\mathcal{N}}
\def \vO {\mathcal{O}}

\def \vL {\mathcal{L}}
\def \vJ {\mathcal{J}}

\newtheorem{Theorem}{Theorem}[section]
\newtheorem*{MTheorem}{Main Theorem}
\newtheorem*{Conjecture}{Conjecture}
\newtheorem{Definition}[Theorem]{Definition}
\newtheorem{Proposition}[Theorem]{Proposition}
\newtheorem{Lemma}[Theorem]{Lemma}
\newtheorem{Corollary}[Theorem]{Corollary}

\newtheorem{Rem}[Theorem]{Remark}
\newtheorem{Exam}[Theorem]{Example}
\newenvironment{Remark} {\begin{Rem} \rm  }{\end{Rem}}
\newenvironment{Example} {\begin{Exam} \rm   }{\end{Exam}}

\def \[ {\llfloor}
\def \] {\rrfloor}

\def \cm {\mathbf{m}}
\def \supp {{\mathrm{supp}\, }}

\def \gdc {\mathrm{gdc}\, }
\def \conv {\mathrm{conv}\, }

\def \< {\langle}
\def \> {\rangle}

\def \medge {\mathfrak{e}}

\def \linteg {\lfloor}
\def \rinteg {\rfloor}

\def \For {\cR[[\bfx]]}

\def \lists {\mathbb{L}}

\def \linef {L}

\def \mainlist {\Upsilon}

\def \lex {\mathrm{lex}}

\def \vheight {\mathfrak{h}}

\def \Inv {\mathrm{inv}}
\def \Invp {\mathrm{inv}_1}
\def \Invs {\mathrm{inv}_2}

\def \medge {\mathfrak{e}}

\def \vZ {\mathfrak{z}}

\def \IM {\mathbb{I}}

\def \Ze {\mathrm{Ze}}
\def \Elem {\mathrm{Elem}}
\def \NElem {{\mathrm{NElem}}}

\def \Ax {\mathrm{Ax}}

\renewcommand \emptyset {\varnothing}

\def \hG {\widehat{G}}

\def \Blx {\mathrm{Bl}_x}
\def \Bly {\mathrm{Bl}_y}
\def \Blz {\mathrm{Bl}_z}

\def \Stab {\mathrm{St}\, }

\def \transl {\mathrm{tr}}

\def \Loc {\mathrm{Loc}}

\def \Bad {\mathrm{Bad}}

\def \Badch {\mathcal{B}}

\def \TreeBd {\mathrm{Tr}\mathcal{B}}

\def \Mult {\mathrm{Mult}}
\def \Loc {\mathrm{Loc}}

\def \dist {\mathrm{dist}}

\def \Eq {\mathrm{Eq}}

\def \gG {\mathcal{G}}
\def \lG {\mathfrak{G}}

\def \ffG {\widehat{\mathrm{G}}}

\def \PrG {\mathrm{Pr}\gG}

\def \Nmaps {\mathrm{New}}

\begin{document}
\title[Resolution of Singularities of Vector Fields in Dimension Three]{Resolution of Singularities of Real Analytic Vector Fields in Dimension Three} 
\author{Daniel Panazzolo}
\thanks{This work has been partially supported by the CNPq/Brasil Grant 205904/2003-5 and Fapesp
Grant 02/03769-9.}
\begin{abstract}
Let $\chi$ be an analytic vector field defined in a real analytic manifold of
dimension three.
We prove that all the singularities of $\chi$ can be made elementary 
by a finite number of blowing-ups in the ambient space.
\end{abstract}
\maketitle
\tableofcontents
\section{Introduction}
\subsection{Main Result}
Let $\chi$ be an analytic vector
field defined on a real analytic manifold $M$. We shall say that $\chi$ is
{\em elementary} at a point $p \in M$ if one of the following
conditions holds:
\begin{itemize}
\item[(i)] (Non-singular case) $\chi ( p ) \ne 0$, or
\item[(ii)] (Singular case) $\chi(p) = 0$ and the {\em Jacobian map}
\begin{eqnarray*}
D\chi(p) : \cm_p/\cm_p^2 & \longrightarrow & \cm_p/\cm_p^2,  \quad (\cm_p
\subset \vO_{M,p} \mbox{ is the maximal ideal} ) \\
 {[} {g} {]} \;\;\; & \longmapsto & [ \chi(g) ]
\end{eqnarray*}
has at least one nonzero eigenvalue (here $\chi(\cdot)$ denotes the
action of $\chi$ as a derivation in $\vO_p$).
\end{itemize}
If we fix a local coordinate system $(x_1,\ldots,x_n)$ for
  $M$ at $p$ and write
$$\chi = a_1
  \frac{\partial}{\partial x_1} + \cdots + a_n
  \frac{\partial}{\partial x_n}$$
then the Jacobian map is given by the real matrix
$$ D\chi( p ) =
\left(\begin{array}{c}
\frac{\partial a_i}{\partial x_j} ( 0 )
   \end{array}\right)_{i, j = 1, \ldots, n }
$$
We say that $\chi$ is {\em reduced} if $\mathrm{gdc}(a_1,\ldots,a_n) =
1$ at each point $p \in M$ (this implies that the set $\Ze(\chi) = \{q \in M \mid
\chi(q) = 0\}$ has
codimension strictly greater than one).

Let us enunciate our main result.  We briefly define the necessary concepts
and postpone the details to the next section.

A singularly foliated manifold is a 4-uple $\cM =
(M,\mainlist,\vD,\linef)$ where
\begin{itemize}
\item[(i)] $M$ is a real analytic three-dimensional manifold with corners;
\item[(ii)] $\mainlist \in \lists$ is an ordered list of natural numbers;
\item[(iii)] $\vD = \vD_\mainlist$ is a $\mainlist$-tagged divisor on
  $M$ with normal crossings;
\item[(iv)] $\linef$ is a singular orientable analytic line field on $(M,\vD)$
  which is $\vD$-preserving.
\end{itemize}
At each point $p \in M$, the line field
$\linef$ is locally generated by an analytic vector field $\chi_p$ which is
tangent to the divisor $\vD$.  Such local generator is uniquely defined up to multiplication by a
strictly positive analytic function.

We say that the singularly foliated manifold 
$\cM$ is {\em elementary} at a point $p \in M$ if the local generator
$\chi_p$ is an elementary vector field at $p$. The complement of the set of
elementary points in $M$ will be denoted by $\NElem(\cM)$.

A singularly foliated manifold 
$\cM$ is said to be {\em elementary} if $\NElem(\cM) = \emptyset$.
\begin{MTheorem}
Let $\chi$ be a reduced analytic vector field defined in a real analytic 
three-dimensional manifold $M$ without boundary.
Then, for each relatively compact set $U \subset M$, 
there exists a finite sequence of weighted blowing-ups
\begin{equation}\label{blowupseq}
(U,\emptyset,\emptyset,\linef_\chi|_U) =: \cM_0 \stackrel{\Phi_1}{\longleftarrow} \cM_1 
\stackrel{\Phi_2}{\longleftarrow}
\cdots \stackrel{\Phi_n}{\longleftarrow} \cM_n 
\end{equation}
such that the resulting singularly foliated manifold $\cM_n$ is elementary.
Moreover, the center $Y_i$ of the blowing-up $\Phi_i$ is a smooth analytic  
subset of $\NElem(\cM_i)$, for each $i = 0,\ldots,n-1$.
\end{MTheorem}
In the above enunciate, $\linef_\chi$ denotes the singular orientable line 
field which is associated to the vector field $\chi$.
\subsection{Previous Works}
The theorem of resolution of singularities for vector fields in dimension two 
was present in the work of Bendixson \cite{Be}.  The first complete proof of this result has been given by  Seidenberg in \cite{S}.

In \cite{Pe}, Pelletier gives an alternative proof of such result
through the use of the weighted blowing-ups.

In the book \cite{C1}, Cano proves a result of local reduction of singularities in the formal context for complex three-dimensional vector fields.

The paper \cite{Sa} studies generic equireduction of singularities for vector
fields in arbitrary dimension.

In a recent paper \cite{CMR}, the authors prove a local uniformization theorem for 
analytic vector fields in dimension three.  Their proof is based on the analysis of 
valuations defined by non-oscillating subanalytic integral curves.

The literature on the Newton Polyhedron and its applications is extensive.
For some results related to the use of Newton Polyhedron in resolution
of singularities, we refer the reader to \cite{H1}, \cite{H2} and
\cite{Y}.

In the book \cite{B}, Bruno uses the Newton Polyhedron and Normal Form Theory to describe many explicit algorithms 
for studying the asymptotic behavior of integral curves of vector fields near elementary and nonelementary singular points. 
\subsection{Overview of the Paper}\label{subsect-overview}
The proof of the Main Theorem consists of two parts: the description of the local strategy for resolution of singularities given section~\ref{sect-LocalTheory} and the proof that such local strategy can be {\em globalized}, which will be detailed in section~\ref{sect-GlobalTheory}.

Let us briefly describe the ingredients used in the central result of the paper: the Theorem on Local Resolution of Singularities.

For definiteness, we assume here that $M = (\cR^3_{(x,y,z)},0)$ and that the origin is contained in a divisor with normal crossings $\vD$ which is given either by $\{x = 0\}$ or by $\{xy = 0\}$.  We further
assume that the reduced vector field $\chi$ defined in $M$ is tangent to the divisor $\vD$ and that the origin is a nonelementary singular point.  Finally, we assume that the vertical axis $\{x = y = 0\}$ is not entirely contained in the set $\NElem$ of nonelementary points. 

Using the logarithmic basis $\{x \frac{\partial}{\partial x}, y \frac{\partial}{\partial y}, 
z \frac{\partial}{\partial z}\}$ we can write
$$
\chi = f x \frac{\partial}{\partial x} + g y \frac{\partial}{\partial y} + g z \frac{\partial}{\partial z}
$$
where $f,gy$ and $hz$ are germs in $\cR\{x,y,z\}$. The {\em Newton polyhedron} of $\chi$ (with respect
to the coordinates $(x,y,z)$) is the convex polyhedron
$$
\vN = \conv(\supp(f,g,h)) + \cR^3_{\ge 0}
$$
where $\conv(\cdot)$ denotes the operation of convex closure, $\supp(f,g,h) \subset \cZ^3$ is the set of integer points $\bfv = (v_1,v_2,v_3) \in \cZ^3$ such that the monomial 
$x^{v_1}y^{v_2}z^{v_3}$ has a nonzero coefficient in the Laurent expansion of either $f$, $g$ or $h$; and the $+$ sign is the usual Minkowski sum of convex polyhedrons.

\begin{figure}[htb]
\psfrag{h}{\small $\bfh$}
\psfrag{n}{\small $\bfn$}
\psfrag{v2}{\small $v_2$}
\psfrag{v3}{\small $v_3$}
\psfrag{-1,0}{\tiny $(-1,0)$}
\psfrag{0,-1}{\tiny $(0,-1)$}
\psfrag{(a)}{\small $(a)$}
\psfrag{(b)}{\small $(b)$}
\psfrag{(c)}{\small $(c)$}
\begin{center}
\includegraphics[height=5cm]{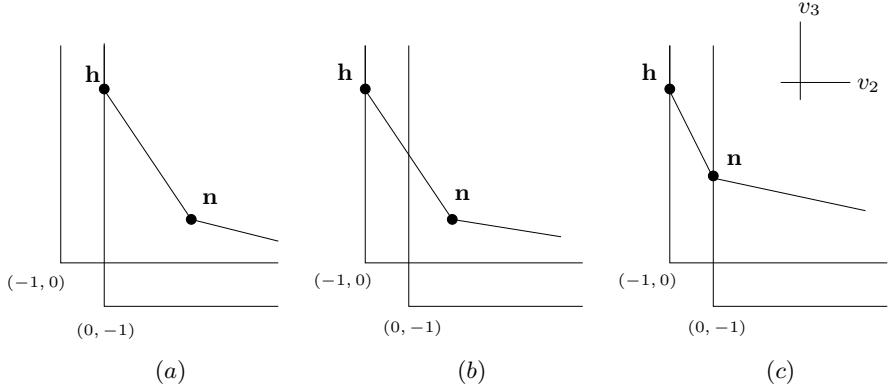}
\end{center}
\caption{Regular and nilpotent configurations.}
\label{fig-planofpaper}
\end{figure}

The {\em higher vertex} of $\vN$ is the vertex $\bfh \in \vN$ which is minimal with respect to the lexicographical ordering in $\cR^3$.  By the hypothesis, it follows that such vertex has the form $\bfh = (0,h_2,h_3)$ for 
some integers $h_2,h_3 \in \cZ_{\ge -1}$.  Moreover, the intersection of $\vN$ with the plane $\{\bfv \in \cR^3 \mid v_1 = 0\}$ is in one of the situations shown in figure~\ref{fig-planofpaper}.

\begin{figure}[htb]
\psfrag{m}{\small $\bfm$}
\psfrag{v1}{\small $v_1$}
\psfrag{v2}{\small $v_2$}
\psfrag{v3}{\small $v_3$}
\psfrag{ml}{\small $\bfm^\prime$}
\psfrag{Nl}{\small $\vN^\prime$}
\psfrag{1/2}{\small $\frac{1}{2}$}
\begin{center}
\includegraphics[height=5cm]{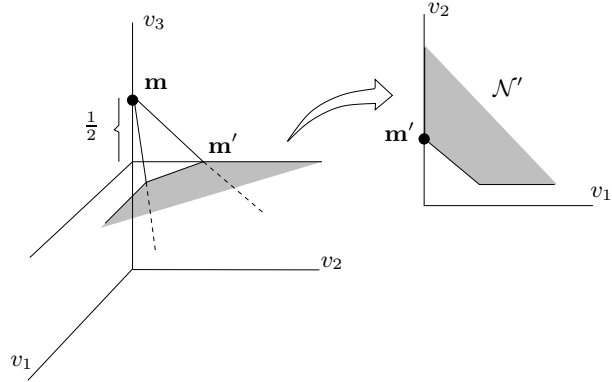}
\end{center}
\caption{The derived polygon.}
\label{fig-planofpaper2}
\end{figure}

Referring to figure~\ref{fig-planofpaper}, the configurations $(a)$ and $(b)$ are called {\em regular} and the configuration $(c)$ is called
{\em nilpotent}.  Like indicated in the figure, we define the {\em main vertex} $\bfm = (m_1,m_2,m_3)$  by $\bfm = \bfh$ in cases $(a)$ and $(b)$ and by $\bfm = \bfn$ in case $(c)$.  Now, we consider the intersection
$$
\vN^\prime = \vN \cap \left\{\bfv \in \cR^3 \mid v_3 = m_3 - \frac{1}{2} \right\}
$$
and call the polygon $\vN^\prime$ the {\em derived polygon} (see figure~\ref{fig-planofpaper2}).

The derived polygon has some similarities with the {\em characteristic polygon} introduced by 
Hironaka~\cite{H2} in his proof of the resolution of singularities for excellent surfaces.  However, there are some essential differences which will be discussed in subsection~\ref{subsect-derivedpolygon}.

Let us denote by $\bfm^\prime = (m_1^\prime,m_2^\prime,m_3 - 1/2)$ the minimal vertex of $\vN^\prime$
(with respect to the lexicographical ordering) and write the displacement vector $\bfm^\prime - \bfm$ 
as $\frac{1}{2}(\Delta_1,\Delta_2,-1)$, for some nonzero rational vector $\Delta = (\Delta_1,\Delta_2) \in \cQ^2$.  

\begin{figure}[htb]
\psfrag{chi}{$\Inv = (3,0,4,0,0,5!\, 2/3)$}
\psfrag{m=(0,-1,4)}{\small $\bfm=(0,-1,4)$}
\psfrag{Delta}{\small $(0,2/3,-1)\; \Rightarrow \;(\Delta_1,\Delta_2) = (0,2/3)$}
\psfrag{m-dd}{\small $m_3  - \frac{1}{\Delta_2}=5/2 \; \Rightarrow \; \vheight = 3$}
\psfrag{(0,1,1)}{\small $(0,1,1)$}
\psfrag{(0,4,-1)}{\small $(0,4,-1)$}
\psfrag{(1,-1,1)}{\small $(1,-1,1)$}
\psfrag{v1}{\small $v_1$}
\psfrag{v2}{\small $v_2$}
\psfrag{v3}{\small $v_3$}
\begin{center}
\includegraphics[height=6.5cm]{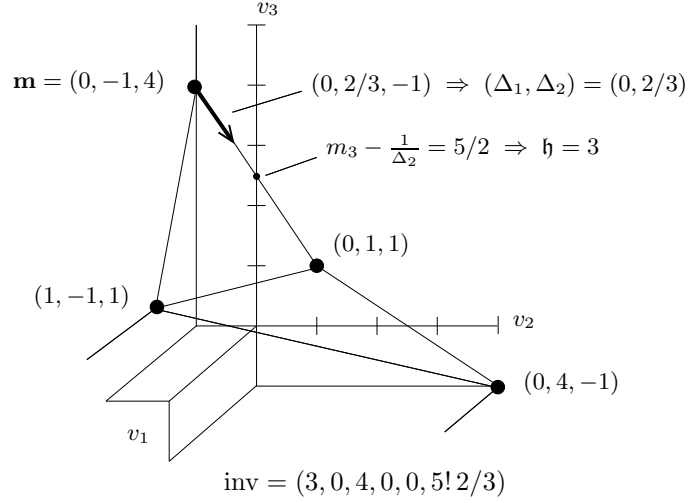}
\end{center}
\caption{Newton polyhedron for $\chi = x^2y\frac{\partial}{\partial x} + 
(z^4 + xz) \frac{\partial}{\partial y} + y^4\frac{\partial}{\partial z}$.}
\label{fig-planofpaper3}
\end{figure}

The {\em main invariant} for the vector field $\chi$ (with respect to the coordinates $(x,y,z)$) is given by the 6-uple of natural numbers
$$
\Inv = (\vheight,m_2 + 1,m_3,\#\iota - 1, \lambda \Delta_1,\lambda \max\{0,\Delta_2\})      
$$
where $\lambda = (m_3 + 1)!$, $\#\iota \in \{1,2\}$ is the number of local irreducible components of the divisor at the origin and the {\em virtual height} $\vheight$ is a natural number defined as follows
$$
\vheight =
\left\{
\begin{matrix}
\lfloor m_3  + 1 - \frac{1}{\Delta_2} \rfloor , &\mbox{ if }m_2 = -1
\mbox{ and }\Delta_1 = 0 \hfill\cr \cr m_3, &\mbox{ if }m_2 =
0\mbox{ or }\Delta_1 > 0\hfill
\end{matrix}
\right.
$$
where $\lfloor \alpha \rfloor := \max \{ c \in \cZ \mid c \le \alpha \}$ (see  figure~\ref{fig-planofpaper3} for an example).

In subsection~\ref{subsect-basicfaceedgeprep} we shall introduce the fundamental notion of {\em stable coordinates}.  Roughly speaking, if we start with a system of local coordinates $(x,y,z)$ as above, we obtain a new system of coordinates $(\widetilde{x},\widetilde{y},\widetilde{z})$ by an analytic change of coordinates of the form
\begin{equation}\label{stab-planpaper}
\widetilde{x} = x, \quad \widetilde{y} = y + G(x), \quad \widetilde{z} = z + F(x,y),
\end{equation}
in such a way that the main invariant $\Inv$, when computed with respect to the new
coordinates $(\widetilde{x},\widetilde{y},\widetilde{z})$, has nice analytic properties such as being an upper-semicontinuous function.  Moreover, we shall see that $\Inv$ is an {\em intrinsic} object attached to the germ of vector field $\chi$, up to fixing an additional geometric structure on the ambient space called an {\em axis} (see subsection~\ref{subsect-axisdef}).  

The local strategy of reduction of singularities will be readout of the Newton polyhedron $\vN$ and main invariant $\Inv$, {\em provided} that these objects are computed with respect to a stable system of coordinates.

The notion of stable coordinates is similar to the notions of {\em well-prepared} and {\em very well-prepared} systems of coordinates, as defined by Hironaka \cite{H2} in the context of function germs.  However, new difficulties appear in the context of vector fields, since the action of the Lie group of coordinates changes given in (\ref{stab-planpaper}) is much harder to study in this situation.  An example of new phenomena is the appearance of the so-called {\em resonant configurations}, described in subsection~\ref{subsect-resonant}.

Let us now briefly introduce a second ingredient of our proof: the notion of weighted blowing-up.  Given a vector of nonzero natural numbers $\bfomega = (\omega_1,\omega_2,\omega_3) \in \cN^3_{>0}$, the {\em $\bfomega$-weighted blowing-up} (with respect to the coordinates $(x,y,z)$) is the proper analytic surjective map given by
\begin{eqnarray*}
\Phi_\bfomega: \cS^{2} \times \cR^+  &\longrightarrow& \cR^3\\
((\overline{x},\overline{y},\overline{z}),\tau) &\longmapsto&
(x,y,z) = (\tau^{\omega_1}\overline{x},\tau^{\omega_2}\overline{y},\tau^{\omega_3}\overline{z})
\end{eqnarray*}
Similarly, for a weight-vector of the form $\bfomega = (\omega_1,0,\omega_3)$ with $\omega_1,\omega_3$ nonzero, we define 
\begin{eqnarray*}
\Phi_\bfomega: \cS^{1} \times \cR^+ \times \cR &\longrightarrow& \cR^3\\
((\overline{x},\overline{z}),\tau,y) &\longmapsto&
(x,y,z) = (\tau^{\omega_1}\overline{x},y,\tau^{\omega_3}\overline{z}).
\end{eqnarray*}
We define analogously the blowing-ups for weight-vectors $\bfomega$ of the form $(0,\omega_2,\omega_3)$
or $(\omega_1,\omega_2,0)$.  Notice that the blowing-up center in these four case is given  
respectively by $\{x=y=z=0\}$, $\{x=z=0\}$, $\{y=z=0\}$ and $\{x=y=0\}$.

Suppose now that $(x,y,z)$ is a stable system of coordinates at the origin, and let
$\Inv \in \cN^6$ be the corresponding main invariant.  Starting from  subsection~\ref{subsect-Newtoninvarlocres}, we prove the following result on local resolution of singularities for vector fields:  There exists a choice of weight-vector $\bfomega$ for which
the strict transform $\widetilde{\chi}$ of the vector field $\chi$ under the $\bfomega$-weighted blowing-up
$$
\Phi_{\bfomega}: \widetilde{M} \rightarrow M
$$
is such that, for each nonelementary point $\widetilde{p} \in \widetilde{M}\cap\Phi^{-1}(0)$, and each choice of a stable system of coordinates $(\widetilde{x},\widetilde{y},\widetilde{z})$ with center at $\widetilde{p}$, the  corresponding main invariant $\widetilde{\Inv}$ is such that
$$
\widetilde{\Inv} <_\lex \Inv
$$
where $<_\lex$ is the usual lexicographical ordering in $\cN^6$.

The second part of the proof of the Main Theorem is given in section~\ref{sect-GlobalTheory}.  There we show that the local strategy for resolution of singularities described in the previous paragraph can be {\em globalized}.  

To prove this, the main ingredient is the fact that both the main invariant $\Inv$ and the choice of weight-vector $\bfomega$ are independent of the given stable system of coordinates.  Moreover, we shall see in subsection~\ref{subsect-uppersemicont} that $\Inv$ is an upper-semicontinous function
on the set $\NElem$ of nonelementary singular points of the vector field.  

Based on these facts, the global strategy of resolution is similar to the one presented by Cano in 
\cite{C2}, based on the notion of generic equireducibility and bad points.  The main distinction is
the fact that our strategy leads to a unique choice of local center and an {\em a priori} absence of cycles, due to a conveniently chosen enumeration of the exceptional divisors.  Moreover, we can guarantee that the blowing-up centers are always contained in the
set of nonelementary points.
\subsection{An Example}\label{example-sanz}
The following example was communicated to me by F.\ Sanz.  
It justifies the use
of \underline{\em weighted} blowing-ups in the resolution of singularities
of vector fields in dimension greater than $2$.
\begin{Example}
Consider the nonelementary germ of vector field
\begin{equation}\label{examplesanz}
\chi = x \left( x \frac{\partial}{\partial x} - 
\alpha y \frac{\partial}{\partial y} - 
\beta z \frac{\partial}{\partial z}\right) + xz \frac{\partial}{\partial y} + 
(y - \lambda x) \frac{\partial}{\partial z}
\end{equation}
for some real constants $\alpha,\beta \ge 0$ and $\lambda > 0$.
We claim that such germ cannot be simplified  by 
any sequence of {\em homogeneous blowing-ups} (i.e.\ blowing-ups with 
weight equal to one) with center contained in the set of singularities.\\
In fact, since the blowing-up center $Y$ is contained in the set of singularities, there are two
possible blowing-up strategies:
\begin{enumerate}[(a)]
\item Blow-up with center at the point $Y = \{x = y = z = 0\}$.
\item  Blow-up with center at the curve $Y = \{x = y = 0\}$.
\end{enumerate}
In the case (a), the $x$-directional blowing-up is given by
$$x = \widetilde{x}, \quad y = \widetilde{x}\,\widetilde{y}, \quad
z = \widetilde{x}\,\widetilde{z}$$
Therefore, the blowing-up of $\chi$ will be given by the following expression
(dropping the tildes to simplify the notation)
$$
\chi = x \left( x \frac{\partial}{\partial x} - 
\alpha^\prime y \frac{\partial}{\partial y} - 
\beta^\prime z \frac{\partial}{\partial z}\right) + xz \frac{\partial}{\partial y} + 
(y - \lambda) \frac{\partial}{\partial z}
$$
where $\alpha^\prime = \alpha + 1$ and $\beta^\prime = \beta + 1$.
If we make the translations $\widetilde{y} = y - \lambda$ and 
$\widetilde{z} = z - \alpha^\prime\lambda$, we obtain 
(dropping again the tildes)
\begin{equation}\label{examplesanz1}
\chi = x \left( x \frac{\partial}{\partial x} - \alpha^\prime 
y \frac{\partial}{\partial y} - \beta^\prime z 
\frac{\partial}{\partial z}\right) + xz \frac{\partial}{\partial y} + 
(y - \lambda^\prime x)\frac{\partial}{\partial z}
\end{equation}
where $\lambda^\prime = \alpha^\prime \beta^\prime \lambda$.  Notice
that the vector field
(\ref{examplesanz1}) can be obtained from (\ref{examplesanz}) 
simply by making the replacement of the constants
$$(\alpha,\beta,\lambda) \rightarrow (\alpha^\prime,\beta^\prime,\lambda^\prime)$$
In the case (b), the $x$-directional blowing-up is given by
$$
x = \widetilde{x}, \quad y = \widetilde{x}\,\widetilde{y}, \quad
z = \widetilde{z}
$$
and we get (dropping the tildes)
$$
\chi = x \left( x \frac{\partial}{\partial x} - 
\alpha^\prime y \frac{\partial}{\partial y} - 
\beta z \frac{\partial}{\partial z} \right) + z \frac{\partial}{\partial y} + 
x(y - \lambda) \frac{\partial}{\partial z}
$$
where $\alpha^\prime = \alpha + 1$.  After the translation 
$\widetilde{y} = y - \lambda$, we obtain
\begin{equation}\label{examplesanz2}
\chi = x \left( x \frac{\partial}{\partial x} - 
\alpha^\prime y \frac{\partial}{\partial y} - 
\beta z \frac{\partial}{\partial z} \right) + 
(z - \overline{\lambda} x)  \frac{\partial}{\partial y} + 
xy \frac{\partial}{\partial z}
\end{equation}
where $\overline{\lambda} = \alpha^\prime\lambda$. Notice
that the vector field
(\ref{examplesanz2}) can be obtained from (\ref{examplesanz}) 
simply by making the replacement of the constants
$$(\alpha,\beta,\lambda) \rightarrow (\alpha^\prime,\beta,\lambda^\prime)$$\\
and interchanging the roles of $y$ and $z$ variables.  This proves that no improvement
has been made neither by the blowing-up (a) nor by the blowing-up (b).
\end{Example}
\begin{Remark}
Up to some additional computations, we can further prove that no improvement can be made if we choose, as blowing-up centers, {\em arbitrary} analytic curves which are left invariant by the vector field.  
\end{Remark}
\subsection{Acknowledgments}
I would like to thank Felipe Cano for his encouragement and suggestions.  I also thank
Joris van der Hoeven for useful discussions on the subject.\\
Some parts of this work
have been done at the Universidad de Valladolid, Universit\' e de Bourgogne and IMPA. 
I thank these institutions for their hospitality.\\
Finally, I would like to thank the anonymous referee for the many valuable comments and suggestions, which greatly improved the presentation of this work.
\subsection{Notations}
\begin{itemize}
\item $\cN = \{n \in \cZ \mid n \ge 0\}$, $\cN_{\ge k} = \{n \in \cZ \mid n
      \ge k\}$.
\item $\cR^+ = \{x \in \cR \mid x \ge 0\}$, $\cR^* = \{x \in \cR \mid x \ne 0\}$
\item $\cR_{> \alpha} = \{x \in \cR \mid
  x > \alpha\}$, $\cR_{\ge \alpha} = \{x \in \cR \mid
  x \ge \alpha\}$.
\item $\cR^n_{s} = (\cR^+)^s \times \cR^{n-s}$.
\item $\overline{\cR} = \cR \cup \{\infty\}$ is the extended  field of real
      numbers, with the usual extended ordering relation (we define similarly
$\overline{\cQ}$ and $\overline{\cZ}$).
\item $\mathrm{Mat}(n,\cR)$ is the set of $n\times n$ real
  matrices.
\item $\cR[[\bfx]]$ is the ring of real formal series in the variables
  $\bfx$.
\item $\lists$ is the set of all {\em reverse ordered} lists of
  natural numbers.  A typical element $\iota \in \lists$ is written
$$
\iota = [i_1,\ldots,i_k],\quad \mbox{ where }i_1,\ldots, i_k \in
\cN\mbox{ and }i_1 > i_2 > \cdots > i_k$$
$\# \iota = k$ denotes the length of the list $\iota$.
\item For two lists
$\iota,\rho \in \lists\,$, we denote by $\iota \cup \rho$,
  $\iota \cap \rho$ and $\iota \setminus \rho$ the new lists which are
  obtained by the usual operations of
  concatenation, intersection and difference (for instance, $[3,2,1]
  \cup [5,3] = [5,3,2,1]$, $[5,3,2] \cap [3,1] = [3]$ and $[5,4,3,2]
  \setminus [5,4,1] = [3,2]$).
\item For $\bfu,\bfv \in \cR^n$, the notation $\bfv <_\lex \bfu$ indicates that
$\bfv$ is {\em lexicographically} smaller than $\bfu$, i.e. the relation
$$
\bigvee_{i=0}^{n-1} \quad [(v_1,\ldots,v_{i}) = (u_1,\ldots,u_{i})\; \wedge \; v_{i+1} < u_{i+1}]
$$
holds.
\item For $\alpha \in \cR$, $\lfloor \alpha \rfloor = \max\{n \in \cZ \mid n \le \alpha\}$ and 
$\lceil \alpha \rceil = \min\{n \in \cZ \mid n \ge \alpha\}$.
\end{itemize}

\section{Blowing-up and Singularly Foliated Manifolds}\label{sect-SFM}
\subsection{Manifolds with Corners}
We shall work on the category of analytic manifolds with corners.  Recall that a
$n$-dimensional manifold with corners $M$ is a paracompact topological space which is locally
modeled by 
$$
\cR^n_s = \{(x_1,\ldots,x_n) \in \cR^n \mid x_1 \ge 0, \ldots x_s \ge 0\}
$$ 
A {\em local chart} (or local coordinate system) at a point $p \in M$ is a pair $(U,\phi)$ such that 
$U \subset M$ is an open neighborhood of $p$ and $\phi: U \rightarrow \cR^n_s$ is a diffeomorphism with
$\phi(p) = 0$.  Note that $\phi(U \cap \partial M)$ is mapped to 
$\partial \cR^n_s := \cup_{i=1}^s \{x_i = 0\}$.  

For notational simplicity, we shall sometimes omit the subscript $s$ when we refer to the space $\cR^n_s$.  

The number $b(p) := s$ is the number of boundary components which meet at
$p$ (it is independent of the choice of local chart).  
Note that $\partial M := \{ p \in M \mid b(p) \ge 1\}$.

In this work, we shall say that a subset $N \subset M$ is a
{\em submanifold} if for each point $p \in N$ there exists a local chart (as defined above)
such that
$$
N = \cR^n_s \cap \{x_{i_1} = \cdots = x_{i_k} = 0\}
$$
for some sublist of indices $[i_1,\ldots,i_k] \subset [n,\ldots,1]$.  

The connected components of $\partial M \setminus \partial \partial M$ will play a role similar to
the irreducible components of the exceptional divisors in the classical results of resolution of singularities.  For this reason, we call an {\em irreducible divisor} (or divisor component) of $M$
to a connected codimension one submanifold which is contained in $\partial M$.
A {\em divisor with normal crossings} (or, shortly, a divisor) is a subset $\vD \subset \partial M$ formed
by some union of irreducible divisors.

We shall denote by $\vO_M$ the sheaf of germs of analytic functions on $M$.  If there is no risk of ambiguity, the stalk
of $\vO_M$ at a point $p \in M$ will be simply denoted by $\vO_p$ .

We refer to \cite{Mi} and \cite{Me} for further details on the theory of manifolds with corners.
\subsection{Singularly Foliated Manifolds}
Let $M$ be a real analytic three-dimensional manifold (with corners) and 
let $\mainlist \in  \lists$ be a list of natural numbers.
\begin{Definition}
A {\em $\mainlist$-tagged divisor} on $M$ is a divisor with normal
crossings $\vD \subset
M$ together with a bijection
$$ \mainlist \longrightarrow  \mbox{Set of irreducible components of
}\vD $$
which associates to each index $i \in \mainlist$ an irreducible component
$D_i \subset \vD$.  We shall shortly write $\vD = \vD_\mainlist$ to
indicate that $\vD$ is a $\mainlist$-tagged divisor.
\end{Definition}
Let $\chi$ be an analytic vector field on $M$. Given a point $p \in M$ and a prime germ
$g \in \cm_p$ (where $\cm_p \subset \vO_p$ is the maximal ideal),
consider the ideal $\vI_\chi g\, \subset \vO_p$ which is generated by
the set
$$\{\chi(h) \mid h \in  (g) \vO_p \}$$
where $\chi(h)$ is the action of $\chi$ (seen as a derivation) on $h
\in \vO_p$.
\begin{Definition}\label{def-nondegvf}
The vector field $\chi$ will be
called {\em nondegenerate with respect to the divisor $\vD$} if for all
point $p \in M$ and all prime $g \in
\cm_p$, one of the following two cases occurs:
\begin{enumerate}[(i)]
\item The set $(g=0)$ is a not a local irreducible component of the divisor and
      $\vI_\chi (\cm_p)$ is not divisible by $g$.
\item The set $(g=0)$ is a local irreducible component of the divisor and
      $\vI_\chi (g)$ is not divisible by $g^2$.
\end{enumerate}
\end{Definition}
Choose some coordinate system $(x_1,\ldots,x_n)$ at $p$ and suppose
that $g = x_1$.   If we write
$$\chi = a_1
  \frac{\partial}{\partial x_1} + \cdots + a_n
  \frac{\partial}{\partial x_n}, \quad a_1,\ldots,a_n \in \vO_p
$$
then the ideal $\vI_\chi (\cm_p)$ is generated by $\{a_1,a_2,\ldots,a_n\}$ and
the ideal $\vI_\chi(g)$ is generated by $\{a_1,a_2 x_1,\ldots,a_n
x_1\}$.  Hence, the conditions $(i)$ and $(ii)$ of the definition can be rewritten as 
\begin{enumerate}[(i)]
\item If $(x_1 = 0) \not\subset \vD$ then $\{a_1,\ldots,a_n\} \not\subset (x_1)\vO_p$.
\item If $(x_1 = 0) \subset \vD$ then there exist no collection of germs
      $\{b_1,\ldots,b_n\} \subset \vO_p$
such that we can write $a_1 = x_1^2 b_1$ and $a_j = x_1 b_j$ for $j \ge 2$.
\end{enumerate}
\begin{Remark}
Suppose that we consider (as in \cite{C1}) the sheaf $\Theta_M[\log \vD]$ of vector fields adapted to $\vD$ (i.e.\ the dual of the sheaf of logarithmic forms with respect to $\vD$).  Then an 
element $\chi_p \in \Theta_M[\log \vD]_p$ is nondegenerate if and only if the adapted coefficients are without a common divisor.
\end{Remark}
Note that a reduced vector field
(as defined in the introduction) is automatically
nondegenerate. However, a nondegenerate vector field can have a set of
singularities $Z$ of codimension one.  For instance,
$$
\chi = x_1 \frac{\partial}{\partial x_1} + 0 \frac{\partial}{\partial x_2} + \cdots + 0
  \frac{\partial}{\partial x_n}
$$
is a nondegenerate vector field in $\cR^n$ (if the divisor $\vD$ contains
$(x_1 = 0)$).  Note that each singular point on
the hypersurface $Z = \{x_1 = 0\}$ is elementary.
\begin{Remark}\label{remark-singularsmooth}
Let $p \in M$ be a singular point of a nondegenerate vector field $\chi$.
Suppose that, for some neighborhood $U \subset M$ of $p$,
we have $\Ze(\chi) \cap U = \{ f = 0 \}$, for some analytic function $f$ such that
$df(p) \ne 0$.  Then,  writing $\chi = f \ \chi_1$ for some analytic
vector field $\chi_1$ defined in $U$, we  obtain the following equivalence
$$
p \mbox{ is an elementary singular point of } \chi \; \Leftrightarrow \; \chi_1(f)(p) \ne 0
$$
(where $\chi_1$ acts as a derivation on $f$).
\end{Remark}
As we shall see, even if we start with a vector field $\chi$ which is reduced, the procedure of resolution of singularities can produce new vector fields which belong to the more general class of
nondegenerate vector fields.  This is due to the occurrence of the so-called {\em dicritical} situations (see example~\ref{exampl-dicrit}).

A {\em singular orientable analytic line field} on $(M,\vD)$ is given by a collection
of pairs $\linef = \{(U_\alpha,\chi_\alpha)\}_{\alpha \in A}$, where
$\{U_\alpha\}$ is an open covering of $M$ and
$$\chi_\alpha : U_\alpha \rightarrow TU_\alpha$$
is an analytic vector field in $U_\alpha$ which is nondegenerate with respect
to $U_\alpha \cap \vD$ (see definition~\ref{def-nondegvf}) and such that for each pair of
indices
$\alpha,\beta \in A$,
$$
\chi_\alpha = h_{\alpha\beta}\cdot \chi_\beta
$$
for some strictly positive analytic function $h_{\alpha\beta}: U_\alpha \cap
U_\beta \rightarrow \cR_{>0}$.

An analytic vector field $\chi$ defined in
a neighborhood $U \subset M$ of a point $p$ is
called {\em local generator} of $\linef$ if the collection
$\{(U_\alpha,\chi_\alpha)\}_{\alpha \in A} \cup \{(U,\chi)\}$ is still
a singular orientable analytic line field.

Let $Y \subset M$ be an analytic subset and $\linef$ be a singular
orientable analytic line
field on $(M,\vD)$. We shall say that $\linef$ is {\em $Y$-preserving}
if for each point $p \in M$ and each local generator $\chi$ of $\linef$ at $p$,
$$\chi(g) \in \vI(Y_p), \quad\mbox{ for all }g \in \vI(Y_p)$$
where $\vI(Y_p)$ is the ideal in the local ring $\vO_p$ which defines
the germ $Y_p$ and $\chi(g)$ is the action of $\chi$ (seen as a
derivation on $\vO_p$)
on $g$.
\begin{Definition}\label{def-singfolmanif}
A {\em singularly foliated manifold} is a 4-uple $\cM =
(M,\mainlist,\vD,\linef)$ where
\begin{itemize}
\item[(i)] $M$ is an analytic three-dimensional manifold with corners;
\item[(ii)] $\mainlist \in \lists$ is list of natural numbers;
\item[(iii)] $\vD = \vD_\mainlist$ is a $\mainlist$-tagged divisor on
  $M$;
\item[(iv)] $\linef$ is a singular orientable analytic line field on $(M,\vD)$
  which is $\vD$-preserving.
\end{itemize}
\end{Definition}
For each point $p \in M$, we define the {\em incidence list at $p$} as the
sublist
\begin{equation}\label{incidencelist}
\iota_p = \{ i \in \mainlist \mid p \in D_i \}, \quad (\mbox{note that
}0 \le \#\iota_p \le n)
\end{equation}
where $D_i$ is $i^{th}$ irreducible component of $\vD$.
We shall say that $p$ is a {\em divisor point} if $\#\iota_p \ge
1$. 

Given a singularly foliated manifold $\cM = (M,\mainlist,\vD,\linef)$, we consider the 
analytic subsets
$$
\Ze(\cM) = \{ p \in M \mid \chi(p) = 0 \}, \quad \Elem(\cM) = \{ p \in
M \mid \chi \mbox{ is elementary at } p\}
$$
where $\chi$ is a local generator for $\linef$ at $p$. The set
$\NElem(\cM) = \Ze(\cM) \setminus \Elem(\cM)$ will be called the
set of {\em nonelementary singular points} of $\linef$.
\begin{Proposition}
The set $\NElem(\cM)$ is a closed analytic subset of $M$ of codimension
strictly greater that one.
\end{Proposition}
\begin{proof}
Given a point $p \in M$, fix some local coordinates $(x_1,\ldots,x_n)$ and a
local generator $\chi = a_1 \partial/\partial x_1 + \cdots + a_n \partial/\partial x_n$ for 
$\linef$ at  $p$. Then
$\NElem(\cM)$ is locally defined by the analytic conditions
$$
\{ \chi = 0, \mathrm{Spec}D\chi = 0\}
$$
Let us prove that such germ of analytic set $\NElem(\cM)_p$ has codimension strictly greater
than one.  If this is not the case, there exists some prime element
$f \in \cm_p$ such that $\{f = 0\}$ is contained in $\NElem(\cM)_p$.  By
the coherence of $\vO_M$, we can suppose (possibly replacing $p$ by some neighboring
point) that $\{df(p) \ne 0\}$.  Using the
Remark~\ref{remark-singularsmooth}, we conclude that $\chi$ is necessarily divisible by $f^2$.  This contradicts
the hypothesis of nondegeneracy for $\chi$.
\end{proof}

\subsection{Multiplicity and weighted Blowing-up in $\cR^n$}\label{subsect-blowup}
A {\em weight-vector} is a nonzero vector of natural numbers
$\bfomega = (\omega_1,\ldots,\omega_n) \in \cN^n$.

The {\em $\bfomega$-multiplicity} of a
monomial $\bfx^\bfv = x_1^{v_1}\cdots x_n^{v_n}$ (with
$\bfv \in \cZ^n$)
is the integer number
$$
\mu_\bfomega(\bfx^\bfv) = \< \bfomega, \bfv \> :=
\omega_1 v_1 + \cdots + \omega_n v_n
$$
More generally, the {\em $\bfomega$-multiplicity} of a formal series $f \in \For$ is
given by
$$
\mu_\bfomega(f) = \min \{d \in \cN \mid f
\mbox{ has a monomial }*\bfx^\bfv\mbox{ with
}\mu_\bfomega(\bfx^\bfv) = d \}
$$
(where $*$ denotes some nonzero real number).
We denote by $H_\bfomega^d$ the subset of all formal
series with $\bfomega$-multiplicity equal to $d$.

Given formal series $a_1,\ldots,a_n \in \cR[[\bfx]]$, 
the corresponding formal $n$-dimensional vector field 
$$\chi = a_1 \frac{\partial}{\partial x_1} + \cdots + 
a_n \frac{\partial}{\partial x_n}$$
is a derivation on the ring $\cR[[\bfx]]$.
The $\bfomega$-multiplicity of $\chi$ is the
integer number
$$
\mu_\bfomega(\chi) = \max \{ k \in \cZ \mid \chi(H_\bfomega^d) \subset
H_\bfomega^{d + k}, \forall d \in \cN\}
$$
where $\chi(H_\bfomega^d)$ denotes the action of $\chi$ (seen as a
derivation) on the subset $H_\bfomega^d$.
\begin{Remark}
Using the expression $\chi = a_1 \frac{\partial}{\partial x_1} +
\cdots + a_n \frac{\partial}{\partial x_n}$, we have
$$
\mu_\bfomega(\chi) = \min\{\mu_\bfomega(a_1) - \omega_1,  \ldots,
\mu_\bfomega(a_n) - \omega_n\}
$$
\end{Remark}
Given a weight-vector $\bfomega \in \cN^n_{>0}$, 
the {\em $\bfomega$-weighted blowing-up} of $\cR^n$ is the real analytic surjective map
\begin{eqnarray*}
\Phi_\bfomega: \cS^{n-1} \times \cR^+ &\longrightarrow& \cR^n\\
(\overline{\bfx},\tau)\;\;\; &\longmapsto& \tau^\bfomega \overline{\bfx} =
(\tau^{\omega_1}
\overline{x}_1,\ldots,
\tau^{\omega_n}\overline{x}_n)
\end{eqnarray*}
where we put $\cS^{n-1} = \{\overline{\bfx} \in \cR^n \mid \overline{x}_1^2
+ \cdots +
\overline{x}_n^2 = 1\}$.

More generally, given an arbitrary weight-vector $\bfomega \in \cN^n$, we can reorder the
coordinates and write
$\bfomega = (\omega_1,\ldots,\omega_k,0,\ldots,0)$, where $\omega_1,\ldots,
\omega_k$ are strictly positive.  The {\em $\bfomega$-weighted blowing-up} is the map
\begin{eqnarray*}
\Phi_\bfomega: \cS^{k-1} \times \cR^+ \times \cR^{n-k} &\longrightarrow& \cR^n\\
(\overline{\bfx},\tau,\bfx^\prime)\;\;\;\;\;\; &\longmapsto&
(\tau^\bfomega \overline{\bfx},
\bfx^\prime)
\end{eqnarray*}
where $\bfx^\prime = (x_{k+1},\ldots,x_n)$. The sets
$$Y = \{ x_1 = \ldots = x_k = 0\},\quad D =
\Phi^{-1}_\bfomega(Y) = \cS^{k-1} \times \{0\} \times \cR^{n-k}$$
will be called respectively, the {\em blowing-up center} and
{\em exceptional divisor}
of the blowing-up.  The set $\widetilde{M} = \cS^{k-1} \times \cR^+ \times \cR^{n-k}$ will be called 
{\em blowed-up space}.

It is obvious that $\Phi_\bfomega$ restricts to a diffeomorphism
between $\widetilde{M} \setminus D$ and $\cR^n
\setminus Y$.  The blowing-up creates a boundary component $\partial {\widetilde{M}} = D$.

The definition of $\bfomega$-weighted blowing-up can be easily extended to the spaces with corners
$\cR^n_s$, thus defining an analytic surjective map 
$$\Phi_\bfomega: \widetilde{M} \rightarrow \cR^n_s$$
It is easy to describe the effect of the blowing-up on a divisor $\vD \subset \cR^n_s$.
If $\vD \subset \cR^n_s$ is divisor then the set
$$
\widetilde{\vD} = \Phi^{-1}_\bfomega(\vD) \cup D
$$
is a divisor in $\widetilde{M}$. The divisor $\widetilde{\vD}$ will be called the {\em total transform} of $\vD$.
\begin{Remark}
It follows from our definition of manifolds with corners that a 
divisor $\vD \subset \cR^n_s$ is always given by
a finite union of coordinate hyperplanes, namely
$$
\vD = \bigcup_{i \in \iota}{\{x_i = 0\}}
$$
for some sublist $\iota \subset [n,\ldots,1]$.  
\end{Remark}
Now, let us fix an analytic (non identically zero) vector field $\chi$  in $\cR^n_s$. Consider the analytic
vector field $\chi^*$ defined
in $\widetilde{M} \setminus D$ as the
pull-back of $\chi$ under the
diffeomorphism
$$\Phi_\bfomega:\widetilde{M} \setminus D
\rightarrow \cR^n_s \setminus Y$$
The following result is obtained by a straightforward computation:
\begin{Proposition}\label{prop-blowupvf}
Let $m = \mu_\bfomega(\chi)$ be the $\bfomega$-multiplicity of $\chi$
(seen as a formal vector field at the origin).  Then, the new vector field
$$\widetilde{\chi} = \tau^{-m}\cdot \chi^*$$
satisfies the following conditions:
\begin{itemize}
\item[(i)] $\widetilde{\chi}$ has an analytic extension to $\widetilde{M}$
(which we still denote by $\widetilde{\chi}$);
\item[(ii)] The exceptional divisor $D$ is an
invariant manifold for $\widetilde{\chi}$.
\item[(iii)] If $\chi$ is nondegenerate with respect to some divisor $\vD \subset \cR^n_s$
then $\widetilde{\chi}$ is nondegenerate with respect to the total transformed divisor $\widetilde{\vD}$.
\end{itemize}
\end{Proposition}
\begin{proof}
See e.g.\ \cite{P}.
\end{proof}
The vector field $\widetilde{\chi}$ will be called the {\em strict transform} of $\chi$.
\begin{Example}\label{exampl-dicrit}
Let us see an example of a typical {\em dicritical} situation, 
where the blowing-up of a reduced vector field
results into a nondegenerate vector field whose set of singularities has codimension one.  Consider
the vector field in $\cR^3_{(x,y,z)}$
$$
\chi = y \frac{\partial}{\partial y} + z \frac{\partial}{\partial z}
$$
and choose the weight-vector $\bfomega = (0,1,1)$.  The $\bfomega$-weighted blowing-up (with center $Y = \{y =  z = 0\}$)
gives the manifold $\widetilde{M} = \cR_{x} \times \cS^1_{\theta} \times \cR^+_\tau$ with the exceptional divisor $D = \{(x,\theta,\tau) \in \widetilde{M} \mid \tau = 0\}$.   Since $\mu_\bfomega(\chi) = 0$,   the strict transform of $\chi$ is 
the {\em radial} vector field $\widetilde{\chi} = \tau \frac{\partial}{\partial \tau}$, which vanishes identically on $D$.
\end{Example}
\subsection{Directional Charts of Blowing-up}\label{subsect-directionalcharts}
Let $\bfomega = (\omega_1,\ldots,\omega_k,0,\ldots,0)$ be a weight-vector as above.  

Given an index $1 \le r \le k$, the {\em $x_r$-directional} $\bfomega$-weighted blowing-up is the pair of 
analytic maps
$$
\Phi_\bfomega^{r,+}: U \rightarrow \cR^n \cap \{x_r \ge 0\},\quad
\Phi_\bfomega^{r,-}: U \rightarrow \cR^n \cap \{x_r \le 0\}
$$
with domain $U := \cR^{r-1} \times \cR^+ \times \cR^{n-r}$ which are defined as follows.  
Write the coordinates in $U$ as
$(\widetilde{x}_1,\ldots,\widetilde{x}_n)$.  Then, for each choice of sign $\eps \in \{+,-\}$,
the map $\bfx = \Phi_\bfomega^{r,\eps}(\widetilde{\bfx})$ is given by
$$
\left\{
\begin{matrix}
x_i = {\widetilde{x}_r}^{\,\omega_i}\ \widetilde{x}_i, \hfill& \quad \mbox{ for } i = 1,\ldots,r-1,r+1,\ldots,k \cr \cr
x_r = \eps\; {\widetilde{x}_r}^{\,\omega_r}  \hfill& \cr \cr
x_j = \widetilde{x}_j, \hfill& \quad \mbox{ for } j = k+1,\ldots,n \hfill
\end{matrix}
\right.
$$

\begin{figure}[htb]
\psfrag{phi}{\small $\phi^{r,\eps}$}
\psfrag{Phir}{\small $\Phi_\bfomega^{r,\eps}$}
\psfrag{Phi}{\small $\Phi_\bfomega$}
\begin{center}
\includegraphics[height=6cm]{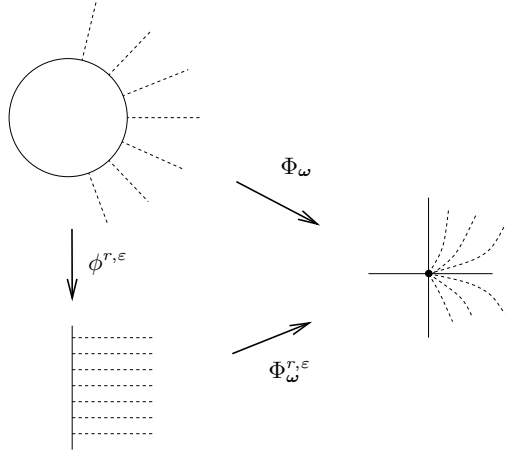}
\end{center}
\caption{The $x_r$-directional chart (dotted lines indicate the fibration $\{d\, \overline{\bfx} = 0\}$).}
\label{fig-directblup}
\end{figure}

\begin{Proposition}\label{prop-directblup}
For $\eps \in \{+,-\}$, there exists an analytic diffeomorphism 
$$
\phi^{r,\eps} : V^{r,\eps} \rightarrow U
$$ 
with domain $V^{r,\eps} := \{(\overline{\bfx},\tau,\bfx^\prime) \in  \cS^{k-1} \times \cR^+ \times \cR^{n-k}
\mid \eps\, \overline{x}_r > 0\}$ which makes the following diagram commutative
\begin{eqnarray*}\label{commutdiagramdirectblowup}
\begin{CD}
V^{r,\eps} @>{\Phi_\bfomega}>> \cR^n \cap \{\eps\, x_r \ge 0\}\\
@V{\phi^{r,\eps}}VV                                   @VV{\mathrm{id}}V \\
U  @>{\Phi_\bfomega^{r,\eps}}>>  \cR^n \cap \{\eps\, x_r \ge 0\}
\end{CD}
\end{eqnarray*}
where $\mathrm{id}$ is the identity map.  
\end{Proposition}
\begin{proof}
See \cite{DR}.
\end{proof}
The pairs $(V^{r,+},\phi^{r,+})$ and $(V^{r,-},\phi^{r,-})$ will be called 
{\em $x_r$-directional charts} of the blowing-up.  

Notice that the exceptional divisor $D$ is mapped 
by $\phi^{r,\eps}$ to the hyperplane $\{\widetilde{x}_r = 0\}$. Moreover, the union of the domains of all 
directional charts, 
$$
V^{1,+} \cup V^{1,-} \cup \cdots \cup V^{k,+} \cup V^{k,-}
$$
gives an open covering of the blowed-up space $\widetilde{\cM} = \cS^{k-1} \times \cR^+ \times \cR^{n-k}$.

\subsection{Weighted Trivializations and Blowing-up on Manifolds}\label{subsect-blupmanif}
Let us fix a weight-vector $\bfomega = (\omega_1,\ldots,\omega_k,0,\ldots,0)$.
We shall say that an analytic map $\phi: U \rightarrow \cR^n$ with domain an open
subset $U \subset \cR^n$ {\em preserves the $\bfomega$-quasihomogeneous structure
on $\cR^n$} if
$$
\phi^*(H_\bfomega^d) \subset H_\bfomega^d
$$
for each natural number $d \in \cN$.

Let $M$ be an $n$-dimensional analytic manifold (with corners) and $Y \subset M$ be a submanifold
of codimension $k$. A {\em trivialization atlas} for $Y \subset M$ is a collection of pairs
$\{(U_\alpha,\phi_\alpha)\}_{\alpha \in A}$ where $\{U_\alpha\}$ is an open
covering of $Y$ and
$$
\phi_\alpha: U_\alpha \rightarrow \cR^n_s
$$
is a local chart such that $\phi_\alpha(Y \cap U_\alpha) = \{0\} \times \cR^{n-k}_{s^\prime}$, for some $s^\prime \le s$.

We shall say that $\{(U_\alpha,\phi_\alpha)\}_{\alpha \in A}$ is {\em $\bfomega$-weighted} trivialization atlas if
for each pair of indices $\alpha,\beta \in A$, the transition map
$$
\phi_{\alpha\beta} := \phi_{\beta} \circ \phi_{\alpha}^{-1}\; : \;
\phi_\alpha(U_\alpha \cap U_\beta) \longrightarrow \phi_\beta(U_\alpha \cap U_\beta)
$$
preserves the $\omega$-quasihomogeneous structure on $\cR^n$.
\begin{Proposition}\label{prop-blowupmanif}
Let $\{(U_\alpha,\phi_\alpha)\}_{\alpha \in A}$ be a $\bfomega$-weighted trivialization
atlas for a submanifold $Y \subset M$.  
Then, there exists an $n$-dimensional analytic
manifold $\widetilde{M}$ and a proper analytic surjective map
$$
\Phi: \widetilde{M} \rightarrow M
$$
such that the following conditions hold:
\begin{enumerate}[(i)]
\item The map $\Phi$ induces a diffeomorphism between $M
      \setminus Y$ and $\widetilde{M} \setminus D$, where
      $D := \Phi^{-1}(Y)$.
\item There exists an a collection of local charts $\{(\widetilde{U}_\alpha,\widetilde{\phi}_\alpha)\}_{\alpha \in A}$  
in $\widetilde{M}$  such that $\{\widetilde{U}_\alpha\}_{\alpha \in A}$ is an open covering of $D$ and
$$
\widetilde{\phi}_\alpha : \widetilde{U}_\alpha \rightarrow \cS^{k-1}_{s^{\prime\prime}} \times \cR^+ \times \cR^{n-k}_{s^\prime}
$$
is an analytic diffeomorphism (where $\cS^{k-1}_{t} = \{\bar \bfx \in \cS^{k-1}\mid \bar x_1 \ge 0, \ldots, \bar x_{t} \ge 0\}$ and $s^\prime + s^{\prime\prime} = s$), such that the following diagram is commutative
\begin{eqnarray*}\label{commutdiagramblowup}
\begin{CD}
\widetilde{U}_\alpha @>{\Phi}>> U_\alpha\\
@V{\widetilde{\phi}_\alpha}VV                                   @VV{\phi_\alpha}V \\
\cS^{k-1}_{s^{\prime\prime}} \times \cR^+ \times \cR^{n-k}_{s^\prime}  @>{\Phi_\bfomega}>>  \cR^n
\end{CD}
\end{eqnarray*}
where $\Phi_\bfomega$ is the $\bfomega$-weighted blowing-up in $\cR^n$.
\end{enumerate}
\end{Proposition}
\begin{proof}
See \cite{DR}, Proposition II.9.
\end{proof}
The map $\Phi: \widetilde{M} \rightarrow M$ will be called a {\em $\bfomega$-weighted blowing-up 
of $M$} with center on $Y$, with respect to the trivialization 
$\{(\widetilde{U}_\alpha,\widetilde{\phi}_\alpha)\}_{\alpha \in A}$.
\begin{Remark}
The existence of a $\omega$-weighted trivialization for a submanifold $Y
\subset M$ can be a strong topological restriction.  This condition can be defined in a more
intrinsic way as the existence of a certain nested sequence of subbundles in the conormal bundle
$N^* Y$ (see e.g.\ \cite{Me}, section 5.15).
\end{Remark}
\subsection{Blowing-up of Singularly Foliated Manifolds}\label{subsect-blupsingfolmanif}
Let $\cM = (M,\mainlist,\vD,\linef)$ be a singularly foliated manifold and $Y \subset M$ be a submanifold which
has a $\bfomega$-weighted trivialization.  The {\em $\bfomega$-weighted} blowing-up of $\cM$ with center $Y$ is
a mapping
$$
\Phi: \widetilde{\cM} \rightarrow \cM,
$$
defined by taking the $4$-uple
$\widetilde{\cM} = (\widetilde{M},\widetilde{\mainlist},\widetilde{\vD},\widetilde{\linef})$, where
\begin{enumerate}[(i)]
\item $\Phi: \widetilde{M} \rightarrow M$ is the $\bfomega$-weighted blowing-up 
of $M$ with center on $Y$;
\item The list $\widetilde{\mainlist}$ is given by $\mainlist \cup [n]$, where 
$n := 1 + \max \{ i \mid i \in \mainlist\}$ if $\mainlist \ne \emptyset$ and $n := 1$ if $\mainlist = \emptyset$;
\item The divisor $\widetilde{\vD}$ is the total transform of $\vD$, with the tagging
$$
\widetilde{\mainlist} \ni i \longrightarrow 
\left\{
\begin{matrix}
D_i^\prime\, , & \mbox{ if }i \in \widetilde{\mainlist} \setminus [n] \hfill \cr\cr
\widetilde{D}\, , & \mbox{ if }i = n \hfill
\end{matrix}
\right.
$$
where  $\widetilde{D} := \Phi^{-1}(Y)$ and 
$D_i^\prime$ is the strict transform of the corresponding divisor $D_i \subset M$ (for each $i \in \mainlist$);
\end{enumerate}
The line field $\widetilde{\linef}$  
is obtained as follows:
Up to some refinement of the coverings, we can suppose that the line field $\linef$ is given by 
a collection 
$\{(U_\beta,\chi_\beta)\}_{\beta \in B}$
where $\chi_\beta$ is a nondegenerate analytic vector field defined in $U_\beta$; and that there exists some subcollection of indices $A \subset B$ such that 
$\{(U_\alpha,\phi_\alpha)\}_{\alpha \in A}$ is the $\bfomega$-weighted trivialization of $Y$.

For each $\alpha \in B$, we can consider the strict transform 
$\widetilde{\chi}_\alpha$ of $\chi_\alpha$ (see Proposition~\ref{prop-blowupvf}) as an analytic vector field
defined in $\widetilde{U}_\alpha = \Phi^{-1}(U_\alpha)$.  

Now, the Proposition~\ref{prop-blowupmanif} implies that
the collection $\{(\widetilde{U}_\alpha,\widetilde{\chi}_\alpha)\}_{\alpha \in B}$
defines a singular line field 
$\widetilde{\linef}$ on $\widetilde{M}$ which satisfies our requirements (see \cite{DR} or \cite{P} for the details).
\subsection{Axis Definition and controllability}\label{subsect-axisdef}
Let $\cM = (M,\mainlist,\vD,\linef)$ be a singularly foliated manifold of dimension three.

As we explained in subsection~\ref{subsect-overview}, the local strategy for the resolution of singularities at a point $p \in \NElem(\cM)$ is based on some invariants attached to the Newton polyhedron. Such Newton polyhedron depends on the vector field $\chi$ which locally generates the line field, but also on a choice of local coordinates $(x,y,z)$ at $p$.  This usually creates difficulties for obtaining a global strategy for the resolution, since the information obtained from the polyhedron is coordinate-dependent.

In order to obtain {\em intrinsic} invariants, we have to restrict the choice of local coordinates and require that they respect to some additional {\em structure} on the ambient space.
We now introduce such structure.
\begin{Definition}\label{def-axis}
An {\em  axis} for $\cM$ is given by a pair $\Ax = (A,\vZ)$, where
$A \subset M$ is an open neighborhood of the set $\NElem(\cM)$
and $\vZ$ is a singular orientable
analytic line field defined on $A$ such that
\begin{enumerate}[(i)]
\item $\vZ$ is $\vD \cap A$-preserving;
\item $\Ze(\vZ) = \emptyset$ (where $\Ze(\vZ)$ is the set of singularities of $\vZ$), and
\item For each point $p \in A \cap \vD$, if we choose a local chart $(U,(x,y,z))$ 
such that $\vZ$ is locally generated by
$\frac{\partial}{\partial z}$ then
$$
\vI_p \not\subset \vJ_p
$$
where $\vI_p \subset \vO_p$ is the ideal which defines the germ of analytic set
$\NElem(\cM)_p$ and $\vJ_p \subset \vO_p$ is the defining ideal of the set $\{x = y = 0\}$ (i.e.\
the leaf of the axis through the point $p$).
\item For each point $p \in A \setminus \vD$, if we choose a local chart $(U,(x,y,z))$
such that $\vZ$ is locally generated by
$\frac{\partial}{\partial z}$ then
$$
\chi(\vJ_p) \not\subset \vJ_p
$$
where $\chi$ is a local generator of $\linef$.
\end{enumerate}
\end{Definition}
The requirement in the item (iii) is equivalent to say that $\{x = y = 0\}$ is not contained in $\NElem(\cM)$.
The (stronger) requirement in the item (iv) is equivalent to say 
that $\{x = y = 0\}$ is not an invariant curve for
the line field $\linef$.

\begin{figure}[htb]
\psfrag{vZ}{\small $\vZ$}
\psfrag{NElem}{\small $\NElem(\cM)$}
\begin{center}
\includegraphics[height=4cm]{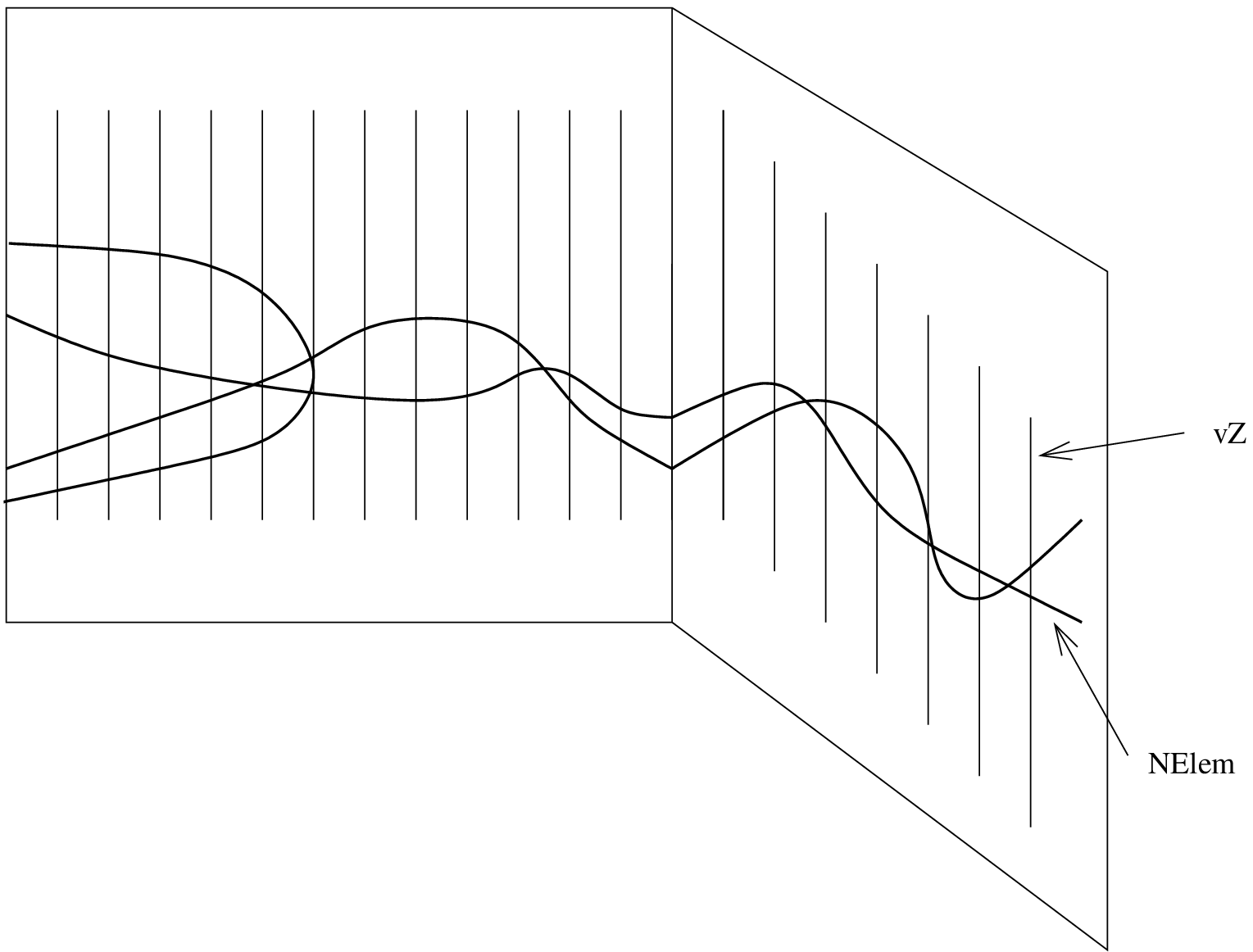}
\end{center}
\caption{Axis.}
\label{fig-Axisdef}
\end{figure}
\begin{Remark}\label{remark-iota=3}
It is not always possible to define an axis for a singularly foliated
line field.  For instance, if there exists a point $p
\in \NElem(\cM)$ such that $\#\iota_p = 3$ then any line field which is
$\vD$-preserving necessarily
vanishes at $p$.  In this case, the requirement in item (ii) of the definition cannot be satisfied.
\end{Remark}
We shall say that the singularly foliated manifold $\cM$ is {\em controllable}  
if there exists an  axis $\Ax$ as defined above.  The pair $(\cM,\Ax)$
will be called a {\em controlled singularly foliated manifold}.

The next Proposition describes a situation where an axis can always be defined: Let $M$ be an analytic manifold of dimension three without boundary and let $\chi$ be a reduced analytic 
vector field defined in $M$.  We consider the singularly foliated manifold $\cM = (M,\mainlist,\vD,\linef)$, where
$\mainlist = \emptyset$, $\vD = \emptyset$ and $\linef = \linef_\chi$ 
is the analytic line field generated by $\chi$.
\begin{Proposition}\label{proposition-definesaxis}
Given a singularly foliated manifold 
$\cM = (M,\emptyset,\emptyset,\linef_\chi)$ as above, there exists an 
axis $\Ax = (A,\vZ)$ for $\cM$.  The pair $(\cM,\Ax)$ is a controlled singularly foliated manifold.
\end{Proposition}
\begin{proof}
The set of nonelementary points of $\cM$ is an one-dimensional analytic subset 
$\NElem \subset M$.  Let $S \subset \NElem$ be the discrete subset of points where $\NElem$ is not locally smooth.

First of all, we are going to define a nonsingular vector field $Z_p$ in an open neighborhood $U_p \subset M$ of 
each point $p \in S$ with the property that
no trajectory of $Z_p$ is a leaf of $\linef \cap U_p$.   
This is easy.  We fix arbitrary local coordinates
$(x,y,z)$ in a neighborhood of $p$ and construct the Newton polyhedron $\vN$ for the vector field $\chi$ with 
respect to these coordinates.  
If the support of $\vN$ contains at least one point in the region 
$$
(\{-1\} \times \{0\} \times \cZ)  \cup (\{0\} \times \{-1\} \times \cZ)
$$
then it suffices to locally define $Z_p$ as the vector field $\frac{\partial}{\partial z}$.  
Otherwise, it follows that $\{x = y = 0\}$ is an invariant curve for the
vector field $\chi$.  In this case, it is immediate to verify that we can choose natural numbers 
numbers $s,t \in \cN_{>0}$ such that the change of coordinates
$$
\widetilde{x} = x + z^s, \quad \widetilde{y} = y + z^t, \quad \widetilde{z} = z
$$
results into a new coordinate system $(\widetilde{x},\widetilde{y},\widetilde{z})$ where the above property holds.

Now, we are going to glue together the collection of vector fields 
$\{Z_p\}_{p \in S}$ in a $C^\infty$-way along the smooth part of $\NElem$: 

\begin{figure}[htb] 
\psfrag{x}{\small $x$}
\psfrag{z}{\small $z$}
\psfrag{p}{\small $p$}
\psfrag{q}{\small $q$}
\psfrag{G}{\small $\Gamma$}
\psfrag{O}{\small $\Strip_\Gamma$}
\begin{center} 
\includegraphics[height=6cm]{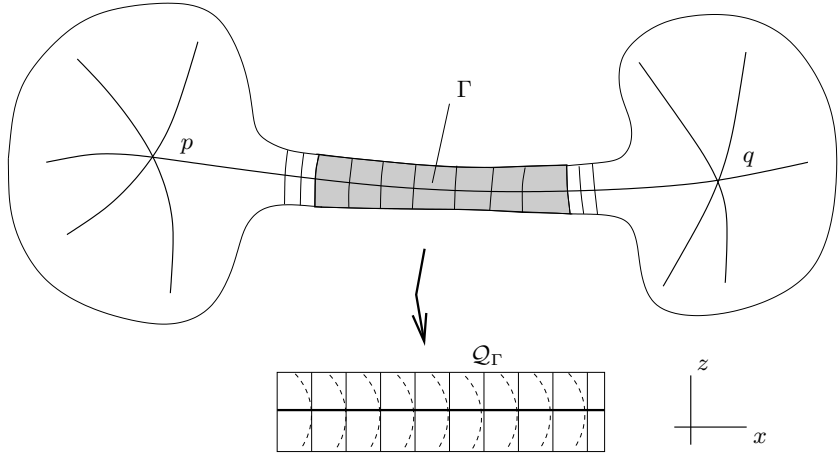}
\end{center}
\caption{The definition of $Z_\Strip$ on the strip $\Strip$.}
\label{fig-defineaxis}
\end{figure}

Let $\Gamma \subset \NElem \setminus S$ be a regular analytic curve connecting two points $p,q \in S$.
Possibly restricting $U_p$ to some smaller neighborhood of $p$, we can assume the $Z_p$ is transversal 
to $\Gamma \cap U_p$ (recall that $\Gamma$ is an analytic arc).  Therefore, in some neighborhood of $\Gamma \cap U_p$, 
we can define a two-dimensional strip $\Strip_p$ formed by the
union of all trajectories of $Z_p$ starting at points of $\Gamma \cap U_p$.
The same argument gives us a two-dimensional strip $\Strip_q$ with base $\Gamma \cap U_q$.

Using the Tubular Neighborhood Theorem, we can glue together these two strips in a 
$C^\infty$ way, as shown in figure~\ref{fig-defineaxis}.  Therefore, we get
a global two-dimensional strip $\Strip$ with base $\Gamma$. Using partitions of unity, it is easy to define a
nonsingular $C^\infty$ vector field $Z_\Strip$ in an open neighborhood of $\Strip$ such that
\begin{enumerate}[(i)]
\item $Z_\Strip$ is tangent to the strip $\Strip$;
\item $Z_\Strip \cap U_p = Z_p$ and $Z_\Strip \cap U_q = Z_q$
\item No trajectory of 
$Z_\Strip$ is invariant left invariant by $\chi$. 
\end{enumerate}
Putting together all such local constructions, we finally obtain a nonsingular 
$C^\infty$ vector field $Z$ defined in an open neighborhood $A \subset M$ of 
$\NElem(\cM)$ which has the property 
$$(P) \quad \mbox{ no trajectory of }Z\mbox{ is invariant by the vector field }\chi$$
Using Grauert's Embedding Theorem~\cite{G}, we can analytically 
embed $M$ in $\cR^k$, for some sufficiently large $k \in \cN$.  
Doing so, the vector field $Z$ can be seen as a map $Z: A \rightarrow \cR^k$ and it is clear that property 
(P) is an open property for the Whitney topology on $C^\infty(\cR^k,\cR^k)$.
Therefore, using Weierstrass Approximation Theorem (in the version of \cite{G}),
we can approximate $Z$ by an 
analytic nonsingular vector field $\widetilde{Z} : A \rightarrow TM$ 
which also has property (P). 
This proves the Proposition. 
\end{proof}
\section{Newton Polyhedron and Adapted Coordinates}\label{sect-NewtonMaps}
\subsection{Adapted Local Charts}\label{subsect-localadaptedchart}
Let $(\cM,\Ax)$ be  a controlled singularly foliated manifold, with
$\cM = (M,\mainlist,\vD,\linef)$ and $\Ax = (A,\vZ)$.

A local chart $(U,(x,y,z))$
centered at a point $p \in A$ will be called an {\em adapted local chart} if
\begin{itemize}
\item $\vZ$ is locally generated by $\frac{\partial}{\partial z}$;
\item If $p \in \vD$ and $\iota_p = [i]$ then $D_i = \{ x = 0 \}$;
\item If $p \in \vD$ and $\iota_p = [i,j]$ (with $i > j$) then $D_i = \{ x = 0\}$ and $D_j = \{y = 0\}$.
\end{itemize}
where $\iota_p$ is the incidence list defined in (\ref{incidencelist}). 

Notice that an adapted local chart can always be defined at a point $p \in A$.  The condition $p \in A \cap \vD$ automatically implies that $\#\iota_p \in \{1,2\}$, by the remark~\ref{remark-iota=3}.

Despite the fact that the definition of adapted local chart is given for all points in $A$, we shall be mostly concentrated (at least until the end of section~\ref{sect-LocalTheory}) on points lying in $A \cap \vD$.

\begin{figure}[htb]
\psfrag{i1}{\small $\iota_p = [i]$}
\psfrag{i2}{\small $\iota_p = [i,j]$}
\psfrag{x}{\small $x$}
\psfrag{y}{\small $y$}
\psfrag{z}{\small $z$}
\psfrag{Di}{\small $D_i$}
\psfrag{Dj}{\small $D_j$}
\psfrag{p}{\small $p$}
\begin{center}
\includegraphics[height=4cm]{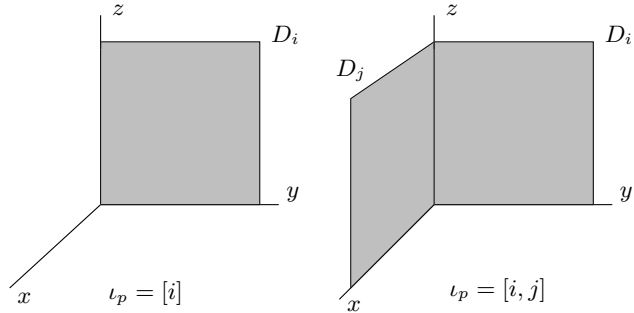}
\end{center}
\caption{Adapted local charts (with $i > j$).}
\label{fig-localpresentation}
\end{figure}

In figure~\ref{fig-localpresentation} we represent the two possible
configurations with the corresponding position of the divisors.

\begin{Proposition}\label{prop-transitionmaps}
Let $(U,(x,y,z))$ and  $(U^\prime,(x^\prime,y^\prime,z^\prime))$ be two adapted local charts at a point
$p \in A$.  Then, the transition map has the form
$$
x^\prime = F(x,y), \quad
y^\prime = G(x,y), \quad
z^\prime = f(x,y) + z w(x,y,z)
$$
where $w$ is a unit and $\frac{\partial(F,G)}{\partial(x,y)}(0,0) \ne 0$.  
More specifically, if $\#\iota_p = 1$ then $F,G$ have the particular form
$$
F(x,y) = x u(x,y), \quad G(x,y) = g(x) + y v(x,y),
$$
where $u,v$ are units and $g(0) = 0$. Similarly, if $\#\iota_p = 2$ then 
$$
F(x,y) = x u(x,y), \quad G(x,y) = y v(x,y),
$$
for some units $u,v$.
\end{Proposition}
\begin{proof}
The change of coordinates should map the vector field
$\frac{\partial}{\partial z}$ into the vector field $U\cdot\frac{\partial}{\partial
z^\prime}$ (for some unit $U \in \vO_p$).  Moreover, if $\#\iota_p \ge 1$, it maps the 
divisor $\{x = 0\}$ into $\{x^\prime = 0\}$. If $\#\iota_p = 2$, the divisor $\{ y = 0
\}$ should also be mapped to $\{y^\prime = 0\}$.
\end{proof}

\subsection{Newton Map and Newton Data}\label{subsection-Newtonmap}
Let $(\cM,\Ax)$ be  a controlled singularly foliated manifold, with
$\cM = (M,\mainlist,\vD,\linef)$ and $\Ax = (A,\vZ)$.

We fix a point $p$ in $A$ and an adapted local chart $(U,(x,y,z))$
centered at $p$.   Our goal is to define the Newton polyhedron of $(\cM,\Ax)$ at $p$ with respect to the 
coordinates $(x,y,z)$.

First of all, we choose an analytic vector field $\chi$ which generates
$\linef$ at $U$. Next, we expand $\chi$ is the {\em logarithmic
basis}:  Consider
the meromorphic functions
$$f := \chi(\ln x) = \frac{\chi(x)}{x}, \quad
g := \chi(\ln y) = \frac{\chi(y)}{y}, \quad
h := \chi(\ln z) = \frac{\chi(z)}{z}$$
where $\chi$ acts as a derivation on $\cR\{x,y,z\}$.  Then, we can write
\begin{equation}\label{logexpansion}
\chi = fx \frac{\partial}{\partial x} + g y \frac{\partial}{\partial y}
+ h z \frac{\partial}{\partial z}.
\end{equation}
\begin{Remark}\label{rem-restrictfs}
If $\chi$ is $\{x = 0\}$-preserving (resp.\ $\{y = 0\}$-preserving, $\{z = 0\}$-preserving)  then
$f$ (resp.\ $g$, $h$) is an analytic germ in $\cR\{x,y,z\}$.
\end{Remark}
We can write the Laurent series expansion of the functions $(f,g,h)$ given in
(\ref{logexpansion}) as
$$
(f,g,h) =
\sum_{\bfv \in \cZ^3}{(f_{\bfv}, g_{\bfv}, h_{\bfv})\cdot
  x^{v_1} y^{v_2} z^{v_3}}$$
where $(f_{\bfv},g_{\bfv},h_{\bfv})$ is a vector in $\cR^3$ for each integer vector $\bfv = (v_1,v_2,v_3) \in \cZ^3$.
The {\em Newton map} for $\chi$ at $p$, relatively to the chart $(U,(x,y,z))$ is the map
\begin{eqnarray*}
\Theta: \cZ^3 &\longrightarrow & \cR^3 \\
          \bfv\hfill &\longmapsto & (f_{\bfv},g_{\bfv}, h_{\bfv}).
\end{eqnarray*}
The {\em support} of $\Theta$ is given by
$$\supp(\Theta) = \{\bfv \in \cZ^n \mid \Theta(\bfv) \ne 0\}.$$
\begin{Remark}\label{remark-newton}
The Newton map $\Theta$ has the following properties:
\begin{itemize}
\item $\supp(\Theta) \subset \cN^3 \cup (\{-1\} \times \cN^2) \cup
  (\cN \times \{-1\} \times \cN) \cup (\cN^2 \times \{-1\})$;
\item $\bfv \in (\{-1\}\times \cN^2) \Rightarrow \Theta(\bfv) \in \cR \times \{0\}\times \{0\}$;
\item $\bfv \in (\cN \times \{-1\} \times \cN) \Rightarrow \Theta(\bfv) \in \{0\}\times \cR \times \{0\}$;
\item $\bfv \in (\cN^2 \times \{-1\}) \Rightarrow \Theta(\bfv) \in \{0\}\times \{0\}\times \cR $.
\end{itemize}
\end{Remark}
The {\em Newton polyhedron} for $(\cM,\Ax)$ at $p$, relatively to the adapted
local chart $(U,(x,y,z))$ is the convex
polyhedron in $\cR^3$ given by
$$
\vN = \conv(\supp(\Theta)) + \cR^3_+
$$
where $\conv(\cdot)$ is the convex closure operation and the $+$ sign
denotes the usual Minkowski sum of convex polyhedrons.
\begin{Lemma}\label{lemma-choicelocgen}
The Newton polyhedron is independent of the choice of the local generator of $\linef$.
\end{Lemma}
\begin{proof}
Indeed, if $\chi,\chi^\prime$ are two local generators, we know that $\chi^\prime = U\, \chi$, for some unit
$U \in \cR\{x,y,z\}$.  Going back to the definition of the Newton polyhedron, 
it is clear that the corresponding polyhedrons 
$\vN$ and $\vN^\prime$ will coincide.
\end{proof}
It is obvious that different choices of local coordinates $(x,y,z)$ lead to different Newton
polyhedrons.  Later on, we shall see that certain essential properties of $\vN$ are preserved 
by the action of the group of coordinate changes defined by Proposition~\ref{prop-transitionmaps}.

From now on, we shall adopt the usual language of the Theory of Convex
Polyhedrons, and refer to the {\em vertices}, {\em edges} and
{\em faces} of $\vN$ (the faces will always be {\em two-dimensional}).

Given a face $F \subset \vN$, there exists a weight-vector $\bfomega \in \cN^3$ and an integer $\mu \in \cZ$ 
such that
$$F = \vN \cap \{\bfv \in \cR^3 \mid \< \bfomega, \bfv \> = \mu\}$$
Notice that if such property is satisfied for a pair $(\bfomega,\mu)$ then it is satisfied on the entire positive ray
$R = \{t\cdot(\bfomega,\mu) \mid t > 0\}$.
The {\em weight-vector} and the {\em multiplicity} associated to $F$ are given by the unique pair 
$(\bfomega,\mu) \in R$ such that $\bfomega = (\omega_1,\omega_2,\omega_3)$ is a nonzero vector of natural numbers satisfying 
$\mathrm{mdc}(\omega_1,\omega_2,\omega_3) = 1$.
\begin{Definition}
The triple $\Omega = ((x,y,z),\iota_p,\Theta)$ will be called a {\em Newton data} for the controlled singularly foliated manifold $(\cM,\Ax)$ (centered) at the point $p$.
\end{Definition}
For notational simplicity, we shall write {\em vertices}, {\em edges}, {\em faces} of $\Omega$ when referring to the corresponding objects of the Newton polyhedron $\vN$.  We shall also refer to the support of the Newton map $\Theta$ simply as $\supp(\Omega)$.
\subsection{Derived Polygon and Displacements}\label{subsect-derivedpolygon}
\begin{figure}[htb]
\psfrag{F}{\small $\mathcal{F}$ (main face)}
\psfrag{xj/y}{\small $v_2$}
\psfrag{xi}{\small $v_1$}
\psfrag{z}{\small $v_3$}
\psfrag{Np}{\small $\vN^\prime$}
\psfrag{f}{\small $f$ (main side)}
\psfrag{e}{\small $\medge$ (main edge)}
\psfrag{n}{\small $\bfm^\prime$}
\psfrag{m}{\small $\bfn$}
\begin{center}
\includegraphics[height=5cm]{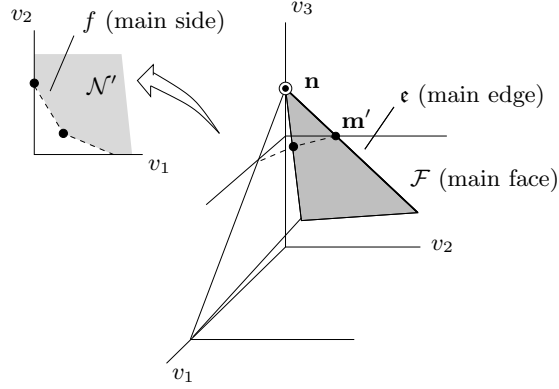}
\end{center}
\caption{The main vertex, the main edge and the derived polygon for $\bfn$.}
\label{fig-invariantdraw2}
\end{figure}

Let us fix a Newton data $\Omega = ((x,y,z),\iota_p,\Theta)$ at a point $p \in A$, and
let $\vN$ be the corresponding Newton polyhedron. The {\em derived polygon} associated to a vertex $\bfn \in \vN$ is given by
$$
\vN^\prime(\bfn) := \vN \cap \{(v_1,v_2,v_3) \in \cR^3 \mid v_3 = n_3 - 1/2\}
$$
Thus, $\vN^\prime(\bfn)$ is a convex polygon contained in the plane $\{\bfv
\mid v_3 = n_3 - 1/2\}$ (see figure~\ref{fig-invariantdraw2}).

\begin{Remark}
The derived polygon has some similarities with the {\em characteristic polygon} introduced by Hironaka is his proof of the resolution of singularities for excellent surfaces (we refer to \cite{H2} for the precise definition of this polygon).  
The following example shows that these two notions are distinct in the context of vector fields:  Consider the germ of vector field
$$
\chi = (z^3 x + xyz^2)\frac{\partial}{\partial x} + xz^3 \frac{\partial}{\partial y} + y^7 \frac{\partial}{\partial z}. 
$$
The associated Newton polyhedron is shown in figure~\ref{fig-charhiron} (left). Let us choose the vertex $\bfn = (0,0,3)$ (this is the minimal vertex of $\vN$ with respect to the lexicographical ordering in $\cR^3$).  Then, the derived polygon and the Hironaka characteristic polyhedron are given respectively by
$$
\vN^\prime = \vN \cap \{\bfv \in \cR^3 \mid v_3 = 5/2\} \quad\mbox{ and }\quad
\vN^{\prime\prime} = \vN \cap \{\bfv \in \cR^3 \mid v_3 = 2\},
$$
The resulting polygons are pictured in the right hand side of figure ~\ref{fig-charhiron}.
\begin{figure}[htb]
\psfrag{v1}{\small $v_1$}
\psfrag{v2}{\small $v_2$}
\psfrag{v3}{\small $v_3$}
\psfrag{Nl}{\small $\vN^\prime$}
\psfrag{Nll}{\small $\vN^{\prime\prime}$}
\psfrag{(0,0,3)}{\small $\bfn = (0,0,3)$}
\psfrag{(0,1,2)}{\small $(0,1,2)$}
\psfrag{(0,6,0)}{\small $(0,0,6)$}
\psfrag{(1,-1,3)}{\small $(1,-1,3)$}
\psfrag{(0,1)}{\small $(0,1)$}
\psfrag{(1/2,0)}{\small $(1/2,0)$}
\psfrag{(0,1/2)}{\small $(0,1/2)$}
\begin{center}
\includegraphics[height=6cm]{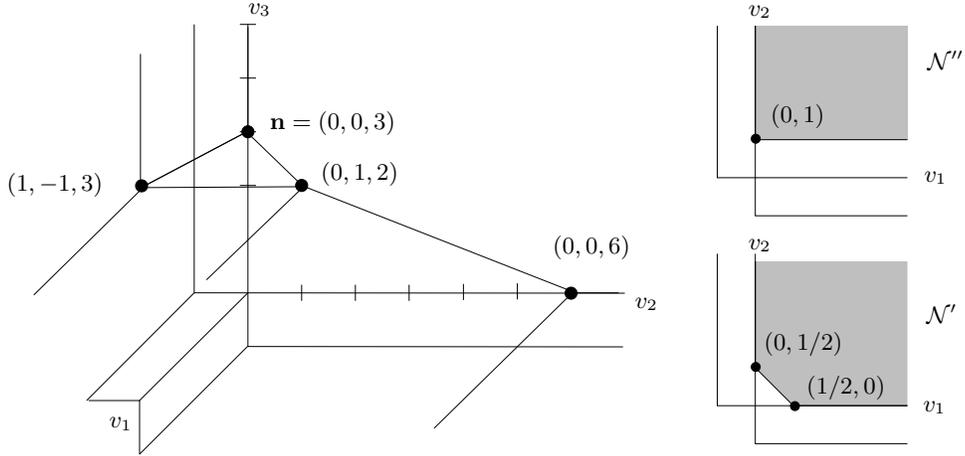}
\end{center}
\caption{The Hironaka characteristic polygon and the derived polygon are distinct.}
\label{fig-charhiron}
\end{figure}
\end{Remark}
\begin{Proposition}\label{prop-derivedpolexists}
Suppose that the vertex $\bfn$ is such that $n_3 \ge 1$.  Then, the
derived polygon $\vN^\prime(\bfn)$ is nonempty.
\end{Proposition}
\begin{proof}
Indeed, suppose by absurd that $\vN^\prime(\bfn) = \emptyset$.  Then, since
$n_3 \ge 1$, the Newton polyhedron $\vN$ should
be contained in the region $\{(v_1,v_2,v_3) \in \cZ^3 \mid v_3 \ge
1\}$.  According to the definition of $\vN$, 
this would imply that the line field $\linef$ is locally generated by a 
vector field $\chi$ which is degenerate (because
the
ideal $\vI_\chi(z)$ would be divisible by $z^2$).  This contradicts our assumptions.
\end{proof}
For the rest of this subsection, let us assume that the derived polygon $\vN^\prime(\bfn)$ is nonempty.

\begin{figure}[htb]
\psfrag{e0}{\small $e_0 $}
\psfrag{e1}{\small $e_1 $}
\psfrag{e2}{\small $e_2 $}
\psfrag{en}{\small $e_n $}
\psfrag{Nl}{\small $\vN^\prime$}
\psfrag{m1l}{\small $m_1^\prime$}
\psfrag{m2l}{\small $m_2^\prime$}
\psfrag{v1}{\small $v_1$}
\psfrag{v2}{\small $v_2$}
\begin{center}
\includegraphics[height=4cm]{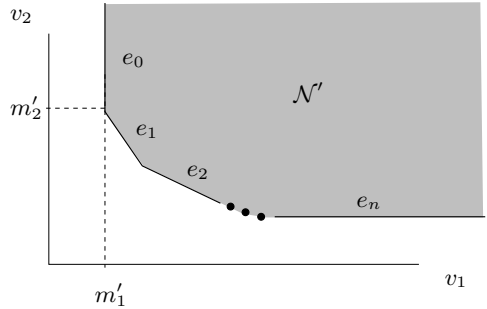}
\end{center}
\caption{The derived polygon.}
\label{fig-derivpolygon}
\end{figure}

The {\em main derived vertex} of $\vN^\prime(\bfn)$ is the minimal vertex
$\bfm^\prime(\bfn)$ of $\vN^\prime(\bfn)$ with
respect to the lexicographical ordering. We write
$$\bfm^\prime(\bfn) = (m_1^\prime(\bfn),m_2^\prime(\bfn),n_3-1/2)$$
The {\em main edge} associated to the vertex $\bfn$ is the unique edge
$\medge(\bfn) \subset \vN$ which contains the segment $\overline{\bfn,\bfm^\prime(\bfn)}$.
\begin{Proposition}\label{prop-grid}
The rational numbers $m_1^\prime(\bfn),m_2^\prime(\bfn)$ always belong to the
finite grid
$$\frac{1}{2(n_3 + 1)!}\cZ$$
\end{Proposition}
\begin{proof}
Indeed, the main edge associated to $\bfn$ has the form
$\medge(\bfn) = \overline{\bfn,\bfv}$, for some vertex 
$\bfv = (v_1,v_2,v_3)$
such that $-1 \le v_3 < n_3$. Then, it is clear that
$$
m_1^\prime(\bfn) = n_1 + \frac{v_1-n_1}{2(n_3-v_3)}\quad\mbox{ and }\quad
m_2^\prime(\bfn) = n_2 + \frac{v_2-n_2}{2(n_3-v_3)}
$$
Now, it suffices to remark that the denominator of such fractions
always lies in the range $\{1,\ldots,2(n_3 + 1)!\}$.
\end{proof}
We picture $\vN^\prime(\bfn)$ in the 2-dimensional plane as in
figure~\ref{fig-derivpolygon},
with the horizontal
 axis corresponding to the $v_1$-coordinate.
Using this representation, we enumerate the sides of
 $\vN^\prime(\bfn)$ from left to right as
$$e_0,e_1,\ldots,e_n$$
with $e_0$ being the infinite vertical side and $e_n$ being the infinite
 horizontal side.

The {\em main side} $f(\bfn)$ of $\vN^\prime(\bfn)$ is defined as follows (see
figure~\ref{fig-derivpolygon}):
\begin{itemize}
\item[(i)] If $m_1^\prime(\bfn) > 0$ then $f(\bfn) := e_0$;
\item[(ii)]  If $m_1^\prime(\bfn) = 0$ then $f(\bfn) = e_1$.
\end{itemize}
By the definition of $\vN^\prime(\bfn)$, to each each side $e \in \vN^\prime(\bfn)$ there
corresponds a unique face $F$ of $\vN$ such that
$$
F \cap \vN^\prime(\bfn) = e
$$
The {\em main face} associated to the vertex $\bfn$ is the unique face $\vF(\bfn) \subset \vN$ such that
$$\vF(\bfn) \cap \vN^\prime(\bfn) = f(\bfn)$$
(see figure~\ref{fig-invariantdraw2}). By construction, the edge $\medge(\bfn)$ can be uniquely written as
$$\medge(\bfn) = \{\bfn + t (\Delta,-1)\mid t \in
I\}$$
where $I \subset \cR$ is a compact interval and
$\Delta \in \cQ^2 \setminus \{(0,0)\}$ is a nonzero vector of rational
numbers.  Using this, the derived vertex $\bfm^\prime(\bfn)$ can be rewritten as $\bfm^\prime(\bfn) = (n_1 + \Delta_1/2,n_2 + \Delta_2/2,n_3 - 1/2)$. 

\begin{figure}[htb]
\psfrag{D}{\small $(\Delta,-1)$}
\psfrag{C}{\small $(C,-1,0)$}
\psfrag{m}{\small $\bfn$}
\psfrag{F}{\small $\vF$}
\psfrag{C = 0}{\small $C = 0$}
\psfrag{0 < C < infty}{\small $0 < C < \infty$}
\psfrag{C = infty}{\small $C = \infty$}
\begin{center}
\includegraphics[height=6cm]{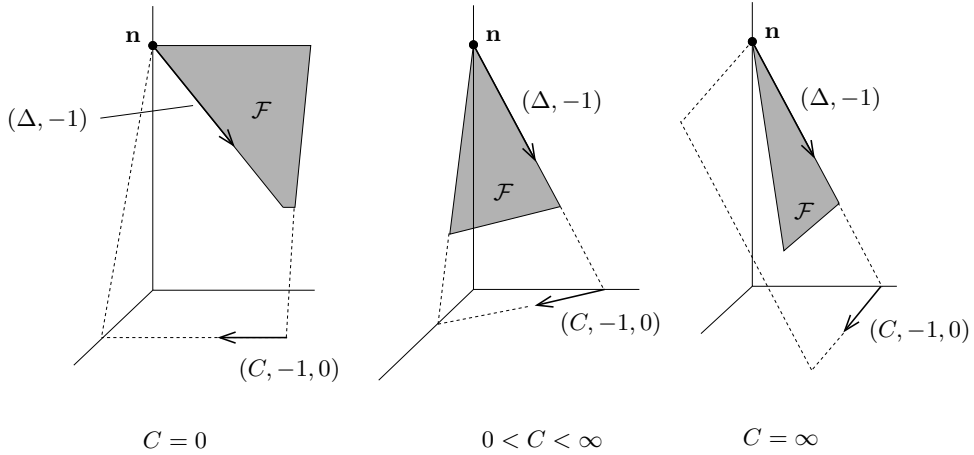}
\end{center}
\caption{The vertical and horizontal displacements.}
\label{fig-DeltaC}
\end{figure}

Similarly, the main side $f(\bfn)$ of $\vN^\prime(\bfn)$ can be uniquely written as
$$f(\bfn) = \{\bfm^\prime(\bfn) + t (C,-1,0)\mid t \in I\}\}$$
where $I \subset \cR$ is an interval and $C$ is a number in
$\overline{\cQ}_{\ge 0} := \cQ_{\ge 0} \cup \{\infty\}$.
\begin{Remark}
Observe that $C = \infty$ and $C = 0$ correspond to the
cases where the main side is
the infinite horizontal and infinite vertical side, respectively (see figure~\ref{fig-DeltaC}).
\end{Remark}
We will call $\Delta(\bfn) := \Delta$ and $C(\bfn) := C$ respectively the {\em vertical displacement vector} and 
{\em horizontal displacement} associated to the vertex $\bfn$.

\subsection{Regular/Nilpotent Configurations and Main Vertex}\label{subsect-regnilpmain}
Let us keep the notations of the previous subsection.  In this subsection, we further assume that the base point $p \in A$ belongs to the divisor $\vD$.

The {\em higher vertex} is the minimal point $\bfh \in \vN$
with respect to the lexicographical ordering in $\cR^3$. It is immediate to see such minimal point always exists 
and that it is a vertex of $\vN$ (see Remark~\ref{remark-newton}).
\begin{Proposition}\label{prop-propbfm}
The higher vertex $\bfh = (h_1,h_2,h_3)$ has the following properties:
\begin{enumerate}[(i)]
\item If $\#\iota_p = 1$ then $h_1 = 0$, $h_2,h_3 \ge -1$;
\item If $\#\iota_p = 2$  then $h_1 = 0$, $h_2 \ge
  0$ and $h_3 \ge -1$;
\end{enumerate}
\end{Proposition}
\begin{proof}
To prove (i), we observe that the surface
$\{x = 0\}$ is preserved by the vector field $\chi$ if and only if
$$\supp(\vN) \cap \{-1\} \times \cN^2 = \emptyset$$
This is equivalent to say that $h_1 \ge 0$ and $h_2,h_3 \ge -1$.  However, if $h_1 \ge 1$, the ideal $\vI_\chi(x)
\in \cR\{x,y,z\}$,
which is generated by $(f x,g xy, h xz)$, would be divisible by
$x^2$.  This would contradict the hypothesis that $\chi$ is a
nondegenerate vector field (see definition~\ref{def-nondegvf}).  The proof of
(ii) is analogous.
\end{proof}
The main edge associated to the higher vertex is given by
\begin{equation}\label{medgebfpbfn}
\medge(\bfh) = \overline{\bfh,\bfn}
\end{equation}
where $\bfn$ is also a vertex of $\vN$.  
It follows from Proposition~\ref{prop-derivedpolexists} that such edge always
exists if $h_3 \ge 1$.  We define $\medge(\bfh) := \emptyset$ if the derived polygon $\vN^\prime(\bfh)$ is empty.
 
We shall say that the Newton data $\Omega$ is in a {\em nilpotent configuration} if the
following three conditions are satisfied:
\begin{enumerate}[(i)]
\item $\#\iota_p = 1$, 
\item $\bfh = (0,-1,h_3)$, for some integer $h_3 \in \cN$, and 
\item $\bfn = (0,0,n_3)$, for some integer $n_3 < h_3$.
\end{enumerate}
If one of these conditions fails, we shall say that $\Omega$ is in
a {\em regular configuration}.

\begin{figure}[htb]
\psfrag{Fbl}{\small $\vF$}
\psfrag{m}{\small $\bfh$}
\psfrag{n}{\small $\bfn$}
\psfrag{medge}{\small $\medge(\bfh)$}
\psfrag{Exceptional configuration}{\small Nilpotent Configuration ($\bfm = \bfn$)}
\psfrag{Regular configuration}{\small Regular Configuration ($\bfm = \bfh$)}
\begin{center}
\includegraphics[height=5cm]{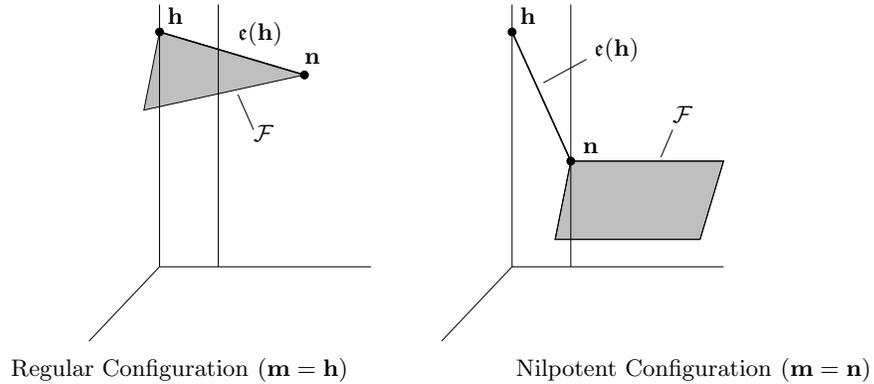}
\end{center}
\caption{The Regular and Nilpotent Configurations.}
\label{fig-regresconf1}
\end{figure}

\begin{Remark}
As we shall see later on section~\ref{sect-LocalTheory}, the treatment of nilpotent configurations constitutes one of the points where the method of resolution of singularities for vector fields differs {\em essentially} from the usual methods of resolution of singularities for functions and analytic sets.  At several points during our proof, we will have the address the delicate issue of the {\em transition} between regular and nilpotent configurations.
\end{Remark}
The {\em main vertex} $\bfm$ of $\Omega$ is chosen as follows:
\begin{enumerate}[(i)]
\item If $\Omega$ is in a regular configuration then $\bfm := \bfh$, and
\item If $\Omega$ is in an nilpotent configuration then $\bfm := \bfn$.
\end{enumerate}
(where $\bfn$ is the vertex defined by (\ref{medgebfpbfn})).
The corresponding vertical and horizontal displacements 
$$\Delta := \Delta(\bfm), \quad C := C(\bfm)$$ 
will be called {\em vertical displacement vector} and {\em horizontal displacement} of $\Omega$. 

The face $\vF := \vF(\bfm)$ and the edge $\medge := \medge(\bfm)$
will be called respectively the {\em main face} and {\em main edge} associated to $\Omega$. The polygon $\vN^\prime = \vN^\prime(\bfm)$ will  be called {\em main derived polygon}.

\subsection{The Class $\Nmaps^{i,\bfm}_{\Delta,C}$}
Let $\Omega = ((x,y,z),\iota_p,\Theta)$ be a Newton data for $(\cM,\Ax)$ at a point $p \in A \cap \vD$.  We shall say that $\Omega$ belongs to the {\em class $\Nmaps^{i,\bfm}_{\Delta,C}$}  if 
\begin{enumerate}[(i)]
\item $\#\iota_p = i$
\item $\bfm$ is the main vertex of $\Omega$;
\item $\Delta$ is the vertical displacement vector of $\Omega$;
\item $C$ is the horizontal displacement of $\Omega$.
\end{enumerate}
The union of all classes of Newton data will be denoted simply
by $\Nmaps$.

Let us consider the Lie group $\gG$ of all polynomial changes of coordinates 
in $\cR^3$ which have the form 
\begin{equation}\label{Gmapdef}
\widetilde{x} = x,\quad \widetilde{y} = y + g(x), \quad \widetilde{z} = z + f(x,y)
\end{equation}
where $f \in \cR [x,y]$ and $g \in \cR[x]$ are real polynomials.  
The group operation is the composition and the
inverse of the map (\ref{Gmapdef}) is simply given by
\begin{equation}\label{invfg}
x = \widetilde{x},\quad y = \widetilde{y} - g(\widetilde{x}), \quad z = \widetilde{z} - f(\widetilde{x},\widetilde{y} - g(\widetilde{x})).
\end{equation}
For shortness, we shall denote the map (\ref{Gmapdef}) by $(f,g) \in \gG$ and its inverse
by $(f,g)^{-1}$.  We define also
the subgroups $\gG^1 = \gG$ and $\gG^2 = \{(f,g) \in \gG \mid g = 0 \}$.
\begin{Remark}
The Lie algebra associated to $\gG$ is the algebra $\lG$ of all 
polynomial vector fields of the form
$$
G(x)\frac{\partial}{\partial y} + F(x,y) \frac{\partial}{\partial z}
$$
with polynomials $G \in \cR[x]$ and $F \in \cR[x,y]$.
\end{Remark}
There is a natural action of the Lie group $\gG$ on the class of Newton data
$\Nmaps$, given as follows:
The action of a map $(f,g) \in \gG$ in the data $\Omega = ((x,y,z),\iota,\Theta)$ is the Newton data given by
$$((\widetilde{x},\widetilde{y},\widetilde{z}),\widetilde{\iota},\widetilde{\Theta}),$$
where $(\widetilde{x},\widetilde{y},\widetilde{z})$ are the coordinates
given by (\ref{Gmapdef}), the list $\widetilde{\iota}$ is the incidence list at the point 
$\widetilde{p} = (\widetilde{x},\widetilde{y},\widetilde{z})^{-1}(0)$ and $\widetilde{\Theta}$ is the Newton map for $(\cM,\Ax)$ at $\widetilde{p}$, relatively to such new adapted local chart.
We denote such action simply by $(f,g)\cdot \Omega$.
\begin{Remark}
In the cases that we will consider more often, we have $f(0) = g(0) = 0$. In this case, $\widetilde{p} = p$, i.e.\ the Newton data $\widetilde{\Omega}$ is centered at the same point as $\Omega$.  If this is not the case, we tacitly assume that the point $\widetilde{p}= (\widetilde{x},\widetilde{y},\widetilde{z})^{-1}(0)$ lies in the domain of definition of the local adapted chart $(U,(x,y,z))$.
\end{Remark}
\subsection{The Subgroups $\gG_{\Delta,C}$}\label{subsect-subgroupsgG}
Recall that the {\em support} of a polynomial $H \in \cR[\bfx]$ is
the subset 
$$\supp(H) = \{\bfv \in \cZ^n \mid \bfx^\bfv \mbox{ is a nonzero monomial 
of }H \}
$$ 
According to the {\em support} of the polynomials
$f$ and $g$ given in (\ref{Gmapdef}), we shall now define several subgroups in $\gG$.

Given $\Delta = (\Delta_1,\Delta_2) \in \cQ^2_{\ge 0}$ and $C \in \overline{\cQ}_{\ge 0}$, we define
$\gG_{\Delta,C}$ as the subgroup of all maps $(f,g) \in \gG$ such that the
following
conditions hold:
\begin{enumerate}[(i)]
\item the support $S_f = \supp(f)$ is contained in the set
\begin{eqnarray*}
\{(a,b) \in \Delta + s (C,-1) \mid s \in \overline{\cQ}_{\ge 0}\} \; \cap\;  \cN^2
\mbox{ if }\Delta_1 = 0\\
\{\Delta\}\; \cap\; \cN^2
\mbox{ if }\Delta_1 > 0
\end{eqnarray*}
\item the support $S_g = \supp(g)$ is contained in the set
\begin{eqnarray*}
\{C\}\; \cap\; \cN \mbox{ if }\Delta_1 = 0\\
\emptyset \mbox{ if }\Delta_1 > 0
\end{eqnarray*}
\end{enumerate}
We further define the subgroups $\gG^1_{\Delta,C}$ and $\gG^2_{\Delta,C}$
as
$$\gG^1_{\Delta,C} = \gG_{\Delta,C} \qquad  
\gG^2_{\Delta,C} = \gG_{\Delta,C} \cap \{(f,g)\in \gG_{\Delta,C} \mid g = 0\}$$
\begin{Remark}
In the above definition, we have the following {\em extreme} cases for 
$\gG_{\Delta,C}$:
\begin{itemize}
\item If $C = \infty$ then $g = 0$ and $f = \xi x^{\Delta_1}y^{\Delta_2}$,
      where the constant $\xi \in \cR$ necessarily vanishes if $\Delta \notin \cN^2$.
\item If $\Delta_1 = 0$ and $C = 0$ then $f \in \cR[y]$ is a polynomial in
$y$ of degree at most $\delta_2$ and $g = \eta$, for some constant $\eta \in \cR$.
\item If $\Delta = (0,0)$ and $C = \infty$ then $g = 0$ and $f = \xi$, for
      some real constant $\xi \in \cR$.
\end{itemize}
In the last two cases, the change of coordinates (\ref{Gmapdef})
correspond to translations $\widetilde{y} = y+\eta, \widetilde{z} = z+f(y)$ and 
$\widetilde{y} = y, \widetilde{z} = z + \xi$, respectively.
\end{Remark}
It will be useful to consider the following decomposition of the group $\gG$:  Define the
subgroup
$$
\gG_{\Delta} := \{(f,g) \in \gG_{\Delta,C}\mid g = 0,f = \xi x^{\Delta_1}y^{\Delta_2}, \; \xi \in \cR\}
$$
(the constant $\xi$ necessarily vanishes if $\Delta \notin \cN^2$), and the normal subgroup $$\gG^{+}_{\Delta,C} = \{(f,g) \in \gG_{\Delta,C} \mid \Delta \notin \supp(f)\}$$
which will be called the subgroup of {\em edge preserving maps}.  It is easy to see that
$$
\gG^{+}_{\Delta,C} \cap \gG_{\Delta} = \{0\}\quad \mbox{ and }\quad
\gG_{\Delta,C} = \gG_{\Delta}\circ \gG^{+}_{\Delta,C} = \gG^{+}_{\Delta,C} \circ\gG_{\Delta}
$$
In other words, $\gG_{\Delta,C}$ is the {\em semi-direct product} of 
$\gG_{\Delta}$ and $\gG^{+}_{\Delta,C}$.  Similar decompositions holds for the subgroups
$\gG^1$ and $\gG^2$.
\begin{Remark}\label{remark-supportedgepreserving}
Later on, we shall need the following remark: For a map $(f,g) \in \gG^{+}_{\Delta,C}$, the support of $f$ is such that
$$
S_f \subset \{(a,b) \in \cN^2 \mid a \ge C \, \tau(\Delta_2) \}
$$
where we define $\tau(\Delta_2) := 1$ if $\Delta_2 \in \cN$, and $\tau(\Delta_2) := \Delta - \linteg
\Delta \rinteg$, otherwise.
\end{Remark}
\subsection{Action of $\gG$ via Adjoint Map}
The action of the group $\gG$ on a Newton data $\Omega$ can be studied via the associated
Lie algebra $\lG$.  Indeed, a map $(f,g) \in \gG$ is the time one map of the flow
associated to the vector field $\Gamma_{f,g} \in \lG$ given by
$$
\Gamma_{f,g} = g(x) \frac{\partial}{\partial y} + f(x,y) \frac{\partial}{\partial z}
$$
Therefore, if $\Omega$ is associated to the vector field $\chi$,
the Newton data $(f,g)\cdot \Omega$ can be obtained from the transformed
vector field 
\begin{equation}\label{BCHformula}
((f,g))_* \chi = \chi + \frac{1}{2} [\Gamma_{f,g},\chi] +
\frac{1}{6}[\Gamma_{f,g},[\Gamma_{f,g},\chi]] + \cdots = \sum_{n=0}^{\infty}{\frac{1}{n!}
(\mathrm{ad}(\Gamma_{f,g}))^n\; \chi}
\end{equation}
because $\mathrm{e}^{\mathrm{ad}(\cdot)} = \mathrm{Ad}(\mathrm{Exp}(\cdot))$.  

Using such remark, we can see how the action of a map in $\gG_{\Delta,C}$ modifies the {\em multiplicity} of a vector field. Let $\bfomega = (\omega_1,\omega_2,\omega_3) \in \cN^3$ be a weight-vector such that $\gdc(\omega_1,\omega_2,\omega_3) = 1$ and 
$$
\< \bfomega, (-1,\Delta_1,\Delta_2) \> = \< \bfomega,(0,-1,C) \> = 0
$$
(with the convention that $\omega_1 = 0$ if $C = \infty$).  Then, $\bfomega$
is the weight-vector associated to the main
face $\vF$ of the Newton polyhedron  $\vN(\Omega)$.

Recall from subsection~\ref{subsect-blowup} that each analytic vector field 
$\chi$ has an associated $\bfomega$-multiplicity
$$
\mu_\bfomega(\chi) := \max \{ k \in \cZ \mid \chi(H_\bfomega^d) \subset
H_\bfomega^{d + k}, \forall d \in \cN\}
$$
\begin{Lemma}\label{lemma-BCH}
If $(f,g) \in \gG_{\Delta,C}$ is nonzero 
then $\mu_\bfomega(\Gamma_{f,g}) = 0$.
\end{Lemma}
\begin{proof}
This follows directly from the definition of the group $\gG_{\Delta,C}$.
\end{proof}
As a consequence, if $\mu_\bfomega(\chi) = m$ and $(f,g)$ is nonzero then
$\mu_\bfomega([\Gamma_{f,g},\chi]) = m$. As a consequence, we have the following result:
\begin{Corollary}\label{corollary-BCH}
Choose $(f,g) \in \gG_{\Delta,C}$.  Then, the vector field
$\widetilde{\chi} = (f,g)_* \chi$ is such that
$$\mu_\bfomega(\widetilde{\chi}) = \mu_\bfomega(\chi).$$
\end{Corollary}
\begin{proof}
It suffices to use the Lemma~\ref{lemma-BCH} and the formula 
(\ref{BCHformula}).
\end{proof}
In the same way, we can prove that the coordinate change associated to 
$(f,g)$ always preserves the {\em $\bfomega$-quasi-homogeneous structure} in $\cR^3$
(see subsection~\ref{subsect-blupmanif}).

Let us now study the action of $\gG_{\Delta,C}$ on a single {\em differential monomial} given by
$$m = x^{v_1} y^{v_2} z^{v_3} \left( \alpha x \frac{\partial}{\partial x} + 
\beta y \frac{\partial}{\partial y} + \gamma z \frac{\partial}{\partial z} \right)$$
with constants $\bfv = (v_1,v_2,v_3) \in \cZ^3$ and $\alpha,\beta,\gamma \in \cR$.  The corresponding Newton map $\Theta$ is such that $\supp(\Theta) = \{\bfv\}$. 

\begin{figure}[htb]
\psfrag{v}{\small $\bfv$}
\psfrag{v1}{\small $v_1$}
\psfrag{v2}{\small $v_2$}
\psfrag{v3}{\small $v_3$}
\psfrag{Delta}{\small $(\Delta_1,\Delta_2,-1)$}
\psfrag{C}{\small $(C,-1,0)$}
\begin{center}
\includegraphics[height=5cm]{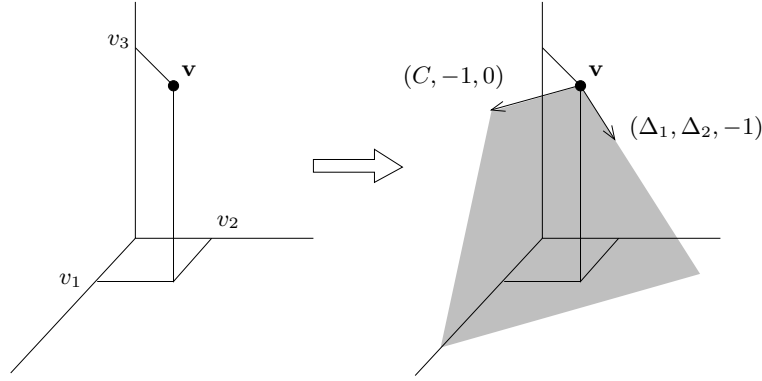}
\end{center}
\caption{The support of $(f,g)_* m$ for a differential monomial $m$.}
\label{fig-actionNewtonmap}
\end{figure}

For the particular case of a map $(f,g) \in \gG_{\Delta,C}$ of the form 
$f = \xi x^{\delta_1} y^{\delta_2}$ and $g = 0$, the coordinate change 
$\widetilde{y} = y + g(x)$ and 
$\widetilde{z} = z + f(x,y)$ maps $m$ to the vector field 
$$\widetilde{m} = x^{v_1} {y}^{v_2} ({z} - \xi x^{\delta_1} {y}^{\delta_2})^{v_3} 
\left( \alpha x \frac{\partial}{\partial x} + \beta {y} \frac{\partial}{\partial {y}} + 
(\gamma {z} + \xi ((\alpha \delta_1 + \beta\delta_2) - \gamma)x^{\delta_1} y^{\delta_2})\frac{\partial}{\partial {z}}\right)$$
(where we drop the tildes). In particular, it is easy to see that the Newton map $\widetilde{\Theta}$ associated to $\widetilde{m}$ has
support contained in the set 
$$
\{\bfu \in \cZ^3 \mid \bfu = \bfv + t(\delta_1,\delta_2,-1),\; t \ge 0\}
$$
Similarly, for a map $(f,g)$ of the form $(f,g)  = (0,\eta x^C)$, we get (dropping the tildes)
$$\widetilde{m} = x^{v_1} ({y} -\eta x^C)^{v_2} {z}^{v_3} 
\left( \alpha x \frac{\partial}{\partial x} + (\beta  {y} + \eta(C \alpha - \beta)) \frac{\partial}{\partial {y}} + \gamma {z} \frac{\partial}{\partial {z}}\right),$$
and the the Newton map $\widetilde{\Theta}$ associated to $\widetilde{m}$ has
support contained in the set 
$$
\{\bfu \in \cZ^3 \mid \bfu = \bfv + s(C,-1,0),\; s \ge 0\}.
$$
Now, an arbitrary map $(f,g) \in \gG_{\Delta,C}$ can be written as the composition of a finite number of maps of the form $(\xi x^{\delta_1} y^{\delta_2},0)$ and 
$(0,\eta x^C)$.  Therefore, the above computations give the following result:
\begin{Lemma}
Consider the differential monomial $m$ given above.  Then, for an arbitrary pair
$(\Delta,C)$, and for an arbitrary map $(f,g) \in \gG_{\Delta,C}$, the Newton data for the vector field $(f,g)_* m$ has its support contained in the set
$$
\{\bfu \in \cZ^3 \mid \bfu = \bfv + t(\Delta_1,\Delta_2,-1) + s(C,-1,0),\; t,s \ge 0\}
$$
(see figure~\ref{fig-actionNewtonmap}).  
\end{Lemma}

More generally, we can consider the action of $(f,g)$ on an arbitrary vector field
$\chi$ as follows.  Write the expansion of $\chi$ as 
$$
\chi = \sum_{i \in I} m_i.
$$
where, for each $i \in I$,  $m_i$ is a differential monomial whose Newton data has support at $\bfv_i \in \cZ^3$. For a map $(f,g) \in \gG_{\Delta,C}$, we clearly have
$$
(f,g)_* \chi = \sum_{i \in I} (f,g)_* m_i.
$$
and this gives the following result:
\begin{Corollary}\label{corollary-howittransforms}
Let $\chi$ be as above. For an arbitrary pair $(\Delta,C)$, and for an arbitrary map $(f,g) \in \gG_{\Delta,C}$, the Newton data for the vector field $(f,g)_* \chi$ has its support contained in the set
$$
\bigcup_{i\in I}\{\bfu \in \cZ^3 \mid \bfu = \bfv_i + t(\Delta_1,\Delta_2,-1) + s(C,-1,0),\; t,s \ge 0\}
$$
\end{Corollary}
\begin{Remark}
In \cite{AGV}, the authors use these kind of coordinate changes to study normal
forms of {\em quasi-homogeneous functions}.  In \cite{H2}, Hironaka uses similar transformations in his definition of well and very well preparations of function germs.
\end{Remark}
\section{Local Theory at $\NElem \cap \vD$}\label{sect-LocalTheory}
Let $(\cM,\Ax)$ be a controlled singularly foliated manifold.
Through this section, we fix a divisor point $p \in A \cap \vD$,
an adapted local chart $(U,(x,y,z))$ for 
$(\cM,\Ax)$ at $p$ and let $\Omega \in \Nmaps^{i,\bfm}_{\Delta,C}$ be the
corresponding Newton data.
\begin{Proposition}
The main vertex 
$\bfm = (m_1,m_2,m_3)$ associated to $\Omega$ is such that $m_1 = 0$ and $m_2 \in \{-1,0\}$.
\end{Proposition}
\begin{proof}
Suppose by absurd that the main vertex $\bfm$ does not satisfy the above requirements.  Then, 
one of the following conditions holds:
\begin{enumerate}[(a)]
\item $m_1 = -1$, or
\item $m_1 \ge 1$, or
\item $m_2 \ge 1$.
\end{enumerate}
In the case (a), it follows from the definition of the Newton polyhedron that
the plane $\{x=0\}$ is not invariant.  This contradicts the definition of an
adapted local chart at $p$ and the hypothesis that $p$ belongs to $A \cap \vD$.

In the case (b), choose a nondegenerate vector field $\chi$ which is a
local generator for the line field $L$.  Then, the condition $m_1 \ge 1$ implies
that
the ideal $\vI_\chi(x) \subset \vO_p$ is divisible by $x^2$.  This contradicts
the definition~\ref{def-nondegvf}.

Therefore, in the case (c) we can assume that $m_1 = 0$ and $m_2 \ge 1$.
This clearly implies that
\begin{equation}\label{suppomegaempty}
\supp(\Omega)\; \cap \; \left( \{0\}\times \{-1,0\} \times \cR \right) \; = \; \emptyset.
\end{equation}
If we write the vector field $\chi$ as
$$
\chi = f x\frac{\partial}{\partial x} + g y \frac{\partial}{\partial y} + h z \frac{\partial}{\partial z}
$$
(with $fx,gy,hz \in \cR\{x,y,z\}$), then (\ref{suppomegaempty}) is equivalent
to say that the functions $f$, $g$ and $h$ vanish identically along the
vertical line $l := \{x = y = 0\}$.
A simple computation
shows that is is equivalent to assert that the Jacobian matrix $D\chi|_l$ has
the form
$$
D\chi|_l = \left[
\begin{matrix}
0 & 0 & 0 \cr
* & 0 & 0 \cr
* & * & 0 \cr
\end{matrix}
\right]
$$
where $*$ denotes some arbitrary real numbers.
As a consequence, the line $l$ is contained in the set of nonelementary points 
$\NElem(\cM)$, which contradicts definition~\ref{def-axis}.  Absurd.
\end{proof}

\subsection{Stable Newton Data and Final Situations}
In the next definitions, we consider the action of the transformation group 
$\gG^i_{\Delta,C}$ on the Newton data $\Omega$. The following notions will be essential in the sequel to study the effect of the {\em translations} is the blowing-up chart.

We say that $\Omega$ is {\em stable} if 
$$
(f,g)\cdot\Omega \in \Nmaps^{i,\bfm}_{\Delta,C}.
$$
for all $(f,g) \in \gG^i_{\Delta,C}$ with $f(0) = g(0) = 0$.
In other words, $\Omega$ is stable if
the action of $\gG^i_{\Delta,C}$ preserves the main vertex, the value of
the vertical displacement vector
$\Delta$ and the value of the horizontal displacement $C$. 

A weaker notion of stability will also be useful.
We say that  $\Omega$ is {\em edge
stable} if for each map $(f,g) \in \gG^i_{\Delta,C}$ with $f(0) = g(0) = 0$ there exists a 
constant $\widetilde{C} \in \overline{\cQ}_{\ge 0}$ such that
$$
(f,g)\cdot\Omega \in \Nmaps^{i,\bfm}_{\Delta,\widetilde{C}}.
$$
\begin{Remark}\label{Remark-translations}
Intuitively,  the notions of stable and edge stable Newton data can be seen as weaker versions of the notion of {\em maximal contact} introduced in the work of Hironaka \cite{H3}.  As we said above, the main goal is to take into account the effect of the translations in the blowing-up chart.\\
More precisely, for a stable Newton data (resp.\ edge stable Newton data) one guarantees that the main invariant strictly decreases after a conveniently chosen blowing-up map followed by any translation of the form $(\widetilde{y} = y + \eta, \widetilde{z} = z + \xi)$ (resp.\ 
$(\widetilde{z} = z + \xi)$).  We refer to subsection~\ref{subsect-Newtoninvarlocres} for the precise statements.\\
In the context of vector fields, the usual notion of maximal contact is too strong and often leads {\em divergent} formal objects.  For instance, the computation of the maximal contact variety (in the sense of \cite{H3}) for the {\em Euler vector field}
$$
x^2\frac{\partial}{\partial x} + (y - x) \frac{\partial}{\partial y}
$$
leads to the formal power series $V = \{y - \sum_{n \ge 1}(n-1)! x^n\}$, which has zero radius of convergence.
\end{Remark}
Using such concepts, we can now identify when the Newton data $\Omega$ is centered at an  
{\em elementary point} $p \in \Elem(\cM)$.

\begin{figure}[htb]
\psfrag{i.a}{\small (i.$a$)}
\psfrag{i.b}{\small (i.$b$)}
\psfrag{i.c}{\small (i.$c$)}
\psfrag{ii}{\small (ii)}
\psfrag{iii}{\small (iii)}
\psfrag{iv}{\small (iv)}
\psfrag{v2}{\small $v_2$}
\psfrag{v3}{\small $v_3$}
\psfrag{mainvertex0}{\small $\bfm = (0,0,0)$}
\psfrag{mainvertex1}{\small $\bfm = (0,-1,0)$}
\psfrag{mainvertex2}{\small $\bfm = (0,0,-1)$}
\psfrag{mainvertex3}{\small $\bfn = (0,0,-1)$}
\begin{center}
\includegraphics[height=7.5cm]{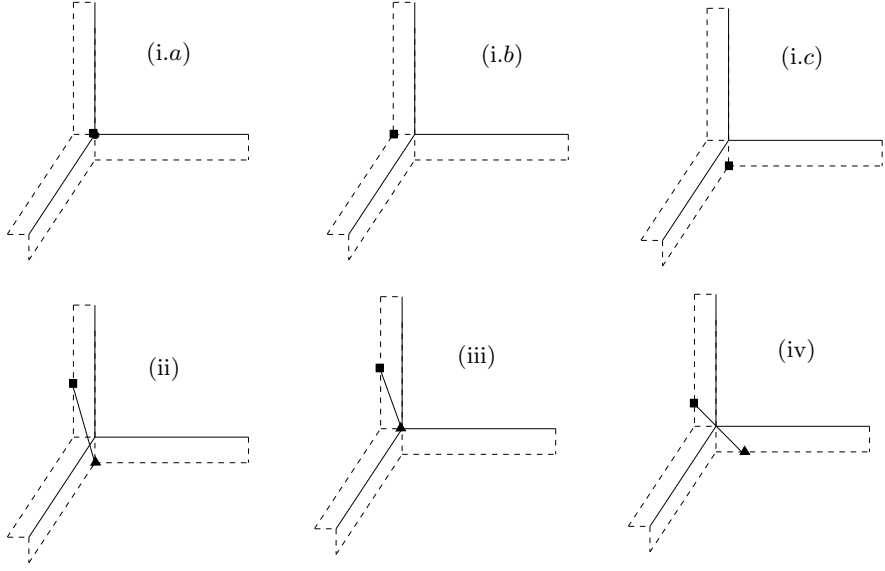}
\end{center}
\caption{The final situations.}
\label{fig-finalsit}
\end{figure}

First of all, we introduce the following notion.
We shall say that $\Omega$ is in a {\em final situation} if the
following holds: If we look at the higher vertex $\bfh  \in \vN$ (see
definition in subsection~\ref{subsect-regnilpmain}) and the associated
edge $\medge(\bfh)$, one of the following conditions is satisfied (see figure~\ref{fig-finalsit}):
\begin{itemize}
\item[(i)] the vertex $\bfh = (h_1,h_2,h_3)$ is such that $h_1 = 0$ and either
$$
(a)\;\; (h_2,h_3) = (0,0), \mbox{ or }\quad (b)\;\; (h_2,h_3) = (-1,0), \mbox{
or }\quad (c)\;\; (h_2,h_3) = (0,-1)
$$
\item[(ii)] The edge $\medge(\bfh)$ is given by $[(0,-1,k),(0,0,-1)]$, for some
  $k \ge 1$.
\item[(iii)] The edge $\medge(\bfh)$ is given by $[(0,-1,k),(0,0,0)]$, for some
  $k \ge 1$.
\item[(iv)] The edge $\medge(\bfh)$ is given by $[(0,-1,1),(0,1,-1)]$ and $\Omega$ is edge stable.
\end{itemize}

The following result justifies the above nomenclature:
\begin{Proposition}\label{prop-finalsituation}
If the Newton data $\Omega$ is in a final situation then it is centered at an elementary point
$p \in \Elem(\cM)$. 
\end{Proposition}
\begin{proof}
Consider a vector field $\chi$ which locally generates the line field
$\linef$, in a neighborhood of $p$. 
If $\chi(p) \ne 0$ we are done. Otherwise, we can write the linear part $D\chi(p)$ as the matrix
$$
\left[
\begin{matrix}
\lambda & 0 & 0 \cr
* & a & b \cr
* & c & d 
\end{matrix}
\right]
$$
where $\lambda, a, b, c, d \in \cR$ and the symbol $*$ denotes 
some arbitrary real constants.  We consider the following cases: 
\begin{enumerate}[(a)]
\item $\lambda \ne 0$;
\item $\lambda = 0$
\end{enumerate}
In case (a), it is clear that $\lambda$ belongs to the spectrum of $D\chi(p)$ and therefore 
$\chi$ is elementary.

In case (b), it clearly suffices to prove the following claim:  {\em The matrix
$$
B = \left[
\begin{matrix}
a & b \cr
c & d 
\end{matrix}
\right]
$$
is not nilpotent (it is obvious that $B \ne 0$).}

To prove such claim, suppose initially that $b = 0$ or $c = 0$.  Then, it
follows from the definition of final situations that 
$(a,d) \ne (0,0)$ and therefore $B$ contains at least one nonzero eigenvalue.

Suppose now that $b \ne 0$ and $c \ne 0$.  It follows that $\Delta = (0,1)$ and $\bfm = (0,-1,1)$.
Assume by absurd that $B$ is nilpotent and consider the $\gG^i_{\Delta,C}$-map 
$(f,g) = ((a/b) y, 0)$, which corresponds
to the coordinate change
$$\widetilde{y} = y, \quad \widetilde{z} = z + (a/b) y$$ 
(see figure~\ref{fig-finalsittrans}). It is easy to see that the Newton data $(f,g)\cdot \Omega$ 
(associated to the local chart 
$(U,(x,\widetilde{y},\widetilde{z})$)
is such that the corresponding matrix $B$ is given by
$$
B = \left[
\begin{matrix}
0 & b \cr
0 & 0 
\end{matrix}
\right]
$$
\begin{figure}[htb]
\psfrag{(f,g)Omega}{\small $(f,g)\cdot \Omega$}
\psfrag{Omega}{\small $\Omega$}
\psfrag{m}{\small $\bfm$}
\begin{center}
\includegraphics[height=4cm]{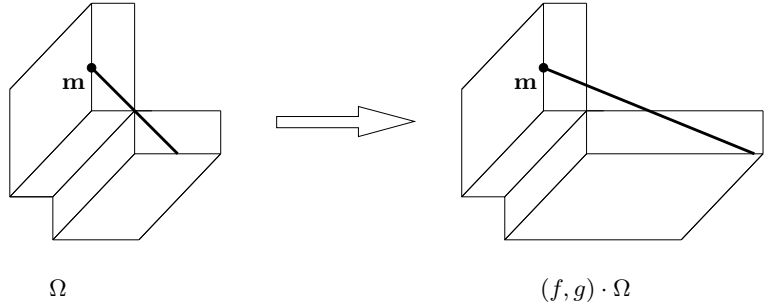}
\end{center}
\caption{The transition from $\Omega$ to $(f,g)\cdot\Omega$.}
\label{fig-finalsittrans}
\end{figure}

This implies that $(f,g)\cdot \Omega$ belongs to the class
$\Nmaps^{i,\bfm}_{\widetilde{\Delta},\widetilde{C}}$, for some
vertical displacement vector $\widetilde{\Delta} >_\lex (0,1)$.  
This contradicts the assumption that $\Omega$ is edge stable.  The claim is proved.
\end{proof}
\begin{Remark}
The above result have the following partial converse (which we will not need in the sequel).  If $\Omega$ is an edge stable Newton data centered at an elementary point $p \in \Elem(\cM)$ then $\Omega$ is necessarily in a final situation.  
\end{Remark}
\subsection{The Local Invariant}\label{subsect-localinvar}
Let us now introduce the main invariant used in the local strategy of the resolution of singularities.  First of all, we prove the following result:
\begin{Lemma}
Suppose that $p$ belongs to $\NElem(\cM)$. Then the main derived polygon $\vN^\prime$ in nonempty.
\end{Lemma}
\begin{proof}
According to Proposition~\ref{prop-derivedpolexists}, it suffices to prove that $m_3 \ge 1$.  But this is a direct consequence of the fact that the Newton data $\Omega$ is not in a final situation.
\end{proof}
Let us suppose that $p$ belongs to $\NElem(\cM)$. Writing the main vertex as $\bfm = (m_1,m_2,m_3)$ and the vertical displacement vector as $\Delta = (\Delta_1,\Delta_2)$, 
we define the {\em virtual height} associated to $\Omega$ as the natural number
$$
\vheight :=
\left\{
\begin{matrix}
\lfloor m_3  + 1 - \frac{1}{\Delta_2} \rfloor , &\mbox{ if }m_2 = -1
\mbox{ and }\Delta_1 = 0 \hfill\cr \cr m_3, &\mbox{ if }m_2 =
0\mbox{ or }\Delta_1 > 0\hfill
\end{matrix}
\right.
$$
where $\lfloor \alpha \rfloor := \max \{ c \in \cZ \mid c \le \alpha \}$.

\begin{figure}[htb]
\psfrag{F}{\small $\mathcal{F}$} 
\psfrag{xj/y}{\small $v_2$}
\psfrag{y}{\small $v_2$} 
\psfrag{xi}{\small $v_1$} 
\psfrag{z}{\small $v_3$} 
\psfrag{Np}{\small $\vN^\prime$} 
\psfrag{f}{\small $f$}
\psfrag{n}{\small $\bfm^\prime$} 
\psfrag{(0,0,m3)}{\small $(0,0,m_3)$} 
\psfrag{(0,-1,m3)}{\small $(0,-1,m_3)$}
\psfrag{(0,0,M)}{\small $(0,0,t)$} 
\psfrag{mu1}{\small $\Inv = (m_3,0,m_3,\ldots)$} 
\psfrag{mu2}{\small $\Inv = (\linteg t + 1 \rinteg, -1,m_3,\ldots)$}
\begin{center}
\includegraphics[height=5cm]{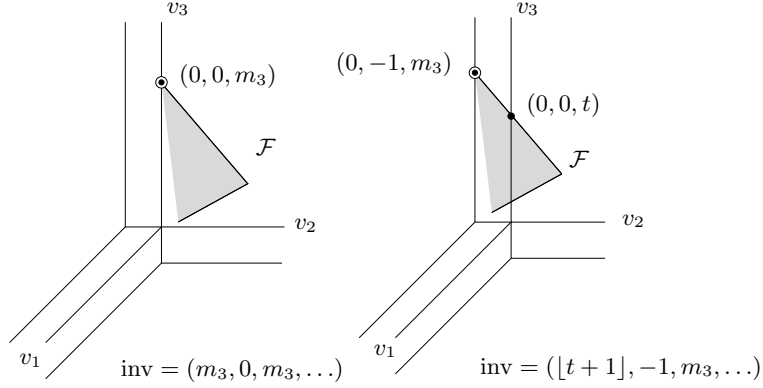}
\end{center}
\caption{The Invariant $\Inv(\Omega)$.}
\label{fig-invariantdraw}
\end{figure}

For $m_2 = -1$ and $\Delta_1 = 0$, the virtual height $\vheight$ is the smallest integer which is strictly
greater than the height of the point of intersection between
the main edge and the vertical plane $\{\bfv = (v_1,v_2) \in \cR^2  \mid v_2 = 0\}$ (as shown
in figure~\ref{fig-invariantdraw}).  

We refer the reader to subsection~\ref{subsect-xdirectblup} for an example which motivates the use of the notion of virtual height.
\begin{Definition}\label{def-multiplicity}
The {\em primary invariant} is the vector
$$\Invp := (\ \vheight,m_2 + 1,m_3\ )$$
The {\em secondary invariant} is the vector
$$
\Invs = (\ \#\iota_p - 1, \ \lambda\ \Delta_1, \ \lambda\ \max\{ 0, \ \Delta_2\}\ )
$$
where $\lambda := 2(m_3 + 1)!$.
The {\em invariant} associated to the Newton data $\Omega$
is the pair
$$\Inv(\Omega) := (\Invp, \Invs)$$
\end{Definition}
\begin{Remark}
It follows from the assumption $\#\iota_p \ge 1$, the choice of $\lambda$ and
the Proposition~\ref{prop-grid} that the vector $\Inv(\Omega)$ 
always belongs to $\cN^6$.
\end{Remark}

\subsection{Regular-Nilpotent Transitions}
Let us introduce the following notation.  Given a subset $A \subset \cZ^3$, let 
$\Omega|_A = ((x,y,z),\iota_p,\Theta|_A)$ be the Newton data which is obtained from $\Omega$ by considering the {\em restricted} Newton map  
$$\Theta|_A : \cZ^3  \rightarrow \cR^3$$
defined as follows: $\Theta|_A(\bfv) = \Theta(\bfv)$, if $\bfv \in A$; $\Theta|_A(\bfv) = 0$, if $\bfv \in \cZ^3\setminus A$.

If $\Omega$ is
associated to a vector field $\chi$ (which is a local generator of the line field at $p$), we denote by $\chi|_A$ the vector field associated to $\Omega|_A$.  Notice that $\chi|_A$ is possibly a degenerate vector field.
\begin{Lemma}\label{lemma-detectnilptrans}
The Newton data $\Omega$ is edge stable if and only if for all
$(f,0) \in \gG_{\Delta}$, there exists a 
constant $\widetilde{C} \in \overline{\cQ}_{\ge 0}$ such that
$$
(f,0)\cdot\Omega \in \Nmaps^{i,\bfm}_{\Delta,\widetilde{C}}.
$$
\end{Lemma}
\begin{proof}
It suffices to notice that each map $(f,g) \in \gG^i_{\Delta,C}$
can be uniquely written as a composition
$$
(f,g) = (f_0,0) \circ (\widetilde{f},\widetilde{g}) 
$$
where $(f_0,0)$ belongs to $\gG_{\Delta}$ and $(\widetilde{f},\widetilde{g})$ is a map belonging to the 
normal subgroup
$\gG^{i,+}_{\Delta,C} \lhd \gG^i_{\Delta,C}$.

Now, Corollary~\ref{corollary-howittransforms} implies that, if we denote by $\medge$ the main edge associated to $\Omega$, 
the Newton data $\widetilde{\Omega} := (\widetilde{f},\widetilde{g}) \cdot
\Omega$ 
is such that
$$\widetilde{\Omega}|_\medge = \Omega|_\medge$$
Moreover, $\medge$ is also the main edge of $\widetilde{\Omega}$. This concludes the proof.  
\end{proof}
Recall that the Newton data can be either in a regular or in a nilpotent
configuration (see subsection~\ref{subsect-regnilpmain}). 
Let us say that $\Omega$ is in a {\em potentially nilpotent situation} if
it is in a regular configuration but 
$$\bfm = (0,-1,m_3) \quad \mbox{ and }\quad \Delta = (0,1)$$
for some $m_3 \ge 1$.  The next result describes some basic aspects of the action of a $\gG^i_{\Delta,C}$-map on $\Omega$.
\begin{Lemma}\label{lemma-regnilp}
Given $(f,g) \in \gG^i_{\Delta,C}$, the Newton data 
$\widetilde{\Omega}:= (f,g) \cdot \Omega$ belongs
to the class
$\Nmaps^{i,\widetilde{\bfm}}_{\widetilde{\Delta},\widetilde{C}}$, where
two cases can occur:
\begin{enumerate}[(i)]
\item If $\Omega$ is not in a potentially nilpotent situation then
\begin{equation}\label{preservesm}
\widetilde{\bfm} = \bfm \quad\mbox{ and }\quad (\widetilde{\Delta},\widetilde{C}) \ge_\lex (\Delta,C)
\end{equation}
In particular, if $\Omega$ is in a regular (resp.\ nilpotent) configuration
then $\widetilde{\Omega}$ is also in a regular (resp.\ nilpotent)
configuration.
\item If $\Omega$ in a potentially nilpotent situation then
\begin{enumerate}[({ii}.a)]
\item Either $\widetilde{\bfm} = \bfm$ and $(\widetilde{\Delta},\widetilde{C})
      \ge_\lex (\Delta,C)$,  or
\item $\widetilde{\Omega}$ is in a nilpotent configuration and
$$\widetilde{\bfm} = (0,0,m_3 - 1), \quad \widetilde{\Delta} >_\lex \Delta$$
\end{enumerate}
\end{enumerate}
\end{Lemma}
\begin{proof}
It suffices to use Corollary~\ref{corollary-howittransforms}.
\end{proof}
In the case where item {\em (ii.b)} holds, we shall say that the data $\Omega$ is in a
{\em hidden nilpotent configuration} and that the transformation $\Omega \rightarrow
\widetilde{\Omega}$ is a {\em regular-nilpotent} transition.
\begin{figure}[htb]
\psfrag{m}{\small $\bfm$}
\psfrag{mt}{\small $\widetilde{\bfm}$}
\psfrag{medge}{\small $\medge$}
\psfrag{medget}{\small $\widetilde{\medge}$}
\psfrag{f0}{\small $(f,g)\cdot$}
\begin{center}
\includegraphics[height=5cm]{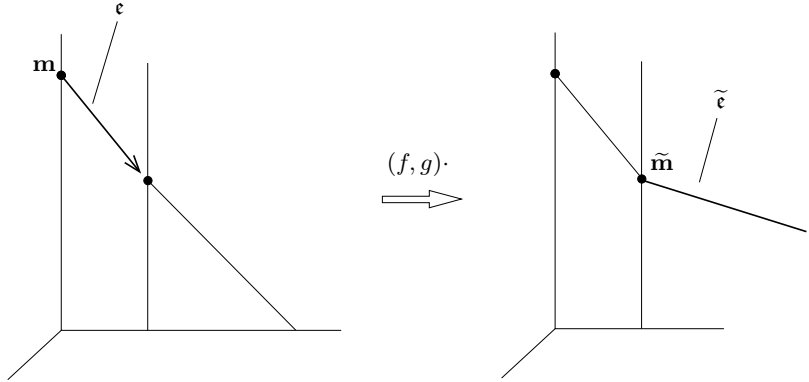}
\end{center}
\caption{The regular-nilpotent transition.}
\label{fig-regnilptrans}
\end{figure}

In view of Lemma~\ref{lemma-detectnilptrans}, a hidden nilpotent configuration may be detected just by the action of the subgroup $\gG_\Delta$.

\subsection{Resonant Configurations}\label{subsect-resonant}
This is a rather technical subsection whose main goal is to characterize those type of Newton data (called {\em resonant configurations}) for which the action of the group $\gG$ is not effective.  This characterization is essential to prove the uniqueness of the local strategy for the resolution of singularities at points $p \in \NElem(\cM) \cap \vD$.

We remark in passing that the occurrence of these resonant configurations have no analog in the theory of resolution of singularities for functions and analytic sets.

Suppose that $\Omega$ is a Newton data associated to an adapted local chart $(U,(x,y,z))$ with center at 
a divisor point
$p \in \vD \cap A$ such that $p \in \NElem(\cM)$.

Let us study how {\em effective} is the action of the group $\gG^i_{\Delta,C}$
on the {\em support} of $\Omega$.  Recall that such support is given by
$$
\supp(\Omega) := \{\bfv \in \cZ^3 \mid \Theta(\bfv) \ne 0\}
$$
where $\Theta:\cZ^3\rightarrow \cR^3$ is the Newton map associated to $\Omega$. To state our next result, we need  the following definition. 
Let $\medge$ be the main edge of 
$\Omega$.  We shall say that $\Omega$ is in {\em $c$-resonant
configuration} (for some $c \in \cN$)
if there exists a map $(f,g) \in \gG^i_{\Delta,c} \setminus \gG_\Delta$ such that
$$
\supp(\, (f,g)\cdot \Omega|_\medge ) \subset \medge
$$
In other words, there is a map 
$(f,g) \in \gG^i_{\Delta,c} \setminus \gG_\Delta$ (i.e.\ {\em not} of the form $(\xi x^{\Delta_1}y^{\Delta_2},0)$)
whose action on the restricted Newton data $\Omega|_\medge$ results into a Newton data which still has the support on
$\medge$.
\begin{Lemma}\label{Lemma-resonant}
Suppose that $\Omega \in \Nmaps^{i,\bfm}_{\Delta,C}$ is 
a Newton data which is in a $c$-resonant
configuration.  Then, $\Omega$ is not edge stable.  Moreover,  $i = 1$ and 
$\Delta = (0,s)$, for some $s >0$.
Considering the associated vector
field $\chi$, one of the following situations occurs:
\begin{enumerate}[(i)]
\item $\Delta = (0,1)$ and the restriction of $\chi$ to the main edge is given by
$$\chi|_\medge = (z + \lambda y)^m \left[
\alpha \left( x \frac{\partial}{\partial x} +  c y\frac{\partial}{\partial y}
+ c z\frac{\partial}{\partial z} \right) + \beta (z + \lambda
  y)\frac{\partial}{\partial y} + \gamma (z + \lambda
  y)\frac{\partial}{\partial z} \right]
$$
for some $m \ge 1$, $\lambda \in \cR$ and $(\alpha,\beta ,\gamma) \in \cR^3$
such that $\beta \ne 0$ and $(\alpha,\gamma + \lambda\beta) \ne (0,0)$.
\item $\Delta = (0,1/\tau)$ for some $\tau \in  \cN_{\ge 2}$, and
$$
\chi|_\medge = z^{\tau m} \left[
\alpha \left( x \frac{\partial}{\partial x} +  c  y\frac{\partial}{\partial y}
+ \frac{c}{\tau} z\frac{\partial}{\partial z} \right) +
\beta z^\tau\frac{\partial}{\partial y} + \gamma z \frac{\partial}{\partial z}
\right]
$$
for some $m \ge 1$ and $(\alpha,\beta,\gamma) \in \cR^3$ such that $\beta \ne
0$ and $(\alpha,\gamma) \ne (0,0)$.
\end{enumerate}
\end{Lemma}
\begin{proof}
If $\Delta = (\Delta_1,\Delta_2)$ for some $\Delta_1 > 0$ then 
$\gG^i_{\Delta,c} = \gG_\Delta$ by definition and nothing has to be proved.

Let us assume that $\Omega$ is edge stable.
Up to a $x$-directional blowing-up with weight-vector $\bfomega = k\cdot (1,c,sc)$ (where $k \in \cN$ is chosen 
in such a way that $\bfomega \in \cN^3$), we can write
$$
\chi|_\medge = F(y,z) x\frac{\partial}{\partial x} + G(y,z)
\frac{\partial}{\partial y} + H(y,z) \frac{\partial}{\partial z}
$$
where $F,G,H$ are $(1,s)$-quasihomogeneous functions of respective degree $M$, $M + 1$
and $M + s$, for some rational number $M \in \cQ$.

After such blowing-up, the map $(f,g)$ is transformed to an
element of the group $\gG^i_{(0,s),0}$.  Keeping the same notation for such
map, we can write
\begin{equation}\label{formfg}
(f,g) = (a_0 + a_1 y + \cdots + a_{k}y^{k}, \eta), \quad \mbox{ with }
a_0,\ldots,a_{k},\eta \in \cR
\end{equation}
where $k := \min\{n \in \cN \mid n \le s\}$ and $(f,g) \ne (a_s y^s,0)$.  
Our 
problem reduces to find conditions on $F,G,H$ such that there is  
one such map for which
\begin{equation}\label{holdssupport}
\supp(\, (f,g) \cdot \Omega|_\medge\, ) \subset \medge
\end{equation}
Let us consider the four possible cases:
\begin{enumerate}
\item $s \in \cN$, $s \ge 2$
\item $s = 1$.
\item $s  = 1/\tau$, $\tau \in \cN$, $\tau \ge 2$.
\item $s \notin \cN \cup 1/\cN$.
\end{enumerate}
In the first case, two possible expressions for $F,G,H$ can appear:
\begin{enumerate}[({1}.a)]
\item If the main vertex has the form
$\bfm = (0,-1,m)$ then 
\begin{equation}\label{generalFGH1}
F(y,z) = F_0(y,z), \; G(y,z) = G_0(y,z), \; H(y,z) = y^{s-1}H_0(y,z)
\end{equation}
where $F_0,G_0,H_0$ are $(1,s)$-quasihomogeneous functions of degree $ms - 1$, $ms$ and $ms$. Moreover, 
$G_0(0,z) = \beta z^m$, for some $\beta \ne 0$.
\item If main vertex has the form
$\bfm = (0,-1,m)$ then
\begin{equation}\label{generalFGH}
F(y,z) = F_1(y,z), \; G(y,z) = G_1(y,z) y, \; H(y,z) = H_1(y,z)z +
H_2(y,z) y^s
\end{equation}
where $F_1,G_1,H_1,H_2$ are $(1,s)$-quasihomogeneous functions of degree $ms$ and
$\mu \in \cR$.
\end{enumerate}
Assuming that condition (\ref{holdssupport}) holds, 
it is easy to see that, in the expression (\ref{generalFGH1}), 
$F_0  \equiv 0$, and $G_0, H_0$ should be a power of a common $(1,s)$-quasihomogeneous form of degree $s$,
$$
G_0(y,z) = \beta (z + \lambda y^s)^m, \quad H_0(y,z) = \gamma (z + \lambda y^s)^m
$$
for some $\beta,\lambda \in \cR^*$ and $\gamma \in \cR$. Looking only to the
function $G_0$, we see that the only 
possible map $(f,g)$ satisfying our requirements is given by
$$
f = -\lambda((y + \eta)^s - y^s),\quad g = \eta, \quad\mbox{ for some }\eta \ne 0
$$ 
In this case, the restricted vector field $\chi|_\medge$ is transformed to (dropping the tildes)
$$
 (z + \lambda y^s)^m \left( \beta \frac{\partial}{\partial y}  + 
\left( -s\lambda\beta y^s + (\gamma +s\lambda\beta) (y - \eta)^{s-1} \right) \frac{\partial}{\partial z}\right)
$$
Therefore, since $s \ge 2$ and (\ref{holdssupport}) is assumed, 
the relation $\gamma +s\lambda\beta = 0$ necessarily holds.  But if we consider the original expression 
of the vector field $\chi|_\medge$ and apply the $\gG_\Delta$-map $(f^\prime,g^\prime)= (\lambda y^s, 0)$, we get
(dropping again the tildes)
$$
z^m \left( \beta \frac{\partial}{\partial y}  + 
\left( \gamma +s\lambda\beta) y^{s-1} \right) \frac{\partial}{\partial z}\right) =  
z^m \left( \beta \frac{\partial}{\partial y} \right)  
$$
(i.e. the support of $\Omega|_\medge$ has a single point). This implies that $\Omega$ is not edge stable. Contradiction.

Similarly for the expression (\ref{generalFGH}), a simple computation shows that if there is a nonzero $(f,g)$ which fixes
the support of $\Omega|_\medge$ then $F_1,G_1,H$ should necessarily be a power of
$(1,s)$-quasihomogeneous form of degree $s$, namely:
$$
F_1(y,z) = \alpha (z + \lambda y^s)^m, \quad G_1(y,z) = \beta (z + \lambda
y^s)^m, \quad H(y,z) = \gamma
(z + \lambda y^s)^{m+1}
$$
for some $\lambda \in \cR$.  Moreover, since $(f,g) \ne (a_s y^s, 0)$, we
conclude from the general expression~(\ref{generalFGH}) that
$G_1 = 0$.  Therefore, $\chi|_\medge$ has the general form
$$
\chi|_\medge = (z + \lambda y^s)^m \left[
\alpha x \frac{\partial}{\partial x} +  \gamma (z + \lambda
  y^s)\frac{\partial}{\partial z} \right]
$$
Notice, however, that this expression implies that $\Omega$ is not edge stable.
In fact, applying the map $(f^\prime,g^\prime)
= (\lambda y^s, 0)$ we get a new Newton data $(f^\prime,g^\prime)\cdot \Omega \in
\Nmaps^{i,\bfm}_{\Delta^\prime,C^\prime}$ with $\Delta^\prime >_\lex \Delta$.  Contradiction.

If we suppose that $s = 1$ then the same reasoning used above leads us to the
general expression
$$
\chi|_\medge = (z + \lambda y)^m \left[
\alpha x \frac{\partial}{\partial x} +  \beta (z + \lambda y)
\frac{\partial}{\partial y} + \gamma (z + \lambda
  y)\frac{\partial}{\partial z} \right]
$$
with $\lambda \in \cR$ and $(\alpha,\beta,\gamma) \in \cR^3 \setminus
\{(0,0,0)\}$. If we apply the map $(f_1,g_1)
= (\lambda y, 0)$ we see that if $\beta = 0$ or $(\alpha,\gamma  +
\lambda\beta = 0)$ then $\Omega$ is not edge stable.  If this is not the case, we are
precisely (up to blowing-up) in the configuration
listed in item (i) of the enunciate.

Suppose now that $s = 1/\tau$, with $\tau \in \cN_{\ge 2}$.  Then, $(f,g) =
(\eta,\xi)$ for some $\xi,\eta \in \cR$ and we obtain
the general expression
$$
\chi|_\medge = z^{\tau m} \left[
\alpha x \frac{\partial}{\partial x} +  \beta z^\tau
\frac{\partial}{\partial y} + \gamma z \frac{\partial}{\partial z} \right]
$$
with $m \ge 1$ and $(\alpha,\beta,\gamma) \in \cR^3$.  Since $\supp(\Omega) \cap \medge$
has at least two points, we see that $\beta \ne
0$ and $(\alpha,\gamma) \ne (0,0)$.  This is
precisely (up to blowing-up) the configuration listed in item (ii) of the enunciate.

It remains to study the case where $s \notin \cN \cup 1/\cN$.  Here, under the
assumption (\ref{holdssupport}), we get
$$
\chi|_\medge = z^{m/s} \left[
\alpha x \frac{\partial}{\partial x} +  \beta y
\frac{\partial}{\partial y} + \gamma z \frac{\partial}{\partial z} \right]
$$
for some $(\alpha,\beta,\gamma) \in \cR^3\setminus\{(0,0,0)\}$.  However, such
expression implies that $\supp(\Omega) \cap \medge$ has a single point.  
Absurd.
\end{proof}
Notice that the configurations (i) and (ii) of the previous Lemma
do not represent edge stable Newton data.  

Indeed, the item (ii) of the Lemma is obviously excluded
because it represents a nilpotent configuration with higher vertex $\bfh =
(0,-1,\tau (m+1))$ and associated edge $\medge(\bfh)$ given by
$$\medge(\bfh) = \overline{\bfh,\bfn}$$
where $\bfn = (0,0,\tau m)$.  In this case, it follows immediately from the definition of the main vertex that the main
vertex associated to $\Omega$ is $\bfn$ and not $\bfh$. The same reasoning can
be used to exclude the item (i) of the Lemma with
$\lambda = 0$.

The item (i) of the Lemma with $\lambda \ne 0$ is also
excluded because it is not edge stable.  Indeed, the coordinate change $\widetilde{z}
= z + \lambda y$ causes a regular-nilpotent transition 
(see Lemma~\ref{lemma-regnilp}).

The following result is an
immediate consequence of Lemma~\ref{Lemma-resonant}.
\begin{Proposition}\label{Proposition-stablenotresonant}
Let $\Omega \in \Nmaps^{i,\bfm}_{\Delta,C}$ be 
an edge stable Newton data. Consider a map $(f,g)
\in \gG^i_{\Delta,c} \setminus \gG_\Delta$ with $c < C$. 
Then, $(f,g)\cdot\Omega$ necessarily belongs to $\Nmaps^{i,\bfm}_{\Delta,c}$.
\end{Proposition}

\begin{figure}[htb]
\psfrag{m}{\small $\bfm$}
\psfrag{(f,g)}{\small $(f,g)\cdot$}
\psfrag{medge}{\small $\medge$}
\psfrag{F}{\small $\vF$}
\psfrag{Ft}{\small $\widetilde{\vF}$}
\begin{center}
\includegraphics[height=5cm]{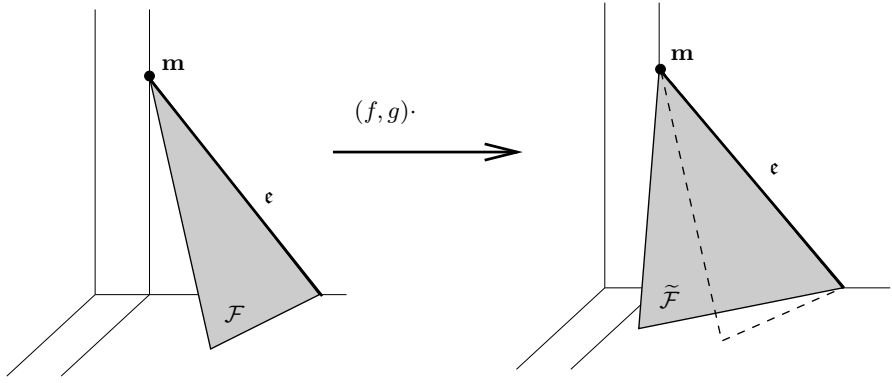}
\end{center}
\caption{The action of a map $(f,g) \in \gG^i_{\Delta,c} \setminus \gG_\Delta$ with $c < C$.}
\label{fig-notresonant}
\end{figure}

In other words, the map $(f,g) \in \gG^i_{\Delta,c} \setminus \gG_\Delta$ with $c < C$ acts
{\em effectively} on $\Nmaps^{i,\bfm}_{\Delta,C}$.

Using the same computations made in the proof of Lemma~\ref{Lemma-resonant}, 
we can immediately prove the following result, which gives a more precise description of the 
action of $\gG^i_{\Delta,c}\setminus \gG_\Delta$ on the main edge of $\Omega$.
\begin{Lemma}\label{Lemma-stablenotresonant}
Suppose that $\Omega \in \Nmaps^{i,\bfm}_{\Delta,C}$ is 
an edge stable Newton data.  Then, for each
$(f,g) \in  \gG^i_{\Delta,c} \setminus \gG_\Delta$, the Newton data 
$$
\widetilde{\Omega} = (f,g) \cdot \left(\, \Omega|_\medge \, \right)$$
(i.e.\ the action of $(f,g)$ in $\Omega|_\medge$)
is such that the following conditions holds:
\begin{enumerate}
\item For each pair $(k,l) \in \supp(f)$, we have  $\supp(\widetilde{\Omega})\; \cap \; \left(\medge + (k,l,
      -1) \right) \ne \emptyset$.
\item For each  $c \in \supp(g)$, we have  $\supp(\widetilde{\Omega}) \; \cap \; 
\left( \medge + (c,-1,0)\right) \ne \emptyset$.
\end{enumerate}
\end{Lemma}
\subsection{Basic Edge Preparation and Basic Face Preparation}\label{subsect-basicfaceedgeprep}
To enunciate the following Lemma, we consider the set 
$$
\Nmaps^{i,\bfm}_{\Delta} := {\displaystyle \bigcup_{C
}{\Nmaps^{i,\bfm}_{\Delta,C}}}. 
$$
of all classes of Newton data with fixed values for $(i,\bfm,\Delta)$.
\begin{Lemma}\label{lemma-uniqueedgemap}
Suppose that the Newton data $\Omega \in \Nmaps^{i,\bfm}_{\Delta,C}$ is centered at a 
point $p \in \NElem(\cM)$. Then, if
$\Omega$ is not edge stable, there exists a unique map
$(f,0) \in \gG_\Delta$ such that
\begin{equation}\label{herefg}
(f,0) \cdot \Omega \not\in \Nmaps^{i,\bfm}_{\Delta}.
\end{equation}
\end{Lemma}
\begin{proof}
Choose an arbitrary map $(f,0) \in \gG_\Delta$ which satisfies
(\ref{herefg}) and define
$\widetilde{\Omega} := (f,0)\cdot \Omega$ (there exists at least one such map by Lemma~\ref{lemma-detectnilptrans}).  Let $\widetilde{\bfm}$ and
$\widetilde{\Delta}$ be the main vertex and the vertical displacement vector associated
to $\widetilde{\Omega}$.

Suppose, first of all, that $\widetilde{\bfm} \ne \bfm$. Then, it follows from
Lemma~\ref{lemma-regnilp} that $\Omega$ is in a hidden nilpotent configuration
and that
$\Omega \rightarrow \widetilde{\Omega}$ is a
regular-nilpotent transition.

In these conditions, we know that $\bfm = (0,-1,m_3)$ and $\Delta = (0,1)$.
Therefore, $m_3 \ge 2$ (because otherwise $\Omega$ would be associated
to an elementary point $p$, by Proposition~\ref{prop-finalsituation}) and moreover $f$ has the
particular form
$$
f = \xi y,
$$
for some constant $\xi \in \cR$.   

Let us suppose that there exists another map
$(f^\prime,0) = (\xi^\prime y, 0)$ such that $\widetilde{\Omega}^\prime
= (f^\prime,0) \cdot \Omega$ does not belong to 
$\Nmaps^{i,\bfm}_{\Delta}$. Then, the composition $(\widetilde{f},0) =
(f^\prime,0) \circ (f,0)^{-1}$ is such that $\widetilde{f}$ is given by 
$\widetilde{f} = (\xi^\prime -
\xi) y$ and it maps
$\widetilde{\Omega}$ to $\widetilde{\Omega}^\prime$.  

We claim that $\widetilde{f} = 0$. Suppose the contrary (i.e.\ $\xi \ne
\xi^\prime$). If $\medge$ is the main edge associated to $\Omega$ then
$$
\chi|_\medge = F(y,z)x \frac{\partial}{\partial x} + G(y,z)
\frac{\partial}{\partial y} + H(y,z) \frac{\partial}{\partial z}
$$
where $F,G,H$ are homogeneous polynomials of degree $m_3-1$, $m_3$ and $m_3$, respectively.  The hypothesis
that $\bfm = (0,-1,m_3)$ 
implies that $G(y,z) = \rho z^{m_3} + \cdots$ for some nonzero constant $\rho \in \cR$.

If we apply the change of coordinates $\widetilde{z} = z + \xi y$ to $\chi|_\medge$, we get a
vector field $\widetilde{\chi}|_\medge = \widetilde{F}x \frac{\partial}{\partial x} + \widetilde{G}
\frac{\partial}{\partial y} + \widetilde{H} \frac{\partial}{\partial
\widetilde{z}}$ with
$$
\widetilde{F} = F, \quad \widetilde{G} = G \quad \mbox{ and }\quad \widetilde{H} = H + \xi G.
$$
The assumption that $\widetilde{\Omega}$ is in a nilpotent configuration implies
that $\widetilde{z} = 0$ is a root of multiplicity $\ge m_3 - 1$ of
$\widetilde{G}(1,\widetilde{z})$.
This is equivalent to say that $z = -\xi$ is a root of multiplicity $\ge m_3 - 1$
of $G(1,z)$.
Let us split the proof in two cases:
\begin{enumerate}[(a)]
\item $\widetilde{\Omega}^\prime$ is in a regular configuration;
\item $\widetilde{\Omega}^\prime$ is in a nilpotent configuration;
\end{enumerate}
In the case (a), the same computations made in the previous paragraph
imply that the polynomial $G(1,z)$ should have $z = -\xi^\prime$ as a root of multiplicity
$\ge m_3$.  This is absurd since $m_3 \ge 2$ and therefore $m_3 + (m_3 - 1)
> m_3$ .

In the case (b), we conclude that $z = -\xi^\prime$ should also be a root of
$G(1,z)$ of multiplicity $m_3 - 1$.  This implies that $2(m_3 - 1) \le m_3$,
i.e.\ $m_3 \le 2$.

Since we assume that $m_3 \ge 2$, it remains to treat the case (b) with $m_3 =
2$.  Here, $\widetilde{\chi}|_\medge$ is necessarily given
(dropping the tildes) by
$$
\widetilde{\chi}|_\medge = \rho z (z + \beta y)\frac{\partial}{\partial y} +
z \left( \alpha x \frac{\partial}{\partial x} + \gamma z\frac{\partial}{\partial z}\right)
$$
for some $(\alpha,\beta,\gamma) \in \cR^3 \setminus \{(0,0,0)\}$ and $\rho \ne
0$.
If we apply the coordinate change $z^\prime = z + \eta y$ (where $\eta := \xi -
\xi^\prime$), we get (dropping the primes)
$$
\rho (z - \eta y)(z + (\beta - \eta) y) \left( \frac{\partial}{\partial y} +
\eta \frac{\partial}{\partial z} \right) +
(z - \eta y) \left( \alpha x \frac{\partial}{\partial x} + \gamma (z - \eta
y)\frac{\partial}{\partial z} \right)
$$
We now use the assumption that $\widetilde{\Omega}^\prime$ is in a nilpotent
configuration.  Looking at the coefficients of $\partial/\partial x$ and
$\partial/\partial y$, this implies that $\alpha = 0$, $\eta = \beta$.
Therefore, the coefficient of $\partial/\partial z$ has the form
$$
(z - \eta y) \left(\rho \eta z + \gamma (z - \eta y) \right)
$$
and, since this expression should be equal to $\gamma^\prime z^2$ (for some
nonzero real constant $\gamma^\prime$), we conclude that necessarily $\eta = 0$. This proves
the claim.

We now prove the Lemma in the simpler case where $\widetilde{\bfm} =
\bfm$. Here, the vertical displacement vector $\widetilde{\Delta}$ is such that
$$
\widetilde{\Delta}>_\lex \Delta.
$$
Suppose that there exists another map $(f^\prime, 0)
\in \gG_\Delta$
such that the Newton data $\widetilde{\Omega}^\prime := (f^\prime,0)\cdot
\Omega$ also has a displacement vector $\widetilde{\Delta}^\prime >_\lex
\Delta$.
We claim that the map
$$(\widetilde{f},0) := (f^\prime,0) \circ (f,0)^{-1} \in
\gG_\Delta$$
which sends $\widetilde{\Omega}$ to $\widetilde{\Omega}^\prime$ is necessarily
the identity.   In fact, if $\widetilde{\chi}$ denotes the
vector field which is associated to $\widetilde{\Omega}$ then
$$\widetilde{\chi}|_\medge = y^{m_2}z^{m_3} \left( \alpha x \frac{\partial}{\partial x} +
\beta y
\frac{\partial}{\partial y} + \gamma z \frac{\partial}{\partial z}\right)$$
for some $(\alpha,\beta,\gamma) \in \cR^3 \setminus
\{(0,0,0)\}$ such that $\alpha = \gamma = 0$ if $m_2 = -1$.

If we write $\widetilde{f} = \xi x^{\Delta_1}y^{\Delta_2}$ then
the map $(\widetilde{f},0)$ transforms $\widetilde{\chi}|_\medge$ to
\begin{equation}\label{mapwidetildef}
y^{m_2}(z - \xi x^{\Delta_1}y^{\Delta_2})^{m_3} \left( \alpha x
\frac{\partial}{\partial x} +
\beta y
\frac{\partial}{\partial y} + (\gamma z - \xi (\gamma - \alpha 
\Delta_1 - \beta \Delta_2)
x^{\Delta_1} y^{\Delta_2})
\frac{\partial}{\partial z}
\right)
\end{equation}
Since $m_3 \ge 1$, the assumption $\widetilde{\Delta}^\prime >_\lex \Delta$
necessarily implies that $\xi = 0$.
\end{proof}
The map $(f,0) \in \gG_{\Delta}$
which is defined by Proposition~\ref{lemma-uniqueedgemap} will be called
{\em basic edge
preparation map} associated to $\Omega$.
\begin{Remark}\label{remark-edgestablefuture}
The expression obtained in (\ref{mapwidetildef}) has the following simple
consequence, which we will need in subsection~\ref{subsect-effectramif}.
Suppose the Newton data $\Omega$ is such that $\Delta_1 > 0$ and 
$$\supp(\Omega|_\medge)\cap \{ \bfv \in \cZ^3 \mid v_3 = m_3 - 1\} = \emptyset.$$
Then, $\Omega$ is edge stable.  Indeed, for all map $(f,0) \in \gG_\Delta$ we know
that $(f,0)\cdot \Omega$ has the same main vertex of $\Omega$ (because no regular-nilpotent transition can occur since $\Delta_1 > 0$).  Moreover, the expression 
(\ref{mapwidetildef}) implies that $\supp((f,0)\cdot \Omega|_\medge)$ contains at least two points.  Therefore, $(f,0)\cdot\Omega$ has the same main edge as $\Omega$.
\end{Remark}
\begin{Lemma}\label{GDeltapreserves}
Suppose that $\Omega \in  \Nmaps^{i,\bfm}_{\Delta,C}$ is an edge stable Newton data.  Then
$$
(f,0) \cdot \Omega
$$
also belongs to $\Nmaps^{i,\bfm}_{\Delta,C}$, for all map $(f,0) \in \gG_\Delta$.
\end{Lemma}
\begin{proof}
If $C = \infty$ or $C = 0$ then nothing has to be proven.  If $0 < C < \infty$
then the main face $\vF$ of the Newton polyhedron $\vN$ associated to $\Omega$
contains at least one vertex $\bfv \in \supp(\Omega)$ which is not in the main
edge $\medge$ and such that
\begin{equation}\label{condbfv}
(\bfv - \Delta) \cap \vN = \emptyset
\end{equation}
Indeed, choose some arbitrary vertex $\bfv^\prime \in \supp(\Omega) \cap \vF \setminus
\medge$ (there exists at least one such vertex because $0 < C < \infty$).  If
such vertex satisfies (\ref{condbfv}) then choose $\bfv = \bfv^\prime$.
Otherwise, there necessarily
exists some $\eps > 0$ such that
segment $\{\bfv - t\Delta \mid t \in [0,-\eps]\}$ is an edge of $\vF$.  Then,
it suffices to choose $\bfv$ to be the other extreme of such edge, i.e.
$\bfv := \bfv^\prime - \eps \Delta$.

Now, if we choose the vertex $\bfv \in \supp(\Omega)$ as above, it is clear
that $\widetilde{\Omega}(\bfv) = \Omega(\bfv)$.  Therefore, since $\bfv$ and
$\medge$ are affinely independent and $\Omega$ is edge stable, we conclude
that $\widetilde{\Omega} \in \Nmaps^{i,\bfm}_{\Delta,C}$.
\end{proof}
In the next Lemma, we consider the action of the subgroup
of edge preserving maps
$\gG^{i,+}_{\Delta,C}$ (see subsection~\ref{subsect-subgroupsgG}).
\begin{Lemma}\label{lemma-uniquestrictmap}
Suppose that $\Omega \in \Nmaps^{i,\bfm}_{\Delta,C}$ is an edge stable Newton data which is not stable.  Then,
there exists a unique edge preserving map $(f,g) \in \gG^{i,+}_{\Delta,C}$ such that
$$
(f,g) \cdot \Omega \mbox{ belongs to }\Nmaps^{i,\bfm}_{\Delta,\widetilde{C}}
$$
for some $\widetilde{C} > C$.
\end{Lemma}
\begin{proof}
To prove the existence part, let $(f,g) \in \gG^i_{\Delta,C}$ be such that
$
(f,g) \cdot \Omega \in \Nmaps^{i,\bfm}_{\Delta,\widetilde{C}}.
$
Then, we can uniquely decompose $(f,g)$ as $(f_2,0) \circ (f_1,g_1)$, where
$(f_2,0)$ belongs to $\gG_\Delta$ and $(f_1,g_1) \in \gG^{i,+}_{\Delta,C}$ is an edge preserving map.

We claim that $\Omega_1 := (f_1,g_1) \cdot \Omega$ belongs to
$\Nmaps^{i,\bfm}_{\Delta,\widetilde{C}}$, for some $\widetilde{C} > C$.  Indeed,
suppose, by absurd, that this is not the case.  Then, $\Omega_1$
is an edge stable Newton data in $\Nmaps^{i,\bfm}_{\Delta,C}$.  Using the
Lemma~\ref{GDeltapreserves}, we conclude that
$(f_2,0) \cdot \Omega_1$ also belongs to $\Nmaps^{i,\bfm}_{\Delta,C}$.  Contradiction.

To prove the uniqueness of $(f,g) \in \gG^{i,+}_{\Delta,C}$, consider
two maps $(f_1,g_1)$ and $(f_2,g_2)$ in $\gG^{i,+}_{\Delta,C}$ such that
$$
\Omega_j := (f_j,g_j)\cdot \Omega \in \Nmaps^{i,\bfm}_{\Delta,C_j}
$$
for some $C_j > C$, $j = 1,2$.  Then, if we define the composed map
$$(\widetilde{f},\widetilde{g}) := (f_2,g_2)\circ(f_1,g_1)^{-1} \in \gG^{i,+}_{\Delta,C}$$
we get $\Omega_2 = (\widetilde{f},\widetilde{g}) \cdot \Omega_1$.  Using
Proposition~\ref{Proposition-stablenotresonant}, we conclude that
$(\widetilde{f},\widetilde{g}) = (0,0)$.
\end{proof}
Given an edge stable Newton data $\Omega$, the map 
$(f,g) \in \gG^{i,+}_{\Delta,C}$
which is defined by Lemma~\ref{lemma-uniquestrictmap} will be called
{\em basic face
preparation map} associated to $\Omega$.
\subsection{Formal Adapted Charts and Invariance of $(\bfm,\Delta,C)$}
For a fixed adapted local chart $(U,(x,y,z))$ at a divisor point $p \in \vD \cap A$, and a choice of local generator for the line field $\linef$, there is an associated
Newton data $\Omega$.  In section~\ref{sect-NewtonMaps}, we have seen how to associate certain 
{\em combinatorial quantities} $(\bfm,\Delta,C)$ to such Newton data.  

A natural question is how such combinatorial quantities depend on the choice of adapted local chart.  Our present
goal is to answer this question.
 
First of all, we need to slightly extend the concept of adapted local chart 
(see subsection~\ref{subsect-localadaptedchart}).

A {\em formal adapted chart} at a divisor point $p \in \vD \cap A$ is a triple $(x,y,z)$ formed by 
elements of the formal completion $\widehat{\vO}_p \supset \vO_p$ (with respect to the Krull topology) such that:
\begin{itemize}
\item The formal functions $x,y,z$ are independent at $p$ (i.e.\ their residue class generate 
$\widehat{\cm}_p/\widehat{\cm}_p^2$).
\item $\vZ$ is locally generated by $\frac{\partial}{\partial z}$;
\item If $\iota_p = [i]$ then $D_i = \{ x = 0 \}$;
\item If $\iota_p = [i,j]$ (with $i > j$) then $D_i = \{ x = 0\}$ and $D_j = \{y = 0\}$.
\end{itemize}
It is immediate to see that the construction of subsection~\ref{subsection-Newtonmap} can 
be carried out in the present setting.  Thus, up to a choice of a local generator $\chi$ for the line field at $p$, there exists a well-defined {\em formal Newton map} 
$$\Theta:\cZ^3 \rightarrow \cR^3$$ 
for $(\cM,\Ax)$ at $p$, relatively to $(x,y,z)$. We call the triple $((x,y,z),\iota_p,\Theta)$ a {\em formal Newton data} for $(\cM,\Ax)$ at $p$.  We define the classes $\Nmaps^{i,\bfm}_{\Delta,C}$ exactly as previously.

Given two formal adapted charts $(x,y,z)$ and $(\widetilde{x},\widetilde{y},\widetilde{z})$ at $p$, 
the transition map is given by
\begin{equation}\label{formalGmaps}
\widetilde{x} = x u(x,y),\quad \widetilde{y} = g(x) + y v(x,y), \quad
\widetilde{z} = f(x,y) + z w(x,y,z)
\end{equation}
where $g \in \cR[[x]]$, $u,v,f \in \cR[[x,y]]$ and $w \in \cR[[x,y,z]]$ are
such that $g(0) = f(0) = 0$, $u,v,w$ are units and $g = 0$ if $\#\iota_p = 2$.  The group of all such 
changes of coordinates forms a group, which we denote by $\ffG$.

\begin{figure}[htb]
\psfrag{Delta}{\small $(\Delta_1,\Delta_2)$}
\psfrag{Sf}{$S_f$}
\psfrag{Sg}{$S_g$}
\psfrag{C}{\small $C$}
\psfrag{(C,-1)}{\small $(C,-1)$}
\begin{center}
\includegraphics[height=3.5cm]{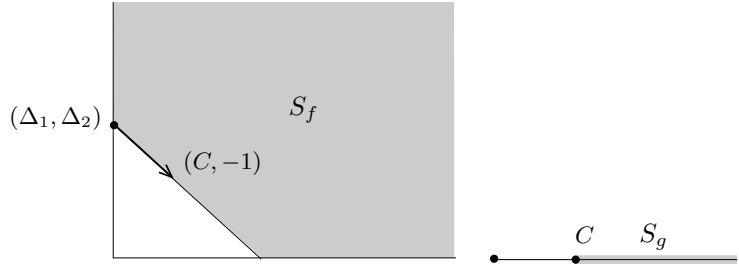}
\end{center}
\caption{The supports of $S_f$ and $S_g$ in the definition of $\ffG^i_{\Delta,C}$.}
\label{fig-supportSfSg}
\end{figure}

An element of $\ffG$ will be shortly denoted by $(f,g,u,v,w)$.  We consider also the subgroups 
$$
\ffG^1 = \ffG \quad \mbox{ and }\quad \ffG^2 = \{(f,g,u,v,w) \in \ffG \mid g = 0 \}.
$$
\begin{Remark}
The Lie Algebra associated to the group $\ffG$ is formed by all formal vector fields having the form
$$
xu(x,y)\frac{\partial}{\partial x} + (g(x) + y v(x,y))\frac{\partial}{\partial y} + 
(f(x,y) + z w(x,y,z))\frac{\partial}{\partial z}
$$
where $u,v,w$ are units and $g(0) = f(0) = 0$.
\end{Remark}
We denote by $\ffG^i_{\Delta,C}$ the subgroup of all maps $(f,g,u,v,w) \in \ffG^i$ such that
the supports of the maps $f$ and $g$ satisfy the following conditions
(see figure~\ref{fig-supportSfSg})
\begin{eqnarray*}
&& S_f \subset \{(a,b) \in \cN^2 \mid \< (1,C), (a,b) - \Delta  \> \,\ge\, 0,\ (a,b)
\ge_\lex \Delta\} \quad \mbox{ and } \cr\cr
&& S_g \subset \{c \in \cN \mid c \ge C\}
\end{eqnarray*}
If $C = \infty$, the former condition is replaced by $S_f \subset \{(a,b) \in \cN^2 \mid b \ge \Delta_2\}$.
We shall say that $\ffG^i_{\Delta,C}$ is the group of {\em $(\Delta,C)$-face maps}.

The following Lemma relates the groups $\ffG^i_{\Delta,C}$ and $\gG^i_{\Delta,C}$.
\begin{Lemma}\label{lemma-normalsubgroupffG}
There exists a normal subgroup $\ffG^{i,+}_{\Delta,C} \lhd \ffG^i_{\Delta,C}$ such that
the quotient $\ffG^i_{\Delta,C}/\ffG^{i,+}_{\Delta,C}$ is naturally isomorphic
to $\gG^i_{\Delta,C}$. We shall say that $\ffG^{i,+}_{\Delta,C}$ is the subgroup of {\em $(\Delta,C)$-face preserving maps}.
\end{Lemma}
\begin{proof}
We define explicitly the subgroup $\ffG^{i,+}_{\Delta,C}$ as follows:
\begin{enumerate}[(a)]
\item If $C \in \{0,\infty\}$ then $(f,g,u,v,w) \in \ffG^i_{\Delta,C}$ belongs to $\ffG^{i,+}_{\Delta,C}$ if and only if
$$
S_f \subset \{(a,b) \in \cN^2 \mid (a,b) >_\lex \Delta\} \quad \mbox{ and } \quad
S_g \subset \{c \in \cN \mid c > C\}
$$
\item If $0 < C < \infty$ then $(f,g,u,v,w) \in \ffG^i_{\Delta,C}$ belongs to $\ffG^{i,+}_{\Delta,C}$ if and only if
$$
S_f \subset \{(a,b) \in \cN^2 \mid \< (1,C), (a,b) - \Delta  \> \,> \, 0\} 
$$
\end{enumerate}
It is immediate to verify that this gives a normal subgroup of $\ffG^i_{\Delta,C}$.  Moreover, 
$$
\gG^i_{\Delta,C} \cap \ffG^{i,+}_{\Delta,C} = \{0\} \quad \mbox{ and }\quad \ffG^i_{\Delta,C} = \ffG^{i,+}_{\Delta,C} \circ \gG^i_{\Delta,C} = \gG^i_{\Delta,C} \circ \ffG^{i,+}_{\Delta,C}
$$
(i.e.\ $\ffG^i_{\Delta,C}$ is the semi-direct product of $\gG^i_{\Delta,C}$ and $\ffG^{i,+}_{\Delta,C}$).
\end{proof}
The group $\ffG^i$ acts in an obvious way on the set of formal Newton data.  Given
$\Omega \in \Nmaps$, we denote by $\ffG^i \cdot \Omega$ 
its orbit under such action. We adopt similar notations for the action of the subgroups $\ffG^i_{\Delta,C}$ and 
$\ffG^{i,+}_{\Delta,C}$.
\begin{Lemma}\label{lemma-normalsubgroupffG2}
Given a formal Newton data $\Omega \in \Nmaps^{i,\bfm}_{\Delta,C}$, the orbit
$$
\ffG^{i,+}_{\Delta,C} \cdot \Omega
$$
lies entirely in the class $\Nmaps^{i,\bfm}_{\Delta,C}$.  If we further assume that $\Omega$ is {\em stable} then
the orbit
$$
\ffG^{i}_{\Delta,C} \cdot \Omega
$$
also lies in the class $\Nmaps^{i,\bfm}_{\Delta,C}$.
\end{Lemma}
\begin{proof}
This is an immediate corollary of Lemma~\ref{lemma-normalsubgroupffG}, Corollary~\ref{corollary-howittransforms} and 
the definition of a stable Newton data.
\end{proof}
As a consequence, we obtain the following result on the invariance of the quantities $(\bfm,\Delta,C)$:
\begin{Proposition}\label{prop-invariancmdeltac}
Let $\Omega \in \Nmaps^{i,\bfm}_{\Delta,C}$ and 
$\widetilde{\Omega} \in \Nmaps^{i,\widetilde{\bfm}}_{\widetilde{\Delta},\widetilde{C}}$ be 
two stable Newton data which lie on the same $\ffG^i$-orbit.  Then 
$(\bfm,\Delta, C) = (\widetilde{\bfm},\widetilde{\Delta},\widetilde{C})$.
\end{Proposition} 
\begin{proof}
Let $(f,g,u,v,w) \in \ffG^i$ be the map such that $(f,g,u,v,w) \cdot \Omega = \widetilde{\Omega}$.
We shall prove that $(f,g,u,v,w)$ belongs to the subgroup $\ffG^{i}_{\Delta,C}$.

We define $P \subset \overline{\cQ}^2 
\times \overline{\cQ}$ as the subset of all pairs $(\Delta_0,C_0)$ such that
$(f,g,u,v,w)$ belongs to the subgroup $\ffG^i_{\Delta_0,C_0}$.  

Since the union $\bigcup_{\Delta,C}{\ffG^i_{\Delta,C}}$ exhausts $\ffG^i$, 
we know that $P$ is nonempty.  Let us fix an element 
$(\overline{\Delta},\overline{C}) \in P$.
Using Lemma~\ref{lemma-normalsubgroupffG},  
we can uniquely write
$$(f,g,u,v,w) = (\overline{f},\overline{g}) \circ (f_1,g_1,u_1,v_1,w_1)$$
with $(\overline{f},\overline{g}) \in \gG^i_{\overline{\Delta},\overline{C}}$ and
$(f_1,g_1,u_1,v_1,w_1) \in \ffG^{i,+}_{\overline{\Delta},\overline{C}}$.
From the discussion of subsection~\ref{subsect-subgroupsgG}, we can further
write the decomposition
\begin{equation}\label{givef1g1}
(\overline{f},\overline{g}) = (f_0,0) \cdot (f_1,g_1)
\end{equation}
with $(f_0,0) \in \gG^i_{\overline{\Delta}}$ and $(f_1,g_1) \in
\gG^{i,+}_{\overline{\Delta},
\overline{C}}$.

First of all, let us assume by absurd that $\bfm \ne \widetilde{\bfm}$.  Then, we immediately see that
either $\widetilde{\Omega}$ or $\Omega$ is in a hidden nilpotent configuration
and that the action of the map $(\overline{f},\overline{g})$ (or its inverse) causes
regular-nilpotent transition.  This contradicts the hypothesis that both $\Omega$
and $\widetilde{\Omega}$ are stable.

Assuming that $\bfm = \widetilde{\bfm}$, let us suppose by absurd that 
$\Delta >_\lex \widetilde{\Delta}$.  Then, the pair $(\widetilde{\Delta},C_0)$
necessarily
lies in the set $P$ (for some constant $C_0$).  Moreover, in the corresponding decomposition (\ref{givef1g1})
for $(\overline{\Delta},\overline{C}) := (\widetilde{\Delta},C_0)$, one has
$$
(f_0,0) = (\xi x^{\widetilde{\Delta_1}}y^{\widetilde{\Delta}_2},0), \quad \mbox{ for some constant }\xi \ne 0.
$$
However, using the above decomposition of $(f,g,u,v,w)$, we immediately see
that
$$
(f_0,0)^{-1}\cdot \widetilde{\Omega} \notin \Nmaps^{i,\widetilde{\bfm}}_{\widetilde{\Delta},\widetilde{C}}
$$
and this contradicts the hypothesis that $\widetilde{\Omega}$ is stable.

Finally, we assume by absurd that $(\bfm,\Delta) = (\widetilde{\bfm},\widetilde{\Delta})$
and $C > \widetilde{C}$. We prove the following:\\
{\em Claim: }There exists a
constant $C_0 < C$ such that the pair $(\Delta,C_0)$ lies in $P$.  Moreover, the decomposition
(\ref{givef1g1}) for $(\overline{\Delta},\overline{C}) := (\Delta,C_0)$ is such that map 
$$
(f_1,g_1) \in \gG^{i,+}_{\Delta,
C_0}
$$
is nonzero.  

Indeed, if the claim is false, the map 
$(f,g,u,v,w)$ should lie in $\ffG^{i}_{\Delta,C}$ and Lemma~\ref{lemma-normalsubgroupffG2} would imply that
$\widetilde{\Omega}$ also lies in $\Nmaps^{i,\bfm}_{\Delta,C}$. This contradicts the assumption that $C > \widetilde{C}$.

Using the above claim and Proposition~\ref{Proposition-stablenotresonant}, we conclude that
$(f_1,g_1)\cdot \Omega$ belongs to $\Nmaps^{i,\bfm}_{\Delta,C_0}$.
Consequently, $\widetilde{\Omega}$ also lies in
$\Nmaps^{i,m}_{\Delta,C_0}$ (i.e.\ $\widetilde{C} = C_0$).  
Taking the inverse map, we see that 
$$(f_1,g_1)^{-1}\cdot\widetilde{\Omega} \notin \Nmaps^{i,\bfm}_{\Delta,C_0}$$  
This contradicts the hypothesis that $\widetilde{\Omega}$ is stable.  The Proposition is proved.
\end{proof}
\subsection{Stabilization of Adapted Charts}
The main goal of this subsection is to prove that one can always find an stable Newton data for $(\cM,\Ax)$ at a nonelementary point $p$ lying on the divisor $\vD$.
\begin{Proposition}
Let $p \in \vD \cap A$ be a divisor point belonging to $\NElem(\cM)$.  
Then, there exists an {\em analytic} adapted local chart
$(U,(x,y,z))$ at $p$ such that the associated Newton data $\Omega = ((x,y,z),\iota_p,\Theta)$ 
is stable. 
\end{Proposition}
This Proposition will be an immediate consequence of the following result:
\begin{Proposition}[Stabilization of Adapted Charts]\label{prop-stabilizmap}
Let $p \in \vD \cap A$ be a divisor point belonging to $\NElem(\cM)$ and let $(U,(x,y,z))$ be an analytic
adapted local chart at $p$.
Then, there exists an analytic change of coordinates 
$$
\widetilde{y} = y + g(x), \quad \widetilde{z} = z + f(x,y), \quad \mbox{ where }f \in \cR\{x,y\}, \quad g \in \cR\{x\}
$$
with $f(0) = g(0) = 0$ such that Newton data associated to the new adapted local chart $(U,(x,\widetilde{y},\widetilde{z}))$ is stable. 
\end{Proposition}
We shall prove this Proposition using two Lemmas which describe the
stabilization of Newton data.

We shall say that a map $(f,g,u,v,w) \in \ffG^i$ is an {\em stabilization map}
for $\Omega$ if $f(0) = g(0) = 0$ and 
$(f,g,u,v,w)\cdot \Omega$ is a stable Newton data.

Similarly, we say that $(f,g,u,v,w) \in \ffG^i$ is an {\em edge stabilization
map} for $\Omega$ if $f(0) = g(0) = 0$ and $(f,g,u,v,w)\cdot \Omega$ is an edge stable 
Newton data.
\begin{Lemma}\label{lemma-edgestable}
Let $\Omega$ the Newton data for a divisor point 
$p \in \vD \cap A$ belonging to $\NElem(\cM)$.  
Then, there exists an edge stabilization map for $\Omega$ which has the form 
$$(f,0,1,1,1) \in \ffG^i,$$
for some $f \in \cR\{x,y\}$. 
\end{Lemma}
\begin{proof}
Define $\Omega_0 := \Omega$ and consider the sequence
\begin{equation}\label{seqomegas}
\Omega_0, \quad \Omega_1 = (f_0,0)\cdot \Omega_0 \quad,\ldots, \quad
\Omega_{n+1} = (f_n,0)\cdot \Omega_n \quad, \ldots
\end{equation}
where each $\Omega_{n+1}$ which is obtained by applying the basic edge preparation map 
$(f_n,0) = (\xi_n x^{a_n}y^{b_n},0)$ to $\Omega_n$ (see subsection~\ref{subsect-basicfaceedgeprep}).  

If there exists a finite natural number $n$ such that $\Omega_n$ is edge
stable, then we are done.  In fact, the
polynomial map
$$
f(x,y) = \sum_{i=0}^{n-1}{\xi_i x^{a_i}y^{b_i}}
$$
is such that $(f,0,1,1,1)\cdot \Omega$ is edge stable.

Otherwise, $\{(a_n,b_n)\}_{n\ge 0}$ is an infinite sequence,
strictly increasing for the lexicographical ordering.  Up to discarding a finite initial segment of the sequence (\ref{seqomegas}), we can assume that all $\Omega_n$ have the same main vertex. In fact, it follows from Lemma~\ref{lemma-regnilp} that a  regular-nilpotent transition in such sequence can occur after {\em at most }$m_3 + 1$ basic edge preparation maps (where $\bfm = (m_1,m_2,m_3)$ is the main vertex of $\Omega_0$).

Therefore, we can assume that for each $n \ge 0$, 
$$
\Omega_n \in \Nmaps_{(a_n,b_n)}^{i,\bfm}
$$
and $(a_{n+1},b_{n+1})>_\lex (a_n,b_n)$. Two cases can occur
(see figure~\ref{fig-casesiiistab}):
\begin{enumerate}[(i)]
\item $a_n \rightarrow \infty$ as $n \rightarrow \infty$;
\item There exist two natural numbers $a,N \in \cN$, with $a \ge \Delta_1$, such that
$$
a_n = a, \mbox{ for all }n \ge N,
$$
and $b_n \rightarrow \infty$ as $n \rightarrow \infty$;
\end{enumerate}

\begin{figure}[htb]
\psfrag{m}{\small $\bfm$}
\psfrag{a}{\small $a$}
\psfrag{D}{\small $\Delta$}
\psfrag{(i)}{(i)}
\psfrag{(ii)}{(ii)}
\begin{center}
\includegraphics[height=6cm]{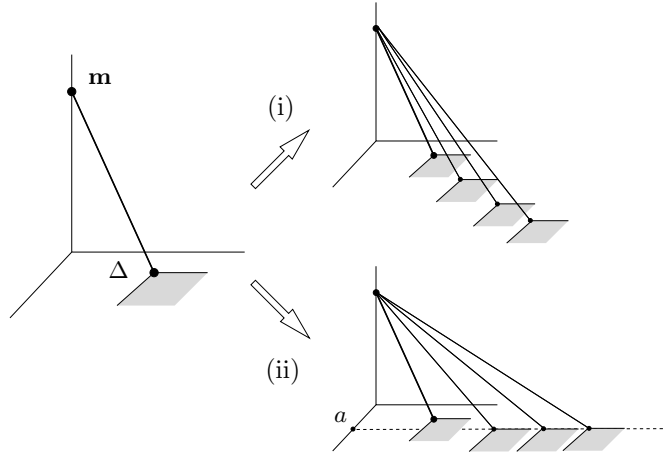}
\end{center}
\caption{Cases (i) and (ii) for the edge stabilization sequence.}
\label{fig-casesiiistab}
\end{figure}

We claim that case $(i)$ cannot occur.  In fact, if we consider the formal series
\begin{equation}\label{formalf}
f(x,y) = \sum_{i=0}^{\infty}{\xi_i x^{a_i}y^{b_i}}
\end{equation}
it is clear that $\widetilde{\Omega} = (f,0,1,1,1)\cdot \Omega$ belongs to
$\Nmaps^{i,\bfm}_{(\infty,\infty)}$.

If we write the vector field associated to $\Omega$ as
$$
\chi = F(x,y,z) \frac{\partial}{\partial x} +
G(x,y,z)\frac{\partial}{\partial y} + H(x,y,z) \frac{\partial}{\partial z}
$$
then the change of coordinates $\widetilde{z} = z + f(x,y)$ gives the formal
vector field
$$
\widetilde{\chi} = \widetilde{F}(x,y,\widetilde{z}) \frac{\partial}{\partial x} +
\widetilde{G}(x,y,\widetilde{z}) \frac{\partial}{\partial y} +
\widetilde{H}(x,y,\widetilde{z})
\frac{\partial}{\partial z}
$$
where $\widetilde{F}(x,y,\widetilde{z}) = F(x,y,\widetilde{z} - f(x,y))$,
$\widetilde{G}(x,y,\widetilde{z}) = G(x,y,\widetilde{z} - f(x,y))$ and
$$
\widetilde{H}(x,y,\widetilde{z}) = H(x,y,\widetilde{z} - f(x,y)) + \frac{\partial
f}{\partial x}(x,y) \widetilde{F}(x,y,\widetilde{z}) +  \frac{\partial
f}{\partial y}(x,y)\widetilde{G}(x,y,\widetilde{z})
$$
If we write $\bfm = (m_1,m_2,m_3)$ then the condition $\widetilde{\Omega} \in
\Nmaps^{i,\bfm}_{(\infty,\infty)}$ implies that there exists a factorization
$$
\left(
\begin{matrix}
\widetilde{F} \cr
\widetilde{G} \cr
\widetilde{H}
\end{matrix}
\right)
=
\widetilde{z}^{\, m_3}
\left(
\begin{matrix}
\widetilde{F}_1 \cr
\widetilde{G}_1 \cr
\widetilde{z}\widetilde{H}_1
\end{matrix}
\right)$$
for some formal germs $F_1,G_1,H_1$.
Going back to the original variables, we get
$$
\left(
\begin{matrix}
F \cr
G \cr
H
\end{matrix}
\right)
=
(z + f(x,y))^{m_3}
\left(
\begin{matrix}
F_1 \cr
G_1 \cr
(z + f(x,y)) H_1
\end{matrix}
\right)
$$
where $(F_1,G_1,H_1)$ is given by
$$
\left(
\begin{matrix}
F_1 \cr
G_1 \cr
H_1
\end{matrix}
\right)
=
\left[
\begin{matrix}
1 & 0 & 0 \cr
0 & 1 & 0 \cr
-\partial f /\partial x & - \partial f /\partial y & 1
\end{matrix}
\right]
\left(
\begin{matrix}
\widetilde{F}_1 \cr
\widetilde{G}_1 \cr
\widetilde{H}_1
\end{matrix}
\right)
$$
Recall now that $m_3 \ge 1$ (because $\Omega$ is centered at a point $p \in \NElem(\cM)$).
Therefore, we can apply the Corollary~\ref{corollary-uniquefactz} of Appendix~A to
conclude that $f(x,y)$ is necessarily an analytic function.

As a consequence, the Newton data $\widetilde{\Omega} = (f,0)\cdot \Omega$ is
analytic.  Notice however that the associated vector field
$\widetilde{\chi}$ violates the condition of being nondegenerate with respect
to the divisor (see definition~\ref{def-nondegvf}).  Indeed,  the
ideal $\vI_{\widetilde{\chi}}(\widetilde{z})$ is generated by
$$
\vI_{\widetilde{\chi}}(\widetilde{z}) = (\widetilde{z}\widetilde{F},\widetilde{z}\widetilde{G},\widetilde{H})
$$
and therefore $\vI_{\widetilde{\chi}}(\widetilde{z})$ is divisible by
$\widetilde{z}$.  Since $(z = 0)$ is not a component of the divisor, we get a contradiction to the assumption that $\widetilde{\chi}$ is nondegenerate.

Suppose now that $(ii)$ holds.  Then, the formal map $f$ given in
(\ref{formalf})
can be written in the form
\begin{equation}\label{faspolyn}
f(x,y) = f_{\delta}(y)x^{\delta} + f_{\delta+1}(y)x^{\delta + 1} + \cdots + f_a(y)x^a, \quad \mbox{ with
}f_\delta,\ldots,f_a \in \cR[[y]]
\end{equation}
where $\delta := \Delta_1$. We claim that $f_\delta,\ldots,f_a$ are analytic germs.

Indeed, let us apply the change of coordinates $\widetilde{z} = z +
f(x,y)$.  Keeping the same notation used above, we get the formal vector field
$$
\widetilde{F}(x,y,\widetilde{z}) \frac{\partial}{\partial x} +
\widetilde{G}(x,y,\widetilde{z}) \frac{\partial}{\partial y} +
\widetilde{H}(x,y,\widetilde{z})
\frac{\partial}{\partial z}
$$
which is associated to the (formal) Newton data $\widetilde{\Omega} = (f,0)\cdot
\Omega$.  From the hypothesis, we know that $\widetilde{\Omega}$ belongs to
$\Nmaps^{i,\bfm}_{\widetilde{\Delta},\widetilde{C}}$, for some
$\widetilde{\Delta} = (\widetilde{a},\widetilde{b})$ such that $\widetilde{a}
> a$.  This is equivalent to say that, if we
consider the homomorphic images $[\widetilde{F}],[\widetilde{G}],[\widetilde{H}]$
in the quotient ring $\widehat{R}_a = \cR[[x,y,z]]/(x^a)\cR[[x,y,z]]$, we get
$$
\left(
\begin{matrix}
[\widetilde{F}] \cr
[\widetilde{G}] \cr
[\widetilde{H}]
\end{matrix}
\right)
=
[\,\widetilde{z}\,]^{m_3}
\left(
\begin{matrix}
[\widetilde{F}_1] \cr
[\widetilde{G}_1] \cr
[\widetilde{z}\widetilde{H}_1]
\end{matrix}
\right)
$$
Going back to the original functions $F,G,H$ and using the same reasoning used
above, we conclude from Corollary~\ref{corollary-uniquefactzquot} of Appendix~A  that
the germ $[z - f(x,y)] \in \widehat{R}_a$ belongs to the
homomorphic image of a convergent germ.  

Therefore, using the expression
(\ref{faspolyn}), we conclude that $f_\delta,\ldots,f_a$ are convergent germs
and the function $f$ given by (\ref{faspolyn}) is analytic.  This proves the claim.

Let us call the map $(f,0) \in \hG^i$ the {\em
extended edge
preparation step} associated to $\Omega$.  If the Newton data 
$$\Omega^{(1)} := (f,0) \cdot \Omega$$ 
is edge prepared, we are done.  Otherwise, we can start all over again, and define a
new edge extended edge preparation step $(f^{(1)},0)$ associated to
$\Omega^{(1)}$ and get $\Omega^{(2)} = (f^{(1)},0)\cdot\Omega^{(1)}$, 
$\Omega^{(3)} = (f^{(2)},0)\cdot\Omega^{(2)}$, and so on.

Suppose that we can iterate such procedure infinitely many times.
Then, we get a sequence of Newton data $\{\Omega^{(n)}\}$,
where each element of the sequence is obtained from its predecessor by an extended edge
preparation step $(f^{(n)},0)$, for some analytic germ $f^{(n)} \in \cR\{x,y\}$ such that
\begin{equation}\label{faspolynn}
f^{(n)} = O(x^{\Delta_1^{(n)}})
\end{equation}
for some strictly increasing sequence $\{\Delta_1^{(n)}\}$ of natural numbers.
We claim that there exists a finite $n \in \cN$ such that $\Omega^{(n)}$ is edge prepared.

Indeed, suppose by absurd that this sequence is infinite.  Then, it follows
from the the expression (\ref{faspolynn}) that the composed map
$$
\mathbf{f}_n = f^{(n)} \circ \cdots \circ f^{(n-1)} \circ \cdots \circ f^{(1)} \circ f
$$
converges in the Krull topology, as $n \rightarrow \infty$, to a formal map $\widetilde{f} \in
\cR[[x,y]]$.
Moreover, the formal Newton data $\widetilde{\Omega} :=
(\widetilde{f},0)\cdot \Omega$ belongs to
$\Nmaps^{i,\bfm}_{(\infty,\infty)}$. Using
the same reasoning used in the proof of item $(i)$, we get a contradiction.
\end{proof}
\begin{Lemma}\label{lemma-facestable}
Let $\Omega$ be an analytic Newton data centered at a divisor point 
$p \in \vD \cap A$ belonging to $\NElem(\cM)$.  Assume that $\Omega$ is edge stable.  Then, there exists
a stabilization map for $\Omega$ which has the form 
$$(f,g,1,1,1) \in \ffG^i$$
for some convergent germs $g \in \cR\{x\}$ and $f \in \cR\{x,y\}$ with $f(0) = g(0) = 0$. 
\end{Lemma}
\begin{proof}
Define $\Omega_0 := \Omega$ and consider the sequence
$$
\Omega_0, \quad \Omega_1 = (f_0,g_0)\cdot \Omega_0 \quad,\ldots, \quad
\Omega_{n+1} = (f_n,g_n)\cdot \Omega_n \quad, \ldots
$$
where $\Omega_{n+1}$ is obtained by applying the basic face preparation
map $(f_n,g_n)$ to $\Omega_n$ (see subsection~\ref{subsect-basicfaceedgeprep}).  
Notice that, for each $n \ge 0$,
$$
\Omega_n \in \Nmaps_{\Delta,C_n}^{i,\bfm}
$$
and $C_0 < C_1 < C_2 < \cdots$ is a strictly increasing sequence of rational
numbers.

\begin{figure}[htb]
\psfrag{D}{\small $\Delta$}
\psfrag{m}{$\bfm$}
\psfrag{C}{\small $C$}
\begin{center}
\includegraphics[height=4cm]{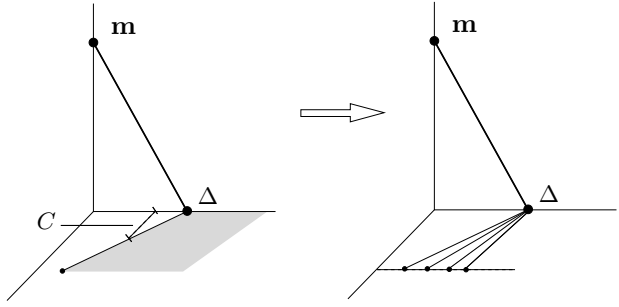}
\end{center}
\caption{The sequence of basic face preparations.}
\label{fig-caseCstab}
\end{figure}

Notice that the rational numbers $C_n$ always belongs to the finite lattice
$\frac{1}{[(m_3 + 1)\Delta]!}\cZ$ (see Remark~\ref{remark-supportedgepreserving}).  
Therefore $C_n \rightarrow \infty$ as $n
\rightarrow \infty$.

If there exists a finite natural number $n$ such that $\Omega_n$ is
stable, then we are done.  In fact, the composed map
$$
(\mathbf{f}_n,\mathbf{g}_n) := (f_n,g_n)\circ \cdots \circ (f_0,g_0)
$$
is a polynomial map and $(\mathbf{f}_n,\mathbf{g}_n,1,1,1) \in \ffG^i$ is such that 
$(\mathbf{f}_n,\mathbf{g}_n,1,1,1)\cdot \Omega$ is stable.

Otherwise, $\{\Omega_n\}$ forms an infinite sequence and the condition that
$(f_n,g_n) \in \gG^i_{\Delta,C_n}$ implies that
the sequence of composed maps $\{(\mathbf{f}_n,\mathbf{g}_n)\}_{n \ge 0}$ converges
(in the Krull topology) to a pair of formal maps $(f,g)$ such that
$$
(f,g,1,1,1) \in \ffG^i
$$
Moreover, $f \in \cR[[x,y]]$ can be written in the form
$$
f(x,y) = f_0(x) + f_1(x) y + \cdots + f_{b-1}(x)y^{b-1}
$$
with each $f_i$ belonging to $\cR[[x]]$ and $b = \lceil\Delta_2\rceil$.
Notice that
$$\widetilde{\Omega} := (f,g,1,1,1) \cdot \Omega$$
is a formal Newton data which belongs to $\Nmaps^{i,\bfm}_{\Delta,\infty}$.

We claim that $(f,g,1,1,1)$ is an analytic map.

Suppose initially that $g = 0$.  Let us write the vector field which is
associated to $\widetilde{\Omega}$ as,
$$
\widetilde{F}(x,\widetilde{y},\widetilde{z}) \frac{\partial}{\partial x} +
\widetilde{G}(x,\widetilde{y},\widetilde{z}) \frac{\partial}{\partial \widetilde{y}} +
\widetilde{H}(x,\widetilde{y},\widetilde{z})
\frac{\partial}{\partial z}
$$
From the hypothesis, we know that the coefficients
$\widetilde{F}(x,\widetilde{y},0)$, $\widetilde{G}(x,\widetilde{y},0)$ and $\widetilde{H}(x,\widetilde{y},0)$ are
such that
\begin{equation}\label{condfgh}
\widetilde{F}(x,\widetilde{y},0),\widetilde{H}(x,\widetilde{y},0) \in (\widetilde{y}^B)\cR[[x,y]], \quad \mbox{
and }\quad \widetilde{G}(x,\widetilde{y},0) \in (\widetilde{y}^{B+1})\cR[[x,y]]
\end{equation}
where $B := \lceil m_2 + \Delta_2 m_3 \rceil$.
Since $\Omega$ is centered at a point $p \in \NElem(\cM)$, we know also that $B \ge 1$ (because
otherwise $\Omega$ would be in a final situation, contradicting Proposition~\ref{prop-finalsituation}).  

\begin{figure}[htb]
\psfrag{e}{\small $\medge$}
\psfrag{m}{$\bfm$}
\psfrag{B}{\small $B$}
\psfrag{ll}{\small $m_2 + \Delta_2 m_3$}
\begin{center}
\includegraphics[height=4cm]{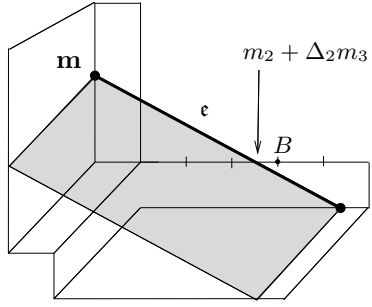}
\end{center}
\caption{The number $B$.}
\label{fig-numberB}
\end{figure}

Hence, we
can use the same reasoning used in the proof of
Lemma~\ref{lemma-edgestable} to show that $f$ is analytic.  

Let us suppose now that $g \ne 0$.
Then, from the definition of $\Nmaps^{i,\bfm}_{\Delta,C}$, we know that $(y = 0)$ is not a
local irreducible component of the divisor at the point $p$. 

Moreover, it follows from the condition (\ref{condfgh}) that the coefficients $\widetilde{F},
\widetilde{G}, \widetilde{H}$ belong to the ideal
$(\widetilde{y},\widetilde{z})$.  This is equivalent to say that
the analytic coefficients $F,G,H$ of the original vector field $\chi$ are contained in the ideal
$$
J = (\; y + g(x), \, z + f(x,-g(x)\; )
$$
Since $(y = 0)$ is not a divisor component and $\chi$ is
nondegenerate, the ideal $J$ is
necessarily the defining ideal of an irreducible one-dimensional component of the germ
$\Ze(\chi)_p$ (the analytic set of zeros of $\chi$).

In other words, the
prime ideal $J$ is an element of the irreducible primary decomposition of
$I := \mathrm{rad}(F,G,H)$.
Therefore, it follows from Lemma~\ref{Lemma-uniquefactgen} that the
functions $g(x)$ and $f(x,-g(x))$ are necessarily analytic.

Now, we can decompose the map $(f,g,1,1,1)$ in a unique way as
$$
(f,g,1,1,1) = (\widetilde{f},0,1,1,1) \circ (f_1,g_1,1,1,1)
$$
where $f_1(x) := f(x,-g(x))$, $g_1(x) := g(x)$  and
$$\widetilde{f} = \widetilde{f}_0(x) + \cdots  + \widetilde{f}_{b-1}(x)y^{b-1}
\in \cR[[x,y]]$$
is a conveniently chosen formal map.  Since
$(f_1,g_1,1,1,1)\cdot \Omega$ is analytic, we conclude as above that
$(\widetilde{f},0,1,1,1)$ is also analytic.  This completes the proof of the Lemma.
\end{proof}
\begin{proof}({\bf of Proposition~\ref{prop-stabilizmap} })  
We define the stabilization map $(f,g,1,1,1) \in \ffG^i$
as the composition
$$
(f,g,1,1,1) := (f_2,g_2,1,1,1) \circ (f_1,0,1,1,1)
$$
where $(f_1,0,1,1,1)$ and $(f_2,g_2,1,1,1)$ are respectively the edge stabilization map given by Lemma~\ref{lemma-edgestable}
and the face stabilization map given by
Proposition~\ref{lemma-facestable}.
\end{proof}
We denote by 
$
\Stab \Omega
$
the {\em stabilized} Newton data defined by the above construction. The transition from $\Omega$ to $\Stab \Omega$ will be called a {\em stabilization} of the Newton data.

We remark that the notion of stable Newton data at a point $p$ is {\em independent} of the choice of the local generator for the line field from which such data is defined (see Lemma~\ref{lemma-choicelocgen}).  

For this reason, and for notational simplicity, we often say that an adapted local chart $(U,(x,y,z))$ at $p$ is {\em stable} (resp.\ {\em edge stable}) whenever the corresponding Newton data $\Omega = ((x,y,z),\iota_p,\Theta)$ is stable (resp.\ {\em edge stable}), where $\Theta$ is defined by fixing some arbitrary choice of local generator for the line field at $p$.
\begin{Remark}
Given a {\em formal} Newton data $\Omega = ((x,y,z),\iota,\Theta)$ at a point $p \in \NElem(\cM)$, it is easy to see that the condition $\Delta_1 > 0$ immediately implies that the {\em formal curve} $\{x = z = 0\}$ lies entirely in the set $\NElem(\cM)$.  \\
Note also that the condition of nondegeneracy for the local generator $\chi$ of the line field guarantees that $\NElem(\cM) \cap \{z = 0\}$ is
an analytic set of dimension at most equal to one, and therefore  $\{x = z = 0\}$ is necessarily an {\em analytic} curve. \\
However, these conditions  {\em do not} imply that
the formal coordinates $(x,y,z)$ are analytic.  This is the reason why we needed extra arguments to prove the analyticity of the stabilization map at the end of the proof of Lemma~\ref{lemma-facestable}.
\end{Remark}
\subsection{Newton Invariant and Local Resolution of Singularities}\label{subsect-Newtoninvarlocres}
Let $(\cM,\Ax)$ be a controlled singularly foliated manifold, and let 
$p \in \vD \cap A$ be a divisor point belonging to $\NElem(\cM)$.  

In Proposition~\ref{prop-stabilizmap},
we have proved that there always exists an 
(analytic) adapted local chart $(U,(x,y,z))$ for $(\cM,\Ax)$ at $p$ such that the associated Newton data $\Omega$ is stable.

The {\em Newton invariant} for $(\cM,\Ax)$ at $p$ is the vector of natural numbers
$$
\Inv(\cM,\Ax,p) = \Inv(\Omega) \in \cN^6 
$$ 
where $\Inv(\Omega)$ is given by definition~\ref{def-multiplicity}.

Let $\vN$ be the Newton polyhedron associated to $\Omega$.
The {\em weight-vector} for $(\cM,\Ax)$ at $p$ is the nonzero vector $\bfomega \in \cN^3$ such that
$\mathrm{mdc}(\omega_1,\omega_2,\omega_3) = 1$ and 
$$
\vF = \vN \cap \{ \bfv \in \cR^3 \mid \< \bfomega, \bfv \> = \mu \}, \mbox{ for some }\mu \in \cZ
$$
where $\vF$ is the main face of $\vN$.  The integer number
$\mu$ given in the formula is called the {\em face
order} for $(\cM,\Ax)$ at $p$.
\begin{Remark}\label{remark-omegadeltaC}
If $\Omega$ belongs to the class $\Nmaps^{i,\bfm}_{\Delta,C}$ then we can
explicitly compute that $\bfomega = k\, \bfalpha$, where $\bfalpha$ is defined as follows:
\begin{enumerate}[(i)]
\item If $C = 0$ then $\bfalpha =  (1,0,\Delta_1)$;
\item If $C = \infty$ then $\bfalpha = (0,1,\Delta_2)$;
\item If $0 < C < \infty$ then $\bfalpha =  (1,C,C\Delta_2)$;
\end{enumerate}
and $k \in \cN$ is the least natural number such that $k \bfalpha$ belongs to $\cN^3$.
\end{Remark}
The {\em local blowing-up center} associated to $(\cM,\Ax)$ at $p$ is the submanifold $Y_p \subset U$ defined as follows:
\begin{enumerate}[(i)]
\item If $\bfomega = (*,*,*)$ then $Y_p = \{x = y = z = 0\}$.
\item If $\bfomega = (*,0,*)$ then $Y_p = \{x = z = 0\}$.
\item If $\bfomega = (0,*,*)$ then $Y_p = \{y = z = 0\}$.
\end{enumerate}
where $*$ denotes nonzero natural numbers.
\begin{Lemma}
The local blowing-up center $Y_p$ lies in $\NElem(\cM)$.
\end{Lemma}
\begin{proof}
Let us consider the case where $Y_p = \{x = z = 0\}$, which corresponds to item (ii) in the above definition  (the reasoning for the item (iii) is analogous). Fixing some local generator
$\chi$ for the line field at $p$, and writing $\bfomega = (\omega_1,0,\omega_3)$ with $\omega_1,\omega_3 \in \cN_{>0}$, we have
$$
\chi = F \frac{\partial}{\partial x} + G \frac{\partial}{\partial y} + H\frac{\partial}{\partial z}
$$  
where $F,G,H$ are analytic germs with $\bfomega$-multiplicity given respectively
by $\mu - \omega_1$, $\mu$ and $\mu - \omega_3$.  Consider now a translation of coordinates
$(\widetilde{x},\widetilde{y},\widetilde{z}) = (x,y - \eta,z)$, for some constant $\eta \in \cR$, and let $\widetilde{\chi} = \widetilde{F} \partial/\partial \widetilde{x} + \widetilde{G}\partial/\partial \widetilde{y} + \widetilde{H}\partial/\partial \widetilde{z}$ be the resulting local generator of the line field.
Since $\omega_2 = 0$, it is obvious that the germs  $\widetilde{F},\widetilde{G},\widetilde{H}$ have $\bfomega$-multiplicity greater or equal to $\mu - \omega_1$, $\mu$ and $\mu - \omega_3$, respectively.
Therefore, the point $\widetilde{p} \in Y_p$ which is the center of the new coordinates
$(\widetilde{x},\widetilde{y},\widetilde{z})$ also belongs to $\NElem(\cM)$.
\end{proof}
\begin{Proposition}\label{prop-preservesstructure}
Let $(U,(x,y,z))$ and $(\widetilde{U},(\widetilde{x},\widetilde{y},\widetilde{z}))$ be two
stable adapted local charts
at $p$.  Then, the corresponding numbers
$$
\Inv(\cM,\Ax,p), \quad \bfomega \quad \mbox{ and }\quad \mu
$$
which are associated to these two charts are equal.  Moreover, the respective local blowing-up centers
$Y_p$ and $\widetilde{Y}_p$ coincide on $U \cap \widetilde{U}$.  Finally, the transition map 
$$
\phi(x,y,z) = (\widetilde{x},\widetilde{y},\widetilde{z}) 
$$
preservers the $\bfomega$-quasihomogeneous structure on $\cR^3$.
\end{Proposition}
\begin{proof}
The first part of the enunciate follows from Proposition~\ref{prop-invariancmdeltac}.

In order to prove the second part of the enunciate, it suffices to remark that the transition map
$\phi$ has the form
$$
\widetilde{x} = x\, u, \quad \widetilde{y} = g(x) + y\, v, \quad \widetilde{z} = f(x,y) + z\, w
$$
and the map $(f,g,u,v,w)$ is a member of the subgroup $\ffG^i_{\Delta,C}$ (by the proof of 
Proposition~\ref{prop-invariancmdeltac}).  Using the explicit 
definition  of such subgroup and the Remark~\ref{remark-omegadeltaC}, we immediately conclude that 
$\phi$ preservers the $\bfomega$-quasihomogeneous structure on $\cR^3$.
\end{proof}

Let $\Omega$ be a stable Newton data for $(\cM,\Ax)$ at $p$, associated to an adapted local chart $(U,(x,y,z))$.

The {\em local blowing-up} for $(\cM,\Ax)$ at $p$ is the $\bfomega$-weighted
blowing-up of
$$
\Phi: \widetilde{\cM} \rightarrow \cM \cap U
$$
with center on $Y_p$, with respect to the trivialization given by $(U,(x,y,z))$.

The apparent arbitrariness in the choice of $(U,(x,y,z))$ can be removed as
follows. Consider two local blowing-ups at $p$, 
$$\Phi_i:\widetilde{\cM}_i \rightarrow \cM \cap U_i, \quad i = 1,2$$ 
which are associated to distinct stable adapted charts
$(U_i,(x_i,y_i,z_i))$ ($i =
1,2$).  

Using Proposition~\ref{prop-preservesstructure}, it follows that 
(up to restricting each $U_i$ to some smaller
neighborhood of $p$), there exists an isomorphism $\Psi:\widetilde{\cM}_1
\rightarrow \widetilde{\cM}_2$ (in the obvious sense of isomorphism between
singularly foliated manifolds) which makes the following diagram commutative: 
\begin{eqnarray*}\label{commutdiagramlocalblowup}
\begin{CD}
\widetilde{\cM}_1 @>{\Psi}>> \widetilde{\cM}_2\\
@V{\Phi_1}VV                                   @VV{\Phi_2}V \\
\cM \cap U  @>{\mathrm{id}}>>  \cM \cap U
\end{CD}
\end{eqnarray*}
where $\mathrm{id}$ is the identity map and $U = U_1 \cap U_2$.

The main Theorem of this section can now be enunciated as follows:
\begin{Theorem}[Local Resolution of Singularities]\label{theorem-localresolution}
Let $(\cM,\Ax)$ be a controlled singularly foliated manifold and let
$p \in \vD \cap A$ be a divisor point in $\NElem(\cM)$.
Consider the {\em local blowing-up} for $(\cM,\Ax)$ at $p$, 
$$
\Phi: \widetilde{\cM} \rightarrow \cM \cap U
$$
with respect to some stable adapted chart $(U,(x,y,z))$.  Then, there exists an axis 
$\widetilde{\Ax} = (\widetilde{A},\widetilde{\vZ}\, )$ for 
$\widetilde{\cM}$ such that each point 
$\widetilde{p} \in \Phi^{-1}(\{p\}) \cap \widetilde{A}$ belonging to
$\NElem(\widetilde{\cM})$ is such that
$$
\Inv(\widetilde{\cM},\widetilde{\Ax},\widetilde{p}\, ) \; <_\lex \; \Inv(\cM,\Ax,p)
$$ 
\end{Theorem}
The proof of this Theorem will be given in subsection~\ref{subsect-localresoltheorem}.
\subsection{Directional Blowing-ups}\label{subsect-directionalexpressions}
Let us fix a stable adapted chart $(U,(x,y,z))$ 
at a divisor point $p \in \NElem(\cM)$
and consider the corresponding $\bfomega$-weighted local blowing-up
$$
\Phi: \widetilde{\cM} \rightarrow \cM \cap U
$$
as defined in the previous subsection.

The Theorem~\ref{theorem-localresolution} will be proved by studying the
effect of such
blowing-up in the $x$, $y$, and $z$-directional charts (see
subsection~\ref{subsect-directionalcharts}).

Let $\Omega = ((x,y,z),\iota_p,\Theta)$ be the Newton data 
associated to the adapted local chart $(U,(x,y,z))$ (for some choice of local generator $\chi$ of $\linef$).
It will be convenient to look at the directional blowing-ups as
transformations on the Newton map $\Theta$.  
For this, we consider the following matrices:
\begin{enumerate}
\item $x$-directional transformation matrices:
$$
B_x = \left[
\begin{matrix}
\omega_1 & \omega_2 & \omega_3 \cr
0 & 1 & 0 \cr
0 & 0 & 1
\end{matrix}
\right]
, \quad
M_x = \left[
\begin{matrix}
{1}/{\omega_1} & 0 & 0 \cr
-{\omega_2}/{\omega_1} & 1 & 0 \cr
-{\omega_3}/{\omega_1} & 0 & 1
\end{matrix}
\right]
$$
\item $y$-directional transformation matrices:
$$
B_y = \left[
\begin{matrix}
1 & 0 & 0 \cr
\omega_1 & \omega_2 & \omega_3 \cr
0 & 0 & 1
\end{matrix}
\right]
, \quad
M_y = \left[
\begin{matrix}
1 & -{\omega_1}/{\omega_2} & 0  \cr
0 & -{1}/{\omega_2} &  0 \cr
0 & -{\omega_3}/{\omega_2} & 1
\end{matrix}
\right]
$$
\item $z$-directional transformation matrices:
$$
B_z = \left[
\begin{matrix}
1 & 0 & 0 \cr
0 & 1 & 0 \cr
\omega_1 & \omega_2 & \omega_3
\end{matrix}
\right]
, \quad
M_z = \left[
\begin{matrix}
1 & 0 & -{\omega_1}/{\omega_3}\cr
0 & 1 & -{\omega_2}/{\omega_3}\cr
0 & 0 & 1/{\omega_3}
\end{matrix}
\right]
$$
\end{enumerate}
We consider also the permutation matrices
$$
I = \left[
\begin{matrix}
0 & 1  & 0  \cr
1 & 0 &  0 \cr
0 & 0 & 1
\end{matrix}
\right]
\quad
J = \left[
\begin{matrix}
0 & 0  & 1  \cr
1 & 0 &  0 \cr
0 & 1 & 0
\end{matrix}
\right]
$$
The {\em directional blowing-ups} of $\Theta$
are the Newton maps $\Blx \Theta$, $\Bly \Theta$ and $\Blz \Theta$
given respectively by
\begin{eqnarray*}
\Blx \Theta (B_x\bfv - \mu\bfe_1) &=& \eps^{v_1}M_x\Theta(\bfv) \quad \mbox{
(defined for }\omega_1 > 0)\\
\Bly \Theta (I\, B_y\bfv - \mu\bfe_1) &=& \eps^{v_2} I \, M_y\Theta(\bfv) \quad \mbox{
(defined for }\omega_2 > 0)\\
\Blz \Theta (J\, B_z\bfv - \mu\bfe_1) &=& \eps^{v_3} J\, M_z\Theta(\bfv)\quad \mbox{
(defined for }\omega_3 > 0)
\end{eqnarray*}
where $\eps \in \{ -1, 1\}$ and $\bfv \in \cZ^3$.  The directional blowing-ups of the Newton data $\Omega$ are defined as follows:
\begin{itemize}
\item {\em $x$-directional blowing-up }: $\Blx \Omega = ((\bar x,\bar y,\bar z),\overline\iota,\Blx \Theta)$;
\item {\em $y$-directional blowing-up }: $\Bly \Omega = ((\bar x,\bar y,\bar z),\overline\iota,\Bly \Theta)$;
\item {\em $z$-directional blowing-up }: $\Blz \Omega = ((\bar x,\bar y,\bar z),\overline\iota,\Blz \Theta)$;
\end{itemize}
where $\overline{\iota} = \iota_p \cup [n]$ (with $n = 1 + \max\{i \mid i \in \mainlist\}$ for $\mainlist \ne \emptyset$ and $n = 1$ for $\mainlist = \emptyset$) and $(\bar x,\bar y,\bar z)$ is a chart respectively defined by the
following singular changes of coordinates
\begin{eqnarray*}
x\mbox{-directional blowing-up}&:&\quad x = \eps\, \overline{x}^{\omega_1}, \quad y = \overline{x}^{\omega_2}\overline{y},
\quad
z = \overline{x}^{\omega_3}\overline{z}, \label{singcoordsx}\\
y\mbox{-directional blowing-up}&:&\quad x = \overline{x}^{\omega_1}\overline{y}, \quad y = \eps\, \overline{x}^{\omega_2},
\quad
z = \overline{x}^{\omega_3}\overline{z},  \label{singcoordsy}\\
z\mbox{-directional blowing-up}&:&\quad x = \overline{x}^{\omega_1}\overline{y}, \quad y = \overline{x}^{\omega_2} \overline{z},
\quad
z = \eps\, \overline{x}^{\omega_3} \label{singcoordsz}
\end{eqnarray*}
followed by a division by $\overline{x}^{\mu}$.
Notice that there exists a cyclic permutation of coordinates in the $y$ and
$z$-directional blowing-ups (corresponding to the permutation matrices $I$ and $J$).

\begin{figure}[htb]
\psfrag{Fbl}{\small $\vF$}
\psfrag{xdir}{\small $x$-dir}
\psfrag{ydir}{\small $y$-dir}
\psfrag{zdir}{\small $z$-dir}
\begin{center}
\includegraphics[height=7cm]{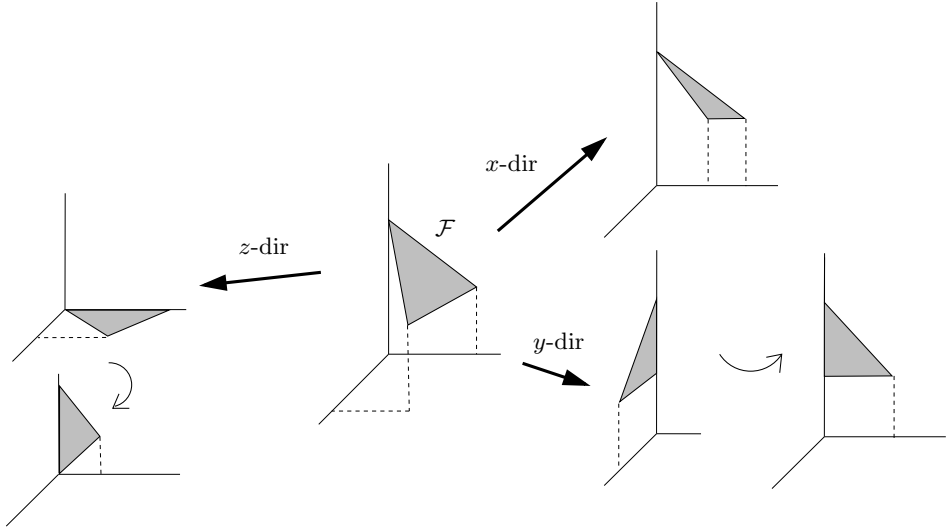}
\end{center}
\caption{The Directional Blowing-Ups.}
\label{fig-dirblowups}
\end{figure}

In figure~\ref{fig-dirblowups}, we give an illustration of the {\em movement}
of the Newton polyhedron which is caused by these maps.

\subsection{$x$-Directional Blowing-up}
Let $(\cM,\Ax)$ be a controlled singularly foliated manifold and let
$p \in \vD \cap A$ be a divisor point in $\NElem(\cM)$.  

Let $\Omega \in \Nmaps^{i,\bfm}_{\Delta,C}$ be a stable Newton data at $p$, with coordinates $(x,y,z)$.  
In this subsection, we assume that the corresponding 
weight-vector $\bfomega = (\omega_1,\omega_2,\omega_3)$ is such that $\omega_1 > 0$.  Geometrically, this assumption means that the local blowing-up center is distinct from the axis
$\{y = z = 0\}$. 

The {\em $x$-directional translation group} is
defined as follows:
\begin{enumerate}[(i)]
\item If $\Delta_1 = 0$ then $\gG^\transl_x := \gG^1_{(0,0),0}$;
\item If $\Delta_1 > 0$ then $\gG^\transl_x := \gG^2_{(0,0),0}$;
\end{enumerate}
In other words, if $\Delta_1 > 0$ then $\gG^\transl_x$ is the group of all translations
$$
\widetilde{z} = z + \xi
$$
for some real constant $\xi \in \cR$.
If $\Delta_1 = 0$ then $\gG^\transl_x$ is the group of all translations
$$
\widetilde{y} = y + \eta, \quad \widetilde{z} = z + \xi
$$
for some real constants $\eta,\xi \in \cR$. In this subsection, we shall prove the following result:
\begin{Proposition}\label{prop-stablexblup}
Given a stable Newton data $\Omega$, let
$\overline{\Omega} = \Blx \Omega$ be its $x$-directional blowing-up.
Then, for each $(\xi,\eta) \in \gG^\transl_x$, either the translated Newton data
$$
\overline{\Omega}_{\xi,\eta} := (\xi,\eta) \cdot \overline{\Omega}
$$
is centered at an elementary point $\widetilde{p} \in \Elem(\widetilde{\cM})$ or 
$$
\Inv(\widetilde{\Omega}_{\xi,\eta}) <_\lex \Inv(\Omega)
$$
where $\widetilde{\Omega}_{\xi,\eta} = \Stab \overline{\Omega}_{\xi,\eta}$ is the stabilization of $\overline{\Omega}_{\xi,\eta}$.
\end{Proposition}
The proof will be given at the end of subsection~\ref{subsect-xdirectDelta=0} and will depend on 
several Lemmas.
  
First of all, let us look at the effect of $\Blx$ on the main face $\vF$.
\begin{Lemma}\label{lemma-directx}
Let $\overline{\Omega} := \Blx \Omega$ be the $x$-directional blowing-up of
$\Omega$.  Then, there exists a
bijective correspondence
\begin{eqnarray*}
\supp(\Omega) \; \cap\;   \vF &\longrightarrow & \supp(\overline{\Omega}) \; \cap\;  \{0\}
\times \cZ^2
\\
\bfv = (v_1,v_2,v_3) &\longmapsto& \pi_x(\bfv) = (0,v_2,v_3)
\end{eqnarray*}
such that the corresponding Newton maps $\overline{\Theta}$ and $\Theta$ satisfy 
$\overline{\Theta}[\pi_x(\bfv)] = M_x \Theta[\bfv]$.
\end{Lemma}
\begin{proof}
This is immediate from the definition of $\Blx$.
\end{proof}

The matrices $M_x$ and $B_x$ which appear in the definition of the 
$x$-directional blowing-up $\Blx$ can be written as products
$B_x = B_x^2\, B_x^1$ and $M_x = M_x^2 \, M_x^1$, where 
$$
B_x^1 = \left[
\begin{matrix}
\omega_1 & 0 & 0 \cr
0 & 1 & 0 \cr
0 & 0 & 1
\end{matrix}
\right]
, \quad
M_x^1 = \left[
\begin{matrix}
{1}/{\omega_1} & 0 & 0 \cr
0 & 1 & 0 \cr
0 & 0 & 1
\end{matrix}
\right]
$$
and 
$$
B_x^2 = \left[
\begin{matrix}
1 & \omega_2 & \omega_3 \cr
0 & 1 & 0 \cr
0 & 0 & 1
\end{matrix}
\right]
, \quad
M_x^2 = \left[
\begin{matrix}
1 & 0 & 0 \cr
-{\omega_2} & 1 & 0 \cr
-{\omega_3} & 0 & 1
\end{matrix}
\right]
$$
Therefore, the map $\Blx$ can be written as the composition $\Blx = \Blx^2 \circ \Blx^1$, where
$$
\Blx^1 \Theta(B_x^1 \bfv) = \eps^{v_1} M_x^1 \Theta(\bfv), \quad \mbox{ and }\quad
\Blx^2 \Theta(M_x^2 \bfv - \mu\bfe_1) =  M_x^2 \Theta(\bfv)
$$
Notice that the maps $\Blx^1$ and $\Blx^2$ correspond respectively 
to the singular changes of coordinates
$$
x = \eps\, \overline{x}^{\omega_1}, \;\; y = \overline{y},
\;\; 
z = \overline{z} \quad \mbox{ and }\quad 
x = \overline{x}, \;\; y = \overline{x}^{\omega_2}\overline{y},
\;\;
z = \overline{x}^{\omega_3}\overline{z}
$$
followed by a division by $\overline{x}^\mu$ (for $\eps \in \{-1,1\}$). 
\subsection{The Effect of a Ramification}\label{subsect-effectramif}
The expressions given in the previous subsection show 
that a $x$-directional blowing-up can always be written 
as a composition of a {\em ramification} $x = \eps \overline{x}^{\omega_1}$ 
followed by a sequence of {\em homogeneous blowing-ups}.
\begin{Example}
For $\bfomega = (2,3,2)$, the $x$-directional blowing-up can be decomposed
as the ramification $(x,y,z) = (\overline{x}^2,\overline{y},\overline{z})$ followed by the sequence of blowing ups
$$
(\overline{x},\overline{y},\overline{z}) = (x_1,x_1y_1,x_1z_1), \;
(x_1,y_1,z_1) = (x_2,x_2 y_2, x_2z_2), \; (x_2,y_2,z_2) = (x_3,x_3 y_3,z_3).
$$
Notice that the last blowing-up has its center on the curve $Y = \{x_2=y_2 = 0\}$.
\end{Example}
The Example from subsection~\ref{example-sanz} shows that the use of ramifications is unavoidable in order to obtain a complete resolution of singularities for vector fields.

Our present goal is to study the effect of a ramification on the Newton data. If $\Omega$ belongs to the class 
$\Nmaps^{i,\bfm}_{\Delta,C}$, it is obvious that 
$\overline{\Omega} = \Blx^1 \Omega$ belongs to the class
$\Nmaps^{i,\bfm}_{\overline{\Delta},\overline{C}}$, where 
$$\overline{\Delta} =
(\omega_1\Delta_1,\Delta_2)\quad\mbox{ and }\quad \overline{C} = \omega_1 C$$
Moreover, we have the following result:
\begin{Lemma}\label{Lemma-commutdiagramBx}
For each map $(f,g) \in \gG^i_{\Delta,C}$ there exists a unique map 
$(\overline{f},\overline{g}) \in  \gG^i_{\overline{\Delta},\overline{C}}$
which makes the following diagram commutative:
\begin{eqnarray*}\label{commutdiagram1}
\begin{CD}
\Omega @>{(f,g)\cdot}>> (f,g)\cdot \Omega\\
@V{\Blx^1}VV                                   @VV{\Blx^1}V \\
\overline{\Omega}  @>{(\overline{f},\overline{g})\cdot}>>    
(\overline{f},\overline{g})\cdot \overline{\Omega}
\end{CD}
\end{eqnarray*}
\end{Lemma}
\begin{proof}
We can explicitly define $\overline{f}(x,y) = f(\eps\,x^{\omega_1},y)$ and
$\overline{g}(x) = g(\eps\,x^{\omega_1})$.
\end{proof}
The next Lemma implies that the stability property is preserved by the
transformation $\Blx^1$.
\begin{Lemma}\label{lemma-stabpreservx}
Suppose that $\Omega$ is stable.  Then 
$\overline{\Omega} = \Blx^1\Omega$ is also a stable Newton data.
\end{Lemma}
\begin{proof}
First of all, let us prove that $\overline{\Omega}$ is edge stable.  For
this, assume
by absurd that there exists a map $(\overline{f},0) \in 
\gG^i_{\overline{\Delta}}$ such that
\begin{equation}\label{cannotoccuredge}
(\overline{f},0)\cdot \overline{\Omega} \notin 
\Nmaps^{i,\bfm}_{\overline{\Delta},\overline{C}}
\end{equation}
We must treat the following two cases:
\begin{enumerate}[(a)]
\item $\Delta_1 = 0$;
\item $\Delta_1 > 0$.
\end{enumerate}
In the case (a), it follows that $\overline{f}$ is a function of $y$ only.
Therefore, the Lemma~\ref{Lemma-commutdiagramBx} implies that 
there exists a map $(f,0) \in 
\gG_\Delta$ such that $(f,0) \cdot \Omega \notin 
\Nmaps^{i,\bfm}_{\Delta,C}$.  This contradicts the fact that $\Omega$ is
stable.

In the case (b), note that $C = 0$.  Write $\Delta_1 = p_1/q_1$ and
$\Delta_2 = p_2/q_2$, with
$\mathrm{mdc}(p_i,q_i) = 1$ ($i=1,2$).  Then, it follows from the Remark~\ref{remark-omegadeltaC}
that the weight-vector is given by $\bfomega = (q_1,0,p_1)$.  Let us split the discussion in two sub-cases:
\begin{enumerate}[({b}.1)]
\item $q_1 = 1$;
\item $q_1 \ge 2$.
\end{enumerate}
In the case $(b.1)$, the map $\Blx^1$ is just the identity map and we are done. 

In the case  $(b.2)$, we notice that the support of the Newton data $\Omega$ is such that
$$
\supp(\Omega|_\medge) \subset \bfm + n\, t \cdot (\frac{p_1}{q_1}, \frac{p_2}{q_2},-1)
\quad\mbox{ for } n \in \cN,$$
where $t$ is the least common multiple of $q_1$ and $q_2$.  In particular, $t
\ge 2$. Therefore, we get
$$
\supp(\Omega|_\medge) \cap \{ \bfv \in \cZ^3 \mid
v_3 = m_3 - 1 \} = \emptyset
$$
and the same property holds for
$\supp(\overline{\Omega}|_{\overline{\medge}})$.  The Remark~\ref{remark-edgestablefuture} now implies that $\overline{\Omega}$ is edge stable.

Now, suppose, by absurd, that there exists an map $(\widehat{f},\widehat{g}) \in 
\gG^i_{\overline{\Delta},\overline{C}} \setminus \gG^i_{\overline{\Delta}}$ such that 
$$
(\widehat{f},\widehat{g}) \cdot \overline{\Omega} \notin 
\Nmaps^{i,\bfm}_{\overline{\Delta},\overline{C}}
$$
Let us look at the action of the inverse map $(\widehat{f},\widehat{g})^{-1}$ 
on the restricted
Newton data $\overline{\Omega}|_\medge$.   
Notice that 
$$(\widehat{f},\widehat{g})^{-1}\cdot\overline{\Omega}|_\medge$$
is precisely the restriction of $\overline{\Omega}$ to the main face $\overline{\vF}$.

Looking at the points on the support of $\overline{\Omega}|_{\overline{\vF}}$ and
using Lemma~\ref{Lemma-stablenotresonant}, we can easily see
that $(\widehat{f},\widehat{g})^{-1}$ should necessarily be of the form
$$
(\widehat{f},\widehat{g})^{-1} = (f(x^{\omega_1},y), g(x^{\omega_1})) 
$$ 
for some $f \in \cR[x,y]$ and $g \in \cR[x]$. Using the
Lemma~\ref{Lemma-commutdiagramBx}, 
this implies that $\Omega$ is not stable.  Absurd.
\end{proof}
\subsection{The $x$-directional Projected Group and the Group $\gG_{\Delta,C}$}
Let us now introduce another subgroup of $\gG$, which will be mainly used for studying the effect of the translations on the $x$-directional blowing-up $\Blx\Omega$.

The {\em $x$-directional projected group} adapted
to $\Nmaps^{i,\bfm}_{\Delta,C}$ is defined as follows:
\begin{enumerate}[(i)]
\item If $\Delta_1 > 0$ then  $\PrG_x := \gG^1_{(0,\Delta_2),\infty}$
\item If $\Delta_1 = 0$ then $\PrG_x := \gG^1_{(0,\Delta_2),0}$
\end{enumerate}
In other words, if $\Delta_1 > 0$ then each $(f,g) \in \PrG_x$ has the form
$$
g = 0, \quad f = \xi y^{\Delta_2}
$$
where the constant $\xi \in \cR$ necessarily vanishes if $\Delta_2 \not\in
\cN$. If $\Delta_1 = 0$ then each $(f,g) \in \PrG_x$ has the form
$$
g(x) = \eta, \quad f(x,y) = a_0 + a_1 y + \cdots + a_b y^b
$$
where $b := \linteg \Delta_2 \rinteg$ and  $\eta, a_0,\ldots,a_b \in \cR$ are constants.
\begin{Lemma}\label{lemma-generaltransl}
Suppose that $\omega_1 = 1$.  Then, given a map $(f,g) \in
\PrG_x$ there exists a unique map $(f_\bfomega,g_\bfomega) \in 
\gG^1_{\Delta,C}$ which
makes the following diagram commutative:
\begin{eqnarray*}\label{commutdiagram}
\begin{CD}
\Nmaps @>{(f_\bfomega,g_\bfomega)\cdot}>> \Nmaps\\
@V{\Blx}VV                                   @VV{\Blx}V \\
\Nmaps  @>{(f,g)\cdot}>>                       \Nmaps
\end{CD}
\end{eqnarray*}
\end{Lemma}
\begin{proof}
Suppose first that $\Delta_1 = 0$ and $(f,g)$ is given by $(f,g) = (\xi y^k,\eta)$, for
some constants $\eta,\xi \in \cR$ and $0 \le k \le b$. The change of coordinates which is associated to
$(f,g)$ is
$$
\widetilde{y} = y + \eta, \quad \widetilde{z} = z + \xi y^{k}
$$
Now, if we apply the blowing-up map $(X,Y,Z) =
(x,x^{\omega_2}y, x^{\omega_3} z)$ on both sides of these equalities
and simplify common powers of $X$,
we get
$$
\widetilde{Y} = Y + \eta X^{\omega_2}, \quad
\widetilde{Z} = Z + \xi X^{\omega_3 - k\omega_2} Y
$$
(notice that $\omega_3 \ge k \omega_2$). Therefore
it suffices to define
$$(f_{\bfomega},g_{\bfomega})  := (\xi X^{\omega_3 -
k\omega_2} Y,\eta X^{\omega_2})$$

Using the same reasoning, we obtain
$f_{\bfomega}$ from an arbitrary polynomial $f$ by making the formal replacement
$$
y^k \rightarrow X^{\omega_3 - k\omega_2}Y^k
$$
and we obtain $g_{\bfomega}$ from $g$ by making the formal replacement
$$
1 \rightarrow X^{\omega_2}
$$
Suppose now that $\Delta_1 > 0$.  Here, the blowing-up map is given by
$(X,Y,Z) = (x,y,x^{\omega_3}z)$ an element $(f,g) \in \PrG_x$ corresponds to a change of
coordinates of the form
$$
\widetilde{z} = z + \xi y^{\Delta_2}
$$
where $\xi = 0$ if $\Delta_2 \notin\cN$.  The corresponding change of coordinates in the $(X,Y,Z)$ variables is given by 
$$
\widetilde{Z} = Z + \xi X^{\omega_3}Y^{\Delta_2}
$$
and, therefore, it suffices to get $(f_{\bfomega},g_{\bfomega}) = (\xi X^{\omega_3}Y^{\Delta_2},0)$.
\end{proof}
\begin{Remark}\label{remark-g=0gw=0}
We remark that if $(f,g) \in \PrG_x$ is such that $g = 0$ then the map 
$(f_{\bfomega},g_{\bfomega}) \in \gG_{\Delta,C}$ given by the Lemma is such that $g_{\bfomega} = 0$.  In particular, for $g = 0$, the map $(f_{\bfomega},g_{\bfomega})$ belongs to the subgroup $\gG^2_{\Delta,C}$.
\end{Remark}
\subsection{$x$-Directional Blowing-up (Case $\Delta_1 = 0$)}\label{subsect-xdirectblup}
In this subsection, we shall study the $x$-directional blowing-up of a
stable Newton data $\Omega$ in the case where $\Delta_1 = 0$. 

Our goal is to prove that the main invariant $\Inv$ is strictly smaller at each nonelementary point $\widetilde{p} \in \Phi^{-1}(\{p\}) \cap \NElem(\widetilde{\cM})$ which lies in the domain of the $x$-directional blowing-up.

\begin{figure}[htb]
\psfrag{v1}{\small $v_1$}
\psfrag{v2}{\small $v_2$}
\psfrag{v3}{\small $v_3$}
\psfrag{(0,0,2)}{\small $\bfm=(0,0,2)$}
\psfrag{t(0,0,2)}{\small $(0,0,2)$}
\psfrag{(1,-1,3)}{\small $(1,-1,3)$}
\psfrag{(0,-1,3)}{\small $\widetilde{m} = (0,-1,3)$}
\psfrag{(0,1,0)}{\small $(0,1,0)$}
\psfrag{Blx}{\small $\Blx$}
\begin{center}
\includegraphics[height=5cm]{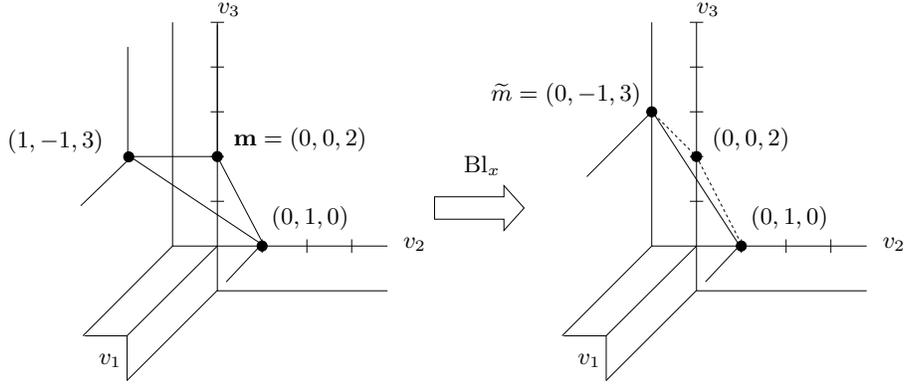}
\end{center}
\caption{The height of the main vertex increase after a $x$-directional blowing-up.}
\label{fig-exemplevirtheight}
\end{figure}

The following Example shows that the height $m_3$ of the main vertex $\bfm$ can increase after a $x$-directional blowing-up.  This is the main reason for introducing the concept of {\em virtual height} $\vheight$ in subsection~\ref{subsect-localinvar}.
\begin{Example}\label{example-virtheigthsaves}
Consider the vector field $\chi = (y^2 + xz^3) \frac{\partial}{\partial y} + z^3 \frac{\partial}{\partial z}$.  The associated Newton polyhedron is pictured in figure~\ref{fig-exemplevirtheight} (left).  The primary invariant is given by $(\vheight,m_2+1,m_3) = (2,1,2)$.  The $x$-directional blowing-up with
weight $\bfomega = (1,2,1)$ results into
the vector field
$$\widetilde{\chi} = (y^2 + z^3)\frac{\partial}{\partial y} + z^3 \frac{\partial}{\partial z}$$
(see figure~\ref{fig-exemplevirtheight} (right)).  Note that $\widetilde{m}_3 = 3 > 2 = m_3$.  However, the primary invariant associated to $\widetilde{\chi}$ is given by $(\widetilde{\vheight},\widetilde{m}_2+1,\widetilde{m}_3) = (2,0,3)$, which is lexicographically smaller than $(2,1,2)$.
\end{Example}

Up to a preliminary
transformation of type $\Blx^1$ (see the previous subsection), we can assume that the
weight-vector $\bfomega$ is such that $\omega_1 = 1$.

\begin{figure}[htb]
\psfrag{n}{\small $n$}
\psfrag{m}{\small $m$}
\psfrag{D}{\small $\Delta$}
\psfrag{Dt}{\small $\widetilde{\Delta}$}
\psfrag{h}{\small $h$}
\psfrag{ht}{\small $\widetilde{h}$}
\psfrag{1}{\small $1$}
\psfrag{0}{\small $0$}
\psfrag{v}{\small $v$}
\begin{center}
\includegraphics[height=5cm]{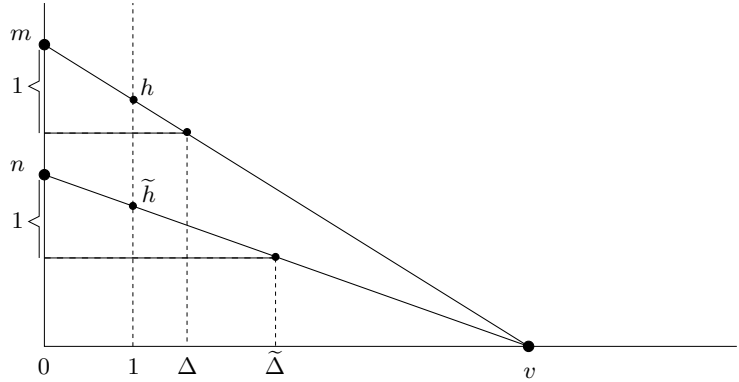}
\end{center}
\caption{Illustration of Lemma~\ref{lemma-simple}.}
\label{fig-lemmasimple}
\end{figure}

The following simple Lemma will be the key to understand the behavior of the virtual height under blowing-up and to prove of Proposition~\ref{prop-stablexblup}.
\begin{Lemma}\label{lemma-simple}
Let us consider three rational points in $\cQ^2$, with coordinates $(0,m)$,
$(0,n)$ and $(v,0)$ such that $v \ge 2$, $1 \le n < m$. Let
$\Delta := \frac{v}{m}$, and $\widetilde{\Delta} := \frac{v}{n}$
be the slope of the lines $\overline{\bfm,\bfv}$ and
$\overline{\bfn,\bfv}$, respectively.  Consider the rational numbers
$$
h := m - \frac{1}{\Delta}, \quad \widetilde{h} := n - \frac{1}{\widetilde{\Delta}}
$$
Then, one necessarily has $\widetilde{h} < h$.  Moreover, one the following situations occurs:
\begin{enumerate}[(i)]
\item $\widetilde{\Delta} < 1$, or
\item $h \ge n-1$.
\end{enumerate}
\end{Lemma}
\begin{proof}
The figure~\ref{fig-lemmasimple} illustrates the statement of the Lemma. The assertion that $\widetilde{h} < h$ is obvious. Suppose now that $v/n = \widetilde{\Delta} \ge 1$.  Then, $\Delta \ge n/m$,
and if we write $m = n + s$ (for some $s > 0$), we get
$$
h = m - \frac{1}{\Delta} \ge \frac{n^2 + n(s-1) - s}{n}
$$
Therefore, the quantity $h - n + 1 \ge s(1 - 1/n)$ is always
positive or zero.
\end{proof}

\begin{Lemma}\label{lemma-tildemnem}
Suppose that $\Delta_1 = 0$.  Suppose further that the $x$-directional blowing-up 
$\Blx \Omega$ is centered at a nonelementary point
$\widetilde{p} \in \NElem(\widetilde{\cM})$. If 
the main vertex $\overline{\bfm}$ of $\Blx \Omega$ is such that
$$
\overline{\bfm} \ne \bfm
$$
then one necessarily has $\Invp(\widetilde{\Omega}) <_\lex \Invp(\Omega)$, where
$\widetilde{\Omega} = \Stab \Blx \Omega$ is the stabilization of $\Blx \Omega$.
\end{Lemma}
\begin{proof}
First of all, let us suppose that the vertex $\overline{\bfm} = (0,\overline{m}_2,\overline{m}_3)$ is such that
$$
(\overline{m}_2,\overline{m}_3) <_\lex (m_2,m_3).
$$
Under such hypothesis, we split the discussion into three cases:
\begin{enumerate}[(a)]
\item $\overline{m}_2 = m_2 = 0$, $\overline{m}_3 < m_3$;
\item $\overline{m}_2 = m_2 = -1$, $\overline{m}_3 < m_3$;
\item $\overline{m}_2 = -1$, $m_2 = 0$.
\end{enumerate}
In the case (a), it is obvious that $\overline{\bfm}$ is also the main vertex of $\widetilde{\Omega}$, because no regular-nilpotent transition can occur in the passage from $\Blx \Omega$ to $\widetilde{\Omega}$.  Hence,
$$
\Invp(\widetilde{\Omega}) =
(\widetilde{\vheight},1, \widetilde{m}_3) <_\lex
(\vheight,1, m_3) = \Invp(\widetilde{\Omega}),
$$
because $\widetilde{\vheight} = \overline{m}_3 < m_3 = \vheight$.

To study the cases (b) and (c), we consider the main vertex $\widetilde{\bfm} = (0,\widetilde{m}_2,\widetilde{m}_3)$ of $\widetilde{\Omega}$.  It is obvious that either $\widetilde{\Omega}$ is in a regular configuration and $\widetilde{\bfm} = \overline{\bfm}$ or
$\widetilde{\Omega}$ is in a nilpotent configuration and there occurs a regular-nilpotent transition in the passage from $\Blx \Omega$ to $\widetilde{\Omega}$.  Notice that, in both cases, we have $\widetilde{m}_3 \le \overline{m}_3$.

To study the case (b), we observe that the main edge $\medge$ of $\Omega$ has the
form $\medge = \overline{\bfm,\bfv}$, for some point $\bfv = (0,v_2,v_3) \in
\supp(\Omega)$ such that $v_2 \ge 1$ and $v_3 < m_3$.
Two cases can occur:
\begin{enumerate}[({b}.i)]
\item $\overline{m}_3 \le v_3$;
\item $\overline{m}_3 > v_3$.
\end{enumerate}
In the case (b.i), we get
$$
\Invp(\widetilde{\Omega}) =
(\widetilde{\vheight},\widetilde{m}_2+1,\widetilde{m}_3) <_\lex
(\vheight,m_2+1, m_3) = \Invp(\widetilde{\Omega}),
$$
because $\widetilde{\vheight} \le  \widetilde{m}_3 \le v_3 < \vheight$.
To treat the case (b.ii), we define the numbers
$$
m := m_3 - v_3, \quad n := \overline{m}_3-v_3, \quad v := v_2 + 1
$$
By the construction, we know that
$$
1 \le n < m, \quad \mbox{ and }\quad v \ge 2
$$
and we can apply the Lemma~\ref{lemma-simple} to the points $(0,m)$, $(0,n)$
and $(v,0)$.  If we denote the displacement
vector of $\Blx \Omega$ by $\overline{\Delta} =
(0,\overline{\Delta}_2)$, and the associated virtual height by
$\overline{\vheight}$, it follows that
$$\overline{\vheight} \le \vheight$$
and one of the following situations occurs:
\begin{enumerate}[({b.ii}.1)]
\item $\overline{\Delta}_2 < 1$, or
\item $\vheight \ge \overline{m}_3$.
\end{enumerate}
In the case (b.ii.1), we know that $\Blx \Omega$ has a stable edge
(because $\gG^i_{\overline{\Delta}} = \{(0,0)\}$).  In particular, no regular-nilpotent transition can occur in the passage from $\Blx \Omega$ to $\widetilde{\Omega}$. Therefore,
$$
\Invp(\widetilde{\Omega}) =
(\overline{\vheight},0, \overline{m}_3) <_\lex
(\vheight,0, m_3) = \Invp(\Omega),
$$
In the case (b.ii.2), if no regular-nilpotent transition occurs in the passage from $\Blx \Omega$ to $\widetilde{\Omega}$ we obtain the same conclusion by the estimate
$\widetilde{\vheight} \le \overline{m}_3 \le \vheight$.  

On the other hand, if there occurs such regular-nilpotent transition then $\widetilde{\vheight} = \widetilde{m}_3 < \overline{m}_3 \le \vheight$ by the definition of nilpotent configurations.

Let us now study the case (c).  Here, keeping the notations of the previous
case, we 
consider the following possibilities:
\begin{enumerate}[({c}.i)]
\item $\overline{m}_3 \le v_3$;
\item $\overline{m}_3 > v_3$.
\end{enumerate}
The case (c.i) is treated exactly as the case (b.i).
To study the case (c.ii), it suffices to consider the points
$$
m := m_3 + \frac{1}{\Delta_2} - v_3, \quad n := \overline{m}_3 - v_3, \quad v
:= v_2 + 1
$$
Since the vertex $\overline{m}$ has the form $\overline{m} =
\pi_x(\bfn)$, for some $\bfn \in \vF \setminus \medge$, it
follows that
$$
0 > \< \bfomega, \bfm - \overline{\bfm} \> = \omega_2(-1) + \omega_3
(\overline{m}_3 - m_3 )
$$
and therefore (since $\Delta_2 = \omega_3/\omega_2$),
$n = \widetilde{m}_3  - v_3 < m_3 + \frac{1}{\Delta_2} - v_3 = m$.  Using the same
arguments of item (b.ii), we conclude that  $\Invp(\widetilde{\Omega}) <_\lex \Invp(\Omega)$.

It remains to study the case where $(\overline{m}_2,\overline{m}_3) >_\lex (m_2,m_3)$.  Here, one necessarily has the conditions 
$$\overline{m}_2 = 0 > -1 = m_2 \quad\mbox{ and }\quad
\overline{m}_3 < \vheight,$$
where the second inequality follows from the fact that $\widetilde{\bfm} = \pi_x (\bfn)$, for some point $\bfn \in \vF$.  Therefore, the Newton data $\Blx \Omega$ is already in a nilpotent configuration and $\overline{\bfm}$ is the main vertex of $\widetilde{\Omega}$.  These conditions imply that 
$\widetilde{\vheight} = \overline{m}_3 < \vheight$.
\end{proof}
\begin{Lemma}\label{lemma-tildemeqm}
Suppose that $\Delta_1 = 0$.  Suppose further that the $x$-directional blowing-up 
$\Blx \Omega$ is centered at a nonelementary point
$\widetilde{p} \in \NElem(\widetilde{\cM})$. If 
the main vertex of $\Blx \Omega$ coincides with that of $\Omega$, then
$$
\Inv(\widetilde{\Omega}) <_\lex \Inv(\Omega)
$$
where $\widetilde{\Omega} = \Stab \Blx \Omega$ is the stabilization of $\Blx\Omega$.
\end{Lemma}
\begin{proof}
We denote by
$$
\overline{\Omega} := (f,0,1,1,1) \cdot \Blx\, \Omega
$$
the analytic edge preparation which is associated $\Blx\Omega$ (see
Lemma~\ref{lemma-edgestable}).

If there exists a regular-nilpotent transition in this preparation (see
Lemma~\ref{lemma-regnilp})
then we are done.  In fact, it is clear that the virtual height 
$\overline{\vheight}$ associated to
$\overline{\Omega}$ is at most equal to $\vheight - 1$ because the main vertex
of  $\overline{\Omega}$ is $\overline{\bfm} = (0,0,m_3
-1)$ and $\overline{\vheight} = m_3 - 1 < m_3  = \vheight$ (a
regular-nilpotent transition can only occur if $m_2 = -1$ and the vertical displacement vector of $\Blx\Omega$ is equal to $\Delta = (0,1)$).

Therefore, let us assume that $\bfm$ is also the main vertex of
$\overline{\Omega}$. Let $\overline{\Delta} = (\overline{\Delta}_1,\overline{\Delta}_2)$ the main
displacement vector which is associated to $\overline{\Omega}$.  Then, by definition
$$
\Invp(\widetilde{\Omega}) = (\overline{\vheight},m_2 + 1,m_3), \quad
\Invs(\widetilde{\Omega}) = (\#\iota_{\widetilde{p}} - 1,M \overline{\Delta}_1, M \max \{0,\overline{\Delta}_2\})
$$
Since $\#\iota_{\widetilde{p}} = \#\iota_p = i$, it clearly suffices to prove the following claim:
$$
\mathrm{Claim: }\quad  \overline{\Delta}_1 = \Delta_1 = 0\quad\mbox{ and }\quad \overline{\Delta}_2 < \Delta_2
$$
To prove such claim, suppose by absurd that either $\overline{\Delta}_1 > 0$ or
$\overline{\Delta}_2 \ge \Delta_2$.  Then, if we write the Taylor series of
the map $f(x,y)$ as
$$
f(x,y) = \sum_{i+j\ge 1}{f_{ij}x^i y^j}
$$
It follows that the polynomial truncation $f_t := \sum_{i =
0}^{b}{f_{i0} x^i}$ (with $b := \linteg \Delta_2 \rinteg$) is such that the
Newton data
$$
\Omega_t := (f_t,0) \cdot \widetilde{\Omega}
$$
has a displacement vector $\Delta_t \ge_\lex \Delta$.  Since the map $(f_t,0)$
belongs to the $x$-directional projected group $\PrG_x$ , it follows from
Lemma~\ref{lemma-generaltransl} that there exists a unique map
$(f_\bfomega,0) \in \gG^i_{\Delta,C}$ such that
$$
\Omega_t = \Blx\ (f_\bfomega,0)\cdot \Omega
$$
We conclude from Lemma~\ref{lemma-directx} that the Newton data
$(f_\bfomega,0)\cdot \Omega$
does not belong to $\Nmaps^{i,\bfm}_{\Delta,C}$.  This contradicts the
hypothesis that $\Omega$ is a stable Newton data.  Absurd.
\end{proof}
\begin{Proposition}\label{prop-reducesinvariant}
Suppose that $\Delta_1 = 0$ and that the $x$-directional blowing-up 
$\Blx \Omega$ is centered at a nonelementary point
$\widetilde{p} \in \NElem(\widetilde{\cM})$.  Then,
the Newton data $\widetilde{\Omega} = \Stab \Blx \Omega$ is such that 
$$
\Inv(\widetilde{\Omega}) <_\lex \Inv(\Omega)
$$
\end{Proposition}
\begin{proof}
This is an immediate consequence of Lemmas~\ref{lemma-tildemnem} and \ref{lemma-tildemeqm}.
\end{proof}
\subsection{Effect of Translations in the $x$-Directional Blowing-up (Case $\Delta_1 = 0$)}
In this subsection, we study the effect of the coordinate translations to the Newton data $\Blx \Omega$.  As we have said in Remark~\ref{Remark-translations}, the notions of stability and edge stability have been introduced precisely to take these effects into account.
\begin{Proposition}\label{prop-reducesinvarianttransl}
Suppose that $\Delta_1 = 0$ and that the Newton data $(\xi,\eta) \cdot \Blx \Omega$ is centered at a nonelementary point
$\widetilde{p} \in \NElem(\widetilde{\cM})$, where $(\xi,\eta) \in \gG^\transl_x$ is a $x$-directional translation.  Then, 
$$
\Inv(\widetilde{\Omega}_{\xi,\eta}) <_\lex \Inv(\Omega)
$$
where $\widetilde{\Omega}_{\xi,\eta} = \Stab (\xi,\eta) \cdot \Blx \Omega$ is the stabilization of $(\xi,\eta) \cdot \Blx \Omega$.
\end{Proposition}
\begin{proof}
Defining $i = \#\iota_p$, we split the proof into two cases:
\begin{enumerate}[(a)]
\item $i = 2$ and $\eta \ne 0$;
\item $i = 1$ or $\eta = 0$.
\end{enumerate}
To treat the case (a), we observe that the $x$-projected face $\widetilde{\vF} :=
\{\pi_x(\bfv) \mid \bfv \in \vF\}$ is equal to the intersection
$\supp(\widetilde{\Omega})\; \cap\; \left( \{0\}\times \cZ^2 \right)$.

In particular, if we denote the main vertex of $\widetilde{\Omega}_{\xi,\eta}$ by
$\widetilde{\bfm}$,
it is immediate to see that
$$
\widetilde{\bfm} \le_\lex \overline{\bfm}
$$
(where $\overline{\bfm}$ is the main vertex of $\Blx\Omega$).  It
follows that
$$
\Invp(\widetilde{\Omega}_{\xi,\eta}) =
(\widetilde{\vheight},\widetilde{m}_2 + 1, \widetilde{m}_3) \le_\lex
(\vheight,m_2 + 1, m_3) = \Invp(\Omega),
$$
because $\widetilde{\vheight} \le m_3 =
\vheight$.  On the other hand, if we have an equality of the
primary multiplicity $\Invp(\cdot)$, then
$$
\Invs(\widetilde{\Omega}_{\xi,\eta}) =
(0,\lambda\ \overline{\Delta}_1, \lambda\ \max\{0,\overline\Delta_2\}) < _\lex
(1,\lambda\ \Delta_1, \lambda\ \max\{0,\Delta_2\}) = \Invs(\Omega),
$$
(because the assumption $\eta \ne 0$ implies that the translated Newton data
$\widetilde{\Omega}_{\xi,\eta}$ is centered at a point $\widetilde{p}$ such that 
$\#\iota_{\widetilde{p}} = 1 < 2 = \#\iota_p$).

We treat now the case (b).  It follows 
from the Lemma~\ref{lemma-generaltransl}, that there exists a unique
map $(f,g) \in \gG^i_{\Delta,C}$ such that
$$
\widetilde{\Omega}_{\xi,\eta} = \Stab \Blx\, (f,g)\cdot \Omega
$$
More explicitly, $(f,g)$ is given by $(\xi x^{\omega_3},\eta
x^{\omega_2})$, where $\eta = 0$ if $i = 2$.

Since $\Omega$ is a stable Newton data, the  Newton data $\Omega_{f,g}
:=(f,g)\cdot \Omega$ is also stable.  Moreover,
$$
\Inv(\Omega_{f,g}) = \Inv(\Omega).
$$
Thus, the result is a direct consequence of applying the
Proposition~\ref{prop-reducesinvariant} to $\Omega_{f,g}$ instead of $\Omega$.
\end{proof}
\subsection{$x$-Directional Blowing-up (Case $\Delta_1 > 0$)}\label{subsect-xdirectDelta=0}
In this subsection, we keep the assumption that $\omega_1 = 1$.  
Recall that this condition can always be obtained up to a preliminary transformation of type $\Blx^1$.
\begin{Lemma}\label{lemma-delta>01}
Suppose that $\Delta_1 > 0$.  Suppose further that the $x$-directional blowing-up 
$\Blx \Omega$ is centered at a nonelementary point
$\widetilde{p} \in \NElem(\widetilde{\cM})$. If 
$$
\Delta_2 \le 0
$$
then $
\Invp(\widetilde{\Omega}) <_\lex \Invp(\Omega)$, where $\widetilde{\Omega} = \Stab \Blx\Omega$ is the stabilization of $\Blx \Omega$.
\end{Lemma}
\begin{proof}
Under the hypothesis of the enunciate, we know that the main edge of $\Omega$ is given by
$\medge = \overline{\bfm,\bfv}$, where $\bfv \in \supp(\Omega)$ is such that
$$
v_2 \le m_2, \quad  v_3 < m_3
$$
Using Lemma~\ref{lemma-directx}, we conclude that the main vertex of
$\Blx\Omega$ is given either by $\overline{\bfm} = (0,v_2,v_3)$ (if $\Blx \Omega$ is in a regular configuration) or by $\overline{\bfm} = (0,0,\overline{m}_3)$ for some $\overline{m}_3 < v_3$ 
(if $\Blx \Omega$ is in a nilpotent configuration).  Therefore, after applying the stabilization map $\Stab$ to $\Blx \Omega$ we get
$$\Invp(\widetilde{\Omega}) = (\widetilde{\vheight},\widetilde{m}_2 + 1, \widetilde{m}_3) <_\lex
(\vheight,m_2 + 1, m_3)$$
because $\widetilde{\vheight} \le v_3 < m_3 = \vheight$.
\end{proof}
\begin{Lemma}\label{lemma-delta>02}
Suppose that $\Delta_1 > 0$.  Suppose further that the $x$-directional blowing-up 
$\Blx \Omega$ is centered at a nonelementary point
$\widetilde{p} \in \NElem(\widetilde{\cM})$. If 
$$
\Delta_2 > 0
$$
then the stabilization $\widetilde{\Omega} = \Stab \Blx\Omega$ of $\Blx \Omega$ is such that 
$$
\Inv(\widetilde{\Omega}) <_\lex  \Inv(\Omega).
$$
\end{Lemma}
\begin{proof}
Under the hypothesis that $\Delta_2 > 0$, we consider separately the
following situations:
\begin{enumerate}[(a)]
\item $\Blx \Omega$ is in a hidden nilpotent configuration, or
\item $\Blx \Omega$ is not in a hidden nilpotent configuration.
\end{enumerate}
In the case (a), let
$$\overline{\Omega} = (f,0)\cdot  \Blx \Omega$$
be the edge preparation of $\Blx \Omega$.  Then,
the virtual height $\overline{\vheight}$ associated to $\overline{\Omega}$
is strictly smaller than $\vheight = m_3$.

Consider now the case $(b)$.  The main vertex of both
$\Blx \Omega$ and $\overline{\Omega}$ is $\bfm$.
Moreover, the displacement vector $\Delta^\prime$ of $\Blx \Omega$ is given by
$$
\Delta_1^\prime = 0, \quad \Delta_2^\prime = \Delta_2.
$$
The argument is now similar to the one used in the proof of
Lemma~\ref{lemma-tildemeqm}.  Let
$$\overline{\Omega} = (f,0)\cdot   \Blx \Omega$$
be the edge preparation of $\Blx \Omega$, and let $\overline{\Delta} =
(\overline{\Delta}_1, \overline{\Delta}_2)$ be the displacement vector
associated to $\overline{\Omega}$. We claim that
$$
\overline{\Delta}_1 = 0, \quad \overline{\Delta}_2 = \Delta_2
$$
Indeed, suppose the contrary.  Then, if we consider the polynomial truncation of $f$
given by $f_t = \xi y^{\Delta_2}$ (with $\xi \in \cR$ equals to zero if $\Delta_2
\notin \cN$), it follows that
$$\Omega_t := (f_t,0)\cdot\Blx \Omega$$
has a displacement vector $\Delta_t >_\lex \Delta$.  Since $(f_t,0)$ belongs to the $x$-projected group
$\PrG_x$, we can use Lemma~\ref{lemma-generaltransl} to conclude that
that there exists a map $(f_\bfomega,0) \in \gG^i_{\Delta,C}$ such that
$$
(f_\bfomega,0)\cdot \Omega
$$
has a vertical displacement vector which (lexicographically) strictly greater than $\Delta$.
But this contradicts the hypothesis that $\Omega$ is stable.  The claim is
proved.

Using such claim, we easily conclude that $\Invp(\widetilde{\Omega}) \le_\lex
\Invp(\Omega)$ and also that $\Invs(\widetilde{\Omega})<_\lex \Invs(\Omega)$.
\end{proof}
\begin{Proposition}\label{prop-reducesinvariantD>0}
Suppose that $\Delta_1 > 0$.  Suppose further that the $x$-directional blowing-up 
$\Blx \Omega$ is centered at a nonelementary point
$\widetilde{p} \in \NElem(\widetilde{\cM})$. Then $\Inv(\widetilde{\Omega}) <_\lex  \Inv(\Omega)$, where $\widetilde{\Omega} = \Stab \Blx\Omega$ is the stabilization of $\Blx \Omega$.
\end{Proposition}
\begin{proof}
It suffices to use the Lemmas~\ref{lemma-delta>01} and \ref{lemma-delta>02}.
\end{proof}
\subsection{Effect of Translations in the $x$-Directional Blowing-up (Case $\Delta_1 > 0$)}
Let us now study the effect of the translations in the $x$-directional blowing-up chart for the case where $\Delta_1 > 0$.
\begin{Proposition}\label{prop-reducesinvariantD>0transl}
Suppose that $\Delta_1 > 0$ and that the Newton data $(\xi,0) \cdot \Blx \Omega$ is centered at a nonelementary point
$\widetilde{p} \in \NElem(\widetilde{\cM})$, where $(\xi,0) \in \gG^\transl_x$ is a $x$-directional translation.  Then, 
$$
\Inv(\widetilde{\Omega}_{\xi,\eta}) <_\lex \Inv(\Omega)
$$
where $\widetilde{\Omega}_{\xi,\eta} = \Stab (\xi,0) \cdot \Blx \Omega$ is the stabilization of $(\xi,\eta) \cdot \Blx \Omega$.
\end{Proposition}
\begin{proof} If $\xi = 0$, this follows from
Lemma~\ref{prop-reducesinvariantD>0}.  

Let us assume that $\xi \ne 0$. We split the proof into three cases:
\begin{enumerate}[(a)]
\item $\Delta_2 < 0$;
\item $\Delta_2 > 0$;
\item $\Delta_2 = 0$.
\end{enumerate}
In the case (a), we know that $\bfm = (0,0,m_3)$.  Moreover, the main vertex
of 
$\Blx\Omega$ is given
by $\overline{\bfm} = (0,-1,\overline{m}_3)$, for some $\overline{m}_3 < m_3$.
It follows that the Newton data $(\xi,0) \cdot \Blx \Omega$ is in a final situation and therefore it is centered at an elementary point $\widetilde{p} \in \Elem(\widetilde{\cM})$. This contradicts the hypothesis on the enunciate.

In the case (b), we have $\overline{\bfm} = \bfm$.
Let $\overline{\chi}$ be the vector field associated to $\Blx \Omega$.
Then,
$$
\overline{\chi}|_{\bfm} = y^{m_2}z^{m_3}\left( \alpha x
\frac{\partial}{\partial x} + \beta y
\frac{\partial}{\partial y} + \gamma z
\frac{\partial}{\partial z}\right), \quad (\alpha,\beta,\gamma) \in \cR^3
\setminus \{(0,0,0)\}
$$
(where $\alpha = \gamma = 0$ if $m_2 = -1$).  Since
$$\pi_x(\vF) \; \cap \; \left(\{0\}\times \{m_2\} \times \cR\right) = \bfm$$
and $m_2 \in \{-1,0\}$, it is clear that after the translation $\widetilde{z}
= z + \xi$ we get a Newton data $(\xi,0)\cdot \Blx \Omega$ which is in a final
situation. Again, this contradicts the hypothesis on the enunciate.

It remains to study the case (c).  Here, we observe that the translation map
$(\xi,0)$ belongs to the $x$-projected group $\PrG_x$. It follows from the
Lemma~\ref{lemma-generaltransl} that there exists a unique
map $(f,g) \in \gG^i_{\Delta,C}$ such that
$$
\widetilde{\Omega}_{\xi} = \Stab \Blx\, (f,g)\cdot \Omega
$$
More explicitly, $(f,g)$ is given by $(\xi x^{\omega_3},0)$.

Since $\Omega$ is a stable Newton data, the same holds for the Newton data $\Omega_{f,g} :=
(f,g)\cdot \Omega$.  Moreover, 
$$
\Inv(\Omega_{f,g}) = \Inv(\Omega).
$$
Thus, the result is a direct consequence of applying the
Proposition~\ref{prop-reducesinvariantD>0} to $\Omega_{f,g}$ instead of $\Omega$.
\end{proof}
We are finally ready to give the proof of 
Proposition~\ref{prop-stablexblup}.
\begin{proof}({\bf of Proposition~\ref{prop-stablexblup}})  In the case $\Delta_1 = 0$, we apply 
Proposition~\ref{prop-reducesinvarianttransl}.
In the case $\Delta_1 > 0$, we apply Proposition~\ref{prop-reducesinvariantD>0transl}.
\end{proof}
\subsection{$y$-Directional Blowing-up}
Let $(\cM,\Ax)$ be a controlled singularly foliated manifold and let
$p \in \vD \cap A$ be a divisor point in $\NElem(\cM)$.  

Let $\Omega \in \Nmaps^{i,\bfm}_{\Delta,C}$ be a stable Newton data, associated to some  adapted chart at $p$ (and some
local generator $\chi$ of $L$).  
In this subsection, we assume that the corresponding 
weight-vector $\bfomega = (\omega_1,\omega_2,\omega_3)$ is such that $\omega_2 > 0$.

The {\em $y$-directional
translation group} is the group $\gG^\transl_y := \gG^1_{(0,0),\infty}$.
In other words, an element of $\gG^\transl_y$ corresponds to a translation
$$
\widetilde{z} = z + \xi
$$
for some constant $\xi \in \cR$.  We denote such element simply by $(\xi,0)$.
\begin{Proposition}\label{prop-stableyblup}
Given a stable Newton data $\Omega$, let
$\overline{\Omega} = \Bly \Omega$ be its $y$-directional blowing-up.
Then, for each $(\xi,0) \in \gG^\transl_y$, either the translated Newton data
$$
\overline{\Omega}_\xi := (\xi,0) \cdot \overline{\Omega}
$$
is centered at an elementary point $\widetilde{p} \in \Elem(\widetilde{\cM})$ or 
$$
\Inv(\widetilde{\Omega}_\xi) <_\lex \Inv(\Omega)
$$
where $\widetilde{\Omega}_\xi = \Stab \overline{\Omega}_\xi$ is the stabilization of $\overline{\Omega}_\xi$.
\end{Proposition}
The proof of the Proposition will be given at the end of subsection~\ref{subsect-ydirectDeltas1}.
First of all, we enunciate the following analog of Lemma~\ref{lemma-directx}:
\begin{Lemma}\label{lemma-directy}
Let $\overline{\Omega} := \Bly \Omega$ be the $y$-directional blowing-up of
$\Omega$.  Then, there exists a
bijective correspondence
\begin{eqnarray*}
\supp(\Omega) \; \cap\;   \vF &\longrightarrow & \supp(\overline{\Omega}) \; \cap\;  \{0\}
\times  \cZ^2
\\
\bfv = (v_1,v_2,v_3) &\longmapsto& \pi_y(\bfv) = (0,v_1,v_3)
\end{eqnarray*}
such that the corresponding Newton maps $\overline{\Theta}$ and $\Theta$ satisfy $\overline{\Theta}[\pi_y(\bfv)] = I\ M_y \Theta[\bfv]$.
\end{Lemma}
\begin{proof}
This is an immediate consequence of the definition of $\Bly$.
\end{proof}
We remark that the Newton data $\Bly \Omega = ((\overline{x},\overline{y},\overline{z}),\overline{\iota},\overline{\Theta})$ is always such that $\#\overline{\iota} = 2$. 

As in the discussion of the $x$-directional blowing-up, we can decompose
the map $\Bly$ in two maps
$\Bly^2$ and $\Bly^1$,  which are respectively associated to the singular changes of coordinates
$$
x = \overline{x}, \;\; y = \eps\, \overline{y}^{\omega_2},
\;\; 
z = \overline{z} \quad \mbox{ and }\quad 
x = \overline{x}^{\omega_1}\overline{y}, \;\; y = \overline{x},
\;\;
z = \overline{x}^{\omega_3}\overline{z}
$$
followed by a division by $\overline{x}^\mu$ (for $\eps \in \{-1,1\}$). The first change of coordinates corresponds to a {\em ramification} and the second change of coordinates can always be written as a composition a finite sequence of homogeneous blowing-ups.

The following Lemma is an analogous to Lemma~\ref{Lemma-commutdiagramBx}.
\begin{Lemma}\label{lemma-stabpreservy}
Suppose that $\Omega$ is edge stable.  Then 
$\overline{\Omega} = \Bly^1\Omega$ is also an edge stable Newton data.
\end{Lemma}
\begin{proof}
The proof is very similar to the proof of Lemma~\ref{lemma-stabpreservx}. We omit the details for shortness.
\end{proof}
Using this Lemma, we can assume that $\omega_2 = 1$ without loss of generality in our results.  
\begin{Lemma}\label{lemma-generaltransly}
Suppose that $\omega_2 = 1$.  Then, given a translation map $(\xi,0) \in
\gG^\transl_y$ there exists a unique map $(f_{\bfomega},0) \in \gG_\Delta$ which
makes the following diagram commutative:
\begin{eqnarray*}\label{commutdiagramy}
\begin{CD}
\Nmaps @>{(f_{\bfomega},g_{\bfomega})\cdot}>> \Nmaps\\
@V{\Bly}VV                                   @VV{\Bly}V \\
\Nmaps  @>{(f,g)\cdot}>>                       \Nmaps
\end{CD}
\end{eqnarray*}
\end{Lemma}
\begin{proof}
The proof is analogous to the proof of Lemma~\ref{lemma-generaltransl}.
Consider the change of coordinates
$$
\widetilde{z} = z + \xi
$$
and apply the map $(X,Y,Z) = (x^{\omega_1}y, \eps x, x^{\omega_3} z)$ on
both sides of the equality.  Cancelling out common powers of $x$, we get
$$
\widetilde{Z} = Z + \bar\xi Y^{\omega_3},\quad \mbox{ for }\bar\xi = \eps^{\omega_3} \xi
$$
which corresponds to the map $(\bar\xi Y^{\omega_3}, 0)$ in the group $\gG_\Delta$.
\end{proof}
\subsection{Effect of Translations in the $y$-Directional Blowing-up}\label{subsect-ydirectDeltas1}
Let us keep the notations of the previous subsection. Recall that we can assume, 
without loss of generality, that $\omega_2 = 1$.

The proof of Proposition~\ref{prop-stableyblup} will be given by considering separately the cases $\Delta_2 > 1$, $\Delta_2 = 1$ and $\Delta_2 < 1$.
\begin{proof}({\bf of Proposition~\ref{prop-stableyblup} when $\Delta_2 > 1$})
Suppose initially that $\xi = 0$.  Write the main edge of $\Omega$
as $\medge = \overline{\bfm,\bfv}$, where $\bfv = (0,v_2,v_3) \in \supp(\Omega)$ is such that
$v_2 > m_2$ and $v_3 < m_3$.

Using Lemma~\ref{lemma-directy}, we conclude that the main vertex of
$\widetilde{\Omega} := \widetilde{\Omega}_0$
is given by $\widetilde{\bfm} = (0,0,v_3)$.  Therefore,
$$
\Invp(\widetilde{\Omega}) = (\widetilde{\vheight},1,v_3) <_\lex (\vheight, m_2
+ 1, m_3) = \Invp(\Omega)
$$
because $\widetilde{\vheight} =
v_3 < m_3 = \vheight$ (the last equality follows from the assumption that
$\Delta_2 > 1$).

Suppose now that $\xi \ne 0$.  We claim that the main vertex
of $\widetilde{\Omega}_\xi$ has the form $\overline{\bfm} =
(0,0,\overline{m}_3)$ for some $\overline{m}_3 \le m_3-1$.

Indeed, if this is not the case then necessarily $\overline{\bfm} = \bfm$ (by
Lemma~\ref{lemma-directy}).  Using Lemma~\ref{lemma-generaltransly}, we
conclude that there exists a map $(f,0) \in \gG_\Delta$ such that
$$
(f,0)\cdot \Omega
$$
has a vertical displacement vector which is strictly greater (lexicographically) than
$\Delta$.  But this contradicts the hypothesis that $\Omega$ is stable (and, in particular, edge stable). The
claim is proved.

Using the claim, we conclude again that $\Invp(\widetilde{\Omega}_\xi) <_\lex
\Invp(\Omega)$.  This proves the Lemma.
\end{proof}
Let us now consider the case $\Delta_2 = 1$.
\begin{proof}({\bf of Proposition~\ref{prop-stableyblup} when $\Delta_2 = 1$})
Let $\bfm = (0,m_2,m_3)$ be the main vertex of $\Omega$ and $\Delta$ be the vertical
displacement vector.  We split the proof into two cases:
\begin{enumerate}[(a)]
\item $m_2 = -1$.
\item $m_2 = 0$.
\end{enumerate}
In the case (a), the primary invariant is given by
$$
\Invp(\Omega) = (m_3,0,m_3)
$$
We claim that $\widetilde{\Omega}_\xi$ has a main vertex $\overline{\bfm} = (0,0,\overline{m}_3)$
such that $\overline{m}_3 \le m_3 - 2$.

Indeed, if there exists a $\xi \in \cR$ such that $\widetilde{\Omega}_\xi$ has
a main vertex with height $\overline{m}_3 = m_3 - 1$ then it follows from the
hypothesis that $\Omega$ should necessarily be in a hidden nilpotent configuration. 
Using Lemma~\ref{lemma-generaltransly}, this
conclusion contradicts the assumption that $\Omega$ is stable and $m_2 = -1$.

As a consequence of the claim, $\Invp(\widetilde{\Omega}_\xi) <_\lex
\Invp(\Omega)$ because $\widetilde{\vheight} = \widetilde{m}_3 \le m_3 - 2 < \vheight$.

The case (b) can be treated as in the proof of the case $\Delta_2 > 1$.
\end{proof}
To conclude the proof of Proposition~\ref{prop-stableyblup}, we treat the case $\Delta_2 < 1$.
\begin{proof}({\bf of Proposition~\ref{prop-stableyblup} when $\Delta_2 < 1$})
We consider separately the following cases:
\begin{enumerate}[(a)]
\item $m_2 = 0$;
\item $m_2 = -1$;
\end{enumerate}
In case (a), we can use exactly the same argument used in the proof of
the case $\Delta_2 > 1$ to conclude that
the main vertex $\overline{\bfm}$ of $\widetilde{\Omega}_\xi$ is such that
$\overline{m}_3 \le m_3 - 1$.  Therefore,
$$
\Invp(\widetilde{\Omega}) = (\overline{\vheight},1,\overline{m}_3) <_\lex (\vheight, 0, m_3) = \Invp(\Omega)
$$
because $\overline{\vheight} = \overline{m}_3 < m_3 = \vheight$.

Let us treat the case (b).  Suppose initially that $\xi = 0$.  Write the main edge of $\Omega$
as $\medge = \overline{\bfm,\bfv}$, where $\bfv = (0,v_2,v_3)$ is such that
$v_2 > 1$ (the strict inequality
follows from the fact that $\Omega$ is not in a nilpotent
configuration). Therefore $v_3 \le  \lfloor m_3 - 1/\Delta_2 \rfloor$.
Using Lemma~\ref{lemma-directy}, we conclude that the main vertex of
$\widetilde{\Omega} := \widetilde{\Omega}_0$
is given by $\widetilde{\bfm} = (0,0,v_3)$.  Therefore
$$
\Invp(\widetilde{\Omega}) = (\widetilde{\vheight},1,v_3) <_\lex (\vheight, 0 , m_3) = \Invp(\Omega)
$$
because, by the definition of the virtual height, $\widetilde{\vheight} =
v_3 \le \lfloor m_3 - 1/\Delta_2 \rfloor < \lfloor m_3 - 1/\Delta_2 + 1
\rfloor = \vheight$.

Suppose now that $\xi \ne 0$. Let $\chi$ be the vector field which is
associated to $\Omega$.  Then, its restriction to the main edge $\medge$ can
be written as
$$
\chi|_\medge = F(y,z) x\frac{\partial}{\partial x} + G(y,z)
\frac{\partial}{\partial y} + H(y,z) \frac{\partial}{\partial z}
$$
where $F,G,H$ are $(\omega_2,\omega_3)$-quasihomogeneous polynomials of degree
$\mu$, $\mu + \omega_2$ and $\mu + \omega_3$, respectively.
The hypothesis $m_2 = -1$ implies that $G(0,z) = \beta z^{m_3}$, for some nonzero
constant $\beta\in \cR$.

Using Lemma~\ref{lemma-directy}, we see that the vector field
$\widetilde{\chi}$ which is associated to $\Bly \Omega$ (before the
translation by $\xi$) is such that its restriction to $\pi_y(\medge)$ has the form
\begin{equation}\label{expreychi}
\frac{1}{\omega_2}G(1,z)
x\frac{\partial}{\partial x} + \left( F(1,z)y - \frac{\omega_1}{\omega_2}
G(1,z)y \right) \frac{\partial}{\partial y} +
\left( H(1,z) - \frac{\omega_3}{\omega_2}
G(1,z)z \right) \frac{\partial}{\partial z}
\end{equation}
Let us consider the polynomial
$$\widetilde{G}(z) := \frac{1}{\omega_2}G(1,z)$$
and denote by $\overline{\vheight}$ the virtual height 
associated to $\widetilde{\Omega}_\xi$.

It follows that $\xi$ is a root of multiplicity $\ge \overline{\vheight}$ of
$\widetilde{G}(z)$.  On the other hand, Corollary~\ref{corollary-pnmult}
of the Appendix~B implies that 
$$
\mu_\xi(\widetilde{G}) \le \lfloor m_3 - \frac{1}{\Delta_2} \rfloor < \lfloor m_3 -
\frac{1}{\Delta_2} + 1 \rfloor = \vheight
$$
This concludes the proof of the Proposition~\ref{prop-stableyblup}.
\end{proof}
\subsection{The $z$-directional Blowing-up}
Let $(\cM,\Ax)$ be a controlled singularly foliated manifold and let
$p \in \vD \cap A$ be a divisor point in $\NElem(\cM)$.  

Let $\Omega \in \Nmaps^{i,\bfm}_{\Delta,C}$ be a stable Newton data at $p$.  
From our constructions, it is clear that 
the associated weight-vector $\bfomega = (\omega_1,\omega_2,\omega_3)$ is always such that $\omega_3 > 0$.
\begin{Lemma}\label{lemma-directz}
Let $\overline{\Omega} := \Blz \Omega$ be the $z$-directional blowing-up of
$\Omega$.  Then, there exists a
bijective correspondence
\begin{eqnarray*}
\supp(\Omega) \; \cap\;   \vF  &\longrightarrow & \supp(\overline{\Omega}) \; \cap\;  \{0\}
\times  \cZ^2
\\
\bfv = (v_1,v_2,v_3) &\longmapsto& \pi_z(\bfv) = (0,v_1,v_2)
\end{eqnarray*}
such that the corresponding Newton maps $\overline{\Theta}$ and $\Theta$ satisfy
$\overline{\Theta}[\pi_z(\bfv)] = J\ M_z \Theta[\bfv]$.
\end{Lemma}
\begin{proof}
This is an immediate consequence of the definition of $\Blz$.
\end{proof}
\begin{Proposition}\label{prop-stablezblup}
Given a stable Newton data $\Omega$, its $z$-directional blowing-up 
$\overline{\Omega} = \Blz \Omega$ is always centered at an elementary point $\widetilde{p} \in \Elem(\widetilde{\cM})$.
\end{Proposition}
\begin{proof}
Using Proposition~\ref{prop-finalsituation}, it is sufficient to prove that $\Blz \Omega$ is in a final situation. 

If we write main vertex of
$\Omega$ as $\bfm = (0,m_2,m_3)$,  it
follows from Lemma~\ref{lemma-directz} that $\widetilde{\Omega}$ contains the
point
$$
\widetilde{\bfm} := \pi_z(\bfm) = (0,0,m_2)
$$
on its support.  It is clear that such point is necessarily the new main vertex of
$\widetilde{\Omega}$.  Moreover, since $m_2 \in \{-1,0\}$, the Newton data
$\widetilde{\Omega}$ is in a final situation.
\end{proof}
\subsection{Proof of the Local Resolution of Singularities}\label{subsect-localresoltheorem}
Let us prove the Theorem~\ref{theorem-localresolution}.
Consider the local blowing-up $\Phi: \widetilde{\cM} \rightarrow \cM \cap U$
defined in the enunciate of the Theorem, and write
$$\widetilde{\cM} = (\widetilde{M},\widetilde{\mainlist},
\widetilde{\vD},\widetilde{\linef})$$
Let $V^x$, $V^y$ and $V^z$ denote the domain of the $x$-directional,
$y$-directional and $z$-directional charts, respectively.

First of all, we define an open subset $\widetilde{A} \subset \widetilde{M}$ and
an analytic line field $\widetilde{\vZ}$ on $\widetilde{A}$ by taking
$$
\widetilde{A} = \Phi^{-1}(A) \cap (V^x \cup V^y) \quad \mbox{ and }
\quad \widetilde{\vZ} = \Phi_*(\vZ) |_{\widetilde{A}}
$$
where $\Phi_*(\vZ)$ denotes the pull-back of $\vZ$.  We remark that 
\begin{enumerate}[(i)]
\item The Proposition~\ref{prop-stablezblup} implies that $\widetilde{A}$ is
      an open neighborhood of $\NElem(\widetilde{\cM})$.
\item On the domain $V^x \cup V^y$, the pull-backed foliation $\Phi_*(\vZ)$ is 
everywhere regular.  Hence $\Ze(\widetilde{\vZ}) = \emptyset$.
\end{enumerate}
It follows that the pair $\widetilde{\Ax} = (\widetilde{A},\widetilde{\vZ})$ satisfies
all the conditions in the definition~\ref{def-axis}.  Hence, $\widetilde{\Ax}$
is an axis for $\widetilde{\cM}$.

Now, let $\widetilde{p} \in \Phi^{-1}(\{p\})$ be a point belonging to
$\NElem(\widetilde{\cM})$.
Then, either $\widetilde{p}$ lies in the domain $V^x$ 
or $\widetilde{p}$ lies in $V^y \setminus V^x$.

Firstly, suppose that $\widetilde{p} \in V^x$ and let
$(\overline{x},\overline{y},\overline{z})$ be the global coordinates of the
$x$-directional chart (given in subsection~\ref{subsect-directionalexpressions}). 
It follows that there exists a unique
pair of constants $(\xi,\eta) \in \cR^2$ such that the coordinates
$$
(\overline{x}, \overline{y} - \eta, \overline{z} - \xi)
$$
define an adapted local chart $\widetilde{p}$. The stabilization of such local
adapted chart corresponds to the stabilization of the Newton data
$(\xi,\eta)\cdot\Blx \Omega$ (where $\Omega$ is the Newton data centered at $p$).
Therefore, it follows from Proposition~\ref{prop-stablexblup} that 
\begin{equation}\label{reducesinvproofloc}
\Inv(\widetilde{\cM},\widetilde{\Ax},\widetilde{p}\, ) \; <_\lex \;
\Inv(\cM,\Ax,p)
\end{equation}
This proves the Theorem in the case where $\widetilde{p} \in V^x$.

Suppose now that $\widetilde{p} \in V^y \setminus V^x$ and let 
$(\overline{x},\overline{y},\overline{z})$ be the global coordinates of the
$y$-directional chart (given in subsection~\ref{subsect-directionalexpressions}).  Then, there exists a unique
constant
$\xi \in \cR$ such that 
$$
(\overline{x}, \overline{y}, \overline{z} - \xi)
$$
defines an adapted local chart at $\widetilde{p}$.  It suffices now to apply the 
Proposition~\ref{prop-stablexblup} to conclude that (\ref{reducesinvproofloc})
also holds.  This proves the Theorem.
\section{Global Theory}\label{sect-GlobalTheory}
\subsection{Upper Semicontinuity of Virtual Height at $\NElem \cap \vD$}
In this subsection, our goal is to prove the upper semicontinuity of the virtual height.
In other words, we will prove that each point $p \in \NElem(\cM)\cap \vD$ has an open neighborhood
$V \subset M$ such that, 
$$
\vheight(\cM,\Ax,q) < \vheight(\cM,\Ax,p)
$$
for each point $q \in \NElem(\cM) \cap \vD \cap V$. For shortness, denote the set of 
nonelementary points simply by $\NElem$, and let 
$$\vheight: \NElem \cap \vD \rightarrow \cN$$
be the virtual height function.
The {\em stratum of virtual height $h$} at $\vD$ is the subset
$$
S_h \cap \vD = \{p \in \NElem(\cM) \cap \vD \mid \vheight(p) = h\}
$$
To enunciate the next result, we introduce the following notion:
Let $D \subset \vD$ be an irreducible component of the divisor and let $p \in \NElem \cap \vD$ be a point in $D$.
We shall say that a local chart $(U,(x,y,z))$ at $p$ is {\em $D$-adapted} if
\begin{itemize}
\item $\vZ$ is locally generated by $\frac{\partial}{\partial z}$;
\item $D = \{x = 0\}$.
\item $\vD \cap U \subset \{xy = 0\}$.
\end{itemize}
We further say that the chart $(U,(x,y,z))$ is {\em $D$-stable} if the corresponding Newton data $\Omega = ((x,y,z),\iota,\Theta)$ is stable (for some choice of local generator for the line field).
\begin{Remark}
Notice that if the point $p$ belongs to the intersection $D \cap D^\prime$ of two divisors then
in a $D$-adapted chart $(U,(x,y,z))$ we necessarily have $D = \{x = 0\}$ and $D^\prime = \{y = 0\}$.
\end{Remark}
\begin{Proposition}\label{prop-delta1D>0=0}
Given an irreducible component of the divisor $D \subset \vD$ and a point $p \in S_h \cap D$, 
 there exists an open neighborhood $V \subset M$ of $p$ such that 
$$\vheight(q) \le h$$ 
for each point $q \in V \cap D \cap \NElem(\cM)$.  Moreover, if we fix a
$D$-stable local 
chart $(U,(x,y,z))$ at $p$,
\begin{itemize}
\item[(i)] If $\Delta_1(D) > 0$ then we locally have $S_h \cap D = \{x = z = 0\}$.
\item[(ii)] If $\Delta_1(D) = 0$ then we locally have $S_h \cap D = \{p\}$.
\end{itemize}
where $\Delta(D) = (\Delta_1(D),\Delta_2(D))$ is the vertical displacement 
vector of the corresponding Newton data $\Omega$.
\end{Proposition}
\begin{proof}
First of all, we consider  the case where $\Delta_1(D) > 0$.
We will show that 
there exists an open neighborhood of the origin $U \subset \cR^2$ such that for each 
$(\xi,\eta) \in U$, the translation map  
\begin{equation}\label{thistransl}
\widetilde{y} = y + \eta, \qquad \widetilde{z} = z + \xi
\end{equation}
is such that one of the following two situations occurs:
\begin{enumerate}[({i}.1)]
\item If $\xi \ne 0$ then the translated Newton data $\widetilde{\Omega} = (\xi,\eta)\cdot \Omega$
is in final situation.
\item If $\xi = 0$ then the virtual height at the translated point $\widetilde{p}$ (i.e.\ 
the point which obtained from $p$ by the local translation (\ref{thistransl})) is equal to $h$.
\end{enumerate}
The item (i.1) is easy.  Indeed, let $\bfm = (0,m_2,m_3)$ be the main vertex of $\Omega$.  Then, 
it is immediate to see that there exists a constant $C > 0$ and a neighborhood
of the 
origin $U \subset \cR^2$ 
such that for each $(\xi,\eta) \in U$, the translated data $\widetilde{\Omega}$ evaluated at 
the point $\widetilde{\bfm} = (0,m_2,0)$ is such that
$$\| \widetilde{\Omega}(\widetilde{\bfm}) \| \ge  C\ |\xi|^{m_3+1}$$
(where $\| \cdot \|$ denotes the Euclidean norm).  Since $m_2 \in \{-1,0\}$,
we 
see that $\widetilde{\Omega}$ is in final situation if $\xi \ne 0$.

Let us prove item (i.2).  If $\xi = 0$, there exists constants $C>0$ and 
$\delta > 0$ such that for each translation $(0,\eta)$ with $|\eta| < \delta$, we have
$$\| \widetilde{\Omega}(\bfm) \| \ge  C,$$
Looking at the restriction of $\Omega$ to the set $\{0\} \times \cZ^2$, the 
translation causes the five possible movements shown in
figure~\ref{fig-possibletrans} 
(we denote by $\widetilde{\bfm}$ the main vertex of $\widetilde{\Omega}$).
In each case, 
it is immediate to see that $\vheight(\widetilde{p}) = h$.

\begin{figure}[htb!] 
\psfrag{m}{\small $\bfm$}
\psfrag{mt}{\small $\widetilde{\bfm}$}
\psfrag{m2=0}{\small $m_2 = 0$}
\psfrag{m2=-1}{\small $m_2 = -1$}
\begin{center} 
\includegraphics[height=8cm]{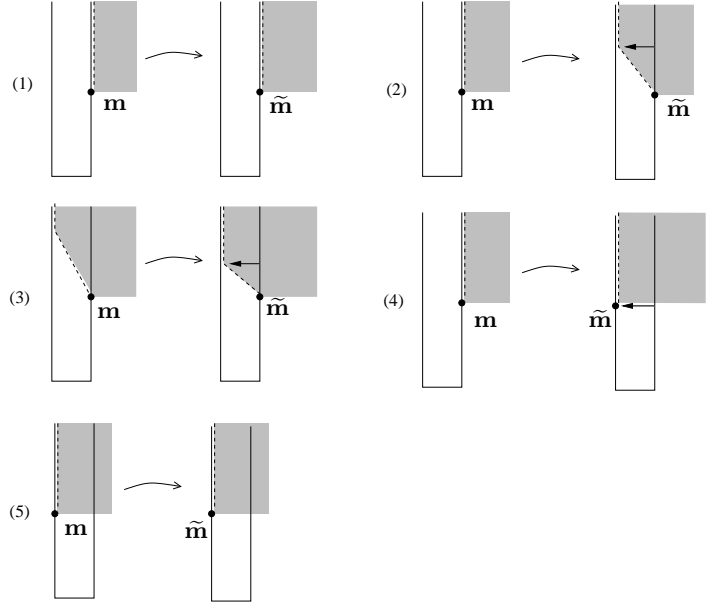}
\end{center}
\caption{The effect of a translation.}
\label{fig-possibletrans}
\end{figure}
 
We proceed now to the proof of the Proposition in the case where $\Delta_1(D) = 0$.  We will show that 
there exists an open neighborhood of the origin $U \subset \cR^2$ such that for each 
$(\xi,\eta) \in U$, the translation map  
$$
\widetilde{y} = y + \eta,  \qquad \widetilde{z} = z + \xi
$$
is such that one of the following two situations occurs:
\begin{enumerate}[({ii}.1)]
\item If $\xi \ne 0$ and $\eta = 0$ then the translated Newton data
      $\widetilde{\Omega} = (\xi,\eta)\cdot \Omega$ 
is in final situation.
\item If $\eta \ne 0$ then the virtual height of the translated point
      $\widetilde{p}$ is 
strictly smaller that $h$.
\end{enumerate}
The proof of item (ii.1) is analogous to the proof of item (i.1).

In order to prove item (ii.2), we consider the following {\em blowing-up} 
in the parameters $(\xi,\eta)$.
\begin{eqnarray*}
\phi: \cR^+\times \cS^1 &\longrightarrow& \cR^2 \\
(r,\theta) &\longmapsto& (\eta, \xi)  = (r \cos(\theta),r^s \sin(\theta))
\end{eqnarray*}
where $s := \Delta_2(D)$. We claim that there exists a neighborhood 
$\overline{U} \subset \cR^+\times \cS^1$ of the set $\{r = 0\}$ such that the 
corresponding neighborhood of the origin
$U := \phi(\overline{U})$ satisfies the conditions stated above.

From item (ii.1), we know that
$\widetilde{\Omega}$ is in final situation for $\theta = \pi/2$ and $r$ sufficiently small.  
Therefore, since this is an open condition,  
there exists an open neighborhood $\overline{V} \subset \cR^+\times \cS^1$ of the point 
$(r,\theta) = (0,\pi/2)$ such that
$\widetilde{\Omega}$ is in final situation for each translation with $(\xi,\eta)$ in $\phi(\overline{V})$.

To complete the study, it suffices to study the collection of all translations $(\xi,\eta)$
which are contained in the image
of the directional  blowing-up
$$
\eta = \overline{\eta}, \quad \xi = \overline{\eta}^s \overline{\xi}
$$
with $\overline{\eta}$  varying in $\cR$ and $\overline{\xi}$ belonging to some compact subset $K \subset \cR$.

For this, we fix some stable local chart $(V,(x,y,z))$ at $p$ and 
consider the local blowing-up for $(\cM,\Ax)$, 
$$
\Phi: \widetilde{\cM} \rightarrow \cM\cap V
$$
(see subsection~\ref{subsect-Newtoninvarlocres}).  Under such blowing-up, the above collection of translations 
can be studied in the domain of the $y$-directional chart, where the blowing-up can be written as 
$$
x = \overline{y}^{\omega_1}\; \overline{x},\quad y = \overline{y}^{\omega_2}, 
\quad z = \overline{y}^{\omega_3} 
\; \overline{z} 
$$
(with $(\omega_2:\omega_3) = (1:s)$). Fixed $\overline{\xi} \in K$, let $\overline{p}$ 
denote the point on the exceptional divisor $\widetilde{D} = \Phi^{-1}(Y_p)$ 
which is obtained by the $y$-directional blowing-up followed by the vertical 
translation
$$
\widetilde{z} = \overline{z} + \overline{\xi}
$$
It follows from the proof of Proposition~\ref{prop-stableyblup} that 
$\vheight(\overline{p}) < \vheight(p) = h$.  
Therefore, using the compactness of $K$, 
it suffices to prove the following:\\
{\em Claim:} Let $(\overline{U}, (\bar x, \bar y, \bar z))$ be a stable local chart at $\overline{p}$.  Then, 
there exists a constant $\delta > 0$ such that for each translation 
\begin{equation}\label{onlyhoriztrans}
\widetilde{y} = \bar y + \bar \eta
\end{equation}
with $|\bar \eta| < \delta$, the corresponding translated point $\widetilde{p}$ is such that 
$\vheight(\widetilde{p}) \le \vheight(\overline{p})$.

The proof of such claim is similar to the proof of item (ii.1).  Let 
$\overline{\Omega}$ be the Newton data at 
$\overline{p}$ and $\overline{\Delta} =
(\overline{\Delta}_1,\overline{\Delta}_2)$ be 
the corresponding vertical 
displacement.  Two cases can appear:
\begin{enumerate}[({ii.2}.(a))] 
\item $\overline{\Delta}_1 > 0$
\item $\overline{\Delta}_1 =  0$
\end{enumerate}
In the case (ii.2.a), the item (i) treated above implies that 
$\vheight(\widetilde{p}) = \vheight(\overline{p})$.

In the case (ii.2.b), let $\bfm = (0,0,m_3)$ be the main vertex associated
to $\overline{\Omega}$ and let $\medge = \overline{\bfm,\bfv}$ be the corresponding main edge.  
Then, if we write $\bfv = (v_1,v_2,v_3)$ (with $v_3 < m_3$), there exists
constant $C,\delta > 0$ such that  
$$
\|\widetilde{\Omega}(\widetilde{\bfv})\| \ge C\ |\overline{\eta}|^{v_2+1},\quad 
\mbox{ with }\widetilde{\bfv} = (0,0,v_3)
$$ 
where $\widetilde{\Omega}$ is the Newton data obtained by the translation 
(\ref{onlyhoriztrans}).  We easily conclude that
$\vheight(\widetilde{p}) \le \vheight(\overline{p})$.
\end{proof}

\subsection{Upper Semicontinuity of the Invariant at $\NElem \cap \vD$}\label{subsect-uppersemicont}
Using the results of the previous subsection, let us prove the upper semicontinuity of the function 
$$
\Inv : \NElem \cap \vD \rightarrow \cN^6
$$
where $\Inv(p)$ is a shorter notation for $\Inv(\cM,\Ax,p)$.

We recall that the invariant 
$\Inv(p) = (\Invp(p),\Invs(p))$ is given by
$$
\Invp = (\vheight,m_2+1,m_3), \quad \Invs = (\# \iota_p - 1,\lambda\ \Delta_1,\lambda \max\{\Delta_2,0\})
$$
where such quantities are computed using some stable local chart $(U,(x,y,z))$ for $(\cM,\Ax)$ at $p$.
The following remark will be useful in the sequel:
\begin{Remark}\label{rem-remarkinv(p)}
The definition of $\Inv$ implies that
\begin{enumerate}[(1)]
\item If $p \in S_h \cap \vD$ is such that $\#\iota_p = 2$ then $\Inv(p) = (h,1,h,1,\cdot,\cdot)$. 
\item If $p \in S_h \cap \vD$ is such that $\#\iota_p = 1$ then either
$$\Inv(p) = (h,0,m_3,0,\cdot,\cdot)\quad \mbox{ or }\quad  \Inv(p) = (h,1,h,0,\cdot,\cdot)$$
\end{enumerate}
where $\cdot$ denotes some arbitrary natural number.
\end{Remark}
\begin{Lemma}\label{lemma-preservesDelta1}
Let $p \in S_h \cap \vD$ be a point such that $\#\iota_p = 1$.  Assume that the displacement vector 
$\Delta = (\Delta_1,\Delta_2)$ satisfies
$$
\Delta_1 > 0
$$
Then, there exists a neighborhood $V \subset M$ of $p$ such that
for each point $q \in (S_h \cap \vD) \cap V$, the corresponding displacement vector
$\widetilde{\Delta} = (\widetilde{\Delta}_1,\widetilde{\Delta}_2)$ is such 
that $\widetilde{\Delta}_1 = \Delta_1$.
\end{Lemma}
\begin{proof}
Let us fix a stable local chart $(U,(x,y,z))$ for $(\cM,\Ax)$ at $p$.  Then, the local blowing-up center is given by
$Y_p = \{x = z = 0\}$.  
Let 
$\Phi: \widetilde{\cM} \rightarrow \cM \cap U$ be the local blowing-up with center 
$Y_p$ and weight-vector $\bfomega_p =  (q_1,0,p_1)$ 
(where $\Delta_1 = p_1/q_1$ is the irreducible rational representation of $\Delta_1)$.

Suppose, by absurd, that there exists a sequence of real numbers 
$\{\eta_k\} \rightarrow 0$, 
such that the corresponding sequence of Newton data $\widetilde{\Omega}_k$ which are 
obtained by the translations $\widetilde{y}_k = y + \eta_k$ have a displacement vector $\Delta^k = (\Delta_1^k,\Delta_2^k)$ such 
that 
$$
\Delta_1^k > \Delta_1
$$
Let $\{q_k\} \rightarrow p$ denote the sequence of points in $Y_p$ which 
are obtained by such sequence of translations.

Using Lemmas~\ref{lemma-directx} and \ref{lemma-generaltransl}, we see that, for each $k$, 
the set $\Phi^{-1}(q_k)$ contains at least one nonelementary point 
$\widetilde{q}_k$ such that 
\begin{equation}\label{qkheighth}
\vheight(\widetilde{q}_k) = \vheight(q_k) = h
\end{equation}
In fact, we can choose such point as the origin in the $x$-directional 
chart of the blowing-up.

On the other hand, the proof of Lemmas~\ref{lemma-delta>01} and \ref{lemma-delta>02} 
imply that each nonelementary point 
$\widetilde{p}$ in $\Phi^{-1}(p)$ satisfies one of the following conditions:
\begin{enumerate}[(a)]
\item Either $\vheight(\widetilde{p}) < h$, or
\item $\vheight(\widetilde{p}) = h$  and $\widetilde{\Delta}_1 = 0$.
\end{enumerate}
where $\widetilde{\Delta} = (\widetilde{\Delta}_1,\widetilde{\Delta}_2)$ is
the vertical 
displacement vector of the 
Newton data at $\widetilde{p}$ (for some fixed stable local chart).  

Using item (ii) of Proposition~\ref{prop-delta1D>0=0} and the compactness of
$\Phi^{-1}(p)$, 
we conclude that there exists some 
neighborhood $\widetilde{V} \subset \widetilde{M}$ of 
$\Phi^{-1}(p)$ such that each nonelementary point $\widetilde{q} \in 
\widetilde{V} \setminus \Phi^{-1}(p)$ is such that $\vheight(\widetilde{q}) <
h$.  
This contradicts (\ref{qkheighth}).
\end{proof}
\begin{Proposition}\label{Prop-invuppsemicont}
The function $\Inv: \NElem\cap \vD \rightarrow \cN^6$ is 
upper semicontinuous (for the lexicographical ordering on $\cN^6$).
\end{Proposition}
\begin{proof}
Given a point $p \in \NElem \cap \vD$,   we have to prove that 
there exists a neighborhood $V \subset M$ of $p$ such that
for each point $q \in \NElem\cap \vD \cap V$,
$$\Inv(q) \le_\lex \Inv(p).$$
The upper semicontinuity of the initial segment of the local invariant, namely
$$(\vheight,m_2+1,m_3,\#\iota_p - 1)$$ 
is obvious by using the Remark~\ref{rem-remarkinv(p)} and Proposition~\ref{prop-delta1D>0=0}.

Let us fix $p \in S_h \in \vD$.  We claim that there exists a 
neighborhood $V \subset M$ of $p$ such that
for each point $q \in (\NElem\cap \vD) \cap V$, we have 
$$
\Invp(q) = \Invp(p), \#\iota_q = \#\iota_p \quad \Rightarrow \quad \Invs(q) \le_\lex \Invs(p)
$$
Indeed, if $\Invp(q) = \Invp(p)$ and $\#\iota_q = \#\iota_p$ then it follows from items (i) and (ii) of 
Proposition~\ref{prop-delta1D>0=0} and from Remark~\ref{rem-remarkinv(p)} that the 
following properties holds:
$$
\#\iota_p = 1, \quad \mbox{ and }\quad \Delta_1 > 0
$$
Therefore, using Lemma~\ref{lemma-preservesDelta1}, we conclude (up to
restricting 
$V$ to some smaller neighborhood of $p$)
that $\Delta_1^q = \Delta_1$ for each point $q \in (\NElem\cap \vD) \cap
V$.  
Moreover, for each fixed stable local chart $(U,(x,y,z))$ at $p$, it 
is clear that the adapted local chart at $q$ which is obtained by the translation
$$
\widetilde{x} = x, \quad \widetilde{y} = y + \eta, \quad \widetilde{z} = z
$$
(for some appropriately chosen constant $\eta \in \cR$) is also stable.  Therefore, we obviously have 
(up to a new restriction of $V$ to 
some smaller neighborhood of $p$) that $\Delta_2^q \le \Delta_2^p$.  
This concludes the proof.
\end{proof}
\subsection{Points at $\NElem \setminus \vD$ and Generic Newton Polygon}
A point $p \in \NElem \setminus \vD$ will be called {\em smooth} if the germ
of analytic set $\NElem_p$ is locally a smooth 
one-dimensional analytic curve.

We shall say that an adapted local chart $(U,(x,y,z))$ for $(\cM,\Ax)$ at $p$ is {\em smoothly adapted} if
$$\NElem = \{y = z = 0\}$$
It follows from Proposition~\ref{prop-transitionmaps} that the transition map between two smoothly adapted
local charts $(U,(x,y,z))$ and $(U^\prime,(x^\prime,y^\prime,z^\prime))$ has the form
\begin{equation}\label{rem-formsmoothly}
x^\prime = f(y) + x u(x,y), \quad y^\prime = y v(x,y), \quad z^\prime = y h(x,y) + z w(x,y,z)
\end{equation}
where $f,u,v,h,w$ are analytic functions such that $f(0) = 0$ and $u,v,w$ are units.

\begin{figure}[htb] 
\psfrag{v1}{\small $v_1$}
\psfrag{v2}{\small $v_2$}
\psfrag{v3}{\small $v_3$}
\psfrag{NG}{\small $\vN_G$}
\begin{center} 
\includegraphics[height=5cm]{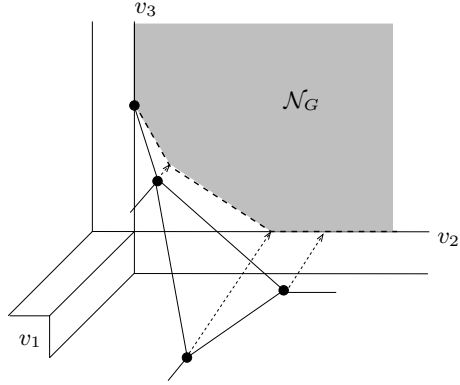}
\end{center}
\caption{The generic Newton polygon.}
\label{fig-genpol}
\end{figure}

Let $\Omega$ be the Newton data for $(\cM,\Ax)$ at the point smooth $p$, relatively to some smoothly adapted
local chart $(U,(x,y,z))$.  
The {\em generic Newton map} associated to $\Omega$ is the map $\Theta_G:\cZ^2 \rightarrow \{0,1\}$ is given by
$$
\Theta_G(\bfv) = \left\{
\begin{matrix}
0,& \quad \mbox{ if }(\cZ \times \{\bfv\}) \cap \supp(\Omega) = \emptyset,\cr\cr
1,& \quad \mbox{ if }(\cZ \times \{\bfv\}) \cap \supp(\Omega) \ne \emptyset 
\end{matrix}
\right.
$$
(see figure~\ref{fig-genpol}). The {\em generic Newton polygon} associated to
$\Omega$ is the convex polygon in $\cR^2$ given by $\vN_G(\Omega) = \supp(\Theta_G) + \cR^2_+$.  

The {\em generic higher vertex} of $\Omega$  is the minimal point $\bfp_G \in 
\vN_G$ for the lexicographical ordering. 
\begin{Remark}
The generic Newton polygon can be equivalently defined as 
$$
\vN_G = \pi(\vN)
$$
where $\pi:\cR^3 \rightarrow \cR^2$ is the linear projection $\pi(v_1,v_2,v_3) = (v_2,v_3)$ and 
$\vN = \vN(\Omega)$ is the Newton polyhedron of $\Omega$.
\end{Remark}
The triple $\Omega_G = ((x,y,z),\iota_p,\Theta_G)$ will be called the {\em generic Newton data}  at $p$.

The {\em generic edge} associated to $\bfp_G$ is the unique edge 
$\medge(\bfp_G) \subset \vN_G$ which intersects the horizontal line
$\{(v_1,v_2) \in \cR^2 \mid v_2 = p_2 - 1/2\}$ (where we write $\bfp_G = (p_1,p_2)$).  We convention that
$\medge(\bfp_G) = \emptyset$ if such intersection is empty.

We shall say that $\Omega_G$ is in a {\em nilpotent configuration}
if the following conditions are satisfied:
\begin{enumerate}[(i)]
\item $\bfp_G = (-1,p_2)$, for some integer $p_2 \in \cZ$;
\item The edge $\medge(\bfp_G)$ has the form $\overline{\bfp_G,\bfn}$, for some vertex $\bfn = (0,n_2)$ with $n_2 \in \cZ$.
\end{enumerate}
If one of these conditions fails, we shall say that $\Omega_G$ is in a {\em regular configuration}.

The {\em generic main vertex}  is a vertex $\bfm_G \in \vN_G$ which is chosen as follows:
\begin{enumerate}[(i)]
\item If $\Omega_G$ is in a regular configuration then $\bfm_G := \bfp_G$.
\item If $\Omega_G$ is in a nilpotent configuration then $\bfm_G := \bfn$ (where  
$\medge(\bfp_G)=\overline{\bfp_G,\bfn}$).
\end{enumerate}
Let us write $\bfm_G = (m_1,m_2)$. The {\em generic main edge} is the edge $\medge_G \subset \vN_G$ which intersects
the horizontal line $\{(v_1,v_2) \in \cR^2 \mid v_2 = m_2 - 1/2\}$.  We convention that 
$\medge_G = \emptyset$ if this intersection is empty.

Note that we can write the generic main edge as
$$\medge_G = \bfm_G + t (\Delta,-1)$$
where $t$ belongs to a real interval  
of the form $[0,L]$ (for some $L > 0$) and $\Delta \in \cQ_{>0}$ is a positive rational number.  
In this setting, we shall shortly say that the generic Newton data $\Omega_G$ belongs to the
class $\Nmaps_{\Delta}^{\bfm_G}$.

The point $p \in \NElem \setminus \vD$ will be called {\em generic} with respect
to an adapted local chart $(U,(x,y,z))$ if the following conditions hold:
\begin{enumerate}[(i)]
\item $p$ is a smooth point.
\item $(U,(x,y,z))$ is a smoothly adapted local chart $p$;
\item All the vertices of the corresponding Newton polyhedron $\vN(\Omega)$
      belong to the region $\{-1,0\} \times \cZ^2$. 
\end{enumerate}
\begin{Remark}
Suppose that $p \in \NElem \setminus \vD$ is generic with respect to
an adapted local chart $(U,(x,y,z))$.  Then, the generic main edge $\medge_G$ defines a face
of the Newton polygon $\vN = \vN(\Omega)$.  More precisely, the Minkowski sum 
$$
\vF := \medge_G + \{t\cdot (1,0,0) \mid t \in \cR_+\}
$$
is a face of $\vN$.
\end{Remark}

\begin{figure}[htb] 
\psfrag{v1}{\small $v_1$}
\psfrag{v2}{\small $v_2$}
\psfrag{v3}{\small $v_3$}
\psfrag{NG}{\small $\vN_G$}
\psfrag{mg}{\small $\medge_G$}
\psfrag{Mg}{\small $\bfm_G$}
\psfrag{v}{\small $\bfv$}
\begin{center} 
\includegraphics[height=5.5cm]{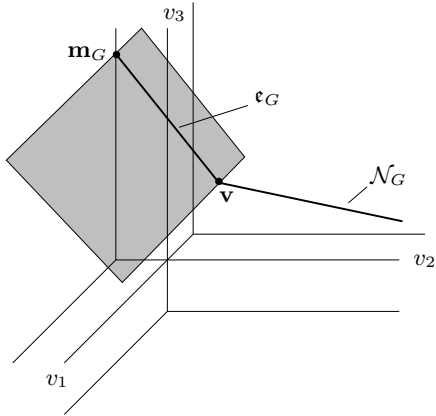}
\end{center}
\caption{Newton polyhedron at a generic point.}
\end{figure}

\subsection{Generic Edge Stability and Equireducible Points}
Let $p \in \NElem \setminus \vD$ be a smooth point and let $(U,(x,y,z))$ be a smoothly adapted local chart
for $(\cM,\Ax)$ at $p$.  We denote respectively by $\Omega$ and $\Omega_G$ the associated
Newton data and generic Newton data (for some choice of local generator for the line field)

Given a rational number $\delta \in \cQ_{\ge 0}$, the group of 
$\gG_{\delta}$-maps is the group of all analytic change of coordinates of the form
$$
\widetilde{z} = z + g(y),
$$
where $g(y)$ is given by $g(y) = \xi y^\delta$ (for some constant 
$\xi \in \cR$) if $\delta \in \cN$ and $g \equiv 0$ otherwise.

The group $\gG_{\delta}$ acts naturally on the class of Newton data via the coordinate change 
$(x,y,z) \rightarrow (x,y,\widetilde{z})$.  
Given a map $g \in \gG_\delta$, we denote the such action on $\Omega$ simply
by $g\cdot \Omega$.

Suppose that $\Omega_G$ belongs to the class $\Nmaps_{\Delta}^{\bfm_G}$.  
We say that $\Omega$ is {\em generic edge stable} if 
$$
(\; g \cdot \Omega \; )_G \in \Nmaps_{\Delta}^{\bfm_G}
$$
for all map $g \in \gG_{\Delta}$ with $g(0) = 0$.  In other words, the generic Newton data associated to $g\cdot \Omega$ lies in the
class $\Nmaps_{\Delta}^{\bfm_G}$ for all $g \in \gG_{\Delta}$.

The local chart $(U,(x,y,z))$ will be called {\em generic edge stable} if $\Omega$ is generic edge stable.

\begin{Lemma}\label{lemma-edgestablegeneric}
Suppose that the smoothly adapted local chart $(U,(x,y,z))$ is {\em not} generic edge stable.
Then, there exists a unique map $g \in \gG_\Delta$ such that
the transformed generic Newton data $( g\cdot \Omega)_G$ does not belongs to 
$\Nmaps_{\Delta}^{\bfm_G}$.
\end{Lemma}
\begin{proof}
The result can be proved by straightforward modifications in the proof of Lemma~\ref{lemma-uniqueedgemap}.  
\end{proof}
Now, we are ready to give the main definition of this subsection:

We shall say that a smooth point $p \in \NElem \setminus \vD$ is {\em
equireducible} if there exist a smoothly adapted local chart $(U,(x,y,z))$ for $(\cM,\Ax)$ at $p$ such that
\begin{enumerate}[(i)]
\item $p$ is generic with respect to $(U,(x,y,z))$;
\item The corresponding Newton data $\Omega$ is generic edge stable.
\end{enumerate}
In this case, $(U,(x,y,z))$ will be called an {\em equireduction chart} for $(\cM,\Ax)$ at $p$.
\begin{Lemma}\label{lemma-indepcoordequi}
Let $(U,(x,y,z))$ and $(U^\prime,(x^\prime,y^\prime,z^\prime))$ be two equireduction charts at 
an equireducible point $p$.  Then, the transition map (see (\ref{rem-formsmoothly}))
has necessarily the form 
$$
x^\prime = f(y) + x u(x,y), \quad y^\prime = y v(x,y), \quad z^\prime = yh(x,y) + z w(x,y,z)
$$
For some analytic functions $f,u,v,h,w$ such that $f(0) = 0$ and $u,v,w$ are units.  Moreover,
the support of the function $H(x,y) = yh(x,y)$ satisfies the following property
$$
\supp(H) \subset \{(v_1,v_2) \in \cN^2 \mid v_2 > \Delta\}
$$
\end{Lemma}
\begin{proof}
This is a direct corollary of Lemma~\ref{lemma-edgestablegeneric}.
\end{proof}
As a consequence of the second part of the Lemma, the Newton data $\Omega^\prime$ which is associated to the chart 
$(U^\prime,(x^\prime,y^\prime,z^\prime))$ is such that 
$$\Omega_G \in \Nmaps_{\Delta}^{\bfm_G} \quad \Longleftrightarrow \quad (\Omega^\prime)_G \in \Nmaps_{\Delta}^{\bfm_G}$$
Let us now characterize generic Newton data which are centered at points in $\NElem \setminus \vD$. 

First of all, we introduce the following notion: A generic Newton data $\Omega_G$ is in 
{\em final situation} if one of the following conditions holds:
\begin{enumerate}[(i)]
\item The generic main vertex $\bfm_G = (m_1,m_2)$ is such that $m_2 \in
      \{-1,0\}$, or
\item The main edge is given by $\medge_G = \overline{\bfm_G, \bfv}$, where
$\bfm_G = (-1,1)$ and $\bfv = (1,-1)$.
\end{enumerate} 
As a consequence of the definition of equireducible point, we get the following result:

\begin{figure}[htb]
\psfrag{i.a}{\small $m_2 = 0$}
\psfrag{i.b}{\small $m_2 = -1$}
\psfrag{i.c}{\small $m_2 = 0$}
\psfrag{ii}{\small $\bfv = (0,-1)$}
\psfrag{iii}{\small $\bfv = (0,0)$}
\psfrag{iv}{\small $\bfm = (-1,1),\;\bfv = (1,-1)$}
\psfrag{v2}{\small $v_2$}
\psfrag{v3}{\small $v_3$}
\begin{center}
\includegraphics[height=6cm]{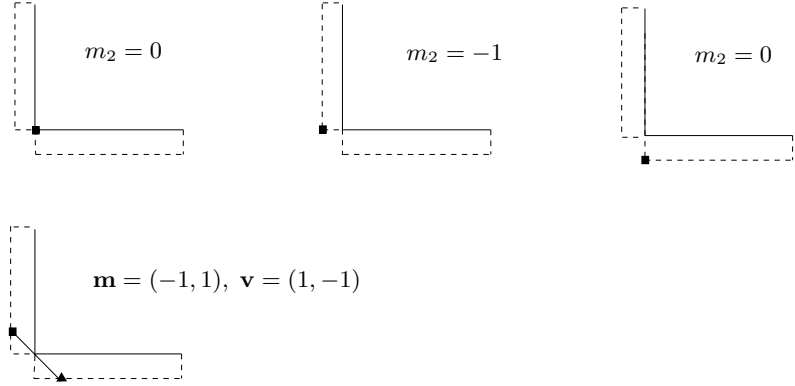}
\end{center}
\caption{The final situations for the generic Newton data.}
\label{fig-finalsit2}
\end{figure}

\begin{Proposition}
Let $p \in \NElem \setminus \vD$ be an equireducible point and let $(U,(x,y,z))$ be an equireduction chart at $p$. 
Then, the associated generic Newton data $\Omega_G$ 
is not in a final situation.
\end{Proposition}
\begin{proof}
Analogous to the proof of Proposition~\ref{prop-finalsituation}.
\end{proof}

\subsection{Local Blowing-up at Equireducible Points}
Let $p \in \NElem(\cM) \setminus \vD$ be an equireducible point.  Let 
$\Omega$ be the Newton data for $(\cM,\Ax)$ at $p$, with respect to some 
equireduction chart $(U,(x,y,z))$.

The {\em generic virtual height} for $(\cM,\Ax)$ at $p$ is defined as
$$
\vheight_G(\cM,\Ax,p) := 
\left\{
\begin{matrix}
\lfloor m_2  + 1 - \frac{1}{\Delta} \rfloor , &\mbox{ if }m_1 = -1 \cr \cr
m_2, &\mbox{ if }m_1 = 0
\end{matrix}
\right.
$$
where $\bfm_G = (m_1,m_2)$ is the main vertex of the generic Newton polygon
of $\vN_G$.

The {\em local blowing-up center} associated to $(\cM,\Ax)$ at $p$ is 
the submanifold  
$$Y_p = \{y = z = 0\}.$$ 
Assume that the generic Newton data $\Omega_G$ belongs to the class $\Nmaps_{\Delta}^{\bfm_G}$.
The {\em weight-vector} associated to $(\cM,\Ax)$ at $p$ is given by
$$\bfomega = (0,q,p)$$
where $\Delta = p/q$ is the irreducible rational representation of $\Delta$.
\begin{Remark}
It follows from the Lemma~\ref{lemma-indepcoordequi} that $\vheight_G(\cM,\Ax,p)$, $Y_p$ and 
$\omega$ are
independent of the choice of the equireduction chart $(U,(x,y,z))$.
\end{Remark}
The {\em local blowing-up} for $(\cM,\Ax)$ at $p$ is the $\bfomega$-weighted
blowing-up 
$$
\Phi: \widetilde{\cM} \rightarrow \cM \cap U
$$
with center on $Y_p$, with respect the trivialization given by $(U,(x,y,z))$.
\begin{Remark}\label{rem-preservesstructure2}
The Lemma~\ref{lemma-indepcoordequi} implies that the transition map between two equireduction charts
always preserves the $\bfomega$-quasihomogeneous structure on $\cR^3$.
\end{Remark}
The following Theorem is a version of the Local Resolution of Singularities for 
equireducible points.
\begin{Theorem}\label{theorem-localresolutiongeneric}
Let $(\cM,\Ax)$ be a controlled singularly foliated manifold and let
$p \in A \setminus \vD$ be an equireducible point in $\NElem(\cM)$.
Consider the {\em local blowing-up} for $(\cM,\Ax)$ at $p$, 
$$
\Phi: \widetilde{\cM} \rightarrow \cM \cap U
$$
with respect to some equireduction chart $(U,(x,y,z))$.  Then, there exists an axis 
$\widetilde{\Ax} = (\widetilde{A},\widetilde{\vZ}\, )$ for 
$\widetilde{\cM}$ such that each point 
$\widetilde{p} \in \Phi^{-1}(\{p\}) \cap \widetilde{A}$ belonging to
$\NElem(\widetilde{\cM})$ is such that
$$\vheight(\widetilde{\cM},\widetilde{\Ax},\widetilde{p}) < \vheight_G(\cM,\Ax,p)$$
\end{Theorem}
\begin{proof}
Analogous to the proof of Theorem~\ref{theorem-localresolution}, using now the definition of equireducible points.
\end{proof}
The {\em local invariant} for $(\cM,\Ax)$ at an equireducible point $p \in \NElem
\setminus \vD$ is the vector of natural numbers
$$
\Inv(\cM,\Ax,p) = (\vheight_G(\cM,\Ax,p),0,0,0,0,0) \in \cN^6
$$

\subsection{Distinguished Vertex Blowing-up}\label{subsection-distingvertblowup}
In this subsection, we describe a procedure which will be used to treat the points $p \in \NElem \setminus \vD$ which are not equireducible.  The basic idea is to {\em include} these points in the divisor $\vD$ by an appropriately chosen weighted blowing-up.

Let $(\cM,\Ax)$ be a controlled singularly foliated manifold.  We fix a point $p \in \NElem$, a local generator
$\chi$ for the line field $\linef$ and a local generator $Z$ for the line field $\vZ$ which defines the axis $\Ax$.

The {\em primitive height} for $(\cM,\Ax)$ at $p$ 
is the minimal integer $h = H(\cM,\Ax,p)$ such that the vector field 
$$
\chi^h := (\vL_Z)^h (\chi) 
$$
is nonzero at $p$.  Here, $(\vL_Z)^h$ is the
$h$-fold composition of the Lie Bracket operator $\vL_Z(\cdot) = [Z,\cdot]$. 

We convention that $H(\cM,\Ax,p) = \infty$ if $\chi^h(p) = 0$  for all $h \in \cN$.  
\begin{Lemma}
For $p \in \NElem \setminus \vD$, the primitive
height $H(\cM,\Ax,p)$ is a well-defined natural number.  Moreover, it is
independent of the choice of the local generators $\chi$ and $Z$. 
\end{Lemma}
\begin{proof}
Let us prove that $H(\cM,\Ax,p)$ is finite.  For this, we fix 
an adapted local chart $(U,(x,y,z))$  at $p$ and write
$$
\chi = F(x,y,z) \frac{\partial}{\partial x} + G(x,y,z) \frac{\partial}{\partial y} + 
H(x,y,z) \frac{\partial}{\partial z}
$$
for some analytic germs $F,G,H$.  We can also choose $Z = \frac{\partial}{\partial z}$.  Therefore,
$$
\chi^h = \frac{\partial^h F}{\partial z^h}(x,y,z) \frac{\partial}{\partial x} + \frac{\partial^h G}{\partial z^h}(x,y,z) 
\frac{\partial}{\partial y} + 
\frac{\partial^h H}{\partial z^h}(x,y,z) \frac{\partial}{\partial z}
$$
If the collection of vector fields
$\{\chi^h\}$ vanishes at the origin for all $h \in \cN$ then the germs $F,G,H$ necessarily belong to the 
ideal $(x,y)\vO_p$.  
This contradicts the fact that $\Ax$ is an axis for $\cM$ (see definition~\ref{def-axis}). 

We now prove that the primitive height is independent of the choice of $\chi$ and $Z$.  For this,  
it suffices to observe that, if we write $\chi^\prime = U\, \chi$ and $Z^\prime = V\,  Z$, for some units $U,V$, then
$$
[Z^\prime,\chi^\prime] = [V Z,U\chi] = UV [Z,\chi] + V Z(U) \chi + U \chi(V) Z
$$
Proceeding by induction, we conclude that $(\vL_{Z^\prime})^h (\chi^\prime)$ vanishes at $p$
if and only if $(\vL_Z)^h (\chi)$ vanishes at $p$.
\end{proof}
An adapted local chart $(U,(x,y,z))$ at a point $p \in \NElem \setminus \vD$ 
will be called {\em strongly adapted} if the associated Newton data $\Omega$ has a polyhedron 
with a vertex of the form $\bfd = (-1,0,d_3)$, where
$$
d_3 = H(\cM,\Ax,p)
$$
The vertex $\bfd$ will be called {\em distinguished vertex}.

The following Lemma shows that we can always construct a strongly adapted local chart.
\begin{Lemma}
Given an adapted local chart $(U,(x,y,z))$ at $p \in \NElem \setminus \vD$, 
there exists a linear change of coordinates of the form 
$$\widetilde{x} = x, \quad \widetilde{y} = y + \xi x, \quad \widetilde{z} = z$$
(for some constant $\xi \in \cR$) such that the resulting local chart  
$(U,(\widetilde{x},\widetilde{y},\widetilde{z}))$ is  strongly adapted.
\end{Lemma}
\begin{proof}
Indeed, since $H = H(\cM,\Ax,p)$ is finite, the Newton data $\Omega$ associated to the chart
$(U,(x,y,z))$ has a Newton polyhedron with at least one vertex of the form $(0,-1,H)$ or one vertex of the form
$(-1,0,H)$.

In the latter case, we are done.  In the former case, it is immediate to see that a change of coordinates as
described in the enunciate leads us to the desired situation.
\end{proof}
Let us fix a strongly adapted local chart $(U,(x,y,z))$ at $p \in \NElem \setminus \vD$.  
Let $\Omega$ be the corresponding Newton data for $(\cM,\Ax)$ and let $\vN(\Omega)$ be its Newton polyhedron. 

The {\em distinguished weight-vector} for $(\cM,\Ax)$ at $p$ (with respect to the chart $(U,(x,y,z))$)
is the weight-vector $\bfomega_\dist \in \cN^3_{>0}$ of minimal norm for which there exists an integer
$\mu \in \cZ$  such that  
$$
\vN \cap \{\bfv \in \cR^3 \mid \< \bfomega_\dist,\bfv \> \; = \; \mu \} = \{\bfd\}
$$
where $\bfd = (-1,0,H(\cM,\Ax,p))$ is the distinguished vertex.
In other words, there exists an integer $\mu$ 
such that the plane $\{ \bfv \mid \< \bfomega_\dist,\bfv \>  =  \mu\}$ intersects
$\vN$ at the single point $\bfd$.

\begin{figure}[htb] 
\psfrag{bfn}{\small $\bfd$}
\psfrag{omegadist}{\small $\bfomega_\dist$}
\psfrag{vN}{\small $\vN$}
\begin{center} 
\includegraphics[height=5cm]{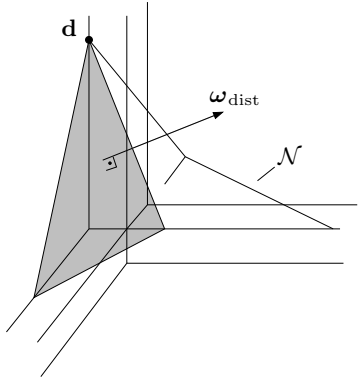}
\end{center}
\caption{Distinguished vertex and $\bfomega_\dist$.}
\end{figure}

The {\em distinguished vertex blowing-up} of $(\cM,\Ax)$ at $p$ 
(with respect to the chart $(U,(x,y,z))$) is the $\bfomega_\dist$-weighted blowing-up
$$
\Phi: \widetilde{\cM} \rightarrow \cM
$$
with center on $\{p\}$, relatively to the local trivialization given by $(U,(x,y,z))$.
\begin{Proposition}\label{prop-distingvertbl}
Let $\Phi: \widetilde{\cM} \rightarrow \cM$ be as above.  
Then, there exists an axis $\widetilde{\Ax}$ for $\widetilde{\cM}$ such that
$$
\vheight(\widetilde{\cM},\widetilde{\Ax},\widetilde{p}) \le H(\cM,\Ax,p).
$$ 
for each point $\widetilde{p} \in \NElem(\widetilde{\cM}) \cap \Phi^{-1}(p)$.
\end{Proposition}
\begin{proof}
Let $\Omega$ be the Newton data associated to the local chart $(U,(x,y,z))$.
We consider separately points lying in the domain of the $z$, $x$ and $y$-directional charts of the blowing-up, and use 
the computations made in the previous section.

In the $z$-directional chart, it follows from Lemma~\ref{lemma-directz}
that $\Blz \Omega$ is in final situation.

In the $x$-directional chart, it follows from Lemma~\ref{lemma-directx}
that the distinguished vertex $\bfd = (0,-1,H(\cM,\Ax,p))$ becomes the higher vertex of
$\Blx \Omega$.  As a 
consequence, 
$\vheight(\widetilde{\cM},\widetilde{\Ax},\widetilde{p}) \le H(\cM,\Ax,p)$ 
for each nonelementary point $\widetilde{p} \in \Phi^{-1}(p)$ which lies 
in the domain of the $x$-directional chart.

In the $y$-directional chart, Lemma~\ref{lemma-directy} implies that the distinguished vertex 
$\bfd = (0,-1,H(\cM,\Ax,p))$ is mapped to the point 
$\widetilde{\bfd} = (0,0,H(\cM,\Ax,p))$ which belongs to the support of 
$\Bly \Omega$. Moreover, 
$$
\supp (\Bly \Omega)\; \cap\; ( \{0\} \times \cZ^2) = \{\widetilde{\bfd}\}
$$
Therefore, we conclude that
$\vheight(\widetilde{\cM},\widetilde{\Ax},\widetilde{p}) \le H(\cM,\Ax,p)$ for each 
point $\widetilde{p} \in \Phi^{-1}(p)$ which lies in the 
domain of the $y$-directional chart.

To finish the proof, we can define an axis $\widetilde{\Ax}$ for 
$\widetilde{\cM}$ exactly as in the proof of Theorem~\ref{theorem-localresolution}.
\end{proof}

\subsection{Nonequireducible Points are Discrete}
Let us now prove that the set of nonequireducible points in $\NElem \setminus \vD$ is finite on each
compact subset of the ambient space.

The first Lemma is an easy result of analytic geometry:
\begin{Lemma}\label{lemma-equiredfinite}
Given an arbitrary point $p \in \NElem$, there exists an open neighborhood $U \subset M$ of $p$
such that each point $q \in (\NElem \setminus \vD) \; \cap \; (U \setminus \{p\})$ is smooth.  
\end{Lemma}
\begin{proof}
Obvious, since the set of nonsmooth points in a Zariski closed subset of the
analytic set $\NElem$.
\end{proof}
\begin{Lemma}\label{lemma-equiredstable}
Let $p \in \NElem \setminus \vD$ an equireducible point.  
Then, there exists an open neighborhood $V \subset M$ of $p$ such that each point  
$q \in \NElem \cap V$ is also an equireducible point.  Moreover, if $(U,(x,y,z))$ is an equireduction
chart at $p$ then the translated coordinates 
$$\widetilde{x} = x + \rho, \quad \widetilde{y} = y, \quad \widetilde{z} = z$$
are equireduction coordinates at $q$ (for some appropriately chosen constant $\rho \in \cR$).
\end{Lemma} 
\begin{proof}
We have to prove that the generic Newton data associated to the translated coordinates
$(\widetilde{x},\widetilde{y},\widetilde{z}) = (x + \rho, y, z)$ 
is edge stable, for all $|\rho|$ sufficiently small.

Suppose, by absurd, that this is not the case. Then, for each $\eps > 0$,
there exists a constant $\rho \in \cR$ with $|\rho| < \eps$ such that the corresponding 
translated point $q$ (with the coordinates $(\widetilde{x},\widetilde{y},\widetilde{z})$) 
satisfies the following:  there exists a $G_\Delta$ map of the form
$$
\overline{z} = \widetilde{z} + \xi \widetilde{y}^\Delta
$$
such that the transformed generic Newton data $\overline{\Omega}_G$ (at the point $q$)
belongs to the class $\Nmaps^{\bfm_G}_{\overline{\Delta}}$ for
some $\overline{\Delta} > \Delta$.

Applying the local blowing-up $\Phi: \widetilde{\cM} \rightarrow \cM \cap U$ 
for $(\cM,\Ax)$ at $p$, we can choose a point $\overline{q} \in \Phi^{-1}(q)$ 
(for instance, the origin in the $y$-directional chart of the blowing-up)
such that its virtual height satisfies
$$\vheight(\overline{q}) \ge \vheight_G(q) \ge \vheight_G(p)$$
But this contradicts the Theorem~\ref{theorem-localresolutiongeneric}
and Proposition~\ref{prop-delta1D>0=0}.
\end{proof}
As an immediate consequence, we get the following result:
\begin{Corollary}\label{corollary-eqopen}
\begin{enumerate}[(i)]
\item The set of equireducible points is an open subset $\Eq$ of $\NElem \setminus
\vD$ (for the topology induced by the topology of $M$).
\item Given two equireducible points $p,q$ on the same connected component of $\Eq$, 
the corresponding generic Newton data necessarily belong to the same class $\Nmaps_{\Delta}^{\bfm_G}$.
\end{enumerate}
\end{Corollary}
Finally, we can state the main result of this subsection:
\begin{Proposition}\label{prop-noneqfinite}
The set of nonequireducible points in $\NElem \setminus \vD$ is finite on each compact
subset $K \subset M$.  
\end{Proposition}
\begin{proof}
By the compactness of $\NElem \cap K$, we just need to prove the following claim:\\
{\em Claim:}  For each point $p \in \NElem\cap K$, there exists an open neighborhood
$U \subset M$ of $p$ such each point in the set $(\NElem \setminus \vD) \cap (U
\setminus \{p\})$ is equireducible.

In order to prove such claim, we consider separately the following cases:
\begin{enumerate}[(1)]
\item $p \in \NElem \setminus \vD$ is an equireducible point;
\item $p \in \NElem \cap \vD$;
\item $p \in \NElem \setminus \vD$ is nonequireducible. 
\end{enumerate}   
In the case (1), the claim is a direct consequence of
Corollary~\ref{corollary-eqopen}.

In the case (2), it suffices to prove that the result holds for each 
irreducible branch of the (possibly singular) germ of analytic set
$\NElem_p$.  Let us fix one such branch, which we denote by $\gamma$. 
Then, two cases can appear:
\begin{enumerate}[({2}.a)]
\item $\gamma = Y_p$.
\item $\gamma \ne Y_p$.
\end{enumerate}
where $Y_p$ is the local blowing-up center for $(\cM,\Ax)$ at $p$.

In the case (2.a), if we fix an arbitrary stable chart $(U,(x,y,z))$ at $p$
then we necessarily have 
$$\gamma = \{y = z = 0\}$$
Using the same reasoning used in the proof of Lemma~\ref{lemma-equiredstable},
we conclude that, for each sufficiently small constant $\rho \in \cR$, the translated
coordinates $(x+\rho,y,z)$ are equireduction coordinates.  Therefore, each
point of $\gamma$ which is sufficiently near $p$ is equireducible.

In the case (2.b), we consider the local blowing-up 
$\Phi: \widetilde{\cM} \rightarrow \cM \cap U$ for $(\cM,\Ax)$ at $p$.
The strict transform of $\gamma$ accumulates at some point 
$$\widetilde{p} \in \NElem(\widetilde{\cM}) \cap \Phi^{-1}(p)$$
We can now repeat the analysis on such point $\widetilde{p}$.  
If we fall in the case (2.b), we 
make another local blowing-up and proceed inductively.

By the Theorem~\ref{theorem-localresolution}, we necessarily fall in case (2.a)
after a finite number of such steps.

Finally, in the case (3), we argue as follows.  Let us fix some strongly
adapted local chart $(U,(x,y,z))$ at $p$ and let 
$\Phi: \widetilde{\cM} \rightarrow \cM \cap U$ be a distinguished
local blowing-up for $(\cM,\Ax)$ at $p$.  

Then,
looking at the strict transform of $\NElem$ and using the compactness of 
$\Phi^{-1}(\{p\})$, the result immediately follows from case (2).
\end{proof}
We shall say that a controlled singularly foliated manifold $(\cM,\Ax)$
is {\em equireducible outside the divisor} if each
point in $\NElem(\cM)  \setminus \vD$ is equireducible.
\begin{Lemma}\label{Lemma-distinguishedprep}
Let $(\cM,\Ax)$ be controlled singularly foliated manifold and $U \subset M$ be a relatively compact subset.
Let $\{p_1,\ldots,p_k\} \subset \NElem \setminus \vD$ be the distinct nonequireducible points of $(\cM,\Ax)$ on 
$U \setminus \vD$.  Then, there exists a blowing-up 
$$
\Phi: \widetilde{\cM} \rightarrow \cM
$$
with center on $p_1$ and an axis $\widetilde{\Ax}$ for $\widetilde{\cM}$ such that the points 
$$\{\Phi^{-1}(p_2),\ldots,\Phi^{-1}(p_k)\}\subset \NElem(\widetilde{\cM})\setminus \widetilde{\vD}$$ 
are the only nonequireducible points for  $(\widetilde{\cM},
\widetilde{\Ax})$ on the relatively compact subset
$\Phi^{-1}(U) \setminus \widetilde{\vD}$.
\end{Lemma}
\begin{proof}
We fix a strongly adapted local chart $(U,(x,y,z))$ at the point $p_1$ and let
$$\Phi : \widetilde{\cM} \rightarrow \cM$$
be a distinguished vertex blowing-up at $p_1$, as 
defined in subsection~\ref{subsection-distingvertblowup}.
The result follows immediately from Proposition~\ref{prop-distingvertbl}.
\end{proof}
\begin{Corollary}
Let $(\cM,\Ax)$ be controlled singularly foliated manifold and let $U \subset M$ be a relatively compact subset.
Then, there exists a finite sequence of blowing-ups
$$
\cM = \cM_0 \buildrel\Phi_1\over \longleftarrow   \cM_1 \longleftarrow \cdots \buildrel  \Phi_k \over \longleftarrow\cM_k
$$
and an axis $\Ax_k$ for $\cM_k$ such that $(\cM_k,\Ax_k)$ is equireducible outside the divisor, when restricted
to $(\Phi_k \circ \cdots \circ \Phi_1)^{-1}(U)$.
\end{Corollary}
Let $(\cM,\Ax)$ be a singularly foliated manifold which is equireducible outside the
divisor. Then, each connected component $Y$ of $\NElem \setminus \vD$ is a
smooth one dimensional analytic curve.

In this case, we define the {\em generic virtual height} for $(\cM,\Ax)$ along $Y$ as the natural number
$$
\vheight(\cM,\Ax,Y) := \vheight_G(\cM,\Ax,p),
$$
where $p$ is an arbitrary point on $Y$. 
Using Corollary~\ref{corollary-eqopen},  one concludes that $\vheight(\cM,\Ax,Y)$ is
independent of the choice of the particular point $p \in Y$.

\subsection{Extending the Invariant to $\NElem \setminus \vD$}
In this section, let us assume that $(\cM,\Ax)$ is equireducible outside the divisor.
In particular, the virtual height function 
$\vheight(\cM,\Ax,\cdot): \NElem \cap \vD \rightarrow \cN$
can be extended to  all the set $\NElem$ by setting 
$$\vheight(\cM,\Ax,\cdot) := \vheight_G(\cM,\Ax,\cdot),\quad\mbox{ on }\NElem \setminus \vD$$ 
We denote such function shortly by $\vheight(p)$.  

The {\em stratum of virtual height h} is the subset 
$$
S_h = \{p \in \NElem \mid \vheight(p) = h\}
$$
\begin{Lemma}\label{lemma-closureYheight}
Given a connected equireducible curve $Y \subset \NElem \setminus \vD$, let
$\overline{Y} \subset \NElem$ the smallest closed analytic subset which contains $Y$.
Then, for each point $p \in \overline{Y} \cap \vD$, we have $\vheight(Y) \le \vheight(p)$.
\end{Lemma}
\begin{proof}
Let $(U,(x,y,z))$ be a stable adapted chart at $p$. Firstly, suppose that the point $p$ is such that 
$Y_p = \overline{Y} \cap U$ (i.e.\ $\overline{Y}$ locally coincides with the local blowing-up at $p$).  
Then, it follows from the same argument used in the proof of 
Lemma~\ref{lemma-preservesDelta1}.
that $\vheight(Y) = \vheight(p)$.

Suppose now that $Y_p \ne \overline{Y}\cap U$ and assume, by absurd, that 
$$
\vheight(Y) > \vheight(p)
$$
We make the local blowing-up $\Phi:\widetilde{\cM}\rightarrow \cM\cap U$ for
$(\cM,\Ax)$ at $p$ and look at the
strict 
transform $Y^\prime$ of the curve 
$Y$.  The closure of such curve $\overline{Y}^\prime$ necessarily intersects the exceptional divisor
$\widetilde{D} = \Phi^{-1}(Y_p)$ in at least one nonelementary point $\widetilde{p}$.  Moreover,
$$
\vheight(Y^\prime) = \vheight(Y) > \vheight(p) \ge \vheight(\widetilde{p})
$$ 
as a consequence of Theorem~\ref{theorem-localresolution}. 

Let us now set $p := \widetilde{p}$, $Y := Y^\prime$ and iterate the process. The
Theorem~\ref{theorem-localresolution} implies that after some finite number of iterations, 
we fall into a situation where $\widetilde{D}$ has no nonelementary points.  
This is a contradiction.
\end{proof}
\begin{Proposition}\label{prop-describesSh}
The function $\vheight: \NElem \rightarrow \cN$ is upper semicontinuous.  
\end{Proposition}
\begin{proof}
This is an immediate consequence of 
Proposition~\ref{prop-delta1D>0=0} and Lemma~\ref{lemma-closureYheight}.
\end{proof}
The Newton invariant
$\Inv(\cM,\Ax,p)$ can also be defined globally on $\NElem$.  We denote it shortly by $\Inv(p)$
and remark that the following relations hold:
\begin{enumerate}[(1)]
\item If $p \in S_h \setminus \vD$ then $\Inv(p) = (h,0,0,0,0,0)$
\item If $p \in S_h \cap \vD$ is such that $\#\iota_p = 2$ then $\Inv(p) = (h,1,h,1,\cdot,\cdot)$. 
\item If $p \in S_h \cap \vD$ is such that $\#\iota_p = 1$ then 
$$\mbox{ Either }\quad\Inv(p) = (h,0,m_3,0,\cdot,\cdot)\quad \mbox{ or }\quad  \Inv(p) = (h,1,h,0,\cdot,\cdot)$$
for some $m_3 \ge h$.
\end{enumerate}
As a consequence of such remark, combined with Propositions~\ref{Prop-invuppsemicont} and
\ref{prop-describesSh}, we conclude that:
\begin{Proposition}\label{Prop-invuppsemicont2}
The function $\Inv: \NElem \rightarrow \cN^6$ is 
upper semicontinuous (for the lexicographical ordering on $\cN^6$).
\end{Proposition}
\subsection{Extended Center, Bad Points and Bad Trees}\label{subsect-extcent}
A controlled singularly foliated manifold $(\cM,\Ax)$ will be called a {\em restriction} if
it is given by the restriction of a controlled singularly foliated manifold $(\cM^\prime,\Ax^\prime)$
to some relatively compact open subset $U$ of the ambient space $M^\prime$.

In this subsection, we shall suppose that
$(\cM,\Ax)$ is a restriction and, moreover, that it is 
equireducible outside the
divisor.  In particular, this implies (by the upper semicontinuity of the height function) that  
$$
\vheight_{\max} := \sup \{\vheight(p) \mid p \in \NElem\}.
$$
is a finite natural number and that the divisor list $\mainlist\in \lists$ has a finite length.

The set $S_{\vheight_{\max}} = \{p \in \NElem \mid \vheight(p) = \vheight_{\max}\}$ 
will be called {\em stratum of maximal height}.
\begin{Lemma}\label{Lemma-describeSvheight}
The stratum of maximal height  $S_{\vheight_{\max}}$ is a closed analytic
subset of $\NElem$.  
Moreover, 
$S_{\vheight_{\max}} \cap \vD$
is an union of isolated points and 
closed analytic curves which have normal crossings with the divisor.
\end{Lemma}
\begin{proof}
The set $S_{\vheight_{\max}}$ is closed by the upper semicontinuity 
of the function $\vheight$.
Moreover, it follows from items (i) and (ii) of Proposition~\ref{prop-delta1D>0=0} that
the set $S_{\vheight_{\max}} \cap \vD$ is locally smooth at each point 
$p \in S_{\vheight_{\max}} \cap \vD$.  
\end{proof}
The {\em extended center} associated to a point $p \in \NElem$ is
the smallest closed analytic subset $\overline{Y}_p \subset \NElem$
which coincides with the local blowing-up center $Y_p$ in a neighborhood of $p$.

\begin{figure}[htb] 
\psfrag{D}{\small $D$}
\psfrag{p}{\small $p$}
\psfrag{Yp}{\small $Y_p$}
\psfrag{Yp=Ybp}{\small $Y_p\equiv \overline{Y}_p$}
\psfrag{Ybp}{\small $\overline{Y}_p$}
\begin{center} 
\includegraphics[height=7.5cm]{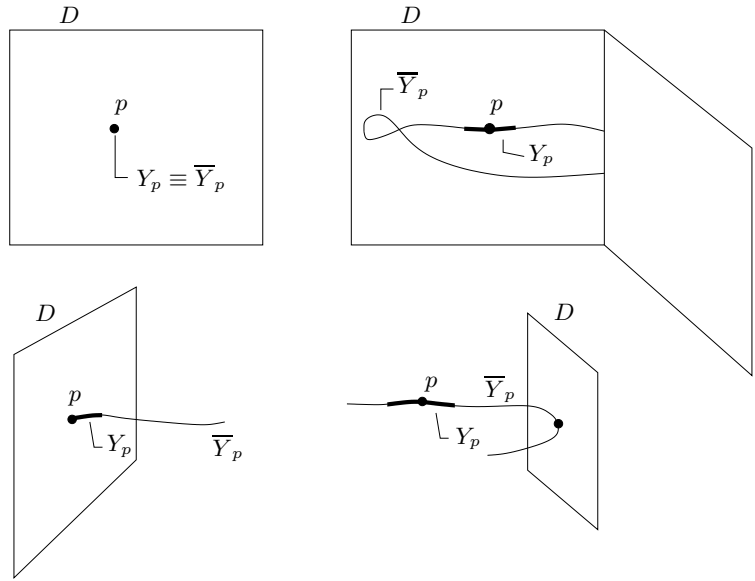}
\end{center}
\caption{The extended centers.}
\end{figure}
\begin{Remark} 
For instance, if $Y_p = \{p\}$ then $\overline{Y}_p = \{p\}$.  On the other 
hand, if $Y_p$ is contained in some irreducible divisor component $D \subset \vD$ then 
$\overline{Y}_p$ is entirely contained in $D$.
\end{Remark}
We say that $\overline{Y}_p$ is a {\em divisorial center} if $\overline{Y}_p \subset \vD$.
\begin{Lemma}
Let $p \in S_{\vheight_{\max}}$ be such that the extended center
$\overline{Y}_p$ is divisorial.  Then, $\overline{Y}_p$ is either
an isolated point or a smooth analytic curve which has normal 
crossings with the divisor $\vD$.
\end{Lemma}
\begin{proof}
This is an obvious consequence of Lemma~\ref{Lemma-describeSvheight}, 
since the extended center $\overline{Y}_p$ is an irreducible closed 
analytic subset of $S_{\vheight_{\max}} \cap \vD$.
\end{proof}
We say that the extended center 
$\overline{Y}_p$ is {\em permissible} at a point $q \in \overline{Y}_p$ if
$$\overline{Y}_q \equiv \overline{Y}_p$$
In other words, $\overline{Y}_p$ is permissible at $q$ if  
the the local blowing-up center for $(\cM,\Ax)$ at $q$ locally coincides with $\overline{Y}_p$. 
A point $q \in \overline{Y}_p$ is called a {\em bad point} if $\overline{Y}_p$ is not permissible at $q$.

We denote by $\Bad(p)$ the set of all bad points in $\overline{Y}_p$.  
We shall say that the extended center $\overline{Y}_p$ is {\em globally permissible} if
$\Bad(p) = \emptyset$.
\begin{Proposition}\label{prop-globallyperm}
Fix a point $p \in S_{\vheight_{\max}}$.
\begin{enumerate}[(i)]
\item If $Y_p = \{p\}$ then $\overline{Y}_p$ is globally permissible.
\item Suppose that $\overline{Y}_p$ is a smooth curve contained in some 
divisor component $D_i \subset \vD$.  
Then each point of $\Bad(p)$ is contained in the intersection $D_i \cap D_j$, for some index $j > i$.
\item If $p \in S_{\vheight_{\max}} \setminus \vD$ is an equireducible point
      then $\Bad(p)$ is a subset of $\overline{Y}_p \cap \vD$.
\end{enumerate}
\end{Proposition}
\begin{proof}
The item (i) is trivial.  To prove item (ii), notice that for each point 
$q \in \overline{Y}_p$, the following three situations can appear
$$
(a)\; \iota_q = [i],\quad (b)\; \iota_q = [i,k],\quad\mbox{ or }\quad (c)\; \iota_q = [j,i]
$$
for some indices $k < i < j$. In the cases (a) and (b), it is cleat that 
the extended center $\overline{Y}_p$ is permissible at $q$ because
$\Delta_1^q = \Delta_1^p > 0$ (by Lemma~\ref{lemma-preservesDelta1}).  
Therefore, a bad point of $\overline{Y}_p$ necessarily lies 
in the intersection of $D_i$ with some divisor $D_j$ of larger index.

The item (iii) is a direct consequence of the assumption that $(\cM,\Ax)$ is 
equireducible outside the divisor $\vD$.
\end{proof}
\begin{Corollary}\label{corollary-badpoint2}
For each point $p \in S_{\vheight_{\max}} \cap D_i$, the following properties hold:
\begin{enumerate}[(1)]
\item If $\#\iota_p = 2$ then the set $\Bad(p)$ has at most one point.  
\item If $\#\iota_p = 1$ and $\overline{Y}_p \subset D_i$ then the $\Bad(p)$ has at most two points.
\end{enumerate}
In both cases each point $q \in \Bad(p)$ is such that  
$\iota_q = [i,j]$, for some index $j > i$.
\end{Corollary}
\begin{proof}
This is a direct consequence of Proposition~\ref{prop-globallyperm} and 
the description of $S_{\vheight_{\max}} \cap \vD$ given by 
Lemma~\ref{Lemma-describeSvheight}.
\end{proof}
\begin{Lemma}\label{lemma-badpoint1}
Let $p \in S_{\vheight_{\max}} \setminus \vD$ be an 
equireducible point.  Then, for each point
$q \in \Bad(p)$, the associated local blowing-up center $Y_q$ is such that
$$
\overline{Y}_q \subset \vD
$$
i.e.\ $\overline{Y}_q$ is necessarily a divisorial center.
\end{Lemma}
\begin{proof}
Indeed, suppose by absurd that $Y_q$ is not a divisorial center and 
$\overline{Y}_q \ne \overline{Y}_p$.  

We fix a stable local chart $(U,(x,y,z))$ at $q$ and let 
$$\Phi: \widetilde{\cM} \rightarrow \cM \cap U$$ be the 
local blowing-up (with center $Y_q$) of $(\cM,\Ax)$ at $q$.  It follows 
from Propositions~\ref{prop-stableyblup} and \ref{prop-stablezblup} that each point 
$\widetilde{p} \in \Phi^{-1}(q)$ is such that either the 
Newton data is in final situation or
$\vheight(\widetilde{p}) < \vheight(q) = \vheight_{\max}$. 

On the other hand, the strict transform of $\overline{Y}_p$ under $\Phi$ 
contains at least one point of  
$\Phi^{-1}(q)$.  This is an absurd, since it
contradicts the fact that $\overline{Y}_p \subset S_{\vheight_{\max}}$.
\end{proof}
A {\em bad chain} is a (possibly infinite)
sequence of points $\{p_n\}_n$ which 
is contained in $S_{\vheight_{\max}}$
and is such that
$$
p_{n+1} \in \Bad(p_n), \; n \ge 0
$$
We shall say that a finite bad chain $\{p_0,\ldots,p_l\}$ is {\em complete} if 
$\Bad(p_l) = \emptyset$.
The number $l$ will be called the {\em length} of the complete bad chain.
\begin{Remark}\label{rem-badchaindivisors}
It follows from Lemma~\ref{lemma-badpoint1} and 
Corollary~\ref{corollary-badpoint2} that for 
a bad chain
$\{p_n\}_n$, we always have
$$\#\iota_{p_1} \ge 1 \quad \mbox{ and }\quad\iota_{p_n} = 2$$ 
for all $n \ge 2$.
\end{Remark}

\begin{figure}[htb] 
\psfrag{Dk}{\small $D_{i_0}$}
\psfrag{Di}{\small $D_{j_0}$}
\psfrag{Dj}{\small $D_{j_1}$}
\psfrag{Ds}{\small $D_{j_2}$}
\psfrag{p}{\small $p_0$}
\psfrag{Yp}{\small $\overline{Y}_{p_0}$}
\psfrag{q}{\small $p_1$}
\psfrag{Yq}{\small $\overline{Y}_{p_1}$}
\psfrag{r=Yr}{\small $p_2 \equiv \overline{Y}_{p_2}$}
\psfrag{k<i<j<s}{\small $i_0 < j_0 < j_1 < j_2$}
\begin{center} 
\includegraphics[height=6cm]{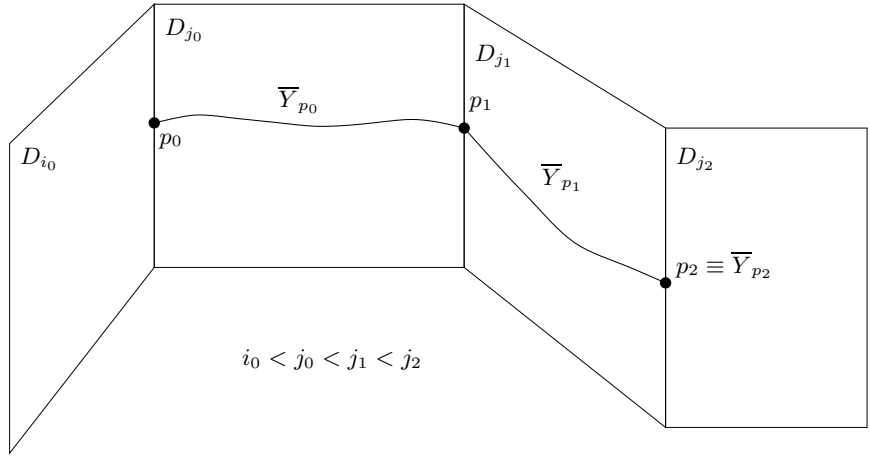}
\end{center}
\caption{The bad chain $\{p_0,p_1,p_2\}$ (here $\Bad(p_2) = \emptyset$).}
\end{figure}

\begin{Lemma}\label{lemma-finiteness}
Each bad chain has a finite number of points.
\end{Lemma}
\begin{proof}
By the Remark~\ref{rem-badchaindivisors}, each bad chain $\{p_n\}_n$ is such that, 
$\#\iota_{p_n} = 2$ for all $n \ge 2$.  
Moreover, if we write $\iota_{p_n} = [j_n,i_n]$ then
$$
i_n < j_n  = i_{n+1} < j_{n+1} = i_{n+2} < j_{n+2} = \cdots
$$
and therefore the indices $\{i_n\} \subset \mainlist$ form a 
strictly increasing sequence (where $\mainlist$ is the list of divisor indices).
Since we supposed that $(\cM,\Ax)$ is a restriction, the list $\mainlist$ is necessarily 
finite. The Lemma is proved.
\end{proof}
Given a point $p \in S_{\vheight_{\max}}$, the {\em bad chains starting at
$p$} is the set $\Badch(p)$ of all complete
bad chains $\{p_n\}_{n\ge 0}$ such that $p_0 = p$.  
\begin{Remark}
It follows from Lemma~\ref{lemma-finiteness} that the 
set $\Badch(p)$ has only a finite number of elements.
\end{Remark}
More generally, given a finite set of points $P  \subset 
S_{\vheight_{\max}}$, we define the
{\em $P$-bad chain} as the union  
$$
\Badch(P):= \bigcup_{p \in P}{\Badch(p)}
$$
of all bad chains starting at points in $P$. Associated to $\Badch(P)$, let us consider 
a directed graph $T = (V,E)$ defined as follows:
\begin{enumerate}
\item The set of vertices $V$ corresponds to the set of points of all bad chains starting at
      $P$ (for simplicity, we identify each element of $V$ with the
      corresponding point in the bad chain).
\item The directed edge $q \rightarrow r$ belongs to the set of edges $E$ if 
there exists a bad chain $\{p_n\}_{n=0}^l$ in $\Badch(P)$ such that
$$
p_i = q, \quad p_{i+1} = r
$$
for some $0 \le i \le l-1$.
\end{enumerate}
\begin{Lemma}
The graph $T = (V,E)$ is a directed tree.
\end{Lemma}
\begin{proof}
We need to prove that $T$ has no cycles.  Let us suppose, 
by absurd, that there exists a cycle in $T$
$$
q_0 \rightarrow q_1 \rightarrow \cdots \rightarrow q_r \rightarrow q_{r+1} = q_0
$$
where $q_{n+1} \in \Bad(q_n)$ for each $0 \le n \le r$.  Let us write $\iota_{q_n} = [i_n]$ 
(if $\#\iota_{q_n} = 1$) and $\iota_{q_n} = [j_n,i_n]$ (if $\#\iota_{q_n} = 2$).

First of all, suppose that the extended center $\overline{Y}_{q_0}$ is 
divisorial (i.e.\ contained in $\vD$).
Then, it follows from Corollary~\ref{corollary-badpoint2} that the sequence 
$i_1<i_2<\cdots<i_n$ is strictly increasing.  No cycle can appear.

Suppose now that $\overline{Y}_{q_0}$ is not divisorial. Then, Lemma~\ref{lemma-badpoint1}
implies that $\overline{Y}_{q_n}$ is divisorial, for all $n\ge 1$.  This 
contradicts the fact that 
$q_{n+1} = q_0$.
\end{proof}
\begin{Definition}
The directed 
tree $T = (V,E)$ defined above will be called {\em Bad tree} associated to 
$P$.  We shall denote it
by $\TreeBd(P)$.
\end{Definition}
From now on, we adopt the usual nomenclature for trees.  Thus, a {\em branch} is any
succession of points and directed edges, 
$$p_0 \rightarrow p_1 \rightarrow \cdots \rightarrow p_k$$
In this case, the number $k$ will be called the {\em length} of the branch.
A point $q \in V$ is a called a {\em descendant} of a point $p$ if there exists a branch of positive length
as above such that $p_0 = p$ and $p_k = q$.
A point $q$ will be called a
{\em terminal} if it has no descendants in the three.
\begin{Remark}
For each terminal point $q \in \TreeBd(P)$,  the extended
center $\overline{Y}_q$ is globally
permissible (because $\Bad(q) = \emptyset$). 
\end{Remark}
The {\em maximal length} of a bad tree is the length $L(\TreeBd(P)) \in \cN$ of 
the longest branch of $\TreeBd(P)$.

Let $F \subset \Badch(P)$ be the set of all terminal points which lie in branches
of maximal length (i.e.\ those branches of $\TreeBd(P)$ which have length $L(\TreeBd(P))$).  
We define the {\em maximal final invariant} of $\TreeBd(P)$ as
$$
\Inv(\TreeBd(P)) := \max{}_\lex \{ \Inv(q) \mid q \in F\},
$$
where the maximum is taken for the lexicographical ordering in $\cN^6$.
The {\em maximal final locus} as the finite set of points
$$
\Loc(\TreeBd(P)) := \{ q \in F \mid \Inv(q) = \Inv(\TreeBd(P))\}
$$
Finally, we define the {\em multiplicity} of the bad tree $\TreeBd(P)$ as the vector 
\begin{equation}\label{def-multbadtree}
\Mult(\TreeBd(P)) := (\; L(\TreeBd(P)) ,\;\Inv(\TreeBd(P)), \;
\#\Loc(\TreeBd(P))\; ) \in \cN^8
\end{equation}
where $\#\Loc(\TreeBd(P))$ is the cardinality of the set $\Loc(\TreeBd(P))$.

\subsection{Maximal Invariant Locus and Global Multiplicity}
In this subsection, we continue to assume that 
$(\cM,\Ax)$ is a a controlled singularly foliated manifold which is a restriction 
(see subsection~\ref{subsect-extcent}) and equireducible outside the
divisor.  

Therefore, the maximal of the invariant $\Inv(p)$,  
$$
\Inv_{\max}(\cM,\Ax) := \sup{}_{\lex} \{\Inv(p) \mid p \in \NElem\}.
$$
is a finite vector in $\cN^6$. If $(\cM,\Ax)$ is clear from the context, we denote such number
simply by $\Inv_{\max}$.  The subset
$$S_{\Inv_{\max}} := \{p \in \NElem \mid \Inv(p) = \Inv_{\max}\} 
\subset S_{\vheight_{\max}}$$
will be called {\em maximal invariant stratum} of $(\cM,\Ax)$.

Consider the subsets $\vD_i := \{p \in \NElem \mid \#\iota_p = i\}$, for $i = 0,1,2$.
We establish the following definitions:
\begin{enumerate}
\item We say that $S_{\Inv_{\max}}$ is of {\em $2$-boundary type} if 
$S_{\Inv_{\max}} \cap \vD_2  \ne \emptyset$
\item We say that $S_{\Inv_{\max}}$ is of {\em $1$-boundary type} if
$S_{\Inv_{\max}} \cap \vD_2  = \emptyset$ and $S_{\Inv_{\max}} 
\cap \vD_1  \ne \emptyset$;
\item We say that $S_{\Inv_{\max}}$ is of {\em $0$-boundary type} if 
$S_{\Inv_{\max}} \cap (\vD_1 \cup \vD_2)  = \emptyset$ and 
$S_{\Inv_{\max}} \cap \vD_0  \ne \emptyset$. 
\end{enumerate}
Using such classification, the following result establishes some properties of $S_{\Inv_{\max}}$:
\begin{Lemma}
The maximal invariant stratum has the following properties:
\begin{enumerate}[(i)]
\item If $S_{\Inv_{\max}}$ is of $2$-boundary type then 
$S_{\Inv_{\max}} \subset \vD_2$.
\item If $S_{\Inv_{\max}}$ is of $1$-boundary type then 
$S_{\Inv_{\max}} \subset \vD_1$.
\item If $S_{\Inv_{\max}}$ is of $0$-boundary type then
$S_{\Inv_{\max}} = S_{\vheight_{\max}} \subset   \vD_0$.
\end{enumerate}
\end{Lemma}
\begin{proof}
The result is a direct consequence of the definition of $\Inv$ and Remark~\ref{rem-remarkinv(p)}.
\end{proof}
In the next Lemmas, we give a more detailed description of $S_{\vheight_{\max}}$:
\begin{Lemma}\label{lemma-2boundary}
A $2$-boundary type
$S_{\Inv_{\max}}$ is  
formed by a union of finite number of distinct points $\{p_1,\ldots,p_m\}$.  
\end{Lemma}
\begin{proof}
It follows immediately from the description of the set $S_{\Inv_{\max}}
\cap \vD$
which is given in Lemma~\ref{Lemma-describeSvheight}.
\end{proof}

\begin{figure}[htb] 
\psfrag{D}{\small $\vD$}
\psfrag{p1}{\small $p_1$}
\psfrag{p2}{\small $p_2$}
\psfrag{p3}{\small $q$}
\psfrag{Ybp1}{\small $\overline{Y}_{p_1}$}
\psfrag{Ybp3}{\small $\overline{Y}_{q}$}
\begin{center} 
\includegraphics[height=5cm]{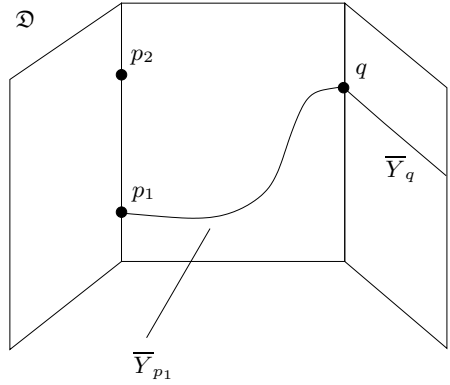}
\end{center}
\caption{The $2$-boundary Maximal Invariant Stratum. Here $\Bad(p_1) = \{q\}$.}
\label{fig-boundarycases1}
\end{figure}

\begin{Lemma}\label{lemma-1boundary}
A $1$-boundary type
$S_{\Inv_{\max}}$ is formed by a
finite union of distinct closed analytic sets
$$
Y_1 \cup \cdots \cup Y_r \quad \cup \quad \{p_1,\ldots,p_m\} \quad \cup \quad \{q_1,\ldots,q_n\}
$$
for some natural numbers $r,m,n \in \cN$,  such that the following conditions holds:
\begin{enumerate}[(i)]
\item $Y_1,\ldots,Y_r$ are globally permissible one-dimensional extended
      centers 
contained in $\vD_1$.
\item Each $p_i$ is
an isolated point of $S_{\Inv_{\max}} \cap \vD_1$ such that $\overline{Y}_{p_i}$  
is a globally permissible extended center contained in $\vD_1$.
\item Each $q_j$ is an isolated point of $S_{\Inv_{\max}} \cap \vD_1$
such that the extended center $\overline{Y}_{q_j}$ is not divisorial.
\end{enumerate}
\end{Lemma}
\begin{proof}
It follows immediately from the description of the set $S_{\vheight_{\max}}
\cap \vD$
which is given in Lemma~\ref{Lemma-describeSvheight} and the
assumption that $S_{\Inv_{\max}} \cap \vD_2 = \emptyset$.
\end{proof}

\begin{figure}[htb] 
\psfrag{D}{\small $\vD$}
\psfrag{p1}{\small $p_1$}
\psfrag{Y1}{\small $Y_1$}
\psfrag{p2}{\small $w$}
\psfrag{q1}{\small $q_1$}
\psfrag{s1}{\small $s_1$}
\psfrag{Ybq1}{\small $\overline{Y}_{q_1}$}
\begin{center} 
\includegraphics[height=5cm]{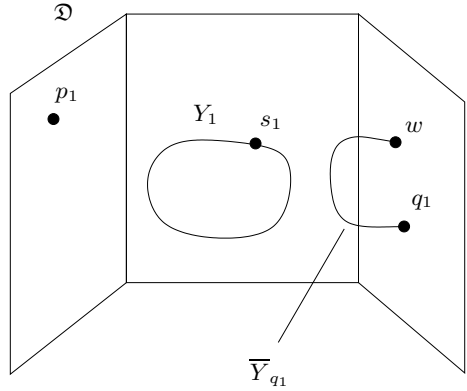}
\end{center}
\caption{The $1$-boundary Maximal Invariant Stratum. Here, $\Bad(q_1) = \{w\}$.}
\label{fig-boundarycases2}
\end{figure}

\begin{Lemma}\label{lemma-0boundary}
A $0$-boundary type
$S_{\Inv_{\max}}$ is formed by a
finite union $Y_1 \cup \cdots \cup Y_r$ of distinct globally permissible 
one-dimensional extended centers which are contained in 
$\vD_0$.
\end{Lemma}
\begin{proof}
It follows immediately from the assumptions that
$S_{\Inv_{\max}} \cap (\vD_2 \cup \vD_1) = \emptyset$ and that $(\cM,\Ax)$ is
equireducible outside the divisor.
\end{proof}

\begin{figure}[htb] 
\psfrag{D}{\small $\vD$}
\psfrag{Y1}{\small $Y_1$}
\psfrag{Y2}{\small $Y_2$}
\psfrag{s1}{\small $s_1$}
\psfrag{s2}{\small $s_2$}
\begin{center} 
\includegraphics[height=5cm]{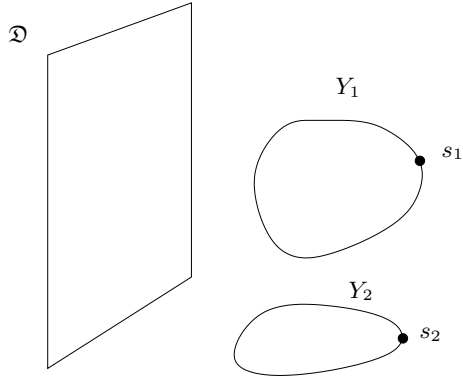}
\end{center}
\caption{The $0$-boundary Maximal Invariant Stratum.}
\label{fig-boundarycases3}
\end{figure}

Based on the above description of $S_{\Inv_{\max}}$, we state the following definition:
A {\em maximal point locus} of $(\cM,\Ax)$ is a finite collection of distinct points 
$$P_{\max} \subset S_{\Inv_{\max}}$$
which is obtained as follows:
\begin{enumerate}[(i)]
\item If $S_{\Inv_{\max}}$ is a $2$-boundary maximal invariant stratum, 
\begin{equation}\label{def-maxptlocus2}
 P_{\max} : =  \{p_1,\ldots,p_m\} = S_{\Inv_{\max}}
\end{equation} 
where $\{p_1,\ldots,p_m\}$ are given by Lemma~\ref{lemma-2boundary} (see figure~\ref{fig-boundarycases1}). 
\item If $S_{\Inv_{\max}}$ is a $1$-boundary maximal invariant stratum, 
\begin{equation}\label{def-maxptlocus1}
P_{\max} = \{s_1,\ldots,s_r, p_1,\ldots,p_m,q_1,\ldots,q_n\} 
\end{equation} 
where each $s_i$ is an arbitrary point contained in the curve $Y_i$, and the set
$\{Y_1,\ldots,Y_r,p_1,\ldots,p_m,q_1,\ldots,q_n\}$ is given by Lemma~\ref{lemma-1boundary}
(see figure~\ref{fig-boundarycases2}). 
\item If $S_{\Inv_{\max}}$ is a $0$-boundary maximal invariant stratum,  
\begin{equation}\label{def-maxptlocus0}
P_{\max} : = \{s_1,\ldots,s_r\}
\end{equation} 
where each $s_i$ is an arbitrary point contained in the curve $Y_i$, and the set
$\{Y_1,\ldots,Y_r\}$ is given by Lemma~\ref{lemma-0boundary}
(see figure~\ref{fig-boundarycases3}). 
\end{enumerate}
The {\em global multiplicity} associated to $(\cM,\Ax)$ is the vector
$$
\Mult(\cM,\Ax) := (\; \Inv_{\max}(\cM,\Ax), \Mult(\TreeBd(P_{\max}))\; )
$$
where $P_{\max}$ is a maximal point locus of $(\cM,\Ax)$ and 
$\Mult(\TreeBd(P_{\max}))$ is the multiplicity of 
the bad tree $\TreeBd(P_{\max})$ (see (\ref{def-multbadtree})).
\begin{Remark}
It is obvious that the value of $\Mult(\cM,\Ax)$ is independent of the 
choice of the points $s_i \in Y_i$ ($i = 1,\ldots,r$) which is 
made in (\ref{def-maxptlocus1}) and (\ref{def-maxptlocus0}).  
\end{Remark}
The  globally permissible extended
center $Y = \overline{Y}_q$ which associated to a terminal 
point $q \in \Loc(\TreeBd(P_{\max}))$ will be called a 
{\em blowing-up center} for $(\cM,\Ax)$.
\begin{Proposition}\label{prop-canblowup}
Let $Y \subset M$ be a blowing-up center for $(\cM,\Ax)$.  Then, there exists a weight-vector $\bfomega \in \cN^3$ 
such that the following properties holds:
\begin{enumerate}[(i)]
\item For each point $p \in Y \setminus \vD$ and each equireduction chart $(U_p,(x_p,y_p,z_p))$ for $(\cM,\Ax)$ at $p$,
$$
\bfomega = \omega_p \quad \mbox{ and }\quad  Y \cap U_p = Y_p
$$
(i.e.\ $\omega$ is the local weight vector at $p$ and $Y \cap U_p$ is the local blowing-up center)
\item For each point $p \in Y \cap \vD$ and each stable adapted chart $(U_p,(x_p,y_p,z_p))$ for $(\cM,\Ax)$ at $p$,
$$
\bfomega = \omega_p \quad \mbox{ and }\quad Y \cap U_p = Y_p
$$
\end{enumerate}
Moreover, the collection of charts $\{(U_p,(x_p,y_p,z_p))\}_{p \in A}$ as defined above
is a $\bfomega$-weighted trivialization atlas for
$Y \subset M$.
\end{Proposition}
\begin{proof}
The items (i) and (ii) follow from the fact that $Y$ is a globally permissible center.

In order to prove the last statement, we need to prove that the transition between two charts in the above trivialization, 
say
$(U_p,(x_p,y_p,z_p))$ and $(U_q,(x_q,y_q,z_q))$, preserves the $\bfomega$-quasihomogeneous structure on $\cR^3$.

If $Y$ is a single point, it suffices to apply Proposition~\ref{prop-preservesstructure}.
If $Y$ is a smooth curve, then we can locally write
$$
Y = (x_p = z_p = 0) \quad\mbox{ or }\quad Y = (y_p = z_p = 0)
$$
In these cases we claim that, for each point $q \in Y$ which is sufficiently near $p$, the respective translated chart 
$$
(\widetilde{x}_q,\widetilde{y}_q,\widetilde{z}_q) = (x_p,y_p - \rho,z_p)  \quad\mbox{ or }\quad
(\widetilde{x}_q,\widetilde{y}_q,\widetilde{z}_q) = (x_p - \rho,y_p,z_p)
$$
is a stable local chart (or an equireduction chart)
for $(\cM,\Ax)$ at $q$ (for some conveniently chosen constant $\rho \in \cR$).

Indeed, such claim can be proved by easy modifications in the proofs of  
Lemmas~\ref{lemma-preservesDelta1} and \ref{lemma-equiredstable}.

Using the claim, combined with Propositions~\ref{prop-preservesstructure} and 
Remark~\ref{rem-preservesstructure2}, we conclude that  
$\{(U_p,(x_p,y_p,z_p))\}_{p \in A}$ is a trivialization of $Y$ which 
preserves the $\bfomega$-quasihomogeneous structure on $\cR^3$.
\end{proof}
It follows from the Proposition~\ref{prop-blowupmanif}  that we can define the 
$\bfomega$-weighted blowing-up of $\cM$ with center on $Y$
\begin{equation}\label{goodblowup}
\Phi : \widetilde{\cM} \rightarrow \cM
\end{equation}
(with respect to the trivialization given by 
Proposition~\ref{prop-canblowup})  The transformed singularly foliated manifold
$\widetilde{\cM}$ is defined according to subsection~\ref{subsect-blupsingfolmanif}.

The map above map will be called a {\em good blowing-up} for $(\cM,\Ax)$.
\subsection{Global Reduction of Singularities}
To state our next result, we recall that a controlled singularly foliated manifold $(\cM,\Ax)$ is 
called a {\em restriction} if it is defined by restriction of a controlled singularly foliated manifold 
$(\cM^{\, \prime},\Ax^\prime)$ to some relatively compact open subset of the ambient space.
\begin{Theorem}\label{theorem-goodblowup}
Let $(\cM,\Ax)$ be a controlled singularly foliated manifold
which is a restriction and equireducible outside the divisor. Let
$$
\Phi:  \widetilde{\cM} \rightarrow \cM
$$
be a good blowing-up for $(\cM,\Ax)$.  
Then, either $\widetilde{\cM}$ is an elementary singularly foliated manifold or there
exists an axis $\widetilde{\Ax} = (\widetilde{A},\widetilde{\vZ})$ for $\widetilde{\cM}$ such that:
\begin{enumerate}[(i)]
\item The singularly foliated manifold $(\widetilde{\cM},\widetilde{\Ax})$ is a restriction;
\item $(\widetilde{\cM},\widetilde{\Ax})$ is equireducible outside the divisor, and 
\item $\Mult(\widetilde{\cM},\widetilde{\Ax}) <_\lex \Mult(\cM,\Ax)$.
\end{enumerate} 
\end{Theorem}
\begin{proof}
Assume that the set $\NElem(\widetilde{\cM})$ of nonelementary points of $\widetilde{\cM}$ in nonempty.  

Let us denote by $Y$ be the blowing-up center and recall the fact
that $\Phi$ locally coincides with the local blowing-up for $(\cM,\Ax)$ at each point $p \in Y$.  From this, 
we conclude from 
Theorems~\ref{theorem-localresolution} and \ref{theorem-localresolutiongeneric} that there exists an axis 
$\widetilde{\Ax} = (\widetilde{A},\widetilde{\vZ})$ for
$\cM$ (obtained by analytic glueing) such that $(\widetilde{\cM},\widetilde{\Ax})$ is a controlled 
singularly foliated manifold such that:
\begin{enumerate}[(i)]
\item $(\widetilde{\cM},\widetilde{\Ax})$ is a restriction, and
\item $(\widetilde{\cM},\widetilde{\Ax})$ is equireducible outside the divisor.
\end{enumerate}  
Moreover, from Theorems~\ref{theorem-localresolution} we conclude that 
$$
\Inv_{\max}(\widetilde{\cM},\widetilde{\Ax}) \le_\lex \Inv_{\max}(\cM,\Ax).
$$
If the inequality is strict, we are done. Otherwise, let us choose a maximal point locus
$\widetilde{P}_{\max}$ for $(\widetilde{\cM},\widetilde{\Ax})$.  

Using again the Theorems~\ref{theorem-localresolution} and \ref{theorem-localresolutiongeneric}, we can write
$$
\widetilde{P}_{\max} = \Phi^{-1}(P_{\max} \setminus Y)
$$
for some maximal point locus $P_{\max}$ of $\cM$.  Indeed, we have 
$\Inv(\widetilde{\cM},\widetilde{\Ax},\widetilde{q}) <_\lex \Inv_{\max}(\cM,\Ax)$ for each 
point $\widetilde{q} \in \Phi^{-1}(Y)$ and therefore $S_{\Inv_{\max}} \cap \Phi^{-1}(Y) = \emptyset$.
  
For shortness, let us write the respective multiplicities of the bad trees $\TreeBd(P_{\max})$ and 
$\TreeBd(\widetilde{P}_{\max})$ simply as
$$\Mult(\TreeBd(P_{\max})) = (L,I,\#) \quad\mbox{ and }\quad 
\Mult(\TreeBd(\widetilde{P}_{\max})) = (\widetilde{L}, \widetilde{I}, \widetilde{\#})$$  
Then, we need to prove that
$$
(L,I,\#) <_\lex (\widetilde{L}, \widetilde{I}, \widetilde{\#})
$$  
First of all, we claim that $L \le \widetilde{L}$.  Indeed, it suffices to study how each branch 
\begin{equation}\label{branchbadthree}
p_0 \rightarrow p_1 \rightarrow \cdots \rightarrow p_l, \qquad (p_0 \in P_{\max})
\end{equation}
of the bad tree $\TreeBd(P_{\max})$ is transformed by the blowing-up.

If $p_i \notin Y$ for all $i = 0,\ldots,l$ then this branch is mapped isomorphically to a branch
of length $l$ of the new bad tree 
$\TreeBd(\widetilde{P}_{\max})$.

Now, suppose that there exists an index $0 \le i \le l$ such that 
$$\{p_0,\ldots,p_{i-1}\} \cap Y = \emptyset \quad \mbox{ and }\quad p_i \in Y$$
If $i = 0$, it is immediate to see that the branch is completely destroyed.
So, we suppose that $i \ge 1$. 
It follows from the Propositions~\ref{prop-stablexblup}, \ref{prop-stableyblup} and \ref{prop-stablezblup} 
that each 
nonelementary point $\widetilde{q} \in \Phi^{-1}(p_i)$ 
satisfies one of the following conditions:
\begin{enumerate}[(1)]
\item $\vheight(\widetilde{\cM},\widetilde{\Ax},\widetilde{q}) < \vheight(\cM,\Ax,p_i)$, or
\item $\vheight(\widetilde{\cM},\widetilde{\Ax},\widetilde{q}) = \vheight(\cM,\Ax,p_i)$ and 
$\Delta_1^{\widetilde{q}} = 0$
\end{enumerate}
where $\Delta^{\widetilde{q}} = (\Delta_1^{\widetilde{q}},\Delta_2^{\widetilde{q}})$ is the vertical displacement 
vector associated to $\widetilde{q}$.  

As a consequence, the points lying in the case (2) are isolated points of 
$\widetilde{S}_{\vheight{\max}} \cap \widetilde{D}$ (where 
$\widetilde{D} := \Phi^{-1}(Y)$ is the exceptional divisor of the blowing-up and $\widetilde{S}_{\vheight{\max}}$
is the stratum of maximal height for $(\widetilde{\cM},\widetilde{\Ax})$).

Now, observe that the strict transform of the extended center
$\overline{Y}_{p_{i-1}}$ intersects the set $\Phi^{-1}(p_i)$ in a unique point $\widetilde{q} \in 
S_{\vheight{\max}}$, which necessarily lies in case (2).  This immediately implies that 
$\Bad(\widetilde{q}) = \emptyset$ and therefore the branch (\ref{branchbadthree}) 
is mapped to a unique branch in 
$\TreeBd(\widetilde{P}_{\max})$, which has one of the following forms 
\begin{eqnarray*}
& \widetilde{p}_0 \rightarrow \widetilde{p}_1 \rightarrow \cdots \rightarrow \widetilde{p}_{i-1} \rightarrow   \widetilde{q},\quad \mbox{ or }\\
& \widetilde{p}_0 \rightarrow \widetilde{p}_1 \rightarrow \cdots \rightarrow \widetilde{p}_{i-1}
\end{eqnarray*}
where $\widetilde{p}_i = \Phi^{-1}(p_i)$ for $i = 0,\ldots,i-1$. In both cases, it is clear that the new branch has a 
length at most equal to the length of the original branch.  We have proved that $\widetilde{L} \le L$.

Let us suppose that $\widetilde{L} = L$.  Then, since the blowing-up creates no new branches of maximal length
$L$, it follows immediately from the Theorem of Local Resolution of
Singularities (Theorem~\ref{theorem-localresolution})  that 
$$\widetilde{I} \le_\lex I$$
It remains to prove that the conditions
$\widetilde{L} = L$ and $\widetilde{I} = I$
imply that $\widetilde{\#} < \#$.  To see this, it suffices to remark that the blowing-up satisfies the 
following properties:
\begin{enumerate}[(1)]
\item The blowing-up $\Phi$ creates no new branches of length $L$;
\item The center $Y$ contains at least one terminal point of a branch which has length exactly equal $L$.
\end{enumerate}
Applying again the Theorems~\ref{theorem-localresolution} and \ref{theorem-localresolutiongeneric}, 
we immediately conclude that $\widetilde{\#} < \#$.  
This completes the proof of the Theorem.
\end{proof}
\subsection{Proof of the Main Theorem}
We are now ready to prove the Main Theorem of this work.
\begin{proof}{\bf (of the Main Theorem)} Let
$\cM = (M,\emptyset,\emptyset,\linef_\chi)$ be the singularly foliated manifold associated to $\chi$ and  
let $\Ax$ be an axis for $\cM$, defined as in 
Proposition~\ref{proposition-definesaxis}. Given a relatively compact subset $U \subset M$,
we denote by $(\cM^\prime, \Ax^\prime)$ the restriction of $(\cM,\Ax)$ to $U$.

Using Lemma~\ref{Lemma-distinguishedprep}, we know that there exists a finite sequence of blowing-ups 
$$
(\cM^\prime, \Ax^\prime) = (\cM_0,\Ax_0) \longrightarrow  (\cM_1,\Ax_1) \longrightarrow \cdots 
\longrightarrow  (\cM_k,\Ax_k)
$$
such that the resulting singularly foliated manifold $(\cM_k,\Ax_k)$ is equireducible outside the divisor.

To finish the proof, it suffices to consider the controlled singularly foliated manifold $(\cM_k,\Ax_k)$ and
apply successively the Theorem~\ref{theorem-goodblowup}.
\end{proof}

\section{Appendix A: Faithfully Flatness of $\cC[[x,y,z]]$}
In the proof of the Stabilization of adapted charts, we need the following simple
consequence of the fact that $\cC\{x,y,z\}$ is a unique factorization domain and that its completion $\cC[[x,y,z]]$ is faithfully flat (see e.g.\ \cite{Ma}, sections 4.C and 24.A).
\begin{Lemma}\label{Lemma-uniquefactgen}
Let $I \subset \cC\{x,y,z\}$ be a nonzero radical ideal and let
$$
I^\prime = J_1^\prime \cap \cdots \cap J_k^\prime
$$
be the irreducible primary decomposition of the ideal $I^\prime = I\cC[[x,y,z]]$ in the ring of formal series
$\cC[[x,y,z]]$.  Then, each $J_i^\prime$
($i = 1,\ldots,k$)
can be written as $J_i^\prime = J_i \cC[[x,y,z]]$, for some prime ideal $J_i \subset \cC\{x,y,z\}$. 
\end{Lemma}
\begin{Corollary}\label{corollary-uniquefactz}
Let $H = (H_1,\ldots,H_r) \in \cR\{x,y,z\}^r$ be a nonzero germ of analytic map.
Suppose that we can write the factorization
$$
\left(
\begin{matrix}
H_1 \cr
H_2 \cr
\vdots \cr
H_r
\end{matrix}
\right)
=
(z - f(x,y))
\left(
\begin{matrix}
S_1 \cr
S_2 \cr
\vdots \cr
S_r
\end{matrix}
\right)
$$
where $f \in \cR[[x,y]]$ and $S_1,\ldots,S_r \in \cR[[x,y,z]]$.  Then,
necessarily $f \in \cR\{x,y\}$ is an analytic germ.
\end{Corollary}
\begin{proof}
Indeed, the hypothesis imply that
the ideal $I = \mathrm{rad}(H_1,\ldots,H_r)$ is contained in the principal
ideal $J = (z - f(x,y))\cC[[x,y,z]]$,  In particular, $J$ is a member of the irreducible
primary decomposition of $I$ in $\cC[[x,y,z]]$.

Therefore, it suffices to apply the previous Lemma
to conclude that $f$ is necessarily an analytic germ.
\end{proof}
Given a nonzero natural number $a \in \cN$, consider now the ideal 
$$\widehat{I}_a = (x^a)\cR[[x,y,z]]
\subset \cR[[x,y,z]].$$  
The elements of the quotient ring $\widehat{R}_a =
\cR[[x,y,z]]/I_a$ are uniquely represented by polynomials in the $x$-variable $\cR[[y,z]]\, [x]$
whose degree is at most $a-1$. We let $R_a$ denote the image of $\cR\{x,y,z\}$ under the quotient map.
\begin{Corollary}\label{corollary-uniquefactzquot}
Let $([H_1],\ldots,[H_r]) \in R_a^r$ be a nonzero germ.
Suppose that we can write the factorization (in $\widehat{R}_a$)
$$
\left(
\begin{matrix}
[H_1] \cr
[H_2] \cr
\vdots \cr
[H_r]
\end{matrix}
\right)
=
(z - [f(x,y)])
\left(
\begin{matrix}
[S_1] \cr
[S_2] \cr
\vdots \cr
[S_r]
\end{matrix}
\right)
$$
where $[f] \in \widehat{R}_a$ and $[S_1],\ldots,[S_r] \in \widehat{R}_a$.  Then,
the germ $[f]$ necessarily lies in $R_a$.
\end{Corollary}
\begin{proof}
It suffices to use the previous Corollary.
\end{proof}
\section{Appendix B: Virtual Height}
Let us start with an elementary version of Descartes Lemma.  
\begin{Lemma}\label{lemma-multpolyn}
Let $Q(z)$ be a polynomial in $\cC[z]$ with $m$ nonzero monomials.  Then, the multiplicity of
$Q$ at a point $\xi \ne 0$ is at most $m - 1$.
\end{Lemma}
\begin{proof}
Given a polynomial $Q \in \cC[z]$, let $\mu(Q)$ be the multiplicity of $Q$ at the origin
(i.e.\ the greatest natural number $k$ such that $z^k$ divides $Q(z)$).

We consider the sequence of polynomials $Q_0(z),Q_1(z),\ldots$ which is inductively defined as follows.
$$
Q_0 = z^{-\mu(Q)} Q\quad \mbox{ and }\quad
Q_{i+1} = z^{-\mu(Q_i^\prime)}Q_i^\prime,\; \mbox{ for }i \ge 0
$$
(where ${ }^\prime = d/dz$).
By induction, we can easily prove that $Q$ has multiplicity $k$ at some point
$\xi \ne 0$ if and only if
$$
Q_0(\xi) = \cdots = Q_{k-1}(\xi) = 0.
$$
However, it follows from the hypothesis and the above construction that
$Q_{m-1}$ is necessarily a nonzero
constant.  Therefore, the maximum multiplicity of $Q$ at a point $\xi \ne 0$ is at most $m-1$.
\end{proof}
For the rest of this section, we shall adopt the following
notation.  Let $P(x_1,\ldots,x_n)$ be a
$n$ variable polynomial whose support is contained in the straight line
$$r : \ t \in \cR^+ \mapsto \bfp + t (\bfDelta,-1)$$
for some $\bfp = (p_1,\ldots,p_n) \in \cN^n$
and some vector $\bfDelta \in \cQ^{n-1}_{\ge 0}$ of the form
$$
\bfDelta = \left( \frac{a_1}{b_1},\cdots,\frac{a_{n-1}}{b_{n-1}}\right),\ \mbox{ with }
a_i \in \cN, b_i \in \cN_*,\ \gcd(a_i,b_i) = 1
$$
Let $c$ be the least common multiple of $b_1,\ldots,b_{n-1}$ and
$Q(z)$ be the one-variable polynomial
$$
Q(z) = P(1,\ldots,1,z) = * z^{p_n} + \cdots
$$
where $*$ denotes some nonzero coefficient.
\begin{Proposition}\label{prop-pnmult}
The multiplicity $\mu_\xi(Q)$ of
the polynomial $Q(z)$ at a point $\xi \ne \cC$ is at most equal to $\linteg p_n/c \rinteg$.
\end{Proposition}
\begin{proof}
Let $\bfp^1,\ldots,\bfp^k$ denote the points of intersection of the straight line
$r(t)$ with the lattice $\cN^n$, ordered according to the last
coordinate
(so that $\bfp^k = \bfp$).
The Lemma~\ref{lemma-multpolyn} implies that the multiplicity
$\mu_\xi(Q)$ is at most equal to $k$.

The result now follows immediately by noticing that each point
$\bfp^{s}$
is necessarily given by
$\bfp^{s} = \bfp + (k - s) c (\bfDelta,-1)$.
\end{proof}
\begin{Corollary}\label{corollary-pnmult}
Suppose that $p_n \ge c+1$ and
that $1 \le a_i < b_i$ for some
$i \in \{1,\ldots,n-1\}$.
Then,
$$
\mu_\xi(Q) \le p_n - \frac{b_i}{a_i},
$$
for all $\xi \ne 0$.
\end{Corollary}
\begin{proof}
Let us prove that the condition $p_n - b_i/a_i \ge
\linteg p_n/c \rinteg$
is satisfied.
Since $b_i \le c$, it is clearly satisfied if
\begin{equation}\label{ineqpn}
p_n \ge \frac{c^2}{c-1}
\end{equation}
Now we use that $p_n \ge c+1$ and $c \ge b_i \ge 2$.
It follows that the inequality
(\ref{ineqpn}) is immediately satisfied when $p_n > c+1$.
For $p_n = c+1$, we compute
$$
\left \linteg \frac{p_n}{c} \right \rinteg = \left \linteg 1 + \frac{1}{c} \right \rinteg = 1 =
p_n - c \le p_n - \frac{b_i}{a_i}
$$
This concludes the proof.
\end{proof}

\section{Appendix C: Comments on Final Models}\label{sect-Comments}
In this appendix, we shall indicate some possible refinements of our Main
Theorem.  First of all, we introduce the notion of {\em strongly elementary vector
field}.

We use the following notation:  Given a matrix $A \in
\mathrm{Mat}(n,\cR)$ and a formal map $R = (R_1,\ldots,R_n) \in \For^n$, the symbol
$$[A \bfx + R]  \frac{\partial}{\partial \bfx}$$
denotes the formal vector field $\displaystyle{\sum_{i=1}^n{\left[ \< A_i, \bfx \> +
    R_i\right] }
\frac{\partial}{\partial x_i}}$, where $A_i$ is the $i^{th}$-row of the matrix
$A$.

Let $\iota \subset [n,\ldots,1]$ be a sublist of indices and
$$\vD = \bigcup_{i \in \iota}{\{x_i = 0\}}$$
be the corresponding
divisor of coordinate hyperplanes in $\cR^n$.

We say that a
formal $n$-dimensional vector field $\eta$ is $\vD$-preserving if it
can be written in the form
$$
\eta = \sum_{i\in \iota}{a_i x_i \frac{\partial}{\partial x_j}} + \sum_{j \in
  [n,\ldots,1] \setminus \iota}{a_j \frac{\partial}{\partial x_j}}$$
where $a_1,\ldots,a_n \in \For$ are formal series.

A formal $n$-dimensional vector field $\eta$ is called a {\em
$\vD$-final model} if $\eta$ is $\vD$-preserving and has one of the
following expressions:
\begin{itemize}
\item[(1)] Non-singular vector field:
$$
\eta = (\lambda + r(\bfx)) \frac{\partial}{\partial x_1}
$$
for some nonzero constant $\lambda \in \cR^*$ and a germ $r \in \For$
with $r(0) = 0$.
\item[(2)] Singular vector field: There exists a decomposition of
  $\cR^n$ into a cartesian product
$$\bfx = (\bfx_+,\bfx_-,\bfx_\IM,\bfx_0) \in \cR^{n_+}\times
\cR^{n_-} \times \cR^{n_\IM} \times \cR^{n_0}$$
with $n_+ + n_- + n_\IM \ge 1$, such that $\eta$
can be written as
$$
\eta = \left[J_+ \bfx_+ + R_+(\bfx)\right] \frac{\partial}{\partial
  \bfx_+} +
\left[J_- \bfx_- + R_-(\bfx)\right] \frac{\partial}{\partial \bfx_-}
+ \left[J_\IM \bfx_\IM + R_\IM(\bfx)\right] \frac{\partial}{\partial \bfx_\IM}
+
R_0(\bfx) \frac{\partial}{\partial \bfx_0}
$$
and the following conditions hold:
\begin{itemize}
\item[(i)] $(J_+,J_-,J_\IM) \in \mathrm{Mat}(n_+,\cR)\times
\mathrm{Mat}(n_-,\cR)\times
\mathrm{Mat}(n_\IM,\cR)$
are matrices whose eigenvalues are all nonzero and have strictly positive real part,
strictly negative real parts and zero real part, respectively.
\item[(ii)] $R_* \in \For^{n_*}$ is a formal germ such that $R_*(0) =
  DR_*(0) = 0$ (for $* \in
  \{+,-,\IM,0\}$) . Moreover,
$$R_+|_{\bfx_+ = 0} = 0, \quad R_-|_{\bfx_- = 0} = 0, \quad
 R_\IM|_{\bfx_\IM = 0} = 0,$$
and
$$
R_0|_{\bfx_- = \bfx_0 = 0} = R_0|_{\bfx_+ = \bfx_0 = 0} = R_0|_{\bfx_\IM
  = \bfx_0 = 0} = 0$$
As a consequence, the eigenspaces $W_+,W_-,W_\IM$ and $W_0$
which correspond respectively to $J_+,J_-,J_\IM$ and the zero matrix
are (formal) invariant manifolds for $\eta$.
\item[(iii)] The zero set $Z = \{\eta = 0\} \subset W_0$ has normal
crossings (i.e.\ it is given by a finite union of intersections of
coordinate hyperplanes).
\item[(iv)] The restricted vector field $\eta_0 = \eta|_{W_0}$ has the form
$$\eta_0 = \mu\ \bfx_0^{\bfalpha} U(\bfx)\, \widetilde{\eta}, \quad\mbox{ with }
\bfx_0^{\bfalpha} = \prod_{i=0}^{n_0}{\bfx_{0,i}^{\, \alpha_i}}$$
where $\mu \in \cR$ is a real constant, $\bfalpha \in \cN^{n_0}$ is a
vector of natural numbers, $U\in \For$ is a
unit and $\widetilde{\eta}$ is a $n_0$-dimensional
$(\vD \cap W_0)$-final model.
\end{itemize}
\end{itemize}
In other words, item (iv) requires that the {\em restriction} of $\eta$ to the manifold $W_0$ is given (up to multiplication by a unit) by a monomial times a vector field $\eta_0$ which is a final model on a space of {\em strictly lower dimension}.

An analytic vector field $\chi$ defined on $M$ is {\em
  $\vD$-strongly elementary} at a point $p \in M$ if there exists a
$\vD$-adapted {\em formal coordinate system} $\bfx = (x_1,\ldots,x_n)$ at $p$
such that $\chi$, written in these coordinates,
is a $\vD$-final model.
\begin{Remark}
We can not replace the words {\em formal coordinate system} by {\em analytic
coordinate system} in the above
definition. It would be too restrictive.  For instance, it would imply that the local center
manifolds are necessarily analytic. 
\end{Remark}
A singularly foliated manifold $\cM = (M,\mainlist,\vD,\linef)$ will be called {\em strongly
elementary} if for each point $p \in M$, the line field $L$ is locally
generated by a vector field $\chi_p$ which is $\vD$-strongly elementary.
\begin{Conjecture}
Let $\chi$ be a reduced analytic vector field defined in a real analytic 
manifold $M$ without boundary.
Then, for each relatively compact set $U \subset M$, 
there exists a finite sequence of weighted blowing-ups
\begin{equation}
(U,\emptyset,\emptyset,\linef_\chi|_U) =: \cM_0 \stackrel{\Phi_1}{\longleftarrow} \cM_1 
\stackrel{\Phi_2}{\longleftarrow}
\cdots \stackrel{\Phi_n}{\longleftarrow} \cM_n 
\end{equation}
such that the resulting singularly foliated manifold $\cM_n$ is strongly
elementary.
\end{Conjecture}
Let us see a few examples of final models in dimensions 1,2 and 3.
\begin{Example}\label{example_finalmod1}
For $n = 1$, the complete list of final models is the following:
\begin{itemize}
\item Non-singular case:
$$
\eta = (\lambda + r(x)) \frac{\partial}{\partial x}
$$
where $r(0) = 0$ and $\lambda \in \cR^*$.
\item Singular case:
$$
\eta = (\lambda x + x r(x)) \frac{\partial}{\partial x},
$$
where $r(0) = 0$ and $\lambda \in \cR^*$.
\end{itemize}
(note that the former case only occurs if $\vD = \emptyset$).
\end{Example}
\begin{Example}\label{example_finalmod2}
For $n = 2$ and $\vD = \emptyset$, the complete list of final models is the following:
\begin{itemize}
\item Nonsingular case:
$$\eta = (\lambda +  r(\bfx))\frac{\partial}{\partial x_1}$$
where $r(0)=0$ and $\lambda \in \cR^*$.
\item Singular case with $n_+ = 1$, $n_- = 1$:
$$\eta = (\lambda_1 x_1 + x_1 r_1(\bfx))\frac{\partial}{\partial x_1} +
(-\lambda_2 x_2 + x_2r_2(\bfx)) \frac{\partial}{\partial x_2}$$
where $\lambda_1,\lambda_2 \in \cR_{>0}$ and $r_i(0) = 0$ for $i = 1,2$.
\item Singular case with $n_\pm = 2$:
$$\eta = (\pm \lambda_1 x_1 + R_1(\bfx))\frac{\partial}{\partial x_1} +
(\pm \lambda_2 x_2 + R_2(\bfx)) \frac{\partial}{\partial x_2}$$
where $\lambda_1,\lambda_2 \in \cR_{>0}$ and $R_i(0) = DR_i(0) = 0$, for
$i = 1,2$.
\item Singular case with $n_{\pm} = 1$, $n_0 = 1$:
$$\eta = (\pm \lambda_1 x_1  +  x_1 r_1(\bfx))
\frac{\partial}{\partial x_1} + \mu x_2^{\alpha} U(\bfx) (\lambda_2  + r_2(\bfx))\frac{\partial}{\partial x_2}
+ $$
where $\lambda_1 \in \cR_{>0}$, $\lambda_2 \in \cR^*$, $\mu \in \cR$,
$\alpha \ge 2$, $U$ is a unit and $r_i(0) = 0$ for $i = 1,2$.
\item Singular case with $n_{\pm} = n_0 = 0$, $n_\IM = 2$:
$$\eta = (\lambda x_2  +  R_1(\bfx))
\frac{\partial}{\partial x_1} + (-\lambda x_1 + R_2(\bfx) )\frac{\partial}{\partial x_2}
$$
where $\lambda \in \cR_{>0}$ and $R_i(0) = DR_i(0) = 0$, for
$i = 1,2$.
\end{itemize}
\end{Example}
\begin{Example}
For $n = 3$, $\vD = \emptyset$, $n_0 = 1$, $n_+ = n_- = 1$, the final model is given by
$$
\eta = (\lambda_1 x_1 + x_1 r_1(\bfx))\frac{\partial}{\partial x_1} +
(-\lambda_2 x_2 + x_2 r_2(\bfx))\frac{\partial}{\partial x_2} + \mu
x_3^{\alpha} U(\bfx) (\lambda_3 + r_3(\bfx)) \frac{\partial}{\partial x_3}
$$
where $\mu \in \cR$, $\lambda_1,\lambda_2 \in \cR_{>0}$, $\lambda_3 \in
\cR^*$, $\alpha \ge 2$, $U$ is a unit and $r_i(0) = 0$ for $i = 1,2,3$.
\end{Example}
\begin{Example}
For $n = 3$, $\vD = \emptyset$, $n_\pm = 0$, $n_0 = 1$ and $n_\IM =
2$,
the final model is given by
$$
\eta = (\lambda x_2 + x_1 r_1(\bfx) + x_2 s_1(\bfx))\frac{\partial}{\partial x_1} +
(-\lambda x_1 + x_1 r_2(\bfx) + x_2 s_2(\bfx))\frac{\partial}{\partial x_2} +
 \mu
x_3^{\alpha} U(\bfx) (\lambda_3 + r_3(\bfx)) \frac{\partial}{\partial x_3}
$$
where $\mu \in \cR$, $\lambda \in \cR_{>0}$, $\lambda_3 \in
\cR^*$, $\alpha \ge 2$, $U$ is a unit and $r_i(0) = s_i(0) = 0$ for $i = 1,2,3$.
\end{Example}
\begin{Example}
For $n = 3$, $\vD = \emptyset$, $n_0 = 2$ and $n_+ = 1$, the final model is given by
$$
\eta = (\lambda_1 x_1 + x_1 r_1(\bfx))\frac{\partial}{\partial x_1} +
\mu x_2^{\alpha}x_3^{\beta} U(\bfx) \left[ \sum_{i = 2}^3{(f_i(x_2,x_3) +
  x_1 r_i(\bfx))\frac{\partial}{\partial x_i}}\right]
$$
where $\mu \in \cR$, $\lambda_1 \in \cR_{>0}$, $\alpha + \beta \ge 1$, $U$
is a unit and the vector field obtained by restriction to the center
manifold $W_0 = \{x_1 = 0\}$, namely
$$\widetilde{\eta} = f_2(x_2,x_3) \frac{\partial}{\partial x_2}
+ f_3(x_2,x_3) \frac{\partial}{\partial x_3}$$
has one of the forms given in the example~\ref{example_finalmod2}.
\end{Example}

\vspace{1cm}

\begin{flushright}
{\sc Daniel Panazzolo}\\
Instituto de Matem\' atica e Estat\' \i stica\\
Universidade de S\~ ao Paulo\\
Rua do Mat\~ ao 1010 - S\~ ao Paulo - SP \\
05508-090 - Brazil\\
{\tt dpanazzo@ime.usp.br}
\end{flushright}
\end{document}